\documentclass[english,11pt]{elsarticle} 

\usepackage[T1]{fontenc}
\usepackage{float}
\usepackage{amsmath}
\usepackage{amsthm}
\usepackage{amssymb}
\usepackage{amsfonts}
\usepackage{bbm}
\usepackage{graphicx}
\usepackage{epstopdf}
\usepackage[pdftex,colorlinks,linkcolor={blue},citecolor={blue}]{hyperref}

 \usepackage{pgfplots}
  \pgfplotsset{compat=newest}
  \usetikzlibrary{plotmarks}
  \usetikzlibrary{arrows.meta}
  \usepgfplotslibrary{patchplots}
  \usepackage{grffile}
  \usepackage{amsmath}
   \usepackage{caption}
\usepackage{subcaption}

\usepackage{enumitem}
\usepackage{epsfig}
\usepackage{amssymb,amsmath,amsfonts,layout,graphicx}
\usepackage{makeidx}
\usepackage{babel}
\usepackage{tikz}
\usepackage{sublabel}
\usepackage{cases}
\usepackage{algorithm}
\usepackage{algorithmic}

\usepackage
[
        letterpaper,
        left=1.91cm,
        right=1.91cm,
        top=1.91cm,
        bottom=2cm,
]
{geometry}

\newtheorem{lemma}{Lemma}
\newtheorem{assumption}{Assumption}
\newtheorem{proposition}{Proposition}

\newtheorem{remark}{Remark}

\newcommand{\x}{{\bf x}}

\usepackage{enumitem}

\newlist{myItem}{itemize}{3}
\setlist[myItem]{nosep,label=\protect\mpbullet,topsep=0pt}

\begin{document}
\begin{frontmatter}
 \title{\textbf{Charging Autonomous Electric Vehicle Fleet for Mobility-on-Demand Services: Plug in or Swap out?}}

\author[]{Jing Gao}
\ead{jgaoax@connect.ust.hk}
\author[]{Sen Li}
\ead{cesli@ust.hk}
\address[]{Department of Civil and Environmental Engineering, The Hong Kong University of Science and Technology}

\begin{abstract}
This paper compares two prevalent charging strategies for electric vehicles, plug-in charging and battery swapping, to investigate which charging strategy is superior for electric autonomous mobility-on-demand (AMoD) systems. To this end, we use a queueing-theoretic model to characterize the vehicle waiting time at charging stations and battery swapping stations, respectively. The model is integrated into an economic analysis of the electric AMoD system operated by a transportation network company (TNC), where the incentives of passengers, the charging/operating shift of TNC vehicles, the operational decisions of the platform, and the planning decisions of the government are captured. Overall, a bi-level optimization framework is proposed for charging infrastructure planning of the electric AMoD system. Based on the proposed framework, we compare the socio-economic performance of plug-in charging and battery swapping, and investigate how this comparison depends on the evolving charging technologies (such as charging speed, battery capacity, and infrastructure cost). At the planning level, we find that when choosing plug-in charging, increased charging speed leads to a transformation of infrastructure from sparsely distributed large stations to densely distributed small stations, while enlarged battery capacity transforms the infrastructure from densely distributed small stations to sparsely distributed large stations. On the other hand, when choosing battery swapping, both increased charging speed and enlarged battery capacity will lead to a smaller number of battery swapping stations. At the operational level, we find that improved charging speed leads to increased TNC profit when choosing plug-in charging, whereas improved charging speed may lead to smaller TNC profit under battery swapping. The above insights are validated through realistic numerical studies.

\end{abstract}


\end{frontmatter}

\section{Introduction}
Automation, electrification, and shared mobility are dominating trends in developing an efficient and sustainable mobility future. Ride-hailing platforms, such as Uber, Lyft, Waymo and DiDi, are working at the forefront to realize this vision. For instance, Waymo has been providing autonomous rides in the East Valley of Phoenix since 2020 and is planning to expand its self-driving ride-hailing services to San Francisco in 2022 \cite{Waymo2022SF}. DiDi, as the world's largest shared electric vehicle (EV) network, has more than one million registered EVs, accounting for 20\% of China's EV mileage \cite{Brittany2020Didi}. With the advancement of autonomous driving and battery technologies, it is envisioned that shared autonomous electric vehicles (SAEVs) will be operated by Transportation Network Companies (TNC) to provide electric autonomous mobility-on-demand (AMoD) services in future cities. The economic and environmental benefits of electric AMoD services can be substantial. It is estimated that: (a) each 200-mile range SAEV could replace 5.5 privately owned cars under Level-2 charging  \cite{chen2016operations}; (b) the cost of electric AMoD services will be \$0.29-\$0.61 per revenue mile \cite{bauer2018cost}, an order of magnitude lower than that of taxis; (c) the adoption of SAEVs can reduce greenhouse emissions by 73\% and energy consumption by 58\% \cite{jones2019contributions}.

Despite the aforementioned benefits, it remains unclear how to plan the charging infrastructure so that the charging needs of the commercial SAEV fleet can be most efficiently accommodated in future cities. There was a heat debate on whether we should prioritize plug-in charging or battery swapping for EVs. Each of the two charging strategies has its advantages and disadvantages. For instance, plug-in charging has a relatively low infrastructure cost, which can be conveniently installed and easily scaled up. However, typical Level-1 and Level-2 chargers require a few hours before charge completion, which incurs significant inconvenience for EV drivers. One way to overcome this limitation is by significantly accelerating the charging speed of EV batteries. As of today, the Tesla Supercharger can recharge up to 200 miles within 15 minutes \cite{Tesla2022Supercharger}, making it ideal for battery top-up in long-distance highway trips. Compared to fast-charging, an alternative approach is to swap the empty battery out of the EV and replace it with a fully-charged battery that is proactively stored in the battery swapping station. Battery swapping can eliminate the need for waiting, making it particularly popular among impatient customers: as one of the most successful Chinese EV companies, Nio has built over 900 battery swapping stations in 59 cities and executed more than 8 million battery swaps for its customers \cite{Kane2022Nio}. However, the major obstacle to the widespread adoption of battery swapping is that it requires standardization of batteries. This can be quite difficult because manufacturers can hardly agree on a uniform size of EV batteries. In addition, battery swapping also requires EV owners to widely accept the concept of battery-as-a-service, where the ownership of vehicles and batteries are separated so that EV drivers can willingly replace their empty batteries. Given these obstacles, many Western countries have chosen to primarily invest in fast-charging technologies instead of battery swapping stations. 

So far, the comparison between plug-in charging and battery swapping is based on the assumption that most vehicles are human-driven and owned by individual drivers. However, with the advancement of self-driving technology, the landscape of urban mobility can be completely different in the future. In the autonomous-driving era, city residents may no longer owe his/her own car, which would otherwise sit at the parking garage for most of the time. Instead, autonomous vehicles will be shared, and the mobility needs of the city residents will be accommodated by electric AMoD services provided by TNC platforms. These changes in urban mobility have great potential to improve the efficiency of vehicle fleet utilization, lower the travel cost of passengers, and reduce energy consumption and carbon emissions from the transport sector. However, these changes also call for a re-investigation of how to choose the optimal charging strategies for future mobility systems with a large fleet of commercial SAEVs.  In this case, the advantages and disadvantages of plug-in charging and battery swapping should be re-evaluated. On the one hand, while it is unarguable that fast charging still remains competitive for energy top-ups of long-distance highway trips, it is not necessarily a dominating charging option when considering large commercial SAEV fleets in the urban context. This is because the commercial SAEV fleet requires super fast charging as any downtime of vehicles in the commercial fleet cuts in the revenue. However, fast-charging infrastructure requires a much higher voltage and current (e.g., 200 A, 400 V) than its Level-1 (e.g., 12 A, 120 V) and Level-2 (e.g., 32 A, 240 V) alternatives \cite{IEC2021charging}. In the urban context, building a large and densely distributed fast-charging network for SAEVs may incur significant voltage drops and cause instability to the already strained distribution power grids, which are prohibitively expensive to upgrade.
On the other hand, battery swapping remains a competitive option for AMoD systems because the commercial SAEV fleet is owed by the TNC platform, thus it is much easier to enforce battery standardization by using the same type of vehicles, and there is no vehicle-battery ownership concerns when one battery is replaced with another. This naturally raises the following questions: (1) which charging strategy is superior for electric AMoD in the urban context, plug-in charging or battery swapping?  (2) how does the comparison between plug-in charging and battery swapping depend on charging technologies, such as charging speed, battery capacity, and infrastructure cost? (3) what is the optimal plan of charging infrastructure that maximizes social welfare?

This paper aims to address the aforementioned questions by comparing the performance of plug-in charging and battery swapping in AMoD systems from the socio-economic perspective. To this end, we need to carefully examine how commercial vehicles are queued up at charging stations/battery swapping stations, how the charging demand of the vehicles relates to the travel demand of the passengers, and how these relations are coordinated by the TNC platform and influenced by the government. A complete understanding of these complex relations requires the development of an equilibrium model, which captures the dynamics of charging stations and battery swapping stations while simultaneously encoding the intimate interactions between distinct stakeholders in the AMoD market. This paper will offer such a model, based on which we will investigate how different charging strategies (e.g., plug-in charging vs battery swapping) affect the planning decisions of the government, how they affect the operational decisions of the TNC platform, and how these impacts vary under different charging technologies (such as charging speed, battery capacity, infrastructure cost, etc). The major contributions of this paper are summarized below:
\begin{itemize}
     \item We compare the socio-economic performance of plug-in charging and battery swapping for electric AMoD systems. To this end, we use a queueing-theoretic model to characterize the waiting time and the blocking probability at charging stations and battery swapping stations and integrate the model into an economic analysis of the AMoD market, where the incentives of passengers, TNC platforms, and the government are captured, and their interactions are characterized. Overall, a bi-level optimization framework is proposed for optimal planning of charging infrastructures. In the upper level, the government determines the infrastructure deployment plan to maximize social welfare. In the lower level, the TNC determines the fleet size and the ride fare to maximize its profit subject to market equilibrium constraints. {\em To our best knowledge, this is the first work that jointly considers the economics of the AMoD market, the planning decisions of the government, and the comparison between plug-in charging and battery swapping for the electric AMoD system}.
     
    \item We evaluate how the planning decisions of the government depend on the choice of charging strategies and the evolution of charging technologies. Based on our model, we find that when choosing plug-in charging, increased charging speed results in a transformation of charging deployment from {\em sparsely distributed large stations} to {\em densely distributed small stations}, while enlarged battery capacity transforms the infrastructure deployment from {\em densely distributed small stations} to {\em sparsely distributed large stations}. On the other hand, when choosing battery swapping, both increased charging speed and enlarged battery capacity lead to a smaller number of battery swapping stations. We identify the reason for this difference and show that the change in planning decisions for plug-in charging arises from the trade-off between the vehicle searching time before finding a charging station and the vehicle waiting time after arriving at the charging station.
    
    \item We evaluate how the operational decisions of the TNC platform depend on the choice of charging strategies and the evolution of charging technologies. We find that when choosing plug-in charging, improved charging speed leads to increased TNC profit. However, when choosing battery swapping, the TNC profit may reduce if the charging speed is relatively high. We point out that this can be attributed to the smaller degree of flexibility when deploying battery swapping stations compared to plug-in charging stations.
\end{itemize}

\section{Related Works}

This paper uses a queueing-theoretic model to characterize the vehicle waiting time at charging stations and battery swapping stations, which is further incorporated into the economic analysis and the charging infrastructure planning of the electric AMoD system. Thus, we review related works from the following three perspectives: (1) queueing models for EV charging, (2) economic analysis of electric AMoD systems, and (3) charging infrastructure planning of electric AMoD systems.

\subsection{Queueing models for EV charging}

Queueing-theoretic models are extensively applied to characterize the vehicle charging process at charging stations. Some works characterize congestion at the charging stations and/or the equilibrium choices of the charging demand. For instance, \cite{liang2014plug} developed a queueing network model to estimate the charging demand of plug-in EVs, in which each charging station is modeled as a service center with multiple chargers as servers, and EVs are modeled as customers requesting charging services. \cite{zhang2016platoon} analyzed the performance of the EV platoons charging at renewable energy-supplied stations. The steady state of the queueing system was derived based on the equilibrium equations. \cite{zhang2018optimal} formulated a dual-mode charging station as a queueing network with multiple servers and heterogeneous service rates. The optimal pricing scheme is designed to guide the charging processes of EVs to minimize the service drop rate. \cite{liu2022electric} adopted the M/D/C queueing model to capture the congestion in charging stations and proposed an EV charging station access equilibrium model to characterize EV user equilibrium choices in charging stations.
Aside from the aforementioned works, some other works also considered the impacts of charging stations on the power grid. For instance, \cite{aveklouris2017electric} established an equilibrium queueing model to characterize the vehicle charging process that considers both the congestion in the distribution grid and the congestion due to the limited number of chargers/charging slots within charging stations. \cite{esmailirad2021extended} considered the public charging stations with both grid-to-vehicle and vehicle-to-grid and developed an extended M/M/K/K queueing model to analyze the charging/discharging process at charging stations with $K$ charge/discharge cord plugs.

Queueing models are also developed to characterize EV charging at battery swapping stations. The mechanisms of battery swapping station is more complicated than that of the charging stations. For instance, \cite{tan2014queueing} proposed a mixed queueing network model to analyze the battery swapping station for electric vehicles, which is comprised of an open queue of EVs and a closed queue of batteries. The equilibrium equations of the queueing systems are formulated, and the steady-state distribution is derived. Based on the novel queueing network model, \cite{sun2017optimal} formulated the charging operation problem of the battery swapping station as a constrained Markov decision process and derived the optimal charging policy that minimizes the charging cost while guaranteeing the quality of service. Likewise, based on the mixed queueing model, \cite{tan2018asymptotic} conducted extensive simulations to evaluate the performance of the battery swapping station and charging station using the blocking probability of electric vehicles as the quality-of-service metric. Different parameters of the battery swapping and charging station are examined, such as the number of parking spaces, swapping islands, chargers, and batteries. \cite{imani2016modeling} built a two-stage priority queueing model to capture the queueing effect of swapping and onboard charging processes. \cite{sloothaak2021complete} introduced a closed Markovian queueing network model to represent the evolution of the battery population within a city and proposed a load-balancing policy to achieve the optimal trade-off between EV users' quality-of-service and operational costs.

{\em Our paper differs from all aforementioned works as we integrate the queueing-theoretic model into the economic analysis of AMoD systems. In our paper, the main focus is on characterizing how the dynamics of the charging station and battery swapping station affect the incentives of passengers, platforms, and the social planner in the AMoD market. These elements are missing in the above-mentioned works. 
}

\subsection{Economic analysis of electric AMoD systems}

A large body of works considered temporal dynamics and focused on the dynamic control of electric AMoD systems. Generally, they considered profit-maximizing AMoD operators and focused on the real-time charging scheduling \cite{zhang2016model,tucker2019online,iacobucci2019optimization,boewing2020vehicle,turan2020dynamic,ding2021integrated,ding2022quality}, routing \cite{zhang2016model,iacobucci2019optimization,boewing2020vehicle,turan2020dynamic}, rebalancing \cite{tucker2019online,ding2021integrated} and pricing \cite{turan2020dynamic,ding2021integrated,ding2022quality,ni2021dynamic} of electric AMoD systems. For instance, \cite{zhang2016model} presented a model predictive control (MPC) algorithm to optimize the vehicle scheduling and routing in the AMoD system with electric vehicle charging constraints. \cite{tucker2019online} investigated an online charge scheduling strategy for the electric AMoD system, in which autonomous electric vehicles are scheduled to charge between AMoD rides and dispatched to the next passenger pick-up locations. \cite{iacobucci2019optimization} optimized the routing, relocation, and charge scheduling of the SAEV fleet at two different time scales by running two MPC optimization algorithms in parallel. \cite{boewing2020vehicle} considered the joint vehicle coordination and charge scheduling problem for electric AMoD systems where a mixed-integer linear program was presented to account for the battery level of vehicles and the energy availability in the power grid along the time horizon. \cite{turan2020dynamic} examined the joint routing, battery charging, and pricing problem for a profit-maximizing electric AMoD operator considering the randomness in trip demand, renewable energy availability, and electricity prices. \cite{ni2021dynamic} designed a pricing mechanism in a bi-level framework for the electric mobility-on-demand system to incentivize passengers to choose electric mobility-on-demand services. \cite{ding2021integrated} proposed a combined operation scheme for the battery swapping station and AMoD system, where an expanded network flow model was developed to determine the dynamic swapping schedules and vehicle rebalancing for profit maximization of the AMoD system. \cite{ding2022quality} considered the interaction between the AMoD fleet and the battery swapping stations. A bi-level optimization problem was formulated to capture the interdependence between the AMoD system and the battery swapping system, wherein at the upper level, the battery swapping operator determines the time-varying and location-varying swapping price to minimize its cost, and at the lower level, the AMoD operator develops fleet scheduling strategies to maximize its profit.

Another strand of studies focused on the steady state and developed economic equilibrium models to investigate the implications of electric AMoD systems. For instance, \cite{turan2019smart} studied the potential benefits of smart charging for electric AMoD services, where a profit-maximizing operator makes decisions on routing, charging, rebalancing, and pricing for AMoD rides based on a network flow model. \cite{rossi2019interaction} studied the interaction between the AMoD system and the electric power network. A joint optimization model was proposed to capture vehicles' charging requirements, time-varying customer demand, battery depreciation, and power transmission constraints. It was proved that the socially optimal solution to the joint problem is a general equilibrium if locational marginal pricing is used for electricity. \cite{turan2021competition} investigated the impacts of competition in electric AMoD systems. A network-flow model was formulated to determine the optimal strategies for profit-maximizing electric AMoD operators in monopoly and duopoly markets. The benefits of introducing competition in the market were demonstrated by distinct metrics, including the prices of rides, aggregate demand served, profits of the firms, and consumer surplus. \cite{bang2021congestion} investigated the problem of routing, rebalancing, and charging for electric AMoD systems concerning traffic congestion at the macroscopic level. A volume-delay function was introduced to characterize traffic congestion, and the routing and rebalancing for vehicles with energy constraints and charging constraints was formulated as an optimization problem.

{\em However, all aforementioned studies considered the electric AMoD system under exogenous charging infrastructures, which neglects charging infrastructure planning and its impacts on the operation of electric AMoD systems. In addition, the majority of these work focused on charging stations, while only a few considered battery swapping. Distinct from these works, we characterize the interactions between the AMoD operator and the planner under both charging stations and battery swapping stations. This enables us to compare the two charging strategies in terms of both planning decisions and operational decisions, which have not been studied before in the aforementioned literature.}

\subsection{Charging infrastructure planning of electric AMoD systems}

Numerous studies have investigated the planning of infrastructure for both charging stations \cite{he2013optimal,he2015deploying,zhu2018charging,kavianipour2021electric,bauer2019electrifying,ma2021optimal} and battery swapping stations \cite{mak2013infrastructure,yang2017optimal,sultana2018placement,sun2019optimal,liang2021configuration}. However, most of these works focused on contemporary transportation systems, where the majority of vehicles are privately-owned and driven by individuals. On the other hand, only a handful of work considers the charging infrastructure planning for electric AMoD systems.
For instance, \cite{loeb2018shared} investigated the operation of SAEVs with charging infrastructure development based on agent-based simulations. Different characteristics of SAEV fleets, including the fleet size, vehicle range, and charge time, were examined, and the simulation results indicated that the number of stations needed for the operation of SAEV fleets depends almost wholly on vehicle range. \cite{vosooghi2020shared} explored the impacts of charging station placement, charging types, and vehicle battery capacity on the efficiency of SAEV services through agent-based simulations.  Their results suggested that the performance of SAEV service can be significantly improved by offering battery swapping infrastructure. \cite{lokhandwala2020siting} proposed a framework to optimize charging infrastructure development for electric vehicle adoption in shared autonomous fleets. Several factors were considered, including the level of autonomous vehicle adoption, ridesharing participation, vehicle queueing at charging stations, and the tradeoff between building new charging stations and expanding existing ones. \cite{luke2021joint} considered the joint optimization of autonomous electric vehicle fleet operations and charging station siting. A linear program was introduced to jointly optimize the charging station siting and macroscopic fleet operations. They showed that small-sized EVs with low procurement costs and high energy efficiencies are the most cost-effective. 


{\em However, the aforementioned studies on AMoD \cite{loeb2018shared,vosooghi2020shared,lokhandwala2020siting,luke2021joint} primarily focus on plug-in charging (instead of battery swapping), and the majority of them \cite{loeb2018shared,vosooghi2020shared,lokhandwala2020siting} deploys agent-based simulation for given passenger demands, neglecting the incentive of passengers, which is a crucial aspect of the ride-hailing market. In contrast, our work differs from these studies as we consider the comparison between plug-in charging and battery swapping, the incentives of decision-makers in the electric AMoD market, and the planning decisions of the government.
}

\section{Economic Equilibrium of the Electrified AMoD Market}

Consider a TNC platform that operates a fleet of autonomous electric vehicles to provide AMoD services. The platform deploys a specific charging strategy and determines the ride fare and the fleet size to maximize its profit subject to the availability of charging infrastructure. In response to the platform's pricing and charging strategy, passengers and the platform make interactive decisions, which constitute the market equilibrium. In this section, we develop an economic equilibrium model to characterize the ride-sourcing market equilibrium with electric autonomous vehicles. The incentives of passengers will be captured, and the dynamic charging process of vehicles will be characterized.



\subsection{Passenger incentives} \label{passenger_incentives}
Passengers decide whether to choose AMoD services based on the total travel cost of the trip. We define the travel cost of AMoD as the weighted sum of {\em average} ride fare $p_f$ and {\em average} waiting time $w^c$:
\begin{equation} \label{generalized_travel_cost}
    c=p_f+\alpha w^c ,
\end{equation}
where $p_f$ is the average trip fare for the AMoD services, $w^c$ is the average passenger waiting time before pick-up, and $\alpha$ is the value of time which indicates the passengers' perception between time and money. 
We assume that the average passenger arrival rate $\lambda$ depends on the generalized travel cost $c$ and define the following passenger demand function:
\begin{equation} \label{passenger_demand_function}
    \lambda = \lambda_0 F_p(p_f+\alpha w^c) ,
\end{equation}
where $\lambda_0$ is the arrival rate of potential passengers (total travel demand in the city); $F_p(\cdot)\in [0,1]$ determines the proportion of potential passengers choosing autonomous ride-hailing services. Note that (\ref{passenger_demand_function}) determines the passenger demand as a function of the average generalized cost, which does not require the cost of distinct travelers to be the same. We further assume that $F_p(\cdot)$ is a strictly decreasing and continuous differentiable function such that a higher travel cost leads to a lower arrival rate of passengers. This includes the well-established logit model as a special case.

The passenger waiting time $w^c$ represents the quality of electric AMoD services, which endogenously depends on the vehicle supply and passenger demand of the ride-hailing market. In ride-hailing services, the passenger waiting time comprises the ride confirmation time (from the ride request being submitted to a vehicle being assigned) and the pickup time (from a vehicle being assigned to passenger pickup). Typically, the ride confirmation time only lasts for several seconds, which is negligible compared to the pickup time ranging from five to eight minutes. Therefore, we neglect the ride confirmation time and approximate the passenger waiting time as the pickup time. Assuming that the platform matches the passenger to the nearest idle vehicle, the pickup time depends on the distance of the nearest idle vehicle to the passenger, which can be further characterized as a monotone function of the average number of idle vehicles $N^I$. We denote the passenger waiting time function as $w^c=F_c(N^I)$ and impose the following assumption:
\begin{assumption}
$F_c(N^I)$ is positive, strictly decreasing with respect to $N^I$, and $\lim_{N^I\to 0} F_c(N^I)=\infty$.
\end{assumption}
The passenger waiting time follows the "square root law", which was well-established and widely applied in street hailing taxi market \cite{douglas1972price}, radio dispatching taxi market \cite{arnott1996taxi}, and online ride-hailing market \cite{li2019regulating}. This leads to the following passenger waiting time function:
\begin{equation} \label{func_passenger_waiting_time}
    w^c=F_c(N^I)=\frac{A}{\sqrt{N^I}} ,
\end{equation}
where $A$ is the scaling parameter capturing possible factors in the matching of idle vehicles and passengers, such as the area of the city, the average traffic speed on the road network, the spatial distribution of passengers/vehicles, etc. The square root law indicates that the average passenger waiting time is inversely proportional to the square root of the number of idle vehicles. The intuition behind (\ref{func_passenger_waiting_time}) is that if both waiting passengers and idle vehicles are randomly distributed across the city, the distance between a passenger and her closest idle vehicle is inversely proportional to the square root of the total number of idle vehicles. A detailed justification of the square root law can be found in \cite{zha2018geometric,li2019regulating}.

\subsection{The operating and charging shifts of TNC vehicles} \label{operating_shifts}

TNC vehicles need to be charged due to the limited battery range. Once the platform decides that a vehicle needs to recharge, it will pause the AMoD services for the vehicle and dispatch it to nearby charging infrastructure for energy top-up (referred to as the "charging shift"). After being fully charged, it resumes to provide AMoD services (referred to as the "operating shift"). Throughout the day, each TNC vehicle regularly switches between the operating shift and the charging shift. On average, the total number of TNC vehicles ($N$) can be decomposed as the sum of two groups: vehicles in the operating shift ($N_1$) and vehicles in the charging shift ($N_2$). This leads to the following relation:
\begin{equation} \label{total_vehicle_conservation}
    N=N_1+N_2 .
\end{equation}
The characterization of $N_1$, $N_2$, and their endogenous relations will be discussed below. 

\subsubsection{The operating shift} 
During the operating shift, autonomous electric vehicles provide ride-hailing services to passengers, in which vehicles experience the following three statuses sequentially and recurrently: (1) status 1: being idle and cruising for passengers; (2) status 2: being assigned with passengers and picking up passengers; (3) status 3: being occupied with passengers and delivering passengers. Let $\tau$ be the average duration of the ride-hailing trips. At the stationary state, the conservation of the total vehicle operating hour yields:
\begin{equation} \label{conservation_operating_vehicles}
    N_1 = N^I +\lambda w^c +\lambda \tau,
\end{equation}
where the first term corresponds to the average number of idle vehicles in the operating shift (status 1), and the second term and third term are derived based on Little's law, which represents the average number of vehicles on the way to pick up passengers (status 2), and the average number of vehicles delivering passengers (status 3), respectively. Let $w^v$ be the average vehicle idle time between rides, i.e., the average duration of status 1. Based on Little's law, the number of idle vehicles $N^I$ relates to $w^v$ through the following equation:
\begin{equation}
    N^I = \lambda w^v
\end{equation}

At the stationary state, the number of vehicles in operating shifts $(N_1)$ and the number of vehicles in charging shifts $(N_2)$ endogenously depend on the charging process of autonomous electric vehicles. Next, we focus on the charging shift of TNC vehicles and derive the relations of $N_1$, $N_2$, and other endogenous variables.

\subsubsection{The charging shift} 
During the charging shift, vehicles are dispatched to nearby charging infrastructures (either charging stations or battery swapping stations) for battery top-up. There may be matching frictions between vehicles and charging infrastructure, so that each vehicle may have to spend some time traveling to the charging infrastructure. Furthermore, upon arrival at the charging infrastructure, each vehicle may also need to wait until a charging device is available for immediate use. Overall, the charging process of autonomous electric TNC vehicles can be divided into three phases: (1) phase 1: searching for nearby charging/battery swapping stations; (2) phase 2: waiting for charging at charging/battery swapping stations; (3) phase 3: being charged at charging/battery swapping stations. Let $N_2^m$, $N_2^w$, and $N_2^s$ be the average number of vehicles in the three phases, and denote $t_m$, $t_w$, and $t_s$ as the average duration of the three phases, respectively. At the stationary state, the conservation of the total number of vehicles in charging shifts yields:
\begin{equation}
    N_2 = N_2^m + N_2^w + N_2^s .
\end{equation}
Besides, based on Little's law, we have the following relations:
\begin{equation} \label{N2m_N2w_N2s}
\begin{cases}
N_2^m = \gamma t_m \\
N_2^w = \gamma t_w \\
N_2^s = \gamma t_s
\end{cases} ,
\end{equation}
where $\gamma$ is the arrival rate of vehicles in charging shifts, which is an endogenous variable representing the charging demand of autonomous electric TNC vehicles. Below we derive the relations among these endogenous variables separately.

To derive $t_m$, we assume that arriving vehicles in the charging shift are dispatched to the nearest charging/battery swapping station that is not full for battery replenishment. Due to the limited capacity, the charging/battery swapping station may be full and have no idle parking spaces for vehicles waiting. In this case, the average vehicle searching time equals the average travel time between the vehicle and the nearest charging/battery swapping station with vacant parking spaces. Let $V$ be the average capacity of each charging/battery swapping station and denote $P_V$ as the probability of a station being full. The average vehicle searching time $t_m$ is inversely proportional to the square root of the number of charging/battery swapping stations that are not full \cite{ahn2015analytical,lee2017optimal,lai2021demand}:
\begin{equation} \label{func_vehicle_searching_time}
    t_m = \frac{B}{\sqrt{K (1-P_V)}} ,
\end{equation}
where $B$ is the scaling parameter, $K$ is the number of charging/battery swapping stations, and $K(1-P_V)$ accounts for the number of charging/battery swapping stations that are not full. It follows the same intuition as the derivation of (\ref{func_passenger_waiting_time}): if both vehicles and charging/battery swapping stations are randomly distributed across the city, the expected distance between a vehicle to the closet available charging/battery swapping station is inversely proportional to the square root of the number of available charging/battery swapping stations. Note that the blocking probability $P_V$ is an endogenous variable that depends on both the charging demand $\gamma$ and the configuration of charging/battery swapping stations.

To derive $P_V$, $t_w$, and $t_s$, we note that they have different relations with other endogenous variables under different charging strategies (i.e., plug-in charging v.s. battery swapping). To reflect these differences, we will characterize the operating mechanism of charging stations and battery swapping stations separately.

\begin{figure}[htbp]
  \centering
  \begin{subfigure}[b]{0.48\textwidth}
    \centering
    \includegraphics[width=\textwidth]{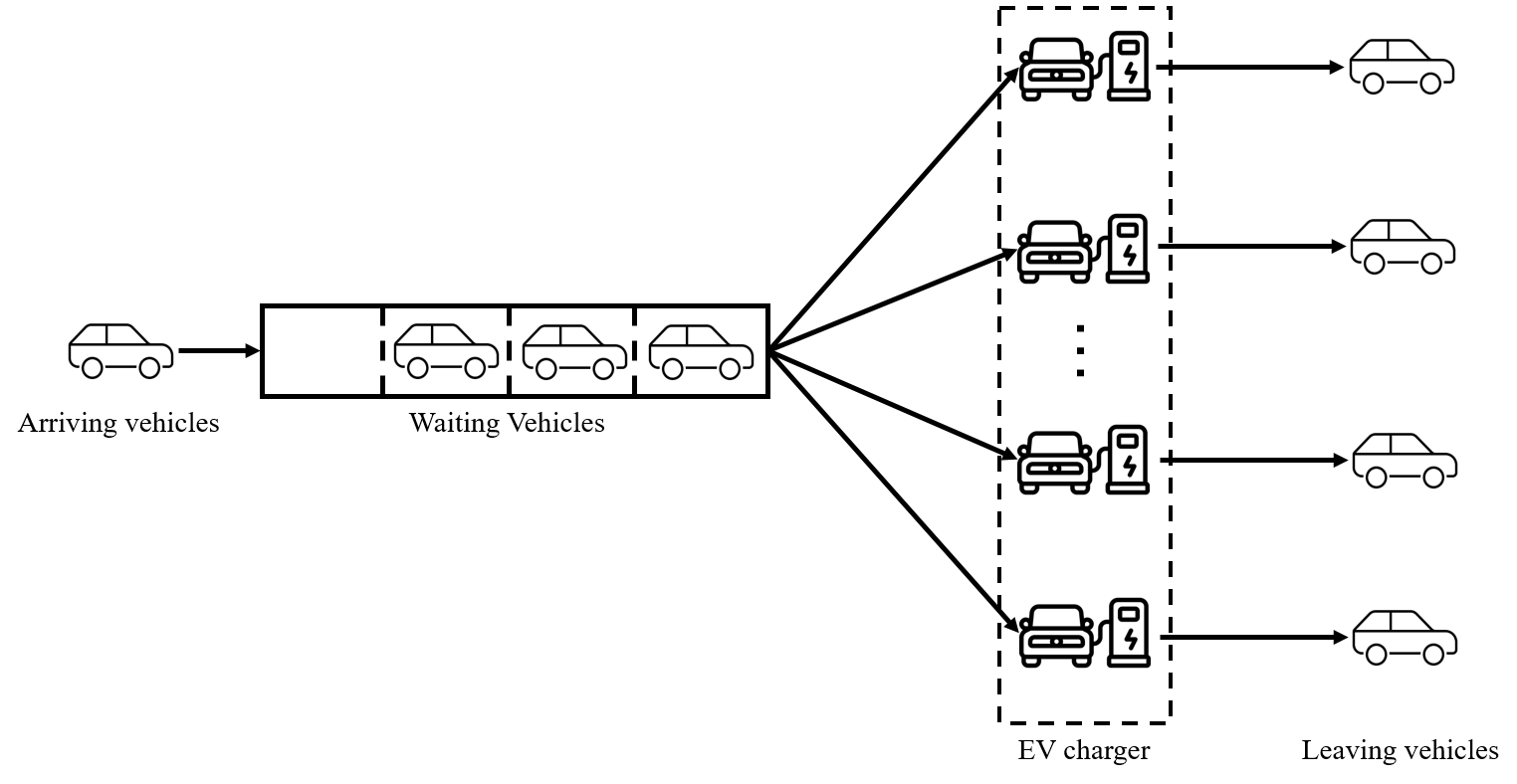}
    \caption{M/M/Q/V queueing model for the charging station.}
    \label{fig:charging_queueing}
  \end{subfigure}
  \hfill
  \begin{subfigure}[b]{0.48\textwidth}
    \centering
    \includegraphics[width=\textwidth]{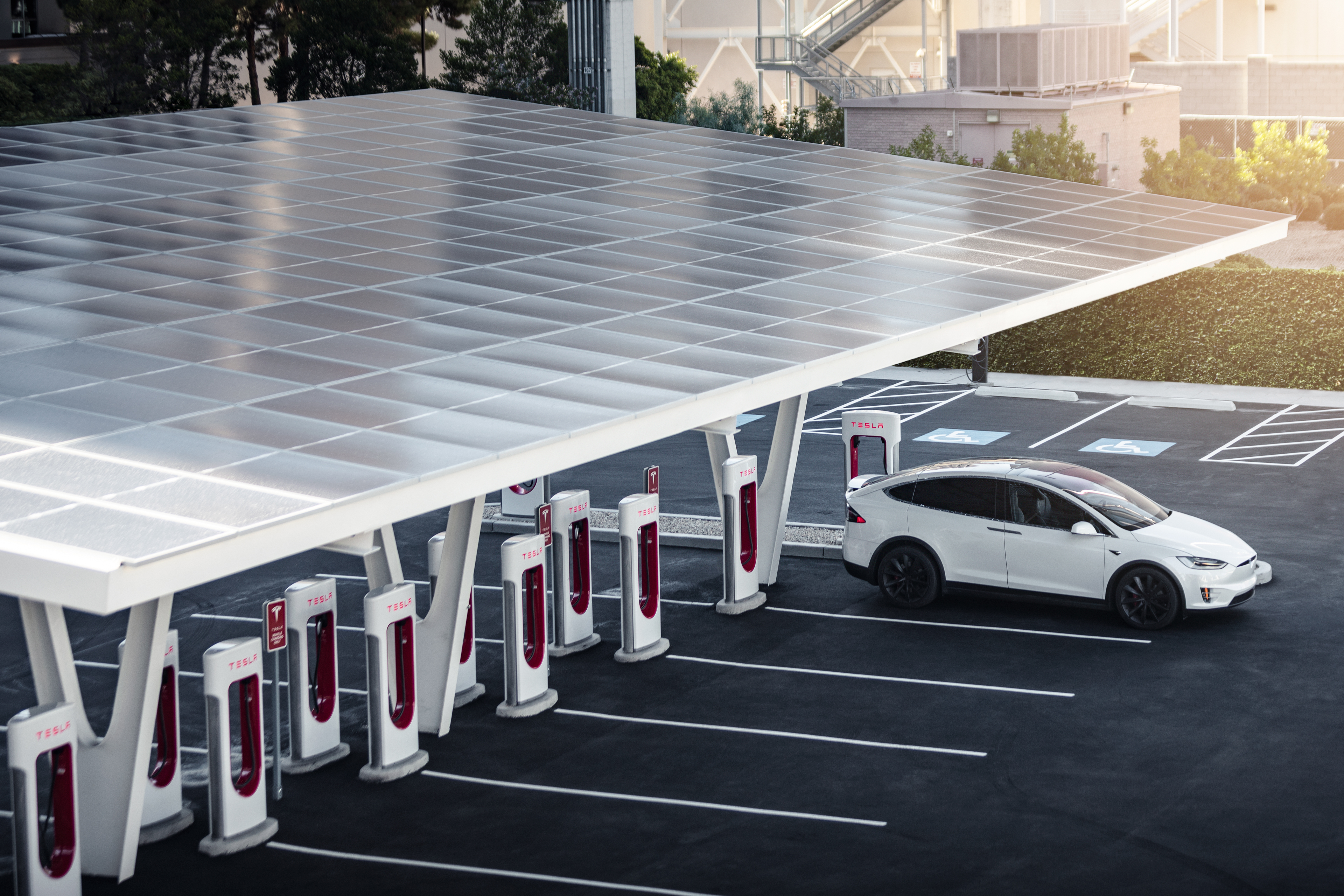}
    \caption[Caption]{The schematic of Tesla's charging station\protect\footnotemark.}
    \label{fig:tesla}
  \end{subfigure}
    \caption{M/M/Q/V model for the charging station and Tesla's charging station.}
  \label{fig:charging_station_illustration}
\end{figure}

{\em Queueing Model for Plug-in Charging}: We use an M/M/Q/V queue model to represent the vehicle charging process at each charging station as shown in Figure \ref{fig:charging_queueing}. Consider a fleet of autonomous electric vehicles with a battery capacity of $C$, which will be dispatched to nearby available charging stations once the battery SoC approaches the minimum SoC. In practice, the platform determines vehicles' charging schedules based on vehicles' real-time operating statuses. Therefore, vehicles arriving at charging/battery swapping stations may have distinct SoCs and different charging times. Without loss of generality, let $\delta$ be the mean battery SoC of vehicles arriving at charging/battery swapping stations\footnote{The mean battery SoC should depends on the charging schedules of the platform, which we regard as exogenous.}. Suppose that there are $K$ charging stations
uniformly distributed across the city, and each station has the capacity $V$ and $Q$ chargers with charging speed $s$ (see Figure \ref{fig:tesla} the schematic of Tesla's charging station). Each charger can be modeled as a server, and arrived vehicles can be regarded as jobs. A charger is `busy' if it is occupied by a vehicle. Otherwise, it is considered as `idle'. We assume that the arrival of TNC vehicles at each charging station follows the Poisson process. The total charging demand of TNC vehicles is $\gamma$, and thus the arrival rate of vehicles at each station is $\gamma/K$. Furthermore, since each vehicle in the charging shift will be charged to full, the average vehicle charging time (service time) is $t_s = (1-\delta) C/s$. If we further assume that the vehicle charging times (service times) are exponentially distributed, then the vehicle charging process at each charging station can be modeled as an M/M/Q/V queue with arrival rate $\gamma/K$ and service rate $1/t_s$. Define the system state as the number of vehicles in the charging station and denote $\rho=\frac{\gamma t_s}{K}$ and $a=\frac{\gamma t_s}{K Q}$. The steady-state distribution of the M/M/Q/V queueing system is given by \cite{allen1990probability}:
\begin{equation} \label{steady_state_probability_mmck}
P_v = \begin{cases}
        \frac{\rho^v}{v!} P_0, \quad & \text{if } 0\leq v \leq Q \\
        \frac{\rho^v}{Q^{v-Q} Q!} P_0, \quad & \text{if } Q \leq v \leq V
    \end{cases}
\end{equation}
where $P_0=\left( \sum_{v=0}^{Q-1} \frac{\rho^v}{v!} + \sum_{v=Q}^V \frac{\rho^v}{Q^{v-Q} Q!} \right)^{-1}$ is the probability of the system having no vehicles, and $P_v$ is the probability of the system holding $v$ vehicles (including those in service). Given the steady-state distribution (\ref{steady_state_probability_mmck}), the blocking probability $P_V$ and the average vehicle waiting time $t_w$ can be characterized as \cite{allen1990probability}:
\begin{equation} \label{blocking_prob_mmck}
    P_V = \frac{\rho^V}{Q^{V-Q} Q!} \left( \sum_{v=0}^{Q-1} \frac{\rho^v}{v!} + \sum_{v=Q}^V \frac{\rho^v}{Q^{v-Q} Q!} \right)^{-1} ,
\end{equation}
and
\begin{equation} \label{tw_mmck}
    t_w = \frac{ a \rho^Q \left[1-a^{V-Q+1}-(1-a)(V-Q+1)a^{V-Q} \right] K }{ \gamma (1-a)^2 (1-P_V) Q!} \left( \sum_{v=0}^{Q-1} \frac{\rho^v}{v!} + \sum_{v=Q}^V \frac{\rho^v}{Q^{v-Q} Q!} \right)^{-1}
\end{equation}
Equation (\ref{blocking_prob_mmck}) and (\ref{tw_mmck}) capture the blocking probability $P_V$ and the average vehicle waiting time $t_w$ as a function of the charging demand $\gamma$, the number of charging stations $K$, the number of chargers at each station $Q$, the capacity of the station $V$, the battery capacity $C$, and the charging speed $s$, respectively. This corresponds to the practice that the vehicle waiting time for charging depends on both the demand, the infrastructure supply, and the battery/charging technologies.

\footnotetext{Retrieved 2023, from https://www.tesla.com/supercharger. Copyright 2023 by Tesla. Reprinted with permission.}

{\em Queueing Model for Battery Swapping}: Compared to plug-in charging, battery swapping has a more complicated operating mechanism since it involves not only battery swapping for electric vehicles but also charging for batteries inside the swapping station. These two processes occur simultaneously and interact with each other. On the one hand, arrived vehicles can only start service when there are fully-charged batteries and unoccupied swapping servers/robots. Otherwise, they have to wait for services, which produces an open queue of electric vehicles. On the other hand, each depleted battery will swap another fully-charged battery out, which implies that there is a closed queue/circulation of batteries inside the station. The strong coupling of the open queue and the closed queue makes it difficult to apply conventional queueing models (e.g., M/M/N queue) to characterize such a system. Inspired by \cite{tan2014queueing}, we use a mixed queueing network to model the battery swapping station\footnote{The mixed queueing model is similar to the model proposed in \cite{tan2014queueing}. However, different from \cite{tan2014queueing}, we primarily focus on characterizing the average vehicle waiting time $t_w$ and integrating it into the economic analysis of the AMoD systems.}. Below we detail the mixed queueing network model.

\begin{figure}[htbp]
  \centering
  \begin{subfigure}[b]{0.48\textwidth}
    \centering
    \includegraphics[width=\textwidth]{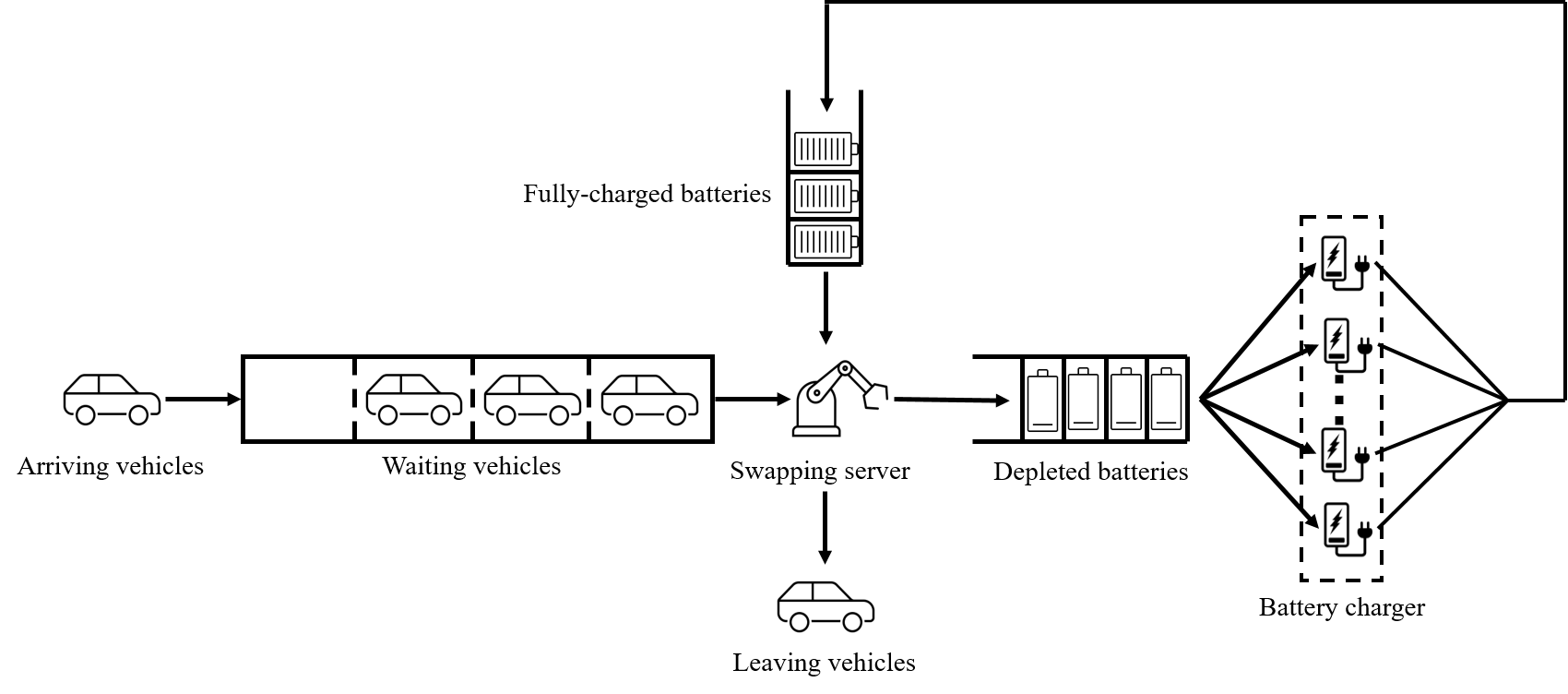}
    \caption{The queueing network model for the battery swapping station.}
    \label{fig:swapping_queueing}
  \end{subfigure}
  \hfill
  \begin{subfigure}[b]{0.48\textwidth}
    \centering
    \includegraphics[width=\textwidth]{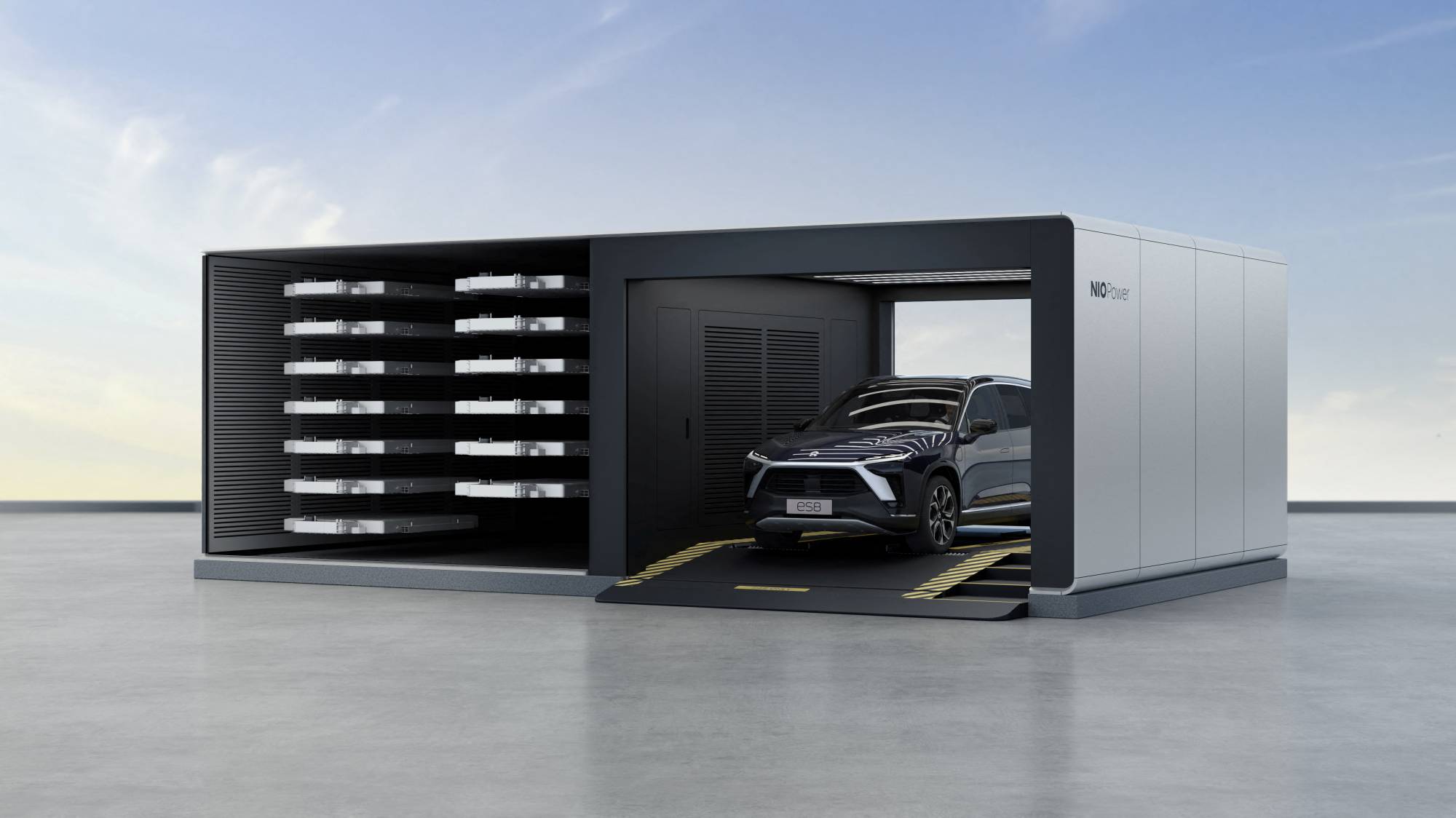}
    \caption[Caption]{The schematic of Nio's battery swap station\protect\footnotemark.}
    \label{fig:nio}
  \end{subfigure}
    \caption{The mixed queueing model for the battery swapping station and Nio's battery swap station.}
  \label{fig:battery_swapping_station_illustration}
\end{figure}

Consider a fleet of electric TNC vehicles with a battery capacity of $C$. There are $K$ battery swapping stations uniformly distributed over the city. The battery swapping station is modeled as a mixed queueing network which consists of an open queue of EVs and a closed queue of EV batteries, as shown in Figure \ref{fig:swapping_queueing}. {Let us take Nio's battery swap station as the prototype (see Figure \ref{fig:nio} the schematic of Nio's battery swap station).} Each swapping station has a single server (i.e., robot/machine) that can execute battery swapping for the customer and $Q$ chargers that can replenish energy for depleted batteries. Batteries swapped out of the vehicles are stored and immediately recharged by the battery chargers. For simplicity, we assume that the number of stored batteries per station is identical to the number of chargers at each station, i.e., $Q$. Due to limited parking spaces, it is assumed that a maximum number of $V$ vehicles can wait for battery swapping services. Therefore, the open queue of the EVs in the mixed queueing network has a capacity/buffer size of $V$. We further assume that the two buffers in the closed queue of EV batteries are large enough to hold all the batteries. Therefore, there will be no deadlock and blocking in the closed queue \cite{onvural1990survey}. The battery swapping server/robot has a constant service time $t_s$ (e.g., 2 minutes) for each battery swapping service, which accounts for the time needed for exchanging batteries. We define the system state of the battery swapping station\footnote{The state is defined for a single battery swapping station and we focus on its stead-state equilibrium. By assuming that the demand is uniformly distributed to distinct stations, we do not need to distinguish one battery swapping station to another.} at any time $t$ as $S_t = (v_t,q_t)$, where $v_t \in[0,V]$ denotes the number of TNC vehicles present in the station, and $q_t \in[0,Q]$ represents the number of fully-charged batteries in the closed queue. Let $\mathcal{S}=\{(v,q)|v\in[0,V],q\in [0,Q]\}$ be the state space. Now we derive the steady-state distribution of the queueing network to characterize the average vehicle waiting time at battery swapping stations.

\footnotetext{Retrieved 2022, from https://www.nio.com/nio-power. Copyright 2022 by NIO. Reprinted with permission.}

We first characterize vehicle arrivals at the queueing network. For any given time $t$, we consider the time interval $(t,t+t_s]$. Let $g(v)$ be the probability of $v$ new vehicle arrivals at the battery swapping station during $(t,t+t_s]$. Assuming that the arrival process of autonomous electric TNC vehicles at each battery swapping station is Poisson with rate $\gamma/K$\footnote{The total charging demand of TNC vehicles is $\gamma$, so the arrival rate of vehicles at each battery swapping station is $\gamma/K$.}, then the probability $g(v)$ can be written as:
\begin{equation} \label{arrival_probability}
    g(v) = \frac{e^{-\gamma t_s / K}(\gamma t_s / K)^v}{v!} \quad v\geq 0,
\end{equation}
which is derived based on the probability density function of Poisson distribution with mean arrival $\gamma t_s / K$.

Next, we consider the dynamic of the closed queue of batteries. Let $f(m,n)$ be the probability of $n$ batteries finishing charging by time $t+t_s$ given that there are $m$ batteries being charged at time $t$. Analogous to vehicle charging service, we assume that battery charging times of different chargers follow independent and identical exponential distributions. Once the battery level approaches the minimum SoC, TNC vehicles go to nearby battery swapping stations to swap depleted batteries with fully-charged batteries. Similarly, let $\delta$ be the average battery SoC of vehicles that arrive at battery swapping stations. The corresponding average battery charging time is $t_c = (1-\delta)C/s$. In this case, $f(m,n)$ can be considered as the probability of getting $n$ batteries fully charged within $(t,t+t_s]$ in $\min\{Q,m\}$ independent Bernoulli trials with success probability $1-e^{-\frac{t_s}{t_c}}$. Therefore, $f(m,n)$ is given by the probability mass function of the binomial distribution:
\begin{equation} \label{binomial_probabiltiy}
    f(m,n) = \begin{pmatrix}
\min\{Q,m\} \\
n
\end{pmatrix} \left(1-e^{-\frac{t_s}{t_c}}\right)^n e^{-\frac{t_s}{t_c}\left(\min\{Q,m\}-n\right)} ,
\end{equation}
where $n\in \left[0,\min\{Q,m\}\right]$, indicating that the number of batteries finishing charging by $(t,t+t_s]$ is bounded above by the number of batteries being charged at time $t$ and the total number of chargers.

Given (\ref{arrival_probability}) and (\ref{binomial_probabiltiy}), we can further establish the dynamics of the queueing system. Without loss of generality, let $S_t = (v_t,q_t)$ be the state of the system at time $t$ holding $v_t$ electric vehicles and $q_t$ fully-charged batteries, and let $S_{t+t_s}=(v_{t+t_s},q_{t+t_s})$ be the system state at time $t+t_s$ with $v_{t+t_s}$ electric vehicles and $q_{t+t_s}$ fully-charged batteries. Denote $P_t(v_t,q_t)=\mathrm{Pr}\left(S_t=(v_t,q_t)\right)$ and $P_{t+t_s}(v_{t+t_s},q_{t+t_s})=\mathrm{Pr}\left(S_{t+t_s}=(v_{t+t_s},q_{t+t_s})\right)$ as the corresponding probability of the system state $S_t=(v_t,q_t)$ and $S_{t+t_s}=(v_{t+t_s},q_{t+t_s})$, respectively. Based on the Bayes' law, the state probability at time $t+t_s$ can be written as:
\begin{equation} \label{Bayes_law}
    P_{t+t_s}\left(S_{t+t_s}\right) = \sum_{S_t \in \mathcal{S}} P_t\left(S_t\right) \cdot \mathrm{Pr}\left(S_{t+t_s}|S_t\right) ,
\end{equation}
where $\mathcal{S}$ is the state space, and the conditional probability $\mathrm{Pr}\left(S_{t+t_s}|S_t\right)$ represents the transition probability from state $S_t$ to $S_{t+t_s}$. To characterize this transition probability, we note that the service time is $t_s$, and thus the vehicle in service at time $t$ will have left the system by time $t+t_s$, and the swapped depleted battery will have started charging at time $t+t_s$. In this regard, all vehicles present at time $t+t_s$ either arrived within $(t,t+t_s]$ or were already waiting for service at time $t$. In addition, all fully-charged batteries present at time $t+t_s$ were either charged during $(t,t+t_s]$ or already fully charged before $t$. Consequently, the transition probability from state $S_t = (v_t,q_t)$ to state $S_{t+t_s}=(v_{t+t_s},q_{t+t_s})$ is given by:
\begin{equation} \label{transition_probability}
\begin{split}
   &\mathrm{Pr}\left(S_{t+t_s}=(v_{t+t_s},q_{t+t_s})|S_t = (v_t,q_t)\right) \\
&=\begin{cases}
    g\left(v_{t+t_s}-v_t+\min\{v_t,q_t,1\}\right) \cdot f\left(Q-q_t,q_{t+t_s}-q_t+\min\{v_t,q_t,1\}\right) \quad v_{t+t_s}\in[0,V-1], q_{t+t_s}\in [0,Q] \\
    \left(1-\sum\limits_{v=0}^{V-1} g\left(v-v_t+\min\{v_t,q_t,1\}\right)\right) \cdot f\left(Q-q_t,q_{t+t_s}-q_t+\min\{v_t,q_t,1\}\right) \quad v_{t+t_s}=V, q_{t+t_s}\in [0,Q] .
\end{cases} 
\end{split}
\end{equation}
The intuition behind (\ref{transition_probability}) is straightforward. The term $\min\{v_t,q_t,1\}$ represents the count of vehicles that finished the battery swapping service within the time interval $(t,t+t_s]$. This count is limited by the number of vehicles $v_t$, the number of fully-charged batteries $q_t$ at time $t$, and the number of swapping servers. Given that each battery swapping station has only one swapping server, either 0 or 1 vehicle will finish the battery swapping service and leave the station during $(t,t+t_s]$, depending on the presence of vehicles and the availability of batteries. Since there are $v_t$ and $v_{t+t_s}$ vehicles in the station at time $t$ and $t+t_s$ respectively, and $\min\{v_t,q_t,1\}$ vehicles have left the station during $(t+t_s]$, when $0\leq v_{t+t_s}\leq V-1$, the number of vehicle arrivals within $(t,t+t_s]$ is $v_{t+t_s}-v_t+\min\{v_t,q_t,1\}$. The corresponding probability is $g\left(v_{t+t_s}-v_t+\min\{v_t,q_t,1\}\right)$ according to (\ref{arrival_probability}). Due to the limited capacity $V$, all cases with new vehicle arrivals not less than $V-1-v_t+\min\{v_t,q_t,1\}$ make state $S_t=(v_t,q_t)$ transit to state $S_{t+t_s}=(V,q_{t+t_s})$. Therefore, the probability of new vehicle arrivals during $(t,t+t_s]$ when $v_{t+t_s}=V$ is given by $1-\sum_{v=0}^{V-1} g\left(v-v_t+\min\{v_t,q_t,1\}\right)$. Furthermore, since there are $q_t$ and $q_{t+t_s}$ fully-charged batteries at time $t$ and $t+t_s$ respectively, and $\min\{v_t,q_t,1\}$ fully-charged batteries have been swapped out during $(t,t_s]$, the number of batteries being charged at time $t$ is $Q-q_t$, and the number of batteries finished charging during $(t,t+t_s]$ is $q_{t+t_s}-q_t+\min\{v_t,q_t,1\}$. According to (\ref{binomial_probabiltiy}), the probability of battery charging is $f\left(Q-q_t,q_{t+t_s}-q_t+\min\{v_t,q_t,1\}\right)$. Finally, the product of vehicle arrival probability and battery charging probability gives the state transition probability (\ref{transition_probability}). 

Next, we derive the steady-state distribution of the mixed queueing network to represent the long-run equilibrium status of the battery swapping station. Denote $P(v,q)=\lim_{t\to \infty} P_t(v_t,q_t)$ as the steady-state probability of the station holding $v$ vehicles and $q$ fully-charged batteries. Plugging (\ref{transition_probability}) into (\ref{Bayes_law}), the steady-state equations can be written as:
\begin{equation} \label{linear_equation_system}
P(v,q) =
\begin{cases}
     \sum\limits_{\hat{v}=0}^V\sum\limits_{\hat{q}=0}^Q P(\hat{v},\hat{q}) \cdot g(v-\hat{v}+\min\{\hat{v},\hat{q},1\}) \cdot f(Q-\hat{q},\hat{q}-q+\min\{\hat{v},\hat{q},1\}) \quad v \in [0,V-1] \\
    \sum\limits_{\hat{v}=0}^V\sum\limits_{\hat{q}=0}^Q P(\hat{v},\hat{q}) \cdot \left(1-\sum\limits_{u=0}^{V-1}g\left(u-\hat{v}+\min\{\hat{v},\hat{q},1\}\right)\right) \cdot f(Q-\hat{q},\hat{q}-q+\min\{\hat{v},\hat{q},1\}) \quad v=V
\end{cases} .
\end{equation}
Equation(s) (\ref{linear_equation_system}) form a finite system of $(V+1)(Q+1)$ linear equations. The solution to (\ref{linear_equation_system}) together with the normalization equation $\sum_{v=0}^V\sum_{q=0}^Q P(v,q)=1$ gives the steady-state distribution of the battery swapping system. 

Finally, based on the steady-state distribution, we derive the blocking probability $P_V$ and the average vehicle waiting time $t_w$ at battery swapping stations. Given $P(v,q)$, the blocking probability $P_V$ of the queueing system, i.e., the probability of the swapping station being full, can be determined as:
\begin{equation} \label{blocking_probability_battery_swapping}
    P_V=\sum_{q=0}^Q P(V,q) .
\end{equation}
The average number of vehicles at each battery swapping station is given by:
\begin{equation} \label{number_of_vehicles_probability}
    \bar{N} = \sum_{v=0}^V \sum_{q=0}^Q v P(v,q) .
\end{equation}
Further, based on Little's law, the average number of vehicles at each station $\bar{N}$ can also be written as:
\begin{equation} \label{number_of_vehicles_littles}
    \bar{N} = \frac{\gamma}{K} \left(1-P_V\right) (t_w+t_s),
\end{equation}
where $\gamma/K$ represents the average arrival rate of vehicles at each station, $t_w$ is the average vehicle waiting time at battery swapping stations, and $t_s$ is the constant service time. Combining equations (\ref{number_of_vehicles_probability}) and (\ref{number_of_vehicles_littles}) derives the average waiting time $t_w$ as:
\begin{equation}
    t_w = \frac{K \sum_{v=0}^V\sum_{q=0}^Q v P(v,q)}{\gamma (1 - P_V)} - t_s .
\end{equation}

\begin{remark}
We assume that the charging time follows the exponential distribution with mean charging time $(1-\delta)C/s$ when deriving the queueing model for both plug-in charging and battery swapping. While the exponential distribution may not fully characterize the exact distribution of vehicles' charging times (should be a truncated distribution, where the truncation occurs at the minimum SoC incurred by the platform), we argue that the independent and identical exponential distributions of charging times capture (i) the average charging time that stems from the platform's charging scheduling and (ii) the independence between different chargers and batteries, which are fundamental aspects of the queueing system. This assumption offers mathematical tractability, and is widely used in EV charging-related literature \cite{tan2014queueing,lai2021demand,zhang2020efficient,zhao2023ev} to serve as reasonable approximations and facilitate the derivation of the analytical model.
\end{remark}

\subsection{The balance of energy} \label{energy_balance}
At the stationary state, the electricity consumed by the autonomous electric vehicle fleet and the electricity charged to the fleet should be balanced regardless of the charging strategies adopted by the TNC platform. Let $l$ be the average hourly electricity consuming rate of autonomous electric TNC vehicles, and then we have the following energy balance equation:
\begin{equation} \label{energy_balance_equation}
    (N-N_2^w-N_2^s) l = \gamma (1-\delta) C,
\end{equation}
where $\gamma$ is the charging demand of TNC vehicles, $N_2^w$ is the number of TNC vehicles waiting at charging stations (or battery swapping stations), and $N_2^s$ is the number of vehicles being charged at the charging station (or serviced at the battery swapping station). The term $N-N_2^w-N_2^s$ accounts for the average number of TNC vehicles traveling on roads, and the left-hand side of (\ref{energy_balance_equation}) represents the hourly electricity consumption by those vehicles cruising on the street. In per-unit time, $\gamma$ vehicles arrive at charging/battery stations with the average battery level $\delta C$ and leave with full battery capacity $C$, and thus the electricity input to the system per hour is $\gamma (1-\delta) C$. The consumed and charged electricity should be balanced at the stationary equilibrium.

\section{Charging Infrastructure Planning: Bilevel Optimization Model}

In this section, we consider the profit maximization decision of the TNC platform and its interaction with the government, which deploys electric vehicle charging infrastructures to promote social welfare. Overall, the problem can be formulated with a bilevel optimization framework: on the upper level, the government designs the infrastructure deployment plan for social welfare maximization; at the lower level, the TNC determines the operational strategy for profit maximization. This section mathematically formulates the bilevel optimization model. The platform's strategic decisions for profit maximization and the government's infrastructure deployment for social welfare maximization will be introduced separately.

\subsection{Profit maximization for the TNC} \label{profit_maximization}
The revenue of the TNC platform is the ride fares collected from passengers, and the cost consists of two parts: (1) electricity cost; and (2) investment and operating cost of the fleet. Let $p_e$ be the average electricity price and let $C_{av}$ be the average hourly operating cost of a TNC vehicle. In each hour, the platform collects $\lambda p_f$ ride fares from passengers, spends $\gamma (1-\delta) C p_e$ electricity costs for charging vehicles and $N C_{av}$ operating costs for the autonomous vehicle fleet. The difference between the revenue and the cost is the platform profit. Given the upper-level government's deployment of charging infrastructures (either charging stations or battery swapping stations), the TNC determines the ride fare and the fleet size to maximize its profit subject to the market equilibrium constraints. The TNC's profit-maximization problem under the strategy of plug-in charging and battery swapping can be formulated as follows, respectively.
\begin{itemize}[leftmargin=0.45cm,topsep=0pt]
    \item TNC's profit maximization under plug-in charging
\end{itemize}
\begin{equation}
\label{optimalpricing_charging}
    \mathop {\max }\limits_{{p_f},{N}} \quad \lambda p_f - \gamma (1-\delta) C p_e - N C_{av}
\end{equation}
\begin{subnumcases}{\label{constraints_charging}}
    \lambda = \lambda_0 F_p(p_f+\alpha w^c) \label{passenger_demand_charging} \\
    w^c = \frac{A}{\sqrt{\lambda w^v}} \label{passenger_waiting_time_charging} \\
    \lambda (w^v+w^c+\tau) = N_1 \label{conservation_N1_charging} \\
    N_1 + N_2 = N \label{conservation_N_charging} \\
    N_2^m + N_2^w + N_2^s = N_2 \label{conservation_N2_charging} \\
    N_2^m = \gamma t_m \label{N2m_charging} \\
    N_2^w = \gamma t_w \label{N2w_charging} \\
    N_2^s = \gamma t_s \label{N2s_charging} \\
    t_m = \frac{B}{\sqrt{K (1-P_V)}} \label{vehicle_searching_time_charging} \\ 
    P_V = \frac{\rho^V}{Q^{V-Q} Q!} \left( \sum_{v=0}^{Q-1} \frac{\rho^v}{v!} + \sum_{v=Q}^V \frac{\rho^v}{Q^{v-Q} Q!} \right)^{-1} \label{blocking_prob_charging} \\
    t_w = \frac{ a \rho^Q \left[1-a^{V-Q+1}-(1-a)(V-Q+1)a^{V-Q} \right] K }{ \gamma (1-a)^2 (1-P_V) Q!} \left( \sum_{v=0}^{Q-1} \frac{\rho^v}{v!} + \sum_{v=Q}^V \frac{\rho^v}{Q^{v-Q} Q!} \right)^{-1} \label{vehicle_waiting_time_charging} \\
    (N-N_2^w-N_2^s)l=\gamma (1-\delta) C \label{energy_balance_charging}
\end{subnumcases}
where $t_s=(1-\delta) C / s$ accounts for the average vehicle charge/service time at charging stations, $\rho = \gamma t_s/K$ and $a=\gamma t_s / KQ$ are auxiliary variables, $K$ represents the number of charging stations, $Q$ denotes the number of chargers per station, $V$ is the capacity of the charging station.

\begin{itemize}[leftmargin=0.45cm,topsep=0pt]
    \item TNC's profit maximization under battery swapping
\end{itemize}
\begin{equation}
\label{optimalpricing_swapping}
    \mathop {\max }\limits_{{p_f},{N}} \quad \lambda p_f - \gamma (1-\delta) C p_e - N C_{av}
\end{equation}
\begin{subnumcases}{\label{constraints_swapping}}
    \lambda = \lambda_0 F_p(p_f+\alpha w^c) \label{passenger_demand_swapping} \\
    w^c = \frac{A}{\sqrt{\lambda w^v}} \label{passenger_waiting_time_swapping} \\
    \lambda (w^v+w^c+\tau) = N_1 \label{conservation_N1_swapping} \\
    N_1 + N_2 = N \label{conservation_N_swapping} \\
    N_2^m + N_2^w + N_2^s = N_2 \label{conservation_N2_swapping} \\
    N_2^m = \gamma t_m \label{N2m_swapping} \\
    N_2^w = \gamma t_w \label{N2w_swapping} \\
    N_2^s = \gamma t_s \label{N2s_swapping} \\
    t_m = \frac{B}{\sqrt{K(1-P_V)}} \label{vehicle_searching_time_swapping} \\
    P_V =\sum_{q=0}^Q P(V,q) \label{blocking_prob_swapping} \\
    t_w = \frac{K \sum_{v=0}^V\sum_{q=0}^Q v P(v,q)}{\gamma (1 - P_V)} - t_s \label{vehicle_waiting_time_swapping} \\ 
    (N-N_2^w-N_2^s)l=\gamma (1-\delta) C \label{energy_balance_swapping}
\end{subnumcases}
where $t_s$ denotes the average vehicle charge/service time at battery swapping stations, $K$ is the number of battery swapping stations, $Q$ is the number of chargers/batteries per station, $V$ is the capacity of the battery swapping station, $P(v,q)$ is the steady-state distribution of the battery swapping station which can be derived from (\ref{linear_equation_system}).

In the optimization (\ref{optimalpricing_charging}) and (\ref{optimalpricing_swapping}), $p_f$ and $N$ are decision variables corresponding to the platform's pricing and fleet sizing strategy; $\lambda$, $w^c$, $w^v$, $N_1$, $N_2$, $N_2^m$, $N_2^w$, $N_2^s$, $\gamma$, $t_m$, $t_w$ and $P_V$ are involved endogenous variables; $\lambda_0$, $\alpha$, $\tau$, $\delta$, $s$, $t_s$, $l$, $p_e$, $A$, $B$, $C$, $C_{av}$, $K$, $Q$ and $V$ are model parameters. Note that there are a few differences between the market equilibrium models for plug-in charging (\ref{constraints_charging}) and battery swapping (\ref{constraints_swapping}). For plug-in charging, the blocking probability $P_V$ (\ref{blocking_prob_charging}) and the vehicle waiting time (\ref{vehicle_waiting_time_charging}) are derived based on the M/M/Q/V queue, whereas for battery swapping, the blocking probability $P_V$ (\ref{blocking_prob_swapping}) and the vehicle waiting time (\ref{vehicle_waiting_time_swapping}) are derived from the mixed queueing model. In addition, the average vehicle charging/service time at charging stations $t_s$ is associated with the charging speed and the battery capacity through $t_s = (1-\delta)C/s$. However, the vehicle service time at battery swapping stations is considered as a constant, which is relatively short. These differences result from the fundamental operating mechanisms of charging stations and battery swapping stations.

\subsection{Social welfare maximization for the government} \label{social_welfare_maximization}

On the upper level, the government deploys charging infrastructures in the city for social welfare maximization. In the ride-sourcing market with autonomous electric vehicles, the social welfare comprises the passengers' surplus, the profit of the TNC, and the cost of constructing and operating charging infrastructures. Let $\phi$ be the per charger investment and operating cost of a charging/battery swapping station in each hour\footnote{Note that the infrastructure cost is pro-rated. In addition, the cost parameter $\phi$ for charging and battery swapping stations may be distinct. Their difference will be reflected in the case study.}. We define social welfare as the difference between the social benefit and the infrastructure cost:
\begin{equation} \label{social_welfare_func}
    \Pi_{SW} = \lambda_0 \int_c^{\infty} F_p(x) dx + \lambda p_f - \gamma (1-\delta) C p_e - N C_{av} - \phi K Q,
\end{equation}
where $c$ denotes the average travel cost of passengers, the term $\lambda_0 \int_c^{\infty} F_p(x) dx$ calculates the passengers' surplus, the term $\lambda p_f - \gamma (1-\delta) C p_e - N C_{av}$ represents the TNC profit, and the last term $\phi K Q$ accounts for the total cost of charging facilities. The government designs the infrastructure deployment scheme to maximize social welfare. In particular, when planning charging stations for vehicle charging, the government determines the number of charging stations $K$ and the number of chargers $Q$ at each station to promote social welfare. Denote $\phi_{vc}$ as the per charger infrastructure cost of charging stations. {The social welfare maximization for the government when planning charging stations can be formulated as the following bi-level optimization:}
\begin{equation}
\label{optimalplanning_charging}
    \mathop {\max }\limits_{{K>0},{Q>0}} \quad \lambda_0 \int_c^{\infty} F_p(x) dx + \lambda p_f - \gamma (1-\delta) C p_e - N C_{av} - \phi_{vc} K Q ,
\end{equation}
\begin{equation*}
    \text{subject to (\ref{optimalpricing_charging}), (\ref{passenger_demand_charging})-(\ref{energy_balance_charging}).}
\end{equation*}
On the other hand, due to the closed structure and technological feasibility, battery swapping stations are usually standardized, and the number of chargers/batteries $Q$ inside the swapping station is pre-designed by the manufacturer for ease of mass production ({see the schematic of Nio's battery swap station in Figure \ref{fig:nio}}). Therefore, when deploying battery swapping stations, the government can only control the total number of battery swapping stations $K$ to maximize social welfare. Let $\phi_{bs}$ be the per charger infrastructure cost of battery swapping stations. The social welfare maximization for the government when deploying battery swapping stations can be cast as:
\begin{equation}
\label{optimalplanning_swapping}
    \mathop {\max }\limits_{{K>0}} \quad \lambda_0 \int_c^{\infty} F_p(x) dx + \lambda p_f - \gamma (1-\delta) C p_e - N C_{av} - \phi_{bs} K Q ,
\end{equation}
\begin{equation*}
    \text{subject to (\ref{optimalpricing_swapping}), (\ref{passenger_demand_swapping})-(\ref{energy_balance_swapping}).}
\end{equation*}

\begin{remark}
In practice, different charging stations may have a distinct number of chargers considering the need for parking spaces and the spatial demand at different locations. This paper neglects these aspects and considers $K$ uniformly distributed homogeneous charging stations with the same number of chargers at the aggregate level. We emphasize that such an aggregate model is sufficient for our purpose since we focus on revealing the interactions among different stakeholders (e.g., passengers, the TNC, and the government) in the electric AMoD market and the long-term charging infrastructure planning at the city level. The spatial imbalance of charging demand can be considered by extending our current model into a network equilibrium model, which is left for future work.
\end{remark}

\subsection{Solution properties}

The optimization (\ref{optimalpricing_charging}) and (\ref{optimalplanning_charging}), and the optimization (\ref{optimalpricing_swapping}) and (\ref{optimalplanning_swapping}) constitute the bilevel optimization framework for the government's infrastructure planning of charging stations and battery swapping stations, respectively. To justify the existence of the market equilibrium and solve the optimization, we first introduce the following three lemmas:

\begin{lemma} \label{lemma_wc_gamma}
Given $A$ and $\tau$, in equilibrium (\ref{constraints_charging}) and (\ref{constraints_swapping}), if $N_1 \geq \left(\sqrt[3]{2}+\sqrt[3]{\frac{1}{4}}\right)\left(A\lambda \right)^{\frac{2}{3}}+\lambda\tau$, there exist strictly positive $w^c$ and $w^v$ satisfying (\ref{passenger_waiting_time_charging})-(\ref{conservation_N1_charging}). Furthermore, $w^c$ can be uniquely determined as a function of $\lambda$ and $\gamma$ from (\ref{passenger_waiting_time_charging})-(\ref{conservation_N1_charging}), i.e., $w^c=w^c(\lambda,N_1)$.
\end{lemma}
The proof of Lemma \ref{lemma_wc_gamma} can be found in Appendix B in our previous work \cite{gao2022shared}. The condition $N_1 \geq \left(\sqrt[3]{2}+\sqrt[3]{\frac{1}{4}}\right)\left(A\lambda \right)^{\frac{2}{3}}+\lambda\tau$ guarantees the existence of strictly positive $w^c$ and $w^v$ in equilibrium (\ref{constraints_charging}) and (\ref{constraints_swapping}). Besides, $w^c=w^c(\lambda,N_1)$ can be uniformly derived under different charging strategies. 

\begin{lemma} \label{lemma_N_gamma}
Given $\delta$, $s$, $t_s$, $l$, $B$, $C$, $K$, $Q$ and $V$, (i) in equilibrium (\ref{constraints_charging}), $P_V$, $N_1$ and $N$ can be uniquely determined as a function of $\gamma$ from (\ref{conservation_N_charging})-(\ref{energy_balance_charging}), denoted as $\hat{P}_V(\gamma)$, $\hat{N}_1(\gamma)$ and $\hat{N}(\gamma)$, where $\hat{N}_1(\gamma)=\frac{\gamma (1-\delta) C}{l}-\frac{\gamma B}{\sqrt{K(1-\hat{P}_V(\gamma))}}$; (ii) in equilibrium (\ref{constraints_swapping}), $P_V$, $N_1$ and $N$ can be uniquely determined as a function of $\gamma$ from (\ref{conservation_N_swapping})-(\ref{energy_balance_swapping}), denoted as $\tilde{P}_V(\gamma)$, $\tilde{N}_1(\gamma)$ and $\tilde{N}(\gamma)$, where $\tilde{N}_1(\gamma)=\frac{\gamma (1-\delta) C}{l}-\frac{\gamma B}{\sqrt{K(1-\tilde{P}_V(\gamma))}}$.
\end{lemma}
The proof of Lemma \ref{lemma_N_gamma} can be found in Appendix A. It states that in equilibrium, given the exogenous model parameters and the upper-level infrastructure deployment plan, the blocking probability $P_V$, the number of operating vehicles $N_1$, and the total fleet size $N$ unilaterally depend on the charging demand $\gamma$. Furthermore, we have the following results on the property of $\hat{N}_1(\gamma)$ and $\tilde{N}_1(\gamma)$.
\begin{lemma} \label{lemma_property_N1}
    Given $\delta$, $s$, $t_s$, $l$, $B$, $C$, $K$, $Q$ and $V$, if $K > \left(\frac{Bl}{(1-\delta)C}\right)^2$, (i) there exist a unique $\hat{\gamma}_0>0$ and $\hat{\gamma}_{*} \in (0,\hat{\gamma}_0)$ such that $\hat{N}_1(\hat{\gamma}_0)=0$, and $\forall \gamma \in (0,\hat{\gamma}_0)$, $0<\hat{N}_1(\gamma)\leq \hat{N}_1(\hat{\gamma}_{*})$; (ii) there exist a unique $\tilde{\gamma}_0>0$ and $\tilde{\gamma}_{*} \in (0,\tilde{\gamma}_0)$ such that $\tilde{N}_1(\tilde{\gamma}_0)=0$, and $\forall \gamma \in (0,\tilde{\gamma}_0)$, $0<\tilde{N}_1(\gamma)\leq \tilde{N}_1(\tilde{\gamma}_{*})$.
\end{lemma}
The proof of Lemma \ref{lemma_property_N1} can be found in Appendix B. It states that under two distinct charging strategies, when the number of charging/swapping stations is higher than $\left(\frac{Bl}{(1-\delta)C}\right)^2$, there exists a positive charging demand $\gamma_0$ at which the number of operating vehicles becomes zero. Furthermore, for any charging demand below $\gamma_0$, there will be a positive number of operating vehicles $N_1$, and the maximum number of operating vehicles is achieved at a $\gamma_{*} \in (0,\gamma_0)$. Without loss of generality, let ${\hat{N}_1}^\text{max}=\hat{N}_1(\hat{\gamma}_{*})$ and ${\tilde{N}_1}^\text{max}=\tilde{N}_1(\tilde{\gamma}_{*})$ be the maximum number of operating vehicles under plug-in charging and battery swapping, respectively.
Based on the results of Lemma \ref{lemma_wc_gamma}-\ref{lemma_property_N1}, we can justify the existence of the market equilibrium (\ref{constraints_charging}) and (\ref{constraints_swapping}):
\begin{proposition} \label{prop_equilibrium_charging_swapping}
Given $\lambda_0$, $\alpha$, $\tau$, $\delta$, $s$, $t_s$, $l$, $A$, $B$, $C$, $K$, $Q$ and $V$, (i) if $K > \left(\frac{Bl}{(1-\delta)C}\right)^{2}$ and $F_p\left(\frac{\alpha A}{\sqrt{ {\hat{N}_1}^\mathrm{max} }}\right)>0$, there exist strictly positive $p_f$, $w^c$, $w^v$, $\lambda$, $\gamma$, $t_m$, $t_w$, $P_V$, $N$, $N_1$, $N_2$, $N_2^m$, $N_2^w$ and $N_2^s$ that constitute a market equilibrium satisfying (\ref{constraints_charging}); (ii) if $K > \left(\frac{Bl}{(1-\delta)C}\right)^{2}$ and $F_p\left(\frac{\alpha A}{\sqrt{{\tilde{N}_1}^\mathrm{max}}}\right)>0$, there exist strictly positive $p_f$, $w^c$, $w^v$, $\lambda$, $\gamma$, $t_m$, $t_w$, $P_V$, $N$, $N_1$, $N_2$, $N_2^m$, $N_2^w$ and $N_2^s$ that constitute a market equilibrium satisfying (\ref{constraints_swapping}).
\end{proposition}


The proof of Proposition \ref{prop_equilibrium_charging_swapping} is deferred to Appendix C. Both (i) and (ii) in Proposition (\ref{prop_equilibrium_charging_swapping}) state the existence of the market equilibrium under distinct conditions. To guarantee the existence of the market equilibrium (\ref{constraints_charging}) and (\ref{constraints_swapping}), Proposition \ref{prop_equilibrium_charging_swapping} states that the supply of charging/battery swapping stations should exceed $\left(\frac{Bl}{(1-\delta)C}\right)^{2}$. Besides, the conditions $F_p\left(\frac{\alpha A}{\sqrt{ {\hat{N}_1}^\mathrm{max} }}\right)>0$ and $F_p\left(\frac{\alpha A}{\sqrt{{\tilde{N}_1}^\mathrm{max}}}\right)>0$ indicate that when the TNC platform sets the ride fare $p_f$ to 0 and has ${\hat{N}_1}^\text{max}$ and ${\tilde{N}_1}^\text{max}$ vehicles in operating under plug-in charging and battery swapping respectively, there will be a positive number of passengers. Overall, these are mild conditions that can be satisfied in reality.

The market equilibrium (\ref{constraints_charging}) and (\ref{constraints_swapping}) involve a bunch of nonlinear constraints, e.g., (\ref{vehicle_waiting_time_charging}) and (\ref{vehicle_waiting_time_swapping}), which makes the optimization problems (\ref{optimalpricing_charging}) and (\ref{optimalpricing_swapping}) non-convex. Besides, due to the complex relations among the endogenous variables and the bi-level structure of the optimization, it is non-trivial to derive the partial derivatives and apply gradient-based algorithms to solve the bi-level optimization. However, since both (\ref{optimalpricing_charging}) and (\ref{optimalpricing_swapping}) are small-scale problems with only two decision variables and several equality constraints, we can equivalently treat $\lambda$ and $\gamma$ as decision variables, absorb all the equality constraints into the objective function, and transform the profit maximization problem into an unconstrained optimization:

\begin{lemma} \label{lemma_unconstrained_charging}
Given the market equilibrium (\ref{constraints_charging}), (\ref{optimalpricing_charging}) is equivalent to the following unconstrained optimization:
\begin{equation} \label{unconstrained_charging}
    \mathop {\max }\limits_{\lambda,\gamma} \quad \tilde{\Pi}_P(\lambda,\gamma) = \lambda \left[F_p^{-1}\left(\frac{\lambda}{\lambda_0}\right)-\alpha w^c\left(\lambda,\hat{N}_1(\gamma)\right)\right]-\gamma (1-\delta) C p_e - \hat{N}(\gamma) C_{av}
\end{equation}
\end{lemma}

\begin{lemma} \label{lemma_unconstrained_swapping}
Given the market equilibrium (\ref{constraints_swapping}), (\ref{optimalpricing_swapping}) is equivalent to the following unconstrained optimization:
\begin{equation} \label{unconstrained_swapping}
    \mathop {\max }\limits_{\lambda,\gamma} \quad \tilde{\Pi}_P(\lambda,\gamma) = \lambda \left[F_p^{-1}\left(\frac{\lambda}{\lambda_0}\right)-\alpha w^c\left(\lambda,\tilde{N}_1(\gamma)\right)\right]-\gamma (1-\delta) C p_e - \tilde{N}(\gamma) C_{av}
\end{equation}
\end{lemma}

Lemma \ref{lemma_unconstrained_charging} and \ref{lemma_unconstrained_swapping} derive the platform profit as a function of $\lambda$ and $\gamma$ under the strategy of charging vehicles and swapping batteries, respectively. By such transformations, for each charging infrastructure plan $(K,Q)$, we can apply the grid search to maximize (\ref{unconstrained_charging}) and (\ref{unconstrained_swapping}), which provides the globally optimal solutions to the lower-level profit maximization problems (\ref{optimalpricing_charging}) and (\ref{optimalpricing_swapping}). After solving the lower-level optimization, we can substitute the optimal solution into (\ref{social_welfare_func}) to compute the corresponding social welfare produced by the infrastructure deployment scheme $(K,Q)$. Finally, we apply a grid search over all possible $K$ and $Q$ to find the optimal infrastructure deployment plan that maximizes the social welfare\footnote{In the planning of battery swapping stations, the government only optimizes the number of battery swapping stations.}.

\section{Case Studies}

This section investigates how the evolvement of charging technologies affects the charging infrastructure planning for the AMoD system. First, we calibrate the model parameters based on real-world TNC data. Second, we present case studies to investigate the impacts of charging speed and battery capacity. Lastly, we compare the TNC market outcomes under two distinct charging strategies.

\subsection{Model parameters}

Consider a case study for New York City. Assume passengers decide whether to choose AMoD services based on a logit model:
\begin{equation} \label{logit_demand_func}
    \lambda = \lambda_0 \frac{\exp(-\epsilon c)}{\exp(-\epsilon c)+\exp(-\epsilon c_0)} ,
\end{equation}
where $c$ is the average generalized travel cost of ride-hailing trips; $\epsilon>0$ and $c_0$ are model parameters. 
Besides, the blocking probability (\ref{blocking_prob_charging}) and the average waiting time (\ref{vehicle_waiting_time_charging}) involve summation and factorial, which restricts the number of chargers per station $Q$ to be integers. For mathematical tractability, we apply analytical continuation to extend these functions to non-integral values of the number of servers, which has been widely used in the queueing theory related literature \cite{smith2003m,dombacher2008stationary}.
\begin{equation}
\begin{cases}
    P_V = \frac{\rho^V (\Gamma(Q+1))^{-1} (Q^{V-Q})^{-1}}{\exp(\rho) \Gamma(Q,\rho) (\Gamma(Q))^{-1} + \rho^Q (1-a^{K-Q+1})/(\Gamma(Q+1)(1-a))} \\
    t_w = \frac{ a \rho^Q \left[1-a^{V-Q+1}-(1-a)(V-Q+1)a^{V-Q} \right] K }{ \gamma (1-a)^2 (1-P_V) \Gamma(Q+1) [\exp(\rho) \Gamma(Q,\rho) (\Gamma(Q))^{-1} + \rho^Q (1-a^{K-Q+1})/(\Gamma(Q+1)(1-a))]}
\end{cases}
\end{equation}
where $\Gamma(\cdot)$ is the gamma function and $\Gamma(\cdot,\cdot)$ is the incomplete gamma function.

With the logit function (\ref{logit_demand_func}), all model parameters involved in the bi-level charging infrastructure planning under plug-in charging and/or battery swapping are:
\begin{equation}
\label{model_parameters}
    \Theta = \{\lambda_0, \alpha, \tau, \epsilon, \delta, c_0, l, s, p_e, A, B, C_{av}, C, Q, V, \phi_{vc}, \phi_{bs}\},
\end{equation}
among which $\lambda_0$, $\epsilon$, $c_0$, $\alpha$, $\tau$ and $A$ are parameters regarding the AMoD service; $B$ is the scaling parameter of the vehicle searching time function; $\delta$, $l$, $s$, $C$, $Q$ and $V$ are parameters that account for vehicles' driving and charging characteristics; $p_e$, $C_{av}$, $\phi_{vc}$ and $\phi_{bs}$ are cost factors that characterize the electricity cost, autonomous vehicle operating cost, and infrastructure costs of charging stations and battery swapping stations, respectively. Below we calibrate these parameters based on real-world data in New York City.


To calibrate the parameters regarding AMoD services, we take the TNC data for the Manhattan Central Business District (CBD) in New York City on a regular weekday. It records all ride-sourcing trips (e.g., origin, destination, pick-up time and drop-off time, etc.) that started from or ended in the Manhattan CBD area, from which we obtain that the average trip length is 2.4 miles and the average trip duration $\tau$ is 16.3 min. For the value of time $\alpha$, the potential passenger demand $\lambda_0$, and the parameters $\epsilon$ and $c_0$ in the logit model, we apply "reverse engineering" to infer their values so that the profit-maximizing solution given them matches the reported TNC data. 
For the scaling parameter $A$ and $B$, we calibrate the passenger waiting time function (\ref{func_passenger_waiting_time}) and vehicle searching time function (\ref{func_vehicle_searching_time}) based on real-world geographic and traffic data of New York City using numerical simulation. The details of the calibration method can be found in Appendix D.

For vehicle charging, the involved model parameters are the average SoC $\delta$, the electricity consuming rate $l$, the average charging speed $s$, the battery capacity $C$, the number of chargers/batteries within the battery swapping station $Q$, and the buffer size of the battery swapping station $V$. Autonomous electric TNC vehicles are assumed with an average electricity consumption of 0.25 kWh/mile when traveling on roads. Thus, the hourly electricity consuming rate $l$ is
\begin{equation}
    \frac{\text{average trip length}}{\text{average trip duration}} \times \text{average energy consumption} \times 60 = \frac{2.4}{16.3} * 0.25 * 60 = 2.21 \text{ kWh/hour} .
\end{equation}
We further assume that the average SoC of vehicles that arrive at charging/battery swapping stations is $\delta=10\%$. For the charging speed $s$ and the battery capacity $C$, since charging technologies and battery technologies evolve over time, we set the nominal value as $s=22 \text{ kW}$ (roughly 1.5 mile/min, Level 2 charging) and $C=25 \text{ kWh}$ (a range of 100 miles). The impacts of different charging speeds and battery capacities will be investigated in Section \ref{impacts_charging_speed} and Section \ref{impacts_battery_capacity}, respectively. We further assume that for battery swapping stations, $Q=6$ and $V=20$.

For cost parameters $p_e$, $C_{av}$, $\phi_{vc}$ and $\phi_{bs}$, we set $p_e=0.12\$/\text{kWh}$ according to the reported average New York commercial electricity rate in February 2020 \cite{Electricity2022nyc}. The investment and operating cost of autonomous vehicles consist of the cost of vehicle automation, coordination and maintenance, etc. We convert all capital costs of the autonomous vehicle to an hourly basis and denote $C_{av}$ as the average hourly operating cost of an autonomous vehicle. We set $C_{av}=15\$/\text{hour}$ as the nominal value such that the hourly operating cost of an autonomous vehicle is considerably lower than the average driver wage of human-driven TNC vehicles \cite{li2019regulating}. 
We transform the capital infrastructure costs into an hourly basis and denote $\phi_{vc}$ and $\phi_{bs}$ as the pro-rated investment cost plus the operating cost of charging stations and battery swapping stations, respectively. {The infrastructure cost for charging stations is relatively low and is set to be $\phi_{vc}=\$8/\text{hour}$, whereas the {\em per charger} cost of battery swapping stations $\phi_{bs}$ is usually several times higher than that of charging stations $\phi_{vc}$, which is set as $\phi_{bs}=\$40/\text{hour}$. In this case, the ratio of per charger cost of battery swapping stations compared to charging stations is $r=\frac{\phi_{bs}}{\phi_{vc}}=5$. }

Overall, the nominal values of model parameters for the numerical example of New York City are summarized below:
\begin{gather*}
    \lambda_0=944/\text{min}, \quad \alpha=2.58, \quad \epsilon=0.155, \quad c_0=15.48, \quad \tau=16.3\text{ min}, \quad \delta=10\%, \\ s=22\text{kW}, \quad t_s=2\text{ min}, \quad l=2.21\text{kW}, \quad p_e=0.12\text{\$/kWh}, \quad \beta_{1,2}=1, \quad A=230, \quad B=230, \\ C=25\text{kWh}, \quad Q=6, \quad V=15, \quad C_{av}=15\text{\$/h}, \quad \phi_{vc}=8\text{\$/hour}, \quad \phi_{bs}=40\text{\$/hour}, \quad r=5.
\end{gather*}

\subsection{Impacts of charging speed} \label{impacts_charging_speed}

This section investigates how charging speed affects the charging infrastructure planning and the market outcomes. Specifically, we vary the value of the model parameter $s$ and solve the bi-level optimizations under distinct charging speeds for two different charging strategies, respectively.

Figure \ref{infrastructure_planning_s} presents the infrastructure planning of charging stations and battery swapping stations under different charging speeds.
Numerical results show that for plug-in charging, as the charging speed improves, the number of charging stations increases (Figure \ref{optimal_K_s}), the number of chargers per station reduces (Figure \ref{optimal_Q_s}), passenger surplus and TNC profit monotonically increases (Figure \ref{PS_s}-\ref{infrastructure_cost_s}), and total social welfare also increases (Figure \ref{SW_s}). For battery swapping, as the charging speed improves, the government prefers to deploy fewer battery swapping stations (Figure \ref{optimal_K_s}). In this case, the passenger surplus and social welfare monotonically increase (Figure \ref{SW_s} and Figure \ref{PS_s}), but the TNC profit initially increases and then drops (Figure \ref{profit_s}).

\begin{figure}[h]
\centering
\begin{subfigure}[b]{0.32\linewidth}
\centering
%
%
\begin{tikzpicture}

\begin{axis}[%
width=1.694in,
height=1.03in,
at={(1.358in,0.0in)},
scale only axis,
xmin=11,
xmax=44,
xtick={11, 22, 33, 44},
xlabel style={font=\color{white!15!black}},
xlabel={Charging speed},
ymin=100,
ymax=400,
ylabel style={font=\color{white!15!black}},
ylabel={Number of stations},
axis background/.style={fill=white},
legend style={at={(1,1)}, anchor=north east, legend cell align=left, align=left, font=\scriptsize, draw=white!12!black}
]
\addplot [color=black, line width=1.0pt]
  table[row sep=crcr]{%
11	138.508034687647\\
11.33	140.129632723886\\
11.66	141.892745374958\\
11.99	143.479077929857\\
12.32	145.19695345023\\
12.65	146.768689197498\\
12.98	148.372187708492\\
13.31	150.024889960815\\
13.64	151.476674483981\\
13.97	152.929353711021\\
14.3	154.543706002467\\
14.63	156.029138692075\\
14.96	157.596991558767\\
15.29	159.020639442869\\
15.62	160.435540978117\\
15.95	161.890531355404\\
16.28	163.414998075556\\
16.61	164.733331023685\\
16.94	166.112915920384\\
17.27	167.417402388822\\
17.6	168.826877438377\\
17.93	170.275489824939\\
18.26	171.585612565519\\
18.59	172.923834815641\\
18.92	174.150587459074\\
19.25	175.548684302122\\
19.58	176.802345675886\\
19.91	178.100794277566\\
20.24	179.476371977128\\
20.57	180.64394194533\\
20.9	181.821349483618\\
21.23	183.314499524822\\
21.56	184.437818400139\\
21.89	185.621465423546\\
22.22	186.878834080674\\
22.55	188.115771406125\\
22.88	189.206972121226\\
23.21	190.604880929504\\
23.54	191.774081837032\\
23.87	192.729258270587\\
24.2	193.940327845875\\
24.53	194.987832344309\\
24.86	196.294997276962\\
25.19	197.382238699844\\
25.52	198.51377140939\\
25.85	199.715737439692\\
26.18	200.747021515305\\
26.51	201.972392941481\\
26.84	203.166303370718\\
27.17	204.185279775528\\
27.5	205.243541053516\\
27.83	206.355565264923\\
28.16	207.35530675413\\
28.49	208.668720860434\\
28.82	209.710230670395\\
29.15	210.667031942559\\
29.48	211.70646077008\\
29.81	212.746859588103\\
30.14	213.800804201247\\
30.47	214.852806809461\\
30.8	215.981939625321\\
31.13	216.890966823787\\
31.46	217.917305641574\\
31.79	218.964994680813\\
32.12	219.97602768276\\
32.45	220.996793218814\\
32.78	222.035449707537\\
33.11	223.04710032917\\
33.44	224.046529075787\\
33.77	225.049457643208\\
34.1	226.012736830541\\
34.43	226.92427633878\\
34.76	227.989054972327\\
35.09	228.919273662566\\
35.42	229.987194245661\\
35.75	230.856448723948\\
36.08	231.901150135283\\
36.41	232.79616689176\\
36.74	233.68921043981\\
37.07	234.588739409696\\
37.4	235.628723422122\\
37.73	236.593623336023\\
38.06	237.473161251633\\
38.39	238.47859060427\\
38.72	239.401506830584\\
39.05	240.325591275289\\
39.38	241.285837566753\\
39.71	242.291290249227\\
40.04	243.174594809868\\
40.37	244.027340080133\\
40.7	245.037719607764\\
41.03	245.864308108554\\
41.36	246.771860493257\\
41.69	247.674352923366\\
42.02	248.57790048158\\
42.35	249.380030196908\\
42.68	250.242060521072\\
43.01	251.081357934728\\
43.34	252.071541994993\\
43.67	252.798157463998\\
44	254.139061525257\\
};
\addlegendentry{plug-in charging}

\addplot [color=blue, dashed, line width=1.0pt]
  table[row sep=crcr]{%
11	340.991973876954\\
11.33	336.523425497398\\
11.66	332.025909423828\\
11.99	327.40478515625\\
12.32	322.662353515625\\
12.65	318.142580986023\\
12.98	313.529825955629\\
13.31	309.08207893371\\
13.64	304.534912132658\\
13.97	300.207521021366\\
14.3	295.825004578001\\
14.63	291.554260346948\\
14.96	287.524414062503\\
15.29	283.28857421875\\
15.62	279.394827684155\\
15.95	275.538659095764\\
16.28	271.703338995576\\
16.61	267.999255657196\\
16.94	264.331436157227\\
17.27	260.772609894048\\
17.6	257.309724390507\\
17.93	253.910832083784\\
18.26	250.633239746094\\
18.59	247.415256500244\\
18.92	244.256782520097\\
19.25	241.114425752312\\
19.58	238.192747533321\\
19.91	235.217285901308\\
20.24	232.373043906409\\
20.57	229.562401771545\\
20.9	226.824951174785\\
21.23	224.147033714667\\
21.56	221.557712531649\\
21.89	219.0071137622\\
22.22	216.503524780273\\
22.55	214.093017577761\\
22.88	211.706161506254\\
23.21	209.329223632811\\
23.54	207.125854492188\\
23.87	204.847717293887\\
24.2	202.69193649292\\
24.53	200.590521143749\\
24.86	198.472226038575\\
25.19	196.4599609375\\
25.52	194.476317986846\\
25.85	192.504880949855\\
26.18	190.557766705751\\
26.51	188.693809531105\\
26.84	186.85607909556\\
27.17	185.046388953924\\
27.5	183.331298828125\\
27.83	181.57958684168\\
28.16	179.883003234858\\
28.49	178.216552734375\\
28.82	176.672363269608\\
29.15	175.022912025452\\
29.48	173.480224632658\\
29.81	171.911621082108\\
30.14	170.398044586182\\
30.47	168.93882743334\\
30.8	167.50488289399\\
31.13	166.0400390625\\
31.46	164.660629630089\\
31.79	163.290405273438\\
32.12	161.956214907696\\
32.45	160.631567216478\\
32.78	159.33226346969\\
33.11	158.056640625\\
33.44	156.817620992661\\
33.77	155.563357099891\\
34.1	154.39453125\\
34.43	153.198242187682\\
34.76	152.006530761719\\
35.09	150.848388668965\\
35.42	149.713230133057\\
35.75	148.62060546875\\
36.08	147.509670257568\\
36.41	146.423336886801\\
36.74	145.386505126952\\
37.07	144.334411714225\\
37.4	143.31655576824\\
37.73	142.315673828125\\
38.06	141.321563720703\\
38.39	140.353405289398\\
38.72	139.404297992581\\
39.05	138.441467290977\\
39.38	137.496936321259\\
39.71	136.572980880737\\
40.04	135.688683390617\\
40.37	134.796523862314\\
40.7	133.917236700655\\
41.03	133.05511623621\\
41.36	132.205009460449\\
41.69	131.379938124155\\
42.02	130.518341064453\\
42.35	129.718011617661\\
42.68	128.906202316284\\
43.01	128.169441223145\\
43.34	127.388763427734\\
43.67	126.587868109345\\
44	125.842666625977\\
};
\addlegendentry{battery swapping}

\end{axis}
\end{tikzpicture}%
\vspace*{-0.3in}
\caption{Optimal number of charging/battery swapping stations $K$ as a function of $s$.}
\label{optimal_K_s}
\end{subfigure}
\hfill
\begin{subfigure}[b]{0.32\linewidth}
\centering
%
%
\begin{tikzpicture}

\begin{axis}[%
width=1.694in,
height=1.03in,
at={(1.358in,0.0in)},
scale only axis,
xmin=11,
xmax=44,
xtick={11, 22, 33, 44},
xlabel style={font=\color{white!15!black}},
xlabel={Charging speed},
ymin=4,
ymax=18,
ylabel style={font=\color{white!15!black}},
ylabel={Number of chargers},
axis background/.style={fill=white},
legend style={at={(1,1)}, anchor=north east, legend cell align=left, align=left, font=\scriptsize, draw=white!12!black}
]
\addplot [color=black, line width=1.0pt]
  table[row sep=crcr]{%
11	17.2274193899362\\
11.33	16.7534633167523\\
11.66	16.2850771763972\\
11.99	15.8523693358067\\
12.32	15.4267672045248\\
12.65	15.03264687129\\
12.98	14.6535252192395\\
13.31	14.2850956974438\\
13.64	13.947999832307\\
13.97	13.6239322793121\\
14.3	13.3008320156937\\
14.63	12.9982674008408\\
14.96	12.7018803685315\\
15.29	12.4275305307434\\
15.62	12.1612424637503\\
15.95	11.9038418424772\\
16.28	11.6516195585039\\
16.61	11.4194071158086\\
16.94	11.1932130189977\\
17.27	10.9779038287364\\
17.6	10.76398352379\\
17.93	10.5547686316118\\
18.26	10.3610430465776\\
18.59	10.1689649080356\\
18.92	9.99250402715252\\
19.25	9.80881002716402\\
19.58	9.64158050555673\\
19.91	9.4749656383328\\
20.24	9.31142912359354\\
20.57	9.16020017374479\\
20.9	9.01444875716319\\
21.23	8.85917457030758\\
21.56	8.72225154749603\\
21.89	8.58674411714099\\
22.22	8.45357501994451\\
22.55	8.32250300219085\\
22.88	8.20165953917934\\
23.21	8.07306986875573\\
23.54	7.95707490491568\\
23.87	7.84893890838586\\
24.2	7.73538807575287\\
24.53	7.63002983600944\\
24.86	7.52030553588775\\
25.19	7.41934973248373\\
25.52	7.32125901135095\\
25.85	7.22163375727831\\
26.18	7.12712526745941\\
26.51	7.03294355952153\\
26.84	6.94084684062281\\
27.17	6.85474645736246\\
27.5	6.76955288178109\\
27.83	6.6856513568195\\
28.16	6.60518018790702\\
28.49	6.52209177397614\\
28.82	6.44289583557309\\
29.15	6.3700649489813\\
29.48	6.2956633000905\\
29.81	6.22351856843474\\
30.14	6.15248235670457\\
30.47	6.08286834771748\\
30.8	6.01223107020054\\
31.13	5.94799224758059\\
31.46	5.88366594472184\\
31.79	5.81914368172481\\
32.12	5.75690940865313\\
32.45	5.69737250823303\\
32.78	5.63614570172789\\
33.11	5.57635560333705\\
33.44	5.51974670226385\\
33.77	5.46368116084446\\
34.1	5.41010292119205\\
34.43	5.35677416130435\\
34.76	5.30316383760893\\
35.09	5.25062851236022\\
35.42	5.19768095797113\\
35.75	5.14998785326962\\
36.08	5.10050306971166\\
36.41	5.05285610578827\\
36.74	5.00653472173504\\
37.07	4.96106898887808\\
37.4	4.91382687659634\\
37.73	4.86920024978083\\
38.06	4.82578607608445\\
38.39	4.78155880892226\\
38.72	4.73990198882683\\
39.05	4.69783356478894\\
39.38	4.6563786964606\\
39.71	4.61591247638267\\
40.04	4.57651032894125\\
40.37	4.53818639706463\\
40.7	4.49913621916358\\
41.03	4.46309880440075\\
41.36	4.42578114088199\\
41.69	4.38997195354828\\
42.02	4.35350524105486\\
42.35	4.31894109617691\\
42.68	4.28567087445104\\
43.01	4.25218389908439\\
43.34	4.21639113981298\\
43.67	4.18665810490005\\
44	4.14727249085844\\
};
\addlegendentry{plug-in charging}

\addplot [color=blue, dashed, line width=1.0pt]
  table[row sep=crcr]{%
11	6\\
11.33	6\\
11.66	6\\
11.99	6\\
12.32	6\\
12.65	6\\
12.98	6\\
13.31	6\\
13.64	6\\
13.97	6\\
14.3	6\\
14.63	6\\
14.96	6\\
15.29	6\\
15.62	6\\
15.95	6\\
16.28	6\\
16.61	6\\
16.94	6\\
17.27	6\\
17.6	6\\
17.93	6\\
18.26	6\\
18.59	6\\
18.92	6\\
19.25	6\\
19.58	6\\
19.91	6\\
20.24	6\\
20.57	6\\
20.9	6\\
21.23	6\\
21.56	6\\
21.89	6\\
22.22	6\\
22.55	6\\
22.88	6\\
23.21	6\\
23.54	6\\
23.87	6\\
24.2	6\\
24.53	6\\
24.86	6\\
25.19	6\\
25.52	6\\
25.85	6\\
26.18	6\\
26.51	6\\
26.84	6\\
27.17	6\\
27.5	6\\
27.83	6\\
28.16	6\\
28.49	6\\
28.82	6\\
29.15	6\\
29.48	6\\
29.81	6\\
30.14	6\\
30.47	6\\
30.8	6\\
31.13	6\\
31.46	6\\
31.79	6\\
32.12	6\\
32.45	6\\
32.78	6\\
33.11	6\\
33.44	6\\
33.77	6\\
34.1	6\\
34.43	6\\
34.76	6\\
35.09	6\\
35.42	6\\
35.75	6\\
36.08	6\\
36.41	6\\
36.74	6\\
37.07	6\\
37.4	6\\
37.73	6\\
38.06	6\\
38.39	6\\
38.72	6\\
39.05	6\\
39.38	6\\
39.71	6\\
40.04	6\\
40.37	6\\
40.7	6\\
41.03	6\\
41.36	6\\
41.69	6\\
42.02	6\\
42.35	6\\
42.68	6\\
43.01	6\\
43.34	6\\
43.67	6\\
44	6\\
};
\addlegendentry{battery swapping}

\end{axis}
\end{tikzpicture}%
\vspace*{-0.3in}
\caption{Optimal number of chargers at each station $Q$ as a function of $s$.}
\label{optimal_Q_s}
\end{subfigure}
\hfill
\begin{subfigure}[b]{0.32\linewidth}
\centering
%
%
\begin{tikzpicture}

\begin{axis}[%
width=1.694in,
height=1.03in,
at={(1.358in,0.0in)},
scale only axis,
xmin=11,
xmax=44,
xtick={11, 22, 33, 44},
xlabel style={font=\color{white!15!black}},
xlabel={Charging speed},
ymin=40000,
ymax=140000,
ylabel style={font=\color{white!15!black}},
ylabel={Social welfare},
axis background/.style={fill=white},
legend style={at={(1,0)}, anchor=south east, legend cell align=left, align=left, font=\scriptsize, draw=white!12!black}
]
\addplot [color=black, line width=1.0pt]
  table[row sep=crcr]{%
11	87087.5685670908\\
11.33	88794.0957715082\\
11.66	90416.6805185042\\
11.99	91961.3852856169\\
12.32	93433.7028524909\\
12.65	94838.6339988819\\
12.98	96180.7210380689\\
13.31	97464.1173659468\\
13.64	98692.6211406563\\
13.97	99869.7046832975\\
14.3	100998.559389164\\
14.63	102082.119561166\\
14.96	103123.082471579\\
15.29	104123.935367852\\
15.62	105086.985249054\\
15.95	106014.359063923\\
16.28	106908.017563864\\
16.61	107769.799648348\\
16.94	108601.402099484\\
17.27	109404.398620075\\
17.6	110180.269076375\\
17.93	110930.378322529\\
18.26	111656.011771368\\
18.59	112358.363004253\\
18.92	113038.556943195\\
19.25	113697.639864856\\
19.58	114336.5978799\\
19.91	114956.349168364\\
20.24	115557.766026309\\
20.57	116141.664409748\\
20.9	116708.809224512\\
21.23	117259.929680107\\
21.56	117795.705165165\\
21.89	118316.77801187\\
22.22	118823.757693873\\
22.55	119317.219871008\\
22.88	119797.70878172\\
23.21	120265.738532984\\
23.54	120721.800095497\\
23.87	121166.357700263\\
24.2	121599.843365221\\
24.53	122022.676924186\\
24.86	122435.257980119\\
25.19	122837.959868114\\
25.52	123231.141392711\\
25.85	123615.144685988\\
26.18	123990.293971344\\
26.51	124356.89977293\\
26.84	124715.255347749\\
27.17	125065.646609413\\
27.5	125408.335439088\\
27.83	125743.584268711\\
28.16	126071.638415785\\
28.49	126392.728289941\\
28.82	126707.085474325\\
29.15	127014.922633355\\
29.48	127316.441733899\\
29.81	127611.843150332\\
30.14	127901.317436912\\
30.47	128185.040255996\\
30.8	128463.19399954\\
31.13	128735.934230773\\
31.46	129003.431298037\\
31.79	129265.836568093\\
32.12	129523.296587259\\
32.45	129775.952248647\\
32.78	130023.945509053\\
33.11	130267.400330096\\
33.44	130506.449046599\\
33.77	130741.212191071\\
34.1	130971.807927386\\
34.43	131198.348548891\\
34.76	131420.94181563\\
35.09	131639.697096635\\
35.42	131854.711338284\\
35.75	132066.083713086\\
36.08	132273.909311412\\
36.41	132478.279911414\\
36.74	132679.280881846\\
37.07	132876.998828399\\
37.4	133071.518411616\\
37.73	133262.916561414\\
38.06	133451.268043428\\
38.39	133636.650649728\\
38.72	133819.136094182\\
39.05	133998.79206248\\
39.38	134175.685542335\\
39.71	134349.881983712\\
40.04	134521.447684038\\
40.37	134690.440252611\\
40.7	134856.922568904\\
41.03	135020.947238982\\
41.36	135182.571602091\\
41.69	135341.855925635\\
42.02	135498.839597398\\
42.35	135653.583379997\\
42.68	135806.134487279\\
43.01	135956.540884646\\
43.34	136104.849492319\\
43.67	136251.102531195\\
44	136395.343732463\\
};
\addlegendentry{plug-in charging}

\addplot [color=blue, dashed, line width=1.0pt]
  table[row sep=crcr]{%
11	41109.6787608299\\
11.33	43790.2004122158\\
11.66	46354.829829337\\
11.99	48809.728526642\\
12.32	51160.7936392518\\
12.65	53413.6288285554\\
12.98	55573.5547751771\\
13.31	57645.5891373623\\
13.64	59634.4865977766\\
13.97	61544.6878990249\\
14.3	63380.3869649874\\
14.63	65145.5033022373\\
14.96	66843.745994834\\
15.29	68478.5832548594\\
15.62	70053.2502605293\\
15.95	71570.7997968808\\
16.28	73034.1021716926\\
16.61	74445.8287031991\\
16.94	75808.5161701507\\
17.27	77124.5414320367\\
17.6	78396.1114629961\\
17.93	79625.3354039503\\
18.26	80814.1932723037\\
18.59	81964.5244473063\\
18.92	83078.0989309607\\
19.25	84156.5517186925\\
19.58	85201.4516145946\\
19.91	86214.2570246323\\
20.24	87196.3553669077\\
20.57	88149.0589220533\\
20.9	89073.6092938793\\
21.23	89971.173434069\\
21.56	90842.8730541279\\
21.89	91689.7490490449\\
22.22	92512.8095601187\\
22.55	93312.9879070476\\
22.88	94091.1884991383\\
23.21	94848.2657049033\\
23.54	95585.0208047205\\
23.87	96302.2356772732\\
24.2	97000.6293991259\\
24.53	97680.9103397364\\
24.86	98343.7329642027\\
25.19	98989.739973278\\
25.52	99619.5212486273\\
25.85	100233.662075348\\
26.18	100832.705528553\\
26.51	101417.174342456\\
26.84	101987.571575611\\
27.17	102544.375265157\\
27.5	103088.039338741\\
27.83	103619.003163111\\
28.16	104137.684081051\\
28.49	104644.478167178\\
28.82	105139.76861909\\
29.15	105623.938413979\\
29.48	106097.321217975\\
29.81	106560.263611526\\
30.14	107013.085695722\\
30.47	107456.097940543\\
30.8	107889.601443984\\
31.13	108313.88196434\\
31.46	108729.213399727\\
31.79	109135.863228788\\
32.12	109534.082331659\\
32.45	109924.117583712\\
32.78	110306.204792841\\
33.11	110680.566841256\\
33.44	111047.429539151\\
33.77	111407.000887183\\
34.1	111759.471174147\\
34.43	112105.055694801\\
34.76	112443.932249666\\
35.09	112776.283112494\\
35.42	113102.28151802\\
35.75	113422.101288161\\
36.08	113735.898340568\\
36.41	114043.837139346\\
36.74	114346.069210776\\
37.07	114642.740754586\\
37.4	114933.994127712\\
37.73	115219.959647779\\
38.06	115500.780225736\\
38.39	115776.581387491\\
38.72	116047.486059948\\
39.05	116313.619057646\\
39.38	116575.086096629\\
39.71	116832.008642169\\
40.04	117084.497186521\\
40.37	117332.653903903\\
40.7	117576.579282339\\
41.03	117816.376522742\\
41.36	118052.140984347\\
41.69	118283.966445964\\
42.02	118511.945606226\\
42.35	118736.156457472\\
42.68	118956.696490625\\
43.01	119173.638258014\\
43.34	119387.07227721\\
43.67	119597.075960468\\
44	119803.716647928\\
};
\addlegendentry{battery swapping}
\end{axis}
\end{tikzpicture}%
\vspace*{-0.3in}
\caption{Social welfare (\$/hour) as a function of $s$.}
\label{SW_s}
\end{subfigure}
\begin{subfigure}[b]{0.32\linewidth}
\centering
%
%
\begin{tikzpicture}

\begin{axis}[%
width=1.694in,
height=1.03in,
at={(1.358in,0.0in)},
scale only axis,
xmin=11,
xmax=44,
xtick={11, 22, 33, 44},
xlabel style={font=\color{white!15!black}},
xlabel={Charging speed},
ymin=75000,
ymax=100000,
ylabel style={font=\color{white!15!black}},
ylabel={Passengers' surplus},
axis background/.style={fill=white},
legend style={at={(1,0)}, anchor=south east, legend cell align=left, align=left, font=\scriptsize, draw=white!12!black}
]
\addplot [color=black, line width=1.0pt]
  table[row sep=crcr]{%
11	79101.1255023089\\
11.33	79824.3558518928\\
11.66	80507.3957548464\\
11.99	81151.1417479823\\
12.32	81762.702064895\\
12.65	82342.9240153728\\
12.98	82896.0464694335\\
13.31	83421.8919149315\\
13.64	83922.6493609173\\
13.97	84400.4239670324\\
14.3	84857.9151905325\\
14.63	85293.573020702\\
14.96	85711.4056920854\\
15.29	86112.5852759544\\
15.62	86494.6625391148\\
15.95	86863.6321046109\\
16.28	87218.663378131\\
16.61	87558.1542620624\\
16.94	87886.9576350899\\
17.27	88202.355657976\\
17.6	88506.6894144033\\
17.93	88799.9959130618\\
18.26	89084.4244804103\\
18.59	89355.6528517347\\
18.92	89622.2896563811\\
19.25	89875.3935437002\\
19.58	90124.7166049968\\
19.91	90363.4485804136\\
20.24	90596.7182605184\\
20.57	90820.5071541205\\
20.9	91039.4764613359\\
21.23	91251.8868125368\\
21.56	91455.9867976726\\
21.89	91654.4036972795\\
22.22	91849.1750854609\\
22.55	92035.6240392951\\
22.88	92218.3243736829\\
23.21	92396.8808297518\\
23.54	92571.5861658896\\
23.87	92738.1493391378\\
24.2	92901.6147834986\\
24.53	93060.4677356359\\
24.86	93217.7043380674\\
25.19	93369.3705543334\\
25.52	93519.7563847498\\
25.85	93664.3965239577\\
26.18	93802.1659282233\\
26.51	93942.0548825076\\
26.84	94077.2679693769\\
27.17	94207.6412553207\\
27.5	94335.2192907602\\
27.83	94461.1773869252\\
28.16	94582.2359566752\\
28.49	94706.386885448\\
28.82	94820.6772343342\\
29.15	94935.7992687326\\
29.48	95047.236500178\\
29.81	95157.3352013807\\
30.14	95265.1058307794\\
30.47	95370.5539158759\\
30.8	95473.1078927628\\
31.13	95573.1306950794\\
31.46	95673.1803404046\\
31.79	95769.931261587\\
32.12	95865.0606817766\\
32.45	95960.7790084327\\
32.78	96051.2671262836\\
33.11	96139.5891998282\\
33.44	96228.810511321\\
33.77	96315.7780450063\\
34.1	96402.0572783012\\
34.43	96484.2202393082\\
34.76	96567.6251066562\\
35.09	96645.8099776473\\
35.42	96724.6186961489\\
35.75	96802.6905088248\\
36.08	96880.542836925\\
36.41	96954.3453635516\\
36.74	97027.6068205595\\
37.07	97099.9364787221\\
37.4	97171.1056759988\\
37.73	97241.7913846048\\
38.06	97309.7815003756\\
38.39	97377.6944039968\\
38.72	97445.145294257\\
39.05	97509.9394771944\\
39.38	97574.5368076539\\
39.71	97639.7501603325\\
40.04	97701.5668439217\\
40.37	97762.3926431286\\
40.7	97824.1341048804\\
41.03	97884.180736801\\
41.36	97942.483733002\\
41.69	98001.2544789531\\
42.02	98057.3459211901\\
42.35	98112.3590367881\\
42.68	98169.1784756205\\
43.01	98223.586253671\\
43.34	98276.3002399879\\
43.67	98330.9734031213\\
44	98382.8058146346\\
};
\addlegendentry{plug-in charging}
\addplot [color=blue, dashed, line width=1.0pt]
  table[row sep=crcr]{%
11	77404.699080199\\
11.33	78578.4381544809\\
11.66	79670.0061860093\\
11.99	80665.9380523207\\
12.32	81569.4796271879\\
12.65	82444.1398215026\\
12.98	83237.5139400057\\
13.31	83998.191766225\\
13.64	84682.1318310505\\
13.97	85348.9662533933\\
14.3	85952.2865861727\\
14.63	86524.4897855641\\
14.96	87091.3966452223\\
15.29	87574.6948876966\\
15.62	88075.5971205959\\
15.95	88541.1806057149\\
16.28	88970.3285573553\\
16.61	89384.9780467864\\
16.94	89769.6081066342\\
17.27	90139.1802959252\\
17.6	90492.896424291\\
17.93	90826.4217814529\\
18.26	91151.7166537915\\
18.59	91458.7841772917\\
18.92	91748.8903067703\\
19.25	92015.2106774245\\
19.58	92296.6710902228\\
19.91	92543.6614637603\\
20.24	92791.4004236938\\
20.57	93022.8042797873\\
20.9	93246.1290139024\\
21.23	93459.7153623729\\
21.56	93669.8067907365\\
21.89	93867.9222371235\\
22.22	94056.3587994599\\
22.55	94244.3928031812\\
22.88	94419.7899161034\\
23.21	94580.6310441527\\
23.54	94758.0259840941\\
23.87	94906.2781103258\\
24.2	95062.6927850306\\
24.53	95215.1375292848\\
24.86	95350.803596433\\
25.19	95493.1404855014\\
25.52	95628.1410878555\\
25.85	95753.1849110625\\
26.18	95870.9514675627\\
26.51	95992.8456877393\\
26.84	96108.6111823749\\
27.17	96218.9550810889\\
27.5	96336.6832264239\\
27.83	96437.5785112342\\
28.16	96539.1350689158\\
28.49	96636.9803988303\\
28.82	96748.663675617\\
29.15	96831.8778284737\\
29.48	96926.5852945466\\
29.81	97008.2087172611\\
30.14	97092.1553637043\\
30.47	97178.5410447488\\
30.8	97262.1541876292\\
31.13	97332.6982251556\\
31.46	97412.1703195508\\
31.79	97486.5022111923\\
32.12	97560.9156120758\\
32.45	97630.6438814665\\
32.78	97698.8296152986\\
33.11	97765.3346575115\\
33.44	97832.7552891661\\
33.77	97891.5511996381\\
34.1	97960.7075879514\\
34.43	98019.272258468\\
34.76	98073.3793019029\\
35.09	98128.5944391685\\
35.42	98183.0720947779\\
35.75	98240.6026009313\\
36.08	98289.9348463124\\
36.41	98339.220271022\\
36.74	98393.2260266742\\
37.07	98439.9708175452\\
37.4	98488.8099290761\\
37.73	98536.6199522307\\
38.06	98581.6115404554\\
38.39	98627.4610226097\\
38.72	98672.9762560266\\
39.05	98712.1456743184\\
39.38	98751.0226633543\\
39.71	98790.1273138008\\
40.04	98833.0983176258\\
40.37	98871.1649440966\\
40.7	98908.2677174732\\
41.03	98945.2872574621\\
41.36	98981.3390804112\\
41.69	99018.917577491\\
42.02	99046.6543324335\\
42.35	99082.782266424\\
42.68	99113.8640277267\\
43.01	99156.0188400441\\
43.34	99187.2333049132\\
43.67	99211.9750804428\\
44	99244.3764865544\\
};
\addlegendentry{battery swapping}
\end{axis}
\end{tikzpicture}%
\vspace*{-0.3in}
\caption{Passengers' surplus (\$/hour) as a function of $s$.} 
\label{PS_s}
\end{subfigure}
\hfill
\begin{subfigure}[b]{0.32\linewidth}
\centering
%
%
\begin{tikzpicture}

\begin{axis}[%
width=1.694in,
height=1.03in,
at={(1.358in,0.0in)},
scale only axis,
xmin=11,
xmax=44,
xtick={11, 22, 33, 44},
xlabel style={font=\color{white!15!black}},
xlabel={Charging speed},
ymin=25000,
ymax=55000,
ylabel style={font=\color{white!15!black}},
ylabel={Platform profit},
axis background/.style={fill=white},
legend style={at={(1,0)}, anchor=south east, legend cell align=left, align=left, font=\scriptsize, draw=white!12!black}
]
\addplot [color=black, line width=1.0pt]
  table[row sep=crcr]{%
11	27075.5310843013\\
11.33	27750.9932110522\\
11.66	28395.1592372742\\
11.99	29006.1102200753\\
12.32	29590.3575850592\\
12.65	30146.2849952539\\
12.98	30678.0793239964\\
13.31	31187.1847317251\\
13.64	31672.3448221468\\
13.97	32137.2739840884\\
14.3	32585.123179604\\
14.63	33013.414276804\\
14.96	33425.9018452527\\
15.29	33821.2209054543\\
15.62	34201.0868190416\\
15.95	34567.6812077071\\
16.28	34921.7492875728\\
16.61	35260.901166307\\
16.94	35589.1424892237\\
17.27	35905.18010355\\
17.6	36211.5774789281\\
17.93	36508.1295993598\\
18.26	36794.0346346751\\
18.59	37070.361416544\\
18.92	37337.8708589398\\
19.25	37597.6358798608\\
19.58	37849.1136701457\\
19.91	38092.8918354683\\
20.24	38330.4999019883\\
20.57	38559.0346027754\\
20.9	38781.4866662027\\
21.23	39000.1640880417\\
21.56	39209.4227431521\\
21.89	39413.4465245165\\
22.22	39612.936556738\\
22.55	39806.3484100083\\
22.88	39993.8737502552\\
23.21	40178.9898717898\\
23.54	40357.8998015971\\
23.87	40529.9697533196\\
24.2	40699.8581771355\\
24.53	40864.3130161175\\
24.86	41027.1404795636\\
25.19	41184.1721929382\\
25.52	41338.3509104273\\
25.85	41488.9394528642\\
26.18	41634.1213983925\\
26.51	41778.5284115341\\
26.84	41919.1569373457\\
27.17	42055.1119395867\\
27.5	42188.3721867741\\
27.83	42319.3778009917\\
28.16	42446.3557713474\\
28.49	42573.9937869725\\
28.82	42695.5376148974\\
29.15	42814.8247732877\\
29.48	42931.8659974186\\
29.81	43046.7801971333\\
30.14	43159.4570117121\\
30.47	43269.8570437965\\
30.8	43378.3527309172\\
31.13	43483.3298495986\\
31.46	43587.4719973832\\
31.79	43689.4154290317\\
32.12	43789.2924130426\\
32.45	43887.9816729547\\
32.78	43984.0715467719\\
33.11	44078.1307120973\\
33.44	44171.0792554362\\
33.77	44262.2220019324\\
34.1	44351.7679911128\\
34.43	44438.7851100971\\
34.76	44525.8232025731\\
35.09	44609.647641559\\
35.42	44693.293122998\\
35.75	44774.6564584796\\
36.08	44855.8666995641\\
36.41	44934.2188141278\\
36.74	45011.4592305808\\
37.07	45087.5497114809\\
37.4	45163.1227680153\\
37.73	45237.2990155639\\
38.06	45309.4439431477\\
38.39	45381.3514908767\\
38.72	45451.9082267607\\
39.05	45520.9296186519\\
39.38	45589.2946051093\\
39.71	45657.2949400216\\
40.04	45723.0091991841\\
40.37	45787.5800515905\\
40.7	45852.4530988121\\
41.03	45915.3000986943\\
41.36	45977.3538392602\\
41.69	46038.869150256\\
42.02	46098.9752166635\\
42.35	46157.6856310751\\
42.68	46216.596894361\\
43.01	46274.1074915372\\
43.34	46331.2069824649\\
43.67	46387.1647668771\\
44	46444.4094275591\\
};
\addlegendentry{plug-in charging}
\addplot [color=blue, dashed, line width=1.0pt]
  table[row sep=crcr]{%
11	46703.6315137884\\
11.33	47061.0125679377\\
11.66	47369.8692954286\\
11.99	47654.7653905325\\
12.32	47905.2199499887\\
12.65	48134.3167380377\\
12.98	48335.0990918421\\
13.31	48524.5963227168\\
13.64	48699.3683074679\\
13.97	48851.1462926\\
14.3	48995.4518017569\\
14.63	49117.5204976131\\
14.96	49242.3426977872\\
15.29	49350.6635695322\\
15.62	49448.6412455633\\
15.95	49539.240677832\\
16.28	49625.145204832\\
16.61	49700.6956827571\\
16.94	49772.1677146407\\
17.27	49839.6791780572\\
17.6	49900.5864121853\\
17.93	49953.9404033927\\
18.26	50005.8540231847\\
18.59	50052.6621394582\\
18.92	50097.3941820545\\
19.25	50135.1495369641\\
19.58	50173.7437525988\\
19.91	50206.4426314491\\
20.24	50237.9474304823\\
20.57	50269.3008541475\\
20.9	50296.3569266533\\
21.23	50319.5120146325\\
21.56	50341.4717161048\\
21.89	50362.3534808597\\
22.22	50380.1896781732\\
22.55	50399.0285998286\\
22.88	50413.1326301562\\
23.21	50426.3707093637\\
23.54	50440.0897930226\\
23.87	50448.696437851\\
24.2	50460.6397746424\\
24.53	50470.9084344728\\
24.86	50478.5585000146\\
25.19	50481.4478954423\\
25.52	50492.3872688913\\
25.85	50495.0856385398\\
26.18	50500.554765665\\
26.51	50502.3412115937\\
26.84	50505.3244975003\\
27.17	50508.4829952085\\
27.5	50508.5692458844\\
27.83	50509.9256537507\\
28.16	50509.6539949471\\
28.49	50509.7940553934\\
28.82	50507.3556785861\\
29.15	50507.1426575574\\
29.48	50504.3593454115\\
29.81	50501.9881590154\\
30.14	50499.4647209661\\
30.47	50496.02401781\\
30.8	50490.8879823051\\
31.13	50487.7649983137\\
31.46	50483.9736562942\\
31.79	50477.1767318664\\
32.12	50473.7073693467\\
32.45	50466.8507936793\\
32.78	50461.7615767075\\
33.11	50457.1371492972\\
33.44	50450.6306042548\\
33.77	50444.6803434066\\
34.1	50438.1356534365\\
34.43	50431.0522388557\\
34.76	50422.8445220151\\
35.09	50416.0564175383\\
35.42	50408.0062524154\\
35.75	50399.9399123725\\
36.08	50392.4979677407\\
36.41	50384.4298797892\\
36.74	50377.2083709562\\
37.07	50368.7418811017\\
37.4	50360.448079264\\
37.73	50351.129442669\\
38.06	50342.3374997444\\
38.39	50332.956710278\\
38.72	50325.6485082977\\
39.05	50314.9820444891\\
39.38	50307.542295039\\
39.71	50296.9953170594\\
40.04	50288.4753930593\\
40.37	50278.9061777964\\
40.7	50268.4358005689\\
41.03	50260.4015460943\\
41.36	50249.7716812274\\
41.69	50239.60759914\\
42.02	50229.0558380074\\
42.35	50221.0913465547\\
42.68	50210.4810485893\\
43.01	50201.056981607\\
43.34	50190.9913287942\\
43.67	50179.9587001958\\
44	50170.0483824211\\
};
\addlegendentry{battery swapping}
\end{axis}
\end{tikzpicture}%
\vspace*{-0.3in}
\caption{Platform profit (\$/hour) as a function of $s$.}
\label{profit_s}
\end{subfigure}
\hfill
\begin{subfigure}[b]{0.32\linewidth}
\centering
%
%
\begin{tikzpicture}

\begin{axis}[%
width=1.694in,
height=1.03in,
at={(1.358in,0.0in)},
scale only axis,
xmin=11,
xmax=44,
xtick={11, 22, 33, 44},
xlabel style={font=\color{white!15!black}},
xlabel={Charging speed},
ymin=8000,
ymax=85000,
ylabel style={font=\color{white!15!black}},
ylabel={Infrastructure cost},
axis background/.style={fill=white},
legend style={at={(1,1)}, anchor=north east, legend cell align=left, align=left, font=\scriptsize, draw=white!12!black}
]
\addplot [color=black, line width=1.0pt]
  table[row sep=crcr]{%
11	19089.0880195194\\
11.33	18781.2532914368\\
11.66	18485.8744736165\\
11.99	18195.8666824407\\
12.32	17919.3567974634\\
12.65	17650.5750117448\\
12.98	17393.404755361\\
13.31	17144.9592807098\\
13.64	16902.3730424079\\
13.97	16667.9932678233\\
14.3	16444.4789809724\\
14.63	16224.8677363397\\
14.96	16014.2250657594\\
15.29	15809.8708135568\\
15.62	15608.7641091027\\
15.95	15416.9542483947\\
16.28	15232.3951018402\\
16.61	15049.2557800218\\
16.94	14874.6980248297\\
17.27	14703.1371414507\\
17.6	14537.9978169568\\
17.93	14377.7471898928\\
18.26	14222.4473437178\\
18.59	14067.6512640256\\
18.92	13921.6035721262\\
19.25	13775.3895587048\\
19.58	13637.2323952426\\
19.91	13499.9912475177\\
20.24	13369.4521361979\\
20.57	13237.8773471484\\
20.9	13112.1539030267\\
21.23	12992.1212204717\\
21.56	12869.7043756592\\
21.89	12751.0722099258\\
22.22	12638.3539483259\\
22.55	12524.7525782954\\
22.88	12414.4893422183\\
23.21	12310.132168558\\
23.54	12207.6858719895\\
23.87	12101.7613921949\\
24.2	12001.6295954126\\
24.53	11902.1038275671\\
24.86	11809.586837512\\
25.19	11715.5828791578\\
25.52	11626.965902466\\
25.85	11538.1912908337\\
26.18	11445.9933552716\\
26.51	11363.6835211115\\
26.84	11281.1695589733\\
27.17	11197.1065854949\\
27.5	11115.2560384463\\
27.83	11036.9709192055\\
28.16	10956.9533122381\\
28.49	10887.6523824797\\
28.82	10809.1293749069\\
29.15	10735.7014086658\\
29.48	10662.660763698\\
29.81	10592.2722481819\\
30.14	10523.2454055793\\
30.47	10455.3707036762\\
30.8	10388.2666241402\\
31.13	10320.5263139052\\
31.46	10257.221039751\\
31.79	10193.5101225261\\
32.12	10131.0565075602\\
32.45	10072.8084327403\\
32.78	10011.3931640028\\
33.11	9950.31958182917\\
33.44	9893.44072015787\\
33.77	9836.78785586768\\
34.1	9782.01734202815\\
34.43	9724.65680051412\\
34.76	9672.50649359903\\
35.09	9615.76052257167\\
35.42	9563.20048086305\\
35.75	9511.26325421836\\
36.08	9462.5002250774\\
36.41	9410.28426626509\\
36.74	9359.78516929404\\
37.07	9310.48736180356\\
37.4	9262.71003239768\\
37.73	9216.17383875452\\
38.06	9167.9574000951\\
38.39	9122.39524514572\\
38.72	9077.9174268354\\
39.05	9032.07703336641\\
39.38	8988.14587042785\\
39.71	8947.16311664211\\
40.04	8903.1283590677\\
40.37	8859.53244210818\\
40.7	8819.66463478833\\
41.03	8778.53359651286\\
41.36	8737.26597017136\\
41.69	8698.26770357435\\
42.02	8657.48154045576\\
42.35	8616.46128786612\\
42.68	8579.64088270218\\
43.01	8541.15286056236\\
43.34	8502.65773013345\\
43.67	8467.03563880358\\
44	8431.87150973023\\
};
\addlegendentry{plug-in charging}
\addplot [color=blue, dashed, line width=1.0pt]
  table[row sep=crcr]{%
11	81838.0737304689\\
11.33	80765.6221193756\\
11.66	79686.2182617188\\
11.99	78577.1484375\\
12.32	77438.96484375\\
12.65	76354.2194366455\\
12.98	75247.158229351\\
13.31	74179.6989440904\\
13.64	73088.3789118379\\
13.97	72049.8050451279\\
14.3	70998.0010987201\\
14.63	69973.0224832674\\
14.96	69005.8593750007\\
15.29	67989.2578125\\
15.62	67054.7586441971\\
15.95	66129.2781829834\\
16.28	65208.8013589382\\
16.61	64319.8213577271\\
16.94	63439.5446777344\\
17.27	62585.4263745714\\
17.6	61754.3338537216\\
17.93	60938.5997001082\\
18.26	60151.9775390625\\
18.59	59379.6615600586\\
18.92	58621.6278048232\\
19.25	57867.4621805549\\
19.58	57166.2594079972\\
19.91	56452.148616314\\
20.24	55769.5305375382\\
20.57	55094.9764251709\\
20.9	54437.9882819485\\
21.23	53795.28809152\\
21.56	53173.8510075957\\
21.89	52561.707302928\\
22.22	51960.8459472656\\
22.55	51382.3242186627\\
22.88	50809.4787615009\\
23.21	50239.0136718747\\
23.54	49710.205078125\\
23.87	49163.452150533\\
24.2	48646.0647583009\\
24.53	48141.7250744998\\
24.86	47633.334249258\\
25.19	47150.390625\\
25.52	46674.316316843\\
25.85	46201.1714279652\\
26.18	45733.8640093804\\
26.51	45286.5142874653\\
26.84	44845.4589829344\\
27.17	44411.1333489418\\
27.5	43999.51171875\\
27.83	43579.1008420032\\
28.16	43171.9207763659\\
28.49	42771.97265625\\
28.82	42401.367184706\\
29.15	42005.4988861084\\
29.48	41635.253911838\\
29.81	41258.789059706\\
30.14	40895.5307006836\\
30.47	40545.3185840015\\
30.8	40201.1718945576\\
31.13	39849.609375\\
31.46	39518.5511112213\\
31.79	39189.697265625\\
32.12	38869.4915778469\\
32.45	38551.5761319548\\
32.78	38239.7432327257\\
33.11	37933.59375\\
33.44	37636.2290382385\\
33.77	37335.2057039738\\
34.1	37054.6875\\
34.43	36767.5781250437\\
34.76	36481.5673828125\\
35.09	36203.6132805515\\
35.42	35931.1752319336\\
35.75	35668.9453125\\
36.08	35402.3208618164\\
36.41	35141.6008528322\\
36.74	34892.7612304685\\
37.07	34640.2588114141\\
37.4	34395.9733843776\\
37.73	34155.76171875\\
38.06	33917.1752929688\\
38.39	33684.8172694555\\
38.72	33457.0315182195\\
39.05	33225.9521498345\\
39.38	32999.2647171021\\
39.71	32777.515411377\\
40.04	32565.2840137481\\
40.37	32351.1657269554\\
40.7	32140.1368081571\\
41.03	31933.2278966904\\
41.36	31729.2022705078\\
41.69	31531.1851497972\\
42.02	31324.4018554688\\
42.35	31132.3227882385\\
42.68	30937.4885559082\\
43.01	30760.6658935547\\
43.34	30573.3032226563\\
43.67	30381.0883462428\\
44	30202.2399902344\\
};
\addlegendentry{battery swapping}
\end{axis}
\end{tikzpicture}%
\vspace*{-0.3in}
\caption{Infrastructure cost (\$/hour) as a function of $s$.}
\label{infrastructure_cost_s}
\end{subfigure}
\caption{Charging infrastructure planning under different charging speeds. Black lines present the results when deploying charging stations; blue dashed lines show the results when deploying battery swapping stations.}
\label{infrastructure_planning_s}
\end{figure}
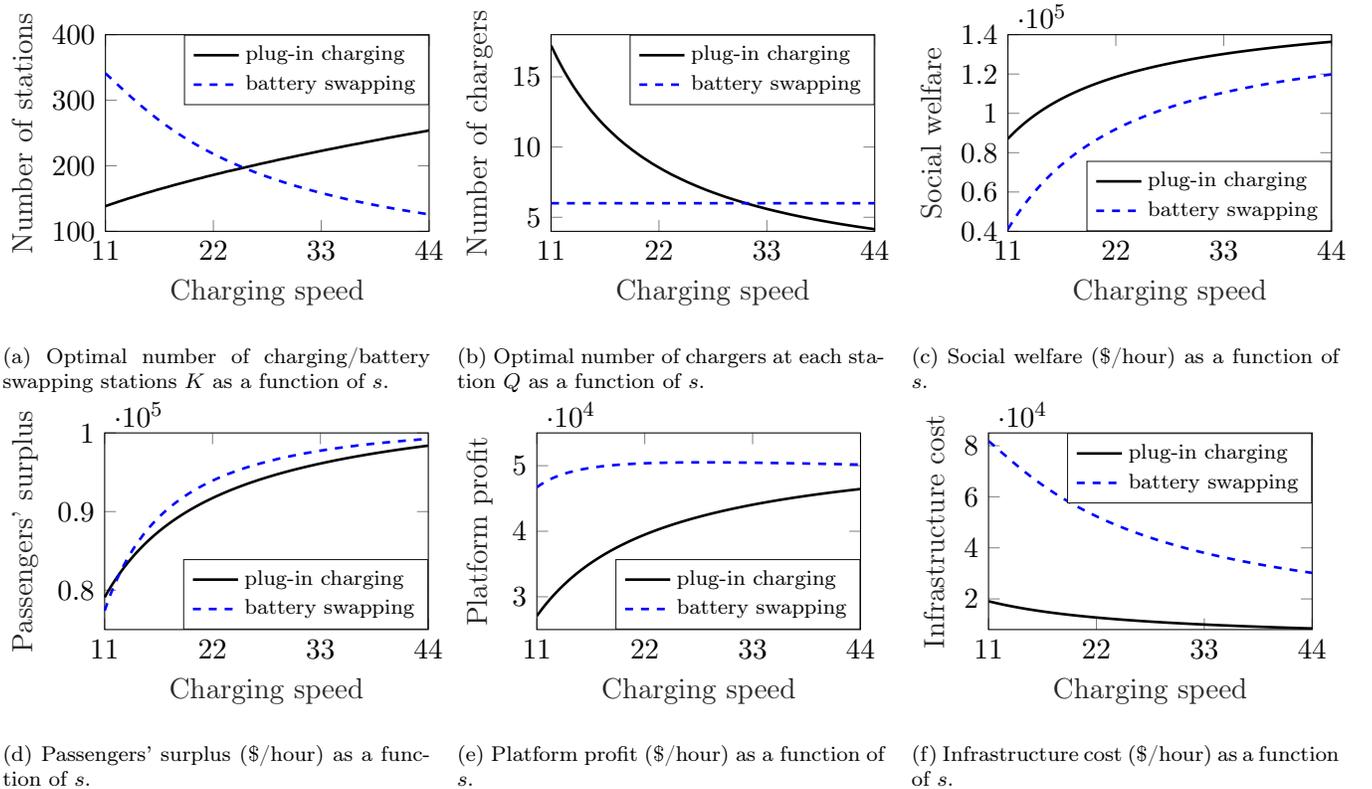

\begin{figure}[th!]
\centering
\begin{subfigure}[b]{0.32\linewidth}
\centering
%
%
\begin{tikzpicture}

\begin{axis}[%
width=1.694in,
height=1.03in,
at={(1.358in,0.0in)},
scale only axis,
xmin=11,
xmax=44,
xtick={11, 22, 33, 44},
xlabel style={font=\color{white!15!black}},
xlabel={Charging speed},
ymin=180,
ymax=225,
ylabel style={font=\color{white!15!black}},
ylabel={Passenger arrivals},
axis background/.style={fill=white},
legend style={at={(1,0)}, anchor=south east, legend cell align=left, align=left, font=\scriptsize, draw=white!12!black}
]
\addplot [color=black, line width=1.0pt]
  table[row sep=crcr]{%
11	183.740778404079\\
11.33	185.243979522721\\
11.66	186.660917537579\\
11.99	187.993919771672\\
12.32	189.258101983792\\
12.65	190.455549042665\\
12.98	191.595300016427\\
13.31	192.67724612295\\
13.64	193.706126475834\\
13.97	194.686471693162\\
14.3	195.62399691941\\
14.63	196.515689125631\\
14.96	197.369899144505\\
15.29	198.189145270147\\
15.62	198.968546921065\\
15.95	199.720437022834\\
16.28	200.443207160278\\
16.61	201.133683803352\\
16.94	201.801812322855\\
17.27	202.442136500731\\
17.6	203.059474186413\\
17.93	203.65395686497\\
18.26	204.229989937235\\
18.59	204.778872110308\\
18.92	205.318065406251\\
19.25	205.829528558203\\
19.58	206.333005311535\\
19.91	206.814772857753\\
20.24	207.285213419369\\
20.57	207.736251688081\\
20.9	208.177308974551\\
21.23	208.604902469837\\
21.56	209.01553266164\\
21.89	209.414509212583\\
22.22	209.805944724356\\
22.55	210.180459086638\\
22.88	210.54725839204\\
23.21	210.905561063393\\
23.54	211.255966419316\\
23.87	211.589885130527\\
24.2	211.917445695898\\
24.53	212.235623153864\\
24.86	212.550426931386\\
25.19	212.85394990844\\
25.52	213.154786145321\\
25.85	213.444011753861\\
26.18	213.719392056707\\
26.51	213.998902809889\\
26.84	214.26896908057\\
27.17	214.529274019766\\
27.5	214.783908037094\\
27.83	215.035221606071\\
28.16	215.27667793923\\
28.49	215.524219042226\\
28.82	215.752025072106\\
29.15	215.981416817767\\
29.48	216.203397414023\\
29.81	216.422645216768\\
30.14	216.637192987379\\
30.47	216.847055844407\\
30.8	217.051100764148\\
31.13	217.250054417464\\
31.46	217.449006991573\\
31.79	217.641348124923\\
32.12	217.830416063269\\
32.45	218.020604770716\\
32.78	218.200355446802\\
33.11	218.375760450279\\
33.44	218.552908273579\\
33.77	218.725539626002\\
34.1	218.896764101914\\
34.43	219.05978209983\\
34.76	219.225226656079\\
35.09	219.380282374063\\
35.42	219.536541732167\\
35.75	219.691306752267\\
36.08	219.845603756195\\
36.41	219.991844032662\\
36.74	220.136982966785\\
37.07	220.280247363106\\
37.4	220.421185542147\\
37.73	220.561139088771\\
38.06	220.695729979178\\
38.39	220.830143025988\\
38.72	220.963616933188\\
39.05	221.091810467094\\
39.38	221.219591908611\\
39.71	221.348569006836\\
40.04	221.470807021825\\
40.37	221.591065452372\\
40.7	221.713113765315\\
41.03	221.831792020892\\
41.36	221.947005435277\\
41.69	222.06312457415\\
42.02	222.173932532847\\
42.35	222.282593750142\\
42.68	222.394805622571\\
43.01	222.502238390087\\
43.34	222.606311366827\\
43.67	222.714236477861\\
44	222.816539020576\\
};
\addlegendentry{plug-in charging}
\addplot [color=blue, dashed, line width=1.0pt]
  table[row sep=crcr]{%
11	180.203138250541\\
11.33	182.652543184675\\
11.66	184.923420488162\\
11.99	186.989428420193\\
12.32	188.858912312405\\
12.65	190.664240908882\\
12.98	192.298058978153\\
13.31	193.861217181837\\
13.64	195.26390736046\\
13.97	196.62899027862\\
14.3	197.861908282923\\
14.63	199.029357410581\\
14.96	200.184199192869\\
15.29	201.16730872053\\
15.62	202.184856490269\\
15.95	203.129406758512\\
16.28	203.998973382238\\
16.61	204.838192815869\\
16.94	205.615804274501\\
17.27	206.362202386424\\
17.6	207.075870516282\\
17.93	207.748168655231\\
18.26	208.403285505041\\
18.59	209.021159136189\\
18.92	209.604426683486\\
19.25	210.139464700187\\
19.58	210.704495422752\\
19.91	211.199969273993\\
20.24	211.696608517789\\
20.57	212.160196947465\\
20.9	212.607321587107\\
21.23	213.034693166086\\
21.56	213.454828021476\\
21.89	213.850792422541\\
22.22	214.227212927109\\
22.55	214.602635843218\\
22.88	214.952654230557\\
23.21	215.273477397885\\
23.54	215.627156049712\\
23.87	215.922599954916\\
24.2	216.234180659508\\
24.53	216.537724886895\\
24.86	216.807753353021\\
25.19	217.090951614645\\
25.52	217.359451444099\\
25.85	217.608059972134\\
26.18	217.842122271293\\
26.51	218.084308914289\\
26.84	218.314243858787\\
27.17	218.533342562243\\
27.5	218.767030478306\\
27.83	218.96724561524\\
28.16	219.168717154486\\
28.49	219.362773273565\\
28.82	219.584210624479\\
29.15	219.749157505953\\
29.48	219.936840706915\\
29.81	220.098556041408\\
30.14	220.264836557636\\
30.47	220.435908409173\\
30.8	220.601451218909\\
31.13	220.741089495354\\
31.46	220.898368111798\\
31.79	221.045443095839\\
32.12	221.192649392672\\
32.45	221.33056025234\\
32.78	221.465394787993\\
33.11	221.596881589071\\
33.44	221.730154173318\\
33.77	221.846357916731\\
34.1	221.983014094243\\
34.43	222.098720331714\\
34.76	222.20560315483\\
35.09	222.314658587556\\
35.42	222.422241271583\\
35.75	222.535835325897\\
36.08	222.633227662643\\
36.41	222.730514436897\\
36.74	222.837103824651\\
37.07	222.929349766922\\
37.4	223.025716027133\\
37.73	223.120039280422\\
38.06	223.208790840362\\
38.39	223.299223462652\\
38.72	223.388985597084\\
39.05	223.466223984041\\
39.38	223.542877541886\\
39.71	223.619971752043\\
40.04	223.704678904328\\
40.37	223.779709918586\\
40.7	223.852833614092\\
41.03	223.925785870503\\
41.36	223.996824000219\\
41.69	224.070862897528\\
42.02	224.125506254436\\
42.35	224.196674604561\\
42.68	224.257896878059\\
43.01	224.34092160812\\
43.34	224.402392937735\\
43.67	224.451113727982\\
44	224.514912653162\\
};
\addlegendentry{battery swapping}
\end{axis}
\end{tikzpicture}%
\vspace*{-0.3in}
\caption{Passenger arrival rate $\lambda$ (per minute) as a function of $s$.}
\label{lambda_s}
\end{subfigure}
\hfill
\begin{subfigure}[b]{0.32\linewidth}
\centering
%
%
\begin{tikzpicture}

\begin{axis}[%
width=1.694in,
height=1.03in,
at={(1.358in,0.0in)},
scale only axis,
xmin=11,
xmax=44,
xtick={11, 22, 33, 44},
xlabel style={font=\color{white!15!black}},
xlabel={Charging speed},
ymin=22.9,
ymax=25,
ylabel style={font=\color{white!15!black}},
ylabel={Passenger cost},
axis background/.style={fill=white},
legend style={at={(1,1)}, anchor=north east, legend cell align=left, align=left, font=\scriptsize, draw=white!12!black}
]
\addplot [color=black, line width=1.0pt]
  table[row sep=crcr]{%
11	24.6421516444848\\
11.33	24.5768162105026\\
11.66	24.5155960435017\\
11.99	24.4583212781313\\
12.32	24.4042847190869\\
12.65	24.3533495483242\\
12.98	24.3050903807265\\
13.31	24.2594763597033\\
13.64	24.216275887449\\
13.97	24.1752713849402\\
14.3	24.1362006185595\\
14.63	24.0991680338527\\
14.96	24.0638080912547\\
15.29	24.0300010856069\\
15.62	23.9979334182383\\
15.95	23.9670848202344\\
16.28	23.9375109562317\\
16.61	23.9093311202117\\
16.94	23.8821304395593\\
17.27	23.8561231885976\\
17.6	23.8311061135136\\
17.93	23.8070673491327\\
18.26	23.7838230981996\\
18.59	23.7617185658137\\
18.92	23.7400458919275\\
19.25	23.7195257684773\\
19.58	23.6993619471538\\
19.91	23.6801007153842\\
20.24	23.661323466733\\
20.57	23.6433493797923\\
20.9	23.6258001110005\\
21.23	23.6088119895644\\
21.56	23.5925212771487\\
21.89	23.5767148162702\\
22.22	23.5612280219277\\
22.55	23.5464300040694\\
22.88	23.5319550550851\\
23.21	23.5178327557052\\
23.54	23.5040382439421\\
23.87	23.4909079115157\\
24.2	23.478041903017\\
24.53	23.4655579685629\\
24.86	23.4532194730757\\
25.19	23.4413353895406\\
25.52	23.4295683500879\\
25.85	23.4182665402365\\
26.18	23.4075158207584\\
26.51	23.3966138657377\\
26.84	23.3860898371714\\
27.17	23.3759550554996\\
27.5	23.3660494627739\\
27.83	23.356281159659\\
28.16	23.3469035757226\\
28.49	23.3372973640312\\
28.82	23.3284638643947\\
29.15	23.3195755063749\\
29.48	23.310980631014\\
29.81	23.302497653648\\
30.14	23.294202369249\\
30.47	23.2860938015856\\
30.8	23.278215302862\\
31.13	23.2705383866615\\
31.46	23.262866439246\\
31.79	23.255454114851\\
32.12	23.2481724009185\\
32.45	23.2408519843593\\
32.78	23.2339374364821\\
33.11	23.2271938948138\\
33.44	23.2203871914702\\
33.77	23.2137577362372\\
34.1	23.2071859172939\\
34.43	23.2009324077187\\
34.76	23.194589132723\\
35.09	23.1886472045246\\
35.42	23.1826621130073\\
35.75	23.1767371839758\\
36.08	23.1708330662202\\
36.41	23.1652398992528\\
36.74	23.159691413267\\
37.07	23.154217084611\\
37.4	23.1488340600619\\
37.73	23.1434910106691\\
38.06	23.1383549155733\\
38.39	23.1332277788843\\
38.72	23.1281386100006\\
39.05	23.1232527834593\\
39.38	23.118384619585\\
39.71	23.1134728824622\\
40.04	23.1088196168928\\
40.37	23.1042434463832\\
40.7	23.0996009258351\\
41.03	23.0950882949732\\
41.36	23.0907090113119\\
41.69	23.0862968937417\\
42.02	23.0820880715866\\
42.35	23.0779622000181\\
42.68	23.0737029745101\\
43.01	23.0696265423322\\
43.34	23.0656788920336\\
43.67	23.0615864718553\\
44	23.0577085196095\\
};
\addlegendentry{plug-in charging}
\addplot [color=blue, dashed, line width=1.0pt]
  table[row sep=crcr]{%
11	24.7975297138852\\
11.33	24.6897042320029\\
11.66	24.5907153000364\\
11.99	24.5014524883518\\
12.32	24.4213184665082\\
12.65	24.3444970548426\\
12.98	24.2754408485797\\
13.31	24.2097787449064\\
13.64	24.1511907128465\\
13.97	24.0944714393822\\
14.3	24.0434925042122\\
14.63	23.9954353107277\\
14.96	23.9480999094457\\
15.29	23.9079606185252\\
15.62	23.8665656085711\\
15.95	23.828275641191\\
16.28	23.7931394330775\\
16.61	23.7593321747417\\
16.94	23.7280959992781\\
17.27	23.6981937168896\\
17.6	23.6696754288744\\
17.93	23.6428748632164\\
18.26	23.6168190397236\\
18.59	23.5922982202038\\
18.92	23.5691983211097\\
19.25	23.5480488790232\\
19.58	23.5257555383334\\
19.91	23.5062416060671\\
20.24	23.4867144486832\\
20.57	23.4685161705013\\
20.9	23.4509909139072\\
21.23	23.4342643070442\\
21.56	23.4178440930214\\
21.89	23.4023894646847\\
22.22	23.3877163915292\\
22.55	23.3731003351313\\
22.88	23.3594895812221\\
23.21	23.3470278285731\\
23.54	23.333305023033\\
23.87	23.3218538797625\\
24.2	23.3097892232321\\
24.53	23.2980474904729\\
24.86	23.2876119306799\\
25.19	23.2766772047005\\
25.52	23.2663192433579\\
25.85	23.2567366247741\\
26.18	23.2477216965271\\
26.51	23.2384009931367\\
26.84	23.2295585026051\\
27.17	23.2211387921911\\
27.5	23.2121649334568\\
27.83	23.2044817875216\\
28.16	23.1967553759883\\
28.49	23.1893180283478\\
28.82	23.1808368762962\\
29.15	23.1745232079358\\
29.48	23.1673432757066\\
29.81	23.1611601802954\\
30.14	23.1548058336633\\
30.47	23.148271869664\\
30.8	23.1419524393645\\
31.13	23.1366244508858\\
31.46	23.1306261872675\\
31.79	23.1250197510599\\
32.12	23.119410901886\\
32.45	23.1141585769279\\
32.78	23.1090256080801\\
33.11	23.104022169729\\
33.44	23.0989528771885\\
33.77	23.0945345555234\\
34.1	23.089340638357\\
34.43	23.0849447004854\\
34.76	23.0808853923868\\
35.09	23.0767449633457\\
35.42	23.072661825925\\
35.75	23.0683520186238\\
36.08	23.0646581226768\\
36.41	23.0609693457219\\
36.74	23.0569291242044\\
37.07	23.0534336621783\\
37.4	23.0497831365385\\
37.73	23.0462110584991\\
38.06	23.0428509358274\\
38.39	23.0394281168918\\
38.72	23.0360316214605\\
39.05	23.0331097638713\\
39.38	23.0302107185524\\
39.71	23.027295699205\\
40.04	23.0240936248193\\
40.37	23.0212580229773\\
40.7	23.0184951337772\\
41.03	23.015739341919\\
41.36	23.0130564511426\\
41.69	23.010260854375\\
42.02	23.0081980119138\\
42.35	23.0055118533324\\
42.68	23.0032015664485\\
43.01	23.0000692344768\\
43.34	22.9977505767207\\
43.67	22.995913171329\\
44	22.9935075396956\\
};
\addlegendentry{battery swapping}
\end{axis}
\end{tikzpicture}%
\vspace*{-0.3in}
\caption{Passenger travel cost $c$ (\$/trip) as a function of $s$.}
\label{c_s}
\end{subfigure}
\hfill
\begin{subfigure}[b]{0.32\linewidth}
\centering
%
%
\begin{tikzpicture}

\begin{axis}[%
width=1.694in,
height=1.03in,
at={(1.358in,0.0in)},
scale only axis,
xmin=11,
xmax=44,
xtick={11, 22, 33, 44},
xlabel style={font=\color{white!15!black}},
xlabel={Charging speed},
ymin=13.9,
ymax=14.5,
ylabel style={font=\color{white!15!black}},
ylabel={Ride fare},
axis background/.style={fill=white},
legend style={at={(1,1)}, anchor=north east, legend cell align=left, align=left, font=\scriptsize, draw=white!12!black}
]
\addplot [color=black, line width=1.0pt]
  table[row sep=crcr]{%
11	14.4735425923376\\
11.33	14.4511616838329\\
11.66	14.4301623380449\\
11.99	14.4104927338708\\
12.32	14.3919155442439\\
12.65	14.3743876382458\\
12.98	14.3577665201027\\
13.31	14.342044189959\\
13.64	14.3271438003091\\
13.97	14.3129917762743\\
14.3	14.2994999459422\\
14.63	14.2867053559781\\
14.96	14.2744829568922\\
15.29	14.2627927908013\\
15.62	14.2516997986938\\
15.95	14.2410251609765\\
16.28	14.2307887045172\\
16.61	14.2210319760309\\
16.94	14.2116120283071\\
17.27	14.2026033601545\\
17.6	14.1939361261688\\
17.93	14.1856062488614\\
18.26	14.1775505580154\\
18.59	14.1698885805913\\
18.92	14.1623755789615\\
19.25	14.1552613464373\\
19.58	14.1482696542828\\
19.91	14.1415906309324\\
20.24	14.1350789863686\\
20.57	14.1288453580039\\
20.9	14.1227588793416\\
21.23	14.1168663610533\\
21.56	14.111215882032\\
21.89	14.1057331266371\\
22.22	14.1003608427197\\
22.55	14.095227906267\\
22.88	14.0902067411635\\
23.21	14.0853076575195\\
23.54	14.0805224288939\\
23.87	14.0759674099649\\
24.2	14.0715043338634\\
24.53	14.0671737818289\\
24.86	14.0628936248162\\
25.19	14.0587716532516\\
25.52	14.0546894903338\\
25.85	14.050769417708\\
26.18	14.0470404331429\\
26.51	14.0432588852918\\
26.84	14.0396085644843\\
27.17	14.0360937204552\\
27.5	14.0326582143996\\
27.83	14.0292702983582\\
28.16	14.0260182296156\\
28.49	14.0226868858753\\
28.82	14.0196237174039\\
29.15	14.0165412798994\\
29.48	14.0135612323342\\
29.81	14.0106198943524\\
30.14	14.0077434795467\\
30.47	14.0049323080871\\
30.8	14.0022007136575\\
31.13	13.9995392023456\\
31.46	13.996879634516\\
31.79	13.9943102323253\\
32.12	13.99178609761\\
32.45	13.9892486570599\\
32.78	13.9868519448164\\
33.11	13.984514901032\\
33.44	13.9821556685393\\
33.77	13.9798581533286\\
34.1	13.9775807018979\\
34.43	13.9754137686563\\
34.76	13.9732157405782\\
35.09	13.9711568460316\\
35.42	13.9690830811672\\
35.75	13.9670303577411\\
36.08	13.9649848909889\\
36.41	13.9630471222288\\
36.74	13.9611251233673\\
37.07	13.9592286636836\\
37.4	13.9573640041576\\
37.73	13.9555133299706\\
38.06	13.9537344977133\\
38.39	13.9519586102177\\
38.72	13.9501961821632\\
39.05	13.9485042473269\\
39.38	13.9468183184017\\
39.71	13.9451177546146\\
40.04	13.943506333703\\
40.37	13.9419219976004\\
40.7	13.9403144459807\\
41.03	13.9387520764309\\
41.36	13.937236026309\\
41.69	13.9357085686952\\
42.02	13.9342515954792\\
42.35	13.932823547223\\
42.68	13.9313492974306\\
43.01	13.9299382436625\\
43.34	13.9285718504267\\
43.67	13.9271553402026\\
44	13.925813612318\\
};
\addlegendentry{plug-in charging}
\addplot [color=blue, dashed, line width=1.0pt]
  table[row sep=crcr]{%
11	14.5266340959873\\
11.33	14.4898111477076\\
11.66	14.4559255092254\\
11.99	14.425306913096\\
12.32	14.3977727880346\\
12.65	14.3713393814487\\
12.98	14.3475471744587\\
13.31	14.3249014757323\\
13.64	14.304676407657\\
13.97	14.2850815156997\\
14.3	14.2674577588248\\
14.63	14.2508359284586\\
14.96	14.2344539842567\\
15.29	14.2205562157773\\
15.62	14.2062202750566\\
15.95	14.1929549922224\\
16.28	14.1807790405797\\
16.61	14.1690613244799\\
16.94	14.1582322571452\\
17.27	14.1478641873371\\
17.6	14.1379750868965\\
17.93	14.1286804084222\\
18.26	14.1196437449872\\
18.59	14.1111381971733\\
18.92	14.1031255118312\\
19.25	14.0957892041387\\
19.58	14.0880555132381\\
19.91	14.0812862592441\\
20.24	14.0745124004935\\
20.57	14.0681997673248\\
20.9	14.0621210500056\\
21.23	14.0563183034249\\
21.56	14.0506225644644\\
21.89	14.0452621551078\\
22.22	14.0401729234053\\
22.55	14.0351036628418\\
22.88	14.0303827280553\\
23.21	14.0260607865078\\
23.54	14.0213024235026\\
23.87	14.0173317873261\\
24.2	14.0131477566247\\
24.53	14.0090765657088\\
24.86	14.005458249448\\
25.19	14.0016676476416\\
25.52	13.9980767163835\\
25.85	13.9947542724002\\
26.18	13.9916294555796\\
26.51	13.9883987651533\\
26.84	13.9853343191965\\
27.17	13.9824159147124\\
27.5	13.9793060181627\\
27.83	13.9766437078212\\
28.16	13.973966042538\\
28.49	13.971389125144\\
28.82	13.9684503531727\\
29.15	13.9662632866629\\
29.48	13.9637759284205\\
29.81	13.961633641637\\
30.14	13.9594325175792\\
30.47	13.9571693283421\\
30.8	13.9549803529636\\
31.13	13.9531348438946\\
31.46	13.9510576740805\\
31.79	13.9491159277442\\
32.12	13.9471736837328\\
32.45	13.945354714148\\
32.78	13.9435775899782\\
33.11	13.941845177324\\
33.44	13.9400898520159\\
33.77	13.9385601944891\\
34.1	13.9367622863105\\
34.43	13.935240276618\\
34.76	13.9338352666174\\
35.09	13.9324019284532\\
35.42	13.9309886129882\\
35.75	13.929496830464\\
36.08	13.9282184660484\\
36.41	13.9269417365052\\
36.74	13.9255435166388\\
37.07	13.9243337626432\\
37.4	13.9230707408676\\
37.73	13.9218343911871\\
38.06	13.9206717752173\\
38.39	13.9194875345442\\
38.72	13.9183124613694\\
39.05	13.9173016384253\\
39.38	13.9162984625194\\
39.71	13.9152900751868\\
40.04	13.9141820886849\\
40.37	13.9132012545629\\
40.7	13.9122455048175\\
41.03	13.9112922633746\\
41.36	13.9103641602069\\
41.69	13.9093975462601\\
42.02	13.9086837404379\\
42.35	13.9077547016782\\
42.68	13.9069556703321\\
43.01	13.9058725589254\\
43.34	13.9050705497006\\
43.67	13.9044352521514\\
44	13.9036031197354\\
};
\addlegendentry{battery swapping}
\end{axis}
\end{tikzpicture}%
\vspace*{-0.3in}
\caption{Per-trip ride fare $p_f$ (\$/trip) as a function of $s$.}
\label{pf_s}
\end{subfigure}
\begin{subfigure}[b]{0.32\linewidth}
\centering
%
%
\begin{tikzpicture}

\begin{axis}[%
width=1.694in,
height=1.03in,
at={(1.358in,0.0in)},
scale only axis,
xmin=11,
xmax=44,
xtick={11, 22, 33, 44},
xlabel style={font=\color{white!15!black}},
xlabel={Charging speed},
ymin=3.5,
ymax=4,
ylabel style={font=\color{white!15!black}},
ylabel={Waiting time},
axis background/.style={fill=white},
legend style={at={(1,1)}, anchor=north east, legend cell align=left, align=left, font=\scriptsize, draw=white!12!black}
]
\addplot [color=black, line width=1.0pt]
  table[row sep=crcr]{%
11	3.94132133804156\\
11.33	3.92467229715878\\
11.66	3.90908283157241\\
11.99	3.89450718769787\\
12.32	3.88076324606316\\
12.65	3.86781469382884\\
12.98	3.85555188396269\\
13.31	3.84396595726524\\
13.64	3.83299693299997\\
13.97	3.82258899560696\\
14.3	3.81267467930901\\
14.63	3.80328010770335\\
14.96	3.79431206758233\\
15.29	3.78573964914948\\
15.62	3.77760993005599\\
15.95	3.76979056560385\\
16.28	3.76229544640097\\
16.61	3.75515470704683\\
16.94	3.74826295009773\\
17.27	3.74167435210971\\
17.6	3.73533720439725\\
17.93	3.72924848847723\\
18.26	3.7233614496838\\
18.59	3.71776356016373\\
18.92	3.71227531510308\\
19.25	3.7070792333488\\
19.58	3.70197375692676\\
19.91	3.69709693195807\\
20.24	3.69234282184668\\
20.57	3.68779225650711\\
20.9	3.68334931459649\\
21.23	3.67904869322136\\
21.56	3.67492457175066\\
21.89	3.67092313551673\\
22.22	3.66700278263876\\
22.55	3.6632566270552\\
22.88	3.65959236973703\\
23.21	3.65601747991695\\
23.54	3.6525255097086\\
23.87	3.64920174478713\\
24.2	3.64594479424555\\
24.53	3.64278456850155\\
24.86	3.63966118149591\\
25.19	3.63665261096471\\
25.52	3.63367397664885\\
25.85	3.63081283818934\\
26.18	3.62809123550987\\
26.51	3.62533138776971\\
26.84	3.62266715995626\\
27.17	3.62010129265286\\
27.5	3.61759350712183\\
27.83	3.6151204888763\\
28.16	3.61274625818101\\
28.49	3.61031413882014\\
28.82	3.60807757635303\\
29.15	3.60582721956415\\
29.48	3.60365092972085\\
29.81	3.60150300747893\\
30.14	3.59940267042724\\
30.47	3.59734941608469\\
30.8	3.59535449193974\\
31.13	3.59341053655655\\
31.46	3.59146775377132\\
31.79	3.58959065214174\\
32.12	3.58774662918935\\
32.45	3.58589276251914\\
32.78	3.58414166343632\\
33.11	3.58243371852009\\
33.44	3.58070989260888\\
33.77	3.57903084608862\\
34.1	3.5773663625566\\
34.43	3.57578241824127\\
34.76	3.57417573338944\\
35.09	3.57267068158643\\
35.42	3.57115466350392\\
35.75	3.56965380861812\\
36.08	3.56815820745396\\
36.41	3.5667413864434\\
36.74	3.56533577127895\\
37.07	3.56394900035945\\
37.4	3.56258529298618\\
37.73	3.56123165918545\\
38.06	3.55993039451939\\
38.39	3.55863146072348\\
38.72	3.55734202629358\\
39.05	3.55610408377224\\
39.38	3.55487065937338\\
39.71	3.55362601854561\\
40.04	3.55244700898831\\
40.37	3.55128738324915\\
40.7	3.55011103870324\\
41.03	3.54896752656681\\
41.36	3.54785774612514\\
41.69	3.54673966087074\\
42.02	3.54567305275481\\
42.35	3.54462738480428\\
42.68	3.5435479368525\\
43.01	3.54251484444563\\
43.34	3.54151435721196\\
43.67	3.54047718281113\\
44	3.539494150113\\
};
\addlegendentry{plug-in charging}
\addplot [color=blue, dashed, line width=1.0pt]
  table[row sep=crcr]{%
11	3.98096729375888\\
11.33	3.95344693189742\\
11.66	3.92821309721357\\
11.99	3.90548278110691\\
12.32	3.88509522421455\\
12.65	3.86556498968757\\
12.98	3.84802080392286\\
13.31	3.83134777874966\\
13.64	3.81647841286415\\
13.97	3.80208911770642\\
14.3	3.78916075402611\\
14.63	3.7769765047555\\
14.96	3.76497904077093\\
15.29	3.75480790804184\\
15.62	3.74431989671103\\
15.95	3.73462040657695\\
16.28	3.72572108236347\\
16.61	3.71715924428752\\
16.94	3.7092495124546\\
17.27	3.70167811222963\\
17.6	3.69445749689065\\
17.93	3.68767226930009\\
18.26	3.68107569563425\\
18.59	3.67486822598081\\
18.92	3.66902046871258\\
19.25	3.66366654065292\\
19.58	3.65802326554082\\
19.91	3.65308346776085\\
20.24	3.64814032875569\\
20.57	3.6435334896033\\
20.9	3.63909684647348\\
21.23	3.6348627921005\\
21.56	3.63070601882054\\
21.89	3.62679353084376\\
22.22	3.62307886361392\\
22.55	3.61937855515096\\
22.88	3.61593288882433\\
23.21	3.61277792328113\\
23.54	3.6093033331513\\
23.87	3.60640391179706\\
24.2	3.6033494056618\\
24.53	3.60037632742794\\
24.86	3.5977339849736\\
25.19	3.59496494459646\\
25.52	3.59234206471873\\
25.85	3.589915640455\\
26.18	3.58763265153003\\
26.51	3.58527218138892\\
26.84	3.58303262922811\\
27.17	3.58090034010803\\
27.5	3.5786274865481\\
27.83	3.57668142624047\\
28.16	3.57472454784895\\
28.49	3.57284066015653\\
28.82	3.57069245082306\\
29.15	3.56909299274145\\
29.48	3.56727416561476\\
29.81	3.56570796072028\\
30.14	3.56409818452872\\
30.47	3.5624428454736\\
30.8	3.5608418939538\\
31.13	3.559492095733\\
31.46	3.55797229193294\\
31.79	3.55655186950221\\
32.12	3.55513070471053\\
32.45	3.55379994681391\\
32.78	3.55249923182244\\
33.11	3.55123139240503\\
33.44	3.5499469089816\\
33.77	3.54882727171874\\
34.1	3.54751098916533\\
34.43	3.5463970635145\\
34.76	3.54536826580207\\
35.09	3.54431900577227\\
35.42	3.54328419106078\\
35.75	3.54219193339527\\
36.08	3.54125568086373\\
36.41	3.54032077876619\\
36.74	3.53929674711846\\
37.07	3.53841081377327\\
37.4	3.53748542467865\\
37.73	3.53658010360932\\
38.06	3.53572835682564\\
38.39	3.53486069083241\\
38.72	3.53399967445392\\
39.05	3.53325896335116\\
39.38	3.53252413024534\\
39.71	3.53178512558847\\
40.04	3.53097346361799\\
40.37	3.53025456140093\\
40.7	3.52955411975183\\
41.03	3.52885545680018\\
41.36	3.52817530656422\\
41.69	3.52746639849413\\
42.02	3.5269435160759\\
42.35	3.52626246188148\\
42.68	3.5256767039211\\
43.01	3.52488243238426\\
43.34	3.52429458411632\\
43.67	3.52382865084405\\
44	3.52321876742642\\
};
\addlegendentry{battery swapping}
\end{axis}
\end{tikzpicture}%
\vspace*{-0.3in}
\caption{Passenger waiting time $w^c$ (minute) as a function of $s$.}
\label{wc_s}
\end{subfigure}
\hfill
\begin{subfigure}[b]{0.32\linewidth}
\centering
%
%
\begin{tikzpicture}

\begin{axis}[%
width=1.694in,
height=1.03in,
at={(1.358in,0.0in)},
scale only axis,
xmin=11,
xmax=44,
xtick={11, 22, 33, 44},
xlabel style={font=\color{white!15!black}},
xlabel={Charging speed},
ymin=17.85,
ymax=18.35,
ylabel style={font=\color{white!15!black}},
ylabel={Vehicle idle time},
axis background/.style={fill=white},
legend style={at={(1,0)}, anchor=south east, legend cell align=left, align=left, font=\scriptsize, draw=white!12!black}
]
\addplot [color=black, line width=1.0pt]
  table[row sep=crcr]{%
11	17.8948324839625\\
11.33	17.9005335728395\\
11.66	17.9066248864233\\
11.99	17.9129888525272\\
12.32	17.9195913247345\\
12.65	17.9263524057146\\
12.98	17.9332467215519\\
13.31	17.9402042076551\\
13.64	17.9471948591544\\
13.97	17.954193166155\\
14.3	17.9611959089742\\
14.63	17.9681359871858\\
14.96	17.9750399776336\\
15.29	17.9818978581161\\
15.62	17.9886357260336\\
15.95	17.9953345149137\\
16.28	18.0019578890874\\
16.61	18.0084528189491\\
16.94	18.0148941975112\\
17.27	18.0212118565215\\
17.6	18.0274373129509\\
17.93	18.0335563394146\\
18.26	18.0396026174399\\
18.59	18.0454700449393\\
18.92	18.0513364033615\\
19.25	18.0569943398917\\
19.58	18.0626515080195\\
19.91	18.0681482990479\\
20.24	18.0735939585127\\
20.57	18.0788870325958\\
20.9	18.0841321830038\\
21.23	18.0892806130543\\
21.56	18.0942863083601\\
21.89	18.0992062017011\\
22.22	18.1040862406489\\
22.55	18.1088075819979\\
22.88	18.1134785343884\\
23.21	18.1180861410427\\
23.54	18.1226364904195\\
23.87	18.1270122325991\\
24.2	18.1313437125139\\
24.53	18.1355871713093\\
24.86	18.1398204116107\\
25.19	18.1439370279244\\
25.52	18.1480460416138\\
25.85	18.1520290734615\\
26.18	18.1558484408884\\
26.51	18.1597518367801\\
26.84	18.1635498379197\\
27.17	18.1672365128287\\
27.5	18.1708652189818\\
27.83	18.1744687228525\\
28.16	18.1779529728043\\
28.49	18.1815462015598\\
28.82	18.1848726624784\\
29.15	18.188239610268\\
29.48	18.1915176422847\\
29.81	18.1947718358447\\
30.14	18.1979718728346\\
30.47	18.2011195028703\\
30.8	18.2041938496415\\
31.13	18.2072061950081\\
31.46	18.2102331689372\\
31.79	18.2131732936864\\
32.12	18.2160758992683\\
32.45	18.2190086875254\\
32.78	18.2217922510672\\
33.11	18.2245209022031\\
33.44	18.2272863373004\\
33.77	18.2299928713588\\
34.1	18.2326879278857\\
34.43	18.2352641029849\\
34.76	18.2378881014761\\
35.09	18.240356201246\\
35.42	18.2428522438583\\
35.75	18.2453335461178\\
36.08	18.2478158041256\\
36.41	18.2501759551001\\
36.74	18.2525268186903\\
37.07	18.2548538950368\\
37.4	18.2571508072117\\
37.73	18.2594390381115\\
38.06	18.2616466292781\\
38.39	18.2638569438304\\
38.72	18.2660592407708\\
39.05	18.2681805205771\\
39.38	18.2703002096006\\
39.71	18.2724472859307\\
40.04	18.2744861174569\\
40.37	18.2764985268259\\
40.7	18.2785450372632\\
41.03	18.2805408079549\\
41.36	18.2824835855962\\
41.69	18.2844461255384\\
42.02	18.2863236727399\\
42.35	18.2881698144894\\
42.68	18.2900803935312\\
43.01	18.2919133892541\\
43.34	18.2936932035933\\
43.67	18.2955428246206\\
44	18.297302007332\\
};
\addlegendentry{plug-in charging}
\addplot [color=blue, dashed, line width=1.0pt]
  table[row sep=crcr]{%
11	17.8845208428724\\
11.33	17.8911941137228\\
11.66	17.8992519736871\\
11.99	17.9081354041697\\
12.32	17.917444107415\\
12.65	17.9275800799855\\
12.98	17.9377164669464\\
13.31	17.9482781823004\\
13.64	17.9584678808386\\
13.97	17.9690346013224\\
14.3	17.9791279868008\\
14.63	17.9891718648897\\
14.96	17.9995628436372\\
15.29	18.0087688902462\\
15.62	18.018654701789\\
15.95	18.0281493102014\\
16.28	18.0371624952355\\
16.61	18.0461103648289\\
16.94	18.05461737044\\
17.27	18.0629808744951\\
17.6	18.0711600967365\\
17.93	18.0790264123382\\
18.26	18.0868451598172\\
18.59	18.0943541016451\\
18.92	18.1015664419354\\
19.25	18.1082874670664\\
19.58	18.1154926930471\\
19.91	18.12190448086\\
20.24	18.1284180166468\\
20.57	18.1345772987334\\
20.9	18.1405911323474\\
21.23	18.1464008267994\\
21.56	18.1521773336054\\
21.89	18.1576796460341\\
22.22	18.1629616402329\\
22.55	18.1682797625939\\
22.88	18.1732810611941\\
23.21	18.177904576573\\
23.54	18.1830468698548\\
23.87	18.1873760453013\\
24.2	18.1919720691589\\
24.53	18.1964854190396\\
24.86	18.2005273621412\\
25.19	18.2047968221016\\
25.52	18.2088694199272\\
25.85	18.2126614532011\\
26.18	18.2162544109602\\
26.51	18.2199925961396\\
26.84	18.2235624740519\\
27.17	18.2269793302426\\
27.5	18.2306443966518\\
27.83	18.2338007210336\\
28.16	18.2369894714611\\
28.49	18.2400764361849\\
28.82	18.2436141738356\\
29.15	18.246263076927\\
29.48	18.2492877322732\\
29.81	18.2519026480865\\
30.14	18.2546028146087\\
30.47	18.2573914052109\\
30.8	18.2600991544407\\
31.13	18.2623907465618\\
31.46	18.2649819542621\\
31.79	18.2674117608126\\
32.12	18.2698525910084\\
32.45	18.2721453761881\\
32.78	18.2743954241642\\
33.11	18.2765951390827\\
33.44	18.2788303992185\\
33.77	18.2807854329529\\
34.1	18.2830916127888\\
34.43	18.2850479959264\\
34.76	18.2868611025137\\
35.09	18.2887141405274\\
35.42	18.2905469476955\\
35.75	18.2924864815513\\
36.08	18.2941538562937\\
36.41	18.2958221555645\\
36.74	18.297654337007\\
37.07	18.2992428943241\\
37.4	18.3009071520388\\
37.73	18.3025373571908\\
38.06	18.3040755739064\\
38.39	18.3056460059818\\
38.72	18.307207835099\\
39.05	18.3085541548069\\
39.38	18.3098913381918\\
39.71	18.3112395304042\\
40.04	18.3127220256687\\
40.37	18.314038636113\\
40.7	18.3153233835383\\
41.03	18.316607176869\\
41.36	18.317858739844\\
41.69	18.3191669540393\\
42.02	18.3201314604738\\
42.35	18.3213910230349\\
42.68	18.3224759906359\\
43.01	18.3239503390794\\
43.34	18.3250424581806\\
43.67	18.3259099817705\\
44	18.3270457528851\\
};
\addlegendentry{battery swapping}
\end{axis}
\end{tikzpicture}%
\vspace*{-0.3in}
\caption{Vehicle idle time $w^v$ (minute) as a function of $s$.}
\label{wv_s}
\end{subfigure}
\hfill
\begin{subfigure}[b]{0.32\linewidth}
\centering
%
%
\begin{tikzpicture}

\begin{axis}[%
width=1.694in,
height=1.03in,
at={(1.358in,0.0in)},
scale only axis,
xmin=11,
xmax=44,
xtick={11, 22, 33, 44},
xlabel style={font=\color{white!15!black}},
xlabel={Charging speed},
ymin=7000,
ymax=9500,
ytick={7000,8000,9000},
yticklabels={{7K},{8K},{9K}},
ylabel style={font=\color{white!15!black}},
ylabel={Number of vehicles},
axis background/.style={fill=white},
legend style={at={(1,0)}, anchor=south east, legend cell align=left, align=left, font=\scriptsize, draw=white!12!black}
]
\addplot [color=black, line width=1.0pt]
  table[row sep=crcr]{%
11	8704.66120882723\\
11.33	8729.08877190581\\
11.66	8751.44409062846\\
11.99	8771.98958401757\\
12.32	8791.01736200589\\
12.65	8808.74894766118\\
12.98	8825.34090926643\\
13.31	8840.70896968939\\
13.64	8855.16501561883\\
13.97	8868.68306270079\\
14.3	8881.32639821377\\
14.63	8893.1297829057\\
14.96	8904.23275720905\\
15.29	8914.84506201011\\
15.62	8924.64119544782\\
15.95	8934.03853955909\\
16.28	8942.92058770026\\
16.61	8951.30573193522\\
16.94	8959.38126679178\\
17.27	8966.99111333559\\
17.6	8974.18629731738\\
17.93	8980.98867789415\\
18.26	8987.64948365243\\
18.59	8993.71396824659\\
18.92	8999.91099758753\\
19.25	9005.37094414443\\
19.58	9011.0068248297\\
19.91	9016.17370401086\\
20.24	9021.23287090625\\
20.57	9026.01769725111\\
20.9	9030.75094485214\\
21.23	9035.12849125114\\
21.56	9039.36608526506\\
21.89	9043.42557400296\\
22.22	9047.44296392645\\
22.55	9051.10727195279\\
22.88	9054.81151137556\\
23.21	9058.29484596918\\
23.54	9061.83215065497\\
23.87	9065.08774670545\\
24.2	9068.20356426195\\
24.53	9071.25914687522\\
24.86	9074.27006684782\\
25.19	9077.16238463256\\
25.52	9080.11993064847\\
25.85	9082.80151639479\\
26.18	9085.22778710814\\
26.51	9087.88333528467\\
26.84	9090.37058290233\\
27.17	9092.73475194285\\
27.5	9095.02533560247\\
27.83	9097.31653282441\\
28.16	9099.43913773979\\
28.49	9101.78289860053\\
28.82	9103.66236776092\\
29.15	9105.75691232037\\
29.48	9107.67453934453\\
29.81	9109.61085242337\\
30.14	9111.48741305898\\
30.47	9113.30814477075\\
30.8	9114.99905535632\\
31.13	9116.69741162063\\
31.46	9118.44975944257\\
31.79	9120.05764294581\\
32.12	9121.65748708092\\
32.45	9123.38298309909\\
32.78	9124.8242544706\\
33.11	9126.20925266816\\
33.44	9127.74164662941\\
33.77	9129.19567597433\\
34.1	9130.69246608185\\
34.43	9131.99815363952\\
34.76	9133.40743915943\\
35.09	9134.57579803975\\
35.42	9135.80792419221\\
35.75	9137.12395443565\\
36.08	9138.42957899566\\
36.41	9139.57247883673\\
36.74	9140.74030624485\\
37.07	9141.90172765897\\
37.4	9142.99635305586\\
37.73	9144.13940068552\\
38.06	9145.18300248489\\
38.39	9146.23424521424\\
38.72	9147.33307785945\\
39.05	9148.30384836847\\
39.38	9149.29948320675\\
39.71	9150.37085822922\\
40.04	9151.29983586592\\
40.37	9152.21820387013\\
40.7	9153.19397413543\\
41.03	9154.15641225834\\
41.36	9155.02125724922\\
41.69	9155.96021239226\\
42.02	9156.76208295633\\
42.35	9157.56252548386\\
42.68	9158.50293960473\\
43.01	9159.32859053012\\
43.34	9160.03648380523\\
43.67	9160.98456888766\\
44	9161.60751942919\\
};
\addlegendentry{plug-in charging}
\addplot [color=blue, dashed, line width=1.0pt]
  table[row sep=crcr]{%
11	7310.77657759538\\
11.33	7395.61059642771\\
11.66	7474.39774377534\\
11.99	7546.24144166679\\
12.32	7611.42274030619\\
12.65	7674.44387772945\\
12.98	7731.64660601229\\
13.31	7786.46292834336\\
13.64	7835.8281470417\\
13.97	7883.92620324345\\
14.3	7927.52975336786\\
14.63	7968.91057272675\\
14.96	8009.868562086\\
15.29	8044.97034348827\\
15.62	8081.25311623832\\
15.95	8115.04553581708\\
16.28	8146.27735267437\\
16.61	8176.47203638859\\
16.94	8204.55903226431\\
17.27	8231.57668422406\\
17.6	8257.47171506088\\
17.93	8281.94726882154\\
18.26	8305.83062580418\\
18.59	8328.43327927246\\
18.92	8349.84592354974\\
19.25	8369.60110581515\\
19.58	8390.3751415492\\
19.91	8408.76092949776\\
20.24	8427.17839545997\\
20.57	8444.45619052426\\
20.9	8461.16197319165\\
21.23	8477.18150106419\\
21.56	8492.94465482312\\
21.89	8507.87119639468\\
22.22	8522.11933252448\\
22.55	8536.32345158469\\
22.88	8549.64797066361\\
23.21	8561.96417596855\\
23.54	8575.39463335573\\
23.87	8586.82311081977\\
24.2	8598.79252758364\\
24.53	8610.47760891551\\
24.86	8621.00651966086\\
25.19	8631.97606990257\\
25.52	8642.43258995084\\
25.85	8652.19906748879\\
26.18	8661.45938409597\\
26.51	8670.98703041345\\
26.84	8680.0867095391\\
27.17	8688.8061377543\\
27.5	8698.01354705974\\
27.83	8706.06945087614\\
28.16	8714.15809944474\\
28.49	8721.98460866293\\
28.82	8730.73394322353\\
29.15	8737.5496440413\\
29.48	8745.1300164425\\
29.81	8751.81667267612\\
30.14	8758.64983352424\\
30.47	8765.63659852247\\
30.8	8772.42582793122\\
31.13	8778.32556215775\\
31.46	8784.81594562819\\
31.79	8790.95044342687\\
32.12	8797.08058092764\\
32.45	8802.88617494611\\
32.78	8808.57906811446\\
33.11	8814.14953105382\\
33.44	8819.77174131621\\
33.77	8824.80754095592\\
34.1	8830.52723272093\\
34.43	8835.52915183721\\
34.76	8840.22485296477\\
35.09	8844.98574932902\\
35.42	8849.68891468539\\
35.75	8854.58693230636\\
36.08	8858.93015897503\\
36.41	8863.26200873398\\
36.74	8867.89970570119\\
37.07	8872.04614113819\\
37.4	8876.32402746976\\
37.73	8880.52487281578\\
38.06	8884.53106823364\\
38.39	8888.58633641942\\
38.72	8892.61156092784\\
39.05	8896.20886885237\\
39.38	8899.77916422778\\
39.71	8903.35749650798\\
40.04	8907.18363093005\\
40.37	8910.67865794518\\
40.7	8914.10273573255\\
41.03	8917.51422053566\\
41.36	8920.85473182928\\
41.69	8924.28908108836\\
42.02	8927.06842483453\\
42.35	8930.3930871394\\
42.68	8933.37941021202\\
43.01	8937.08593090934\\
43.34	8940.06711132336\\
43.67	8942.61780117257\\
44	8945.66377696085\\
};
\addlegendentry{battery swapping}
\end{axis}
\end{tikzpicture}%
\vspace*{-0.3in}
\caption{Total number of vehicles $N$ as a function of $s$.}
\label{N_s}
\end{subfigure}
\begin{subfigure}[b]{0.32\linewidth}
\centering
%
%
\begin{tikzpicture}

\begin{axis}[%
width=1.694in,
height=1.03in,
at={(1.358in,0.0in)},
scale only axis,
xmin=11,
xmax=44,
xtick={11, 22, 33, 44},
xlabel style={font=\color{white!15!black}},
xlabel={Charging speed},
ymin=680,
ymax=900,
ytick={700,800,900},
yticklabels={{0.7K},{0.8K},{0.9K}},
ylabel style={font=\color{white!15!black}},
ylabel={Charging arrivals},
axis background/.style={fill=white},
legend style={at={(1,0)}, anchor=south east, legend cell align=left, align=left, font=\scriptsize, draw=white!12!black}
]
\addplot [color=black, line width=1.0pt]
  table[row sep=crcr]{%
11	710.129523752727\\
11.33	715.598945455357\\
11.66	720.748124977386\\
11.99	725.609967613038\\
12.32	730.214508789056\\
12.65	734.590436103833\\
12.98	738.756324481856\\
13.31	742.70941504842\\
13.64	746.485809177371\\
13.97	750.085990061069\\
14.3	753.518833198691\\
14.63	756.794030186313\\
14.96	759.926847821543\\
15.29	762.942928459771\\
15.62	765.813057526469\\
15.95	768.580306337171\\
16.28	771.236147460897\\
16.61	773.787327454254\\
16.94	776.252707184985\\
17.27	778.620475016619\\
17.6	780.896461567671\\
17.93	783.085488434129\\
18.26	785.216140157802\\
18.59	787.243387959396\\
18.92	789.243519627524\\
19.25	791.12819677566\\
19.58	792.994007769592\\
19.91	794.775501994602\\
20.24	796.510702874459\\
20.57	798.186096146253\\
20.9	799.824255141442\\
21.23	801.392729194631\\
21.56	802.920082611229\\
21.89	804.400241978337\\
22.22	805.84838445665\\
22.55	807.233646999525\\
22.88	808.599092908375\\
23.21	809.914923796997\\
23.54	811.214941803379\\
23.87	812.464329121908\\
24.2	813.675441037587\\
24.53	814.86051731991\\
24.86	816.018943866789\\
25.19	817.146915888204\\
25.52	818.262995891142\\
25.85	819.330876243454\\
26.18	820.354921860267\\
26.51	821.385614391472\\
26.84	822.382045750713\\
27.17	823.350569122051\\
27.5	824.295579793774\\
27.83	825.225578287959\\
28.16	826.123597688577\\
28.49	827.030048568894\\
28.82	827.873827932542\\
29.15	828.728967759003\\
29.48	829.551411266196\\
29.81	830.363754170996\\
30.14	831.15757117754\\
30.47	831.93373450984\\
30.8	832.68374812265\\
31.13	833.425254752634\\
31.46	834.161578274715\\
31.79	834.871344178668\\
32.12	835.570558289502\\
32.45	836.274344654234\\
32.78	836.936515782439\\
33.11	837.583355684268\\
33.44	838.238208412241\\
33.77	838.875543465979\\
34.1	839.509686958684\\
34.43	840.114177895841\\
34.76	840.72165947741\\
35.09	841.294602344009\\
35.42	841.86646764626\\
35.75	842.441951640222\\
36.08	843.008038458207\\
36.41	843.549291616536\\
36.74	844.086660939944\\
37.07	844.616650258439\\
37.4	845.131550566163\\
37.73	845.646285838719\\
38.06	846.14380521742\\
38.39	846.635464493379\\
38.72	847.127379869318\\
39.05	847.598554365972\\
39.38	848.066872596572\\
39.71	848.538430920656\\
40.04	848.988625620497\\
40.37	849.432549809433\\
40.7	849.877277894551\\
41.03	850.316598668951\\
41.36	850.739002607007\\
41.69	851.16553296742\\
42.02	851.571100434694\\
42.35	851.972543353667\\
42.68	852.385983445199\\
43.01	852.78149686891\\
43.34	853.157738642836\\
43.67	853.5602866354\\
44	853.91604057096\\
};
\addlegendentry{plug-in charging}
\addplot [color=blue, dashed, line width=1.0pt]
  table[row sep=crcr]{%
11	688.889282897549\\
11.33	697.965948200354\\
11.66	706.423491705887\\
11.99	714.160418932422\\
12.32	721.200943613032\\
12.65	728.023789196248\\
12.98	734.233411128327\\
13.31	740.196268645701\\
13.64	745.579499620535\\
13.97	750.833159028474\\
14.3	755.606419397032\\
14.63	760.143871645113\\
14.96	764.640426395349\\
15.29	768.503708973611\\
15.62	772.499232290012\\
15.95	776.226389562312\\
16.28	779.676555157831\\
16.61	783.01569490843\\
16.94	786.126123804091\\
17.27	789.121108048334\\
17.6	791.994551502325\\
17.93	794.713506280082\\
18.26	797.368611972188\\
18.59	799.883978921163\\
18.92	802.269244036238\\
19.25	804.472616300878\\
19.58	806.789129135439\\
19.91	808.842739921251\\
20.24	810.900404041025\\
20.57	812.832620227365\\
20.9	814.701991378704\\
21.23	816.495825821781\\
21.56	818.26158970808\\
21.89	819.934977169515\\
22.22	821.533395538942\\
22.55	823.127157933336\\
22.88	824.623530175377\\
23.21	826.008077846401\\
23.54	827.516497215317\\
23.87	828.802613789602\\
24.2	830.148937300294\\
24.53	831.463710966594\\
24.86	832.649945407407\\
25.19	833.885210294582\\
25.52	835.063427353806\\
25.85	836.164820608041\\
26.18	837.209794966346\\
26.51	838.284537990478\\
26.84	839.311556002619\\
27.17	840.29615441941\\
27.5	841.335150924056\\
27.83	842.245659929055\\
28.16	843.159772604885\\
28.49	844.044571852133\\
28.82	845.032366809937\\
29.15	845.804189924441\\
29.48	846.661348244902\\
29.81	847.418602257093\\
30.14	848.192140648827\\
30.47	848.982787108162\\
30.8	849.751291047513\\
31.13	850.420256468024\\
31.46	851.155200910254\\
31.79	851.85026962289\\
32.12	852.544780493237\\
32.45	853.202907261759\\
32.78	853.848337521104\\
33.11	854.479997780171\\
33.44	855.117376547394\\
33.77	855.68900710178\\
34.1	856.337172148254\\
34.43	856.904865990458\\
34.76	857.438168770182\\
35.09	857.978699225612\\
35.42	858.51268530659\\
35.75	859.068439214789\\
36.08	859.561933522001\\
36.41	860.054073371492\\
36.74	860.580439195661\\
37.07	861.051641885802\\
37.4	861.537511720332\\
37.73	862.014691294918\\
38.06	862.469954131156\\
38.39	862.930654473868\\
38.72	863.387919279609\\
39.05	863.797090815197\\
39.38	864.203173631589\\
39.71	864.61008322968\\
40.04	865.044753830337\\
40.37	865.442181910855\\
40.7	865.831586032308\\
41.03	866.219512229271\\
41.36	866.599408328404\\
41.69	866.989777748165\\
42.02	867.30651027082\\
42.35	867.684456872147\\
42.68	868.024320117881\\
43.01	868.445060332903\\
43.34	868.784242387595\\
43.67	869.074978660307\\
44	869.42133091735\\
};
\addlegendentry{battery swapping}
\end{axis}
\end{tikzpicture}%
\vspace*{-0.3in}
\caption{Charging arrival rate $\gamma$ (per hour) as a function of $s$.}
\label{gamma_s}
\end{subfigure}
\hfill
\begin{subfigure}[b]{0.32\linewidth}
\centering
%
%
\begin{tikzpicture}

\begin{axis}[%
width=1.694in,
height=1.03in,
at={(1.358in,0.0in)},
scale only axis,
xmin=11,
xmax=44,
xtick={11, 22, 33, 44},
xlabel style={font=\color{white!15!black}},
xlabel={Charging speed},
ymin=6900,
ymax=8600,
ytick={6500,7500,8500},
yticklabels={{6.5K},{7.5K},{8.5K}},
ylabel style={font=\color{white!15!black}},
ylabel={Number of vehicles},
axis background/.style={fill=white},
legend style={at={(1,0)}, anchor=south east, legend cell align=left, align=left, font=\scriptsize, draw=white!12!black}
]
\addplot [color=black, line width=1.0pt]
  table[row sep=crcr]{%
11	7007.16658859273\\
11.33	7062.46485548149\\
11.66	7114.71297523524\\
11.99	7163.97755328538\\
12.32	7210.80079097865\\
12.65	7255.2455102666\\
12.98	7297.63479608885\\
13.31	7337.95302825465\\
13.64	7376.36644751506\\
13.97	7413.03437250257\\
14.3	7448.16274267235\\
14.63	7481.6305698758\\
14.96	7513.74417366093\\
15.29	7544.59254001878\\
15.62	7573.98558493244\\
15.95	7602.3834164098\\
16.28	7629.7210167266\\
16.61	7655.87360054857\\
16.94	7681.21409511713\\
17.27	7705.53200545951\\
17.6	7729.00697956871\\
17.93	7751.64081257924\\
18.26	7773.59876841003\\
18.59	7794.54604652418\\
18.92	7815.14712040987\\
19.25	7834.71031858486\\
19.58	7853.98852697373\\
19.91	7872.45504622516\\
20.24	7890.50582952502\\
20.57	7907.82926973056\\
20.9	7924.78585760896\\
21.23	7941.24012212527\\
21.56	7957.05639011465\\
21.89	7972.43745081302\\
22.22	7987.54079921593\\
22.55	8002.00393383476\\
22.88	8016.18069742935\\
23.21	8030.04018616426\\
23.54	8043.60514490907\\
23.87	8056.54174168051\\
24.2	8069.24171937826\\
24.53	8081.58695488326\\
24.86	8093.81006994064\\
25.19	8105.60391947995\\
25.52	8117.30088653877\\
25.85	8128.55455661084\\
26.18	8139.27643476412\\
26.51	8150.16602347699\\
26.84	8160.69445231378\\
27.17	8170.84892875553\\
27.5	8180.78801630209\\
27.83	8190.60325702834\\
28.16	8200.03919082888\\
28.49	8209.71845172593\\
28.82	8218.63115499883\\
29.15	8227.61052625783\\
29.48	8236.30487032749\\
29.81	8244.89657449712\\
30.14	8253.3082812411\\
30.47	8261.54081725504\\
30.8	8269.54890615705\\
31.13	8277.36105827136\\
31.46	8285.17703033285\\
31.79	8292.73691245575\\
32.12	8300.17151497486\\
32.45	8307.65365886726\\
32.78	8314.728324777\\
33.11	8321.6351937447\\
33.44	8328.61340452526\\
33.77	8335.41677721881\\
34.1	8342.16756392609\\
34.43	8348.5975465363\\
34.76	8355.12583252937\\
35.09	8361.24659968604\\
35.42	8367.41726799824\\
35.75	8373.53137880806\\
36.08	8379.62931933115\\
36.41	8385.41093480669\\
36.74	8391.1512676983\\
37.07	8396.81933098338\\
37.4	8402.39742375939\\
37.73	8407.93855182376\\
38.06	8413.26926922444\\
38.39	8418.5945668842\\
38.72	8423.88463376582\\
39.05	8428.96710593835\\
39.38	8434.03474108049\\
39.71	8439.15176779557\\
40.04	8444.00264878099\\
40.37	8448.77670315346\\
40.7	8453.62306228931\\
41.03	8458.33716272251\\
41.36	8462.91507480054\\
41.69	8467.53025944741\\
42.02	8471.9356678328\\
42.35	8476.25706840716\\
42.68	8480.72086021906\\
43.01	8484.99564172434\\
43.34	8489.1378883192\\
43.67	8493.43457825958\\
44	8497.50892913508\\
};
\addlegendentry{plug-in charging}
\addplot [color=blue, dashed, line width=1.0pt]
  table[row sep=crcr]{%
11	6877.54073508476\\
11.33	6967.21569584907\\
11.66	7050.66125545389\\
11.99	7126.84367949036\\
12.32	7196.00413451219\\
12.65	7263.00058854113\\
12.98	7323.81335194004\\
13.31	7382.16263863825\\
13.64	7434.66278584805\\
13.97	7485.88661564505\\
14.3	7532.26425532534\\
14.63	7576.28104910633\\
14.96	7619.91983477366\\
15.29	7657.14730478883\\
15.62	7695.75705679049\\
15.95	7731.66783215504\\
16.28	7764.78917379785\\
16.61	7796.81135939458\\
16.94	7826.53260292393\\
17.27	7855.10686158869\\
17.6	7882.47089994023\\
17.93	7908.28693786083\\
18.26	7933.47977862194\\
18.59	7957.27297833365\\
18.92	7979.76354291019\\
19.25	8000.42003527226\\
19.58	8022.26096902435\\
19.91	8041.43628485638\\
20.24	8060.6782657514\\
20.57	8078.65948424092\\
20.9	8096.02046804988\\
21.23	8112.63031042916\\
21.56	8128.97531655403\\
21.89	8144.39476787032\\
22.22	8159.0663085681\\
22.55	8173.71186808391\\
22.88	8187.37763611203\\
23.21	8199.91367820235\\
23.54	8213.74514152388\\
23.87	8225.30801044715\\
24.2	8237.51062606844\\
24.53	8249.40776791066\\
24.86	8259.9984493242\\
25.19	8271.11353822376\\
25.52	8281.65842765824\\
25.85	8291.42788128917\\
26.18	8300.63162447853\\
26.51	8310.16033500331\\
26.84	8319.21249600552\\
27.17	8327.84232232147\\
27.5	8337.05224357806\\
27.83	8344.94730485738\\
28.16	8352.89547017084\\
28.49	8360.55519179894\\
28.82	8369.29993367502\\
29.15	8375.81738434904\\
29.48	8383.23620243106\\
29.81	8389.63105474701\\
30.14	8396.20944536442\\
30.47	8402.98029145673\\
30.8	8409.53491710631\\
31.13	8415.06595222241\\
31.46	8421.29838058522\\
31.79	8427.12843323258\\
32.12	8432.96606323143\\
32.45	8438.43683849012\\
32.78	8443.7877770283\\
33.11	8449.00746114671\\
33.44	8454.29967098734\\
33.77	8458.91570730692\\
34.1	8464.34609481154\\
34.43	8468.94515210088\\
34.76	8473.19502642495\\
35.09	8477.53224482756\\
35.42	8481.81219015424\\
35.75	8486.33251595461\\
36.08	8490.20931349412\\
36.41	8494.08273440252\\
36.74	8498.3277282889\\
37.07	8502.00234278216\\
37.4	8505.84231254631\\
37.73	8509.60138597433\\
38.06	8513.1395182672\\
38.39	8516.7455278707\\
38.72	8520.32565500691\\
39.05	8523.40695341972\\
39.38	8526.46531021258\\
39.71	8529.54209611057\\
40.04	8532.92315165534\\
40.37	8535.91886678917\\
40.7	8538.83891697393\\
41.03	8541.75269763827\\
41.36	8544.59047459967\\
41.69	8547.5490519378\\
42.02	8549.73249124514\\
42.35	8552.57705528535\\
42.68	8555.02449306441\\
43.01	8558.34430222418\\
43.34	8560.8025212796\\
43.67	8562.75132452177\\
44	8565.30330647277\\
};
\addlegendentry{battery swapping}
\end{axis}
\end{tikzpicture}%
\vspace*{-0.3in}
\caption{Number of operating vehicles $N_1$ as a function of $s$.}
\label{N1_s}
\end{subfigure}
\hfill
\begin{subfigure}[b]{0.32\linewidth}
\centering
%
%
\begin{tikzpicture}

\begin{axis}[%
width=1.694in,
height=1.03in,
at={(1.358in,0.0in)},
scale only axis,
xmin=11,
xmax=44,
xtick={11, 22, 33, 44},
xlabel style={font=\color{white!15!black}},
xlabel={Charging speed},
ymin=320,
ymax=1700,
ytick={500,1000,1500},
yticklabels={{0.5K},{1K},{1.5K}},
ylabel style={font=\color{white!15!black}},
ylabel={Number of vehicles},
axis background/.style={fill=white},
legend style={at={(1,1)}, anchor=north east, legend cell align=left, align=left, font=\scriptsize, draw=white!12!black}
]
\addplot [color=black, line width=1.0pt]
  table[row sep=crcr]{%
11	1697.4946202345\\
11.33	1666.62391642432\\
11.66	1636.73111539323\\
11.99	1608.01203073219\\
12.32	1580.21657102725\\
12.65	1553.50343739459\\
12.98	1527.70611317758\\
13.31	1502.75594143474\\
13.64	1478.79856810377\\
13.97	1455.64869019822\\
14.3	1433.16365554142\\
14.63	1411.4992130299\\
14.96	1390.48858354812\\
15.29	1370.25252199134\\
15.62	1350.65561051538\\
15.95	1331.65512314928\\
16.28	1313.19957097366\\
16.61	1295.43213138666\\
16.94	1278.16717167465\\
17.27	1261.45910787608\\
17.6	1245.17931774867\\
17.93	1229.3478653149\\
18.26	1214.0507152424\\
18.59	1199.16792172241\\
18.92	1184.76387717766\\
19.25	1170.66062555957\\
19.58	1157.01829785597\\
19.91	1143.7186577857\\
20.24	1130.72704138124\\
20.57	1118.18842752055\\
20.9	1105.96508724318\\
21.23	1093.88836912587\\
21.56	1082.30969515042\\
21.89	1070.98812318994\\
22.22	1059.90216471053\\
22.55	1049.10333811804\\
22.88	1038.63081394621\\
23.21	1028.25465980491\\
23.54	1018.22700574589\\
23.87	1008.54600502494\\
24.2	998.961844883694\\
24.53	989.672191991965\\
24.86	980.459996907184\\
25.19	971.558465152607\\
25.52	962.8190441097\\
25.85	954.246959783954\\
26.18	945.951352344017\\
26.51	937.717311807684\\
26.84	929.676130588557\\
27.17	921.885823187318\\
27.5	914.237319300379\\
27.83	906.71327579607\\
28.16	899.399946910904\\
28.49	892.064446874598\\
28.82	885.031212762083\\
29.15	878.146386062541\\
29.48	871.369669017033\\
29.81	864.714277926254\\
30.14	858.179131817877\\
30.47	851.767327515718\\
30.8	845.450149199266\\
31.13	839.336353349276\\
31.46	833.272729109718\\
31.79	827.320730490053\\
32.12	821.485972106059\\
32.45	815.729324231831\\
32.78	810.095929693595\\
33.11	804.574058923463\\
33.44	799.128242104159\\
33.77	793.778898755513\\
34.1	788.524902155753\\
34.43	783.400607103217\\
34.76	778.281606630054\\
35.09	773.329198353707\\
35.42	768.390656193972\\
35.75	763.592575627593\\
36.08	758.800259664513\\
36.41	754.161544030037\\
36.74	749.589038546548\\
37.07	745.08239667559\\
37.4	740.598929296472\\
37.73	736.200848861766\\
38.06	731.913733260442\\
38.39	727.639678330048\\
38.72	723.44844409363\\
39.05	719.336742430116\\
39.38	715.264742126264\\
39.71	711.219090433659\\
40.04	707.297187084936\\
40.37	703.44150071667\\
40.7	699.570911846125\\
41.03	695.819249535835\\
41.36	692.106182448675\\
41.69	688.42995294485\\
42.02	684.826415123523\\
42.35	681.305457076695\\
42.68	677.782079385676\\
43.01	674.33294880578\\
43.34	670.898595486028\\
43.67	667.549990628075\\
44	664.098590294112\\
};
\addlegendentry{plug-in charging}
\addplot [color=blue, dashed, line width=1.0pt]
  table[row sep=crcr]{%
11	433.235842510616\\
11.33	428.394900578637\\
11.66	423.736488321451\\
11.99	419.397762176437\\
12.32	415.418605793997\\
12.65	411.443289188318\\
12.98	407.833254072246\\
13.31	404.300289705103\\
13.64	401.165361193655\\
13.97	398.0395875984\\
14.3	395.26549804252\\
14.63	392.629523620424\\
14.96	389.948727312347\\
15.29	387.823038699436\\
15.62	385.496059447832\\
15.95	383.377703662043\\
16.28	381.488178876525\\
16.61	379.660676994004\\
16.94	378.026429340377\\
17.27	376.469822635367\\
17.6	375.000815120644\\
17.93	373.660330960704\\
18.26	372.35084718224\\
18.59	371.160300938811\\
18.92	370.082380639547\\
19.25	369.181070542888\\
19.58	368.114172524853\\
19.91	367.324644641382\\
20.24	366.500129708569\\
20.57	365.79670628333\\
20.9	365.141505141772\\
21.23	364.551190635031\\
21.56	363.969338269092\\
21.89	363.476428524359\\
22.22	363.05302395638\\
22.55	362.611583500781\\
22.88	362.270334551578\\
23.21	362.050497766204\\
23.54	361.649491831852\\
23.87	361.515100372614\\
24.2	361.2819015152\\
24.53	361.069841004849\\
24.86	361.008070336666\\
25.19	360.862531678809\\
25.52	360.774162292599\\
25.85	360.771186199626\\
26.18	360.827759617445\\
26.51	360.826695410138\\
26.84	360.874213533576\\
27.17	360.96381543283\\
27.5	360.961303481683\\
27.83	361.12214601876\\
28.16	361.262629273901\\
28.49	361.429416863993\\
28.82	361.434009548511\\
29.15	361.732259692256\\
29.48	361.893814011438\\
29.81	362.185617929111\\
30.14	362.44038815982\\
30.47	362.656307065735\\
30.8	362.89091082491\\
31.13	363.259609935338\\
31.46	363.517565042966\\
31.79	363.822010194292\\
32.12	364.11451769621\\
32.45	364.449336455991\\
32.78	364.791291086164\\
33.11	365.142069907109\\
33.44	365.472070328873\\
33.77	365.891833648995\\
34.1	366.181137909394\\
34.43	366.583999736328\\
34.76	367.029826539824\\
35.09	367.453504501458\\
35.42	367.876724531144\\
35.75	368.254416351746\\
36.08	368.720845480911\\
36.41	369.179274331465\\
36.74	369.571977412294\\
37.07	370.043798356025\\
37.4	370.481714923456\\
37.73	370.923486841451\\
38.06	371.391549966446\\
38.39	371.840808548722\\
38.72	372.285905920936\\
39.05	372.801915432648\\
39.38	373.313854015203\\
39.71	373.81540039741\\
40.04	374.260479274711\\
40.37	374.759791156007\\
40.7	375.263818758624\\
41.03	375.761522897385\\
41.36	376.264257229613\\
41.69	376.740029150559\\
42.02	377.335933589389\\
42.35	377.816031854058\\
42.68	378.35491714761\\
43.01	378.741628685164\\
43.34	379.264590043756\\
43.67	379.866476650805\\
44	380.360470488084\\
};
\addlegendentry{battery swapping}
\end{axis}
\end{tikzpicture}%
\vspace*{-0.3in}
\caption{Number of non-operating vehicles $N_2$ as a function of $s$.}
\label{N2_s}
\end{subfigure}
\begin{subfigure}[b]{0.32\linewidth}
\centering
%
%
\begin{tikzpicture}

\begin{axis}[%
width=1.694in,
height=1.03in,
at={(1.358in,0.0in)},
scale only axis,
xmin=11,
xmax=44,
xtick={11, 22, 33, 44},
xlabel style={font=\color{white!15!black}},
xlabel={Charging speed},
ymin=12,
ymax=21,
ylabel style={font=\color{white!15!black}},
ylabel={Searching time},
axis background/.style={fill=white},
legend style={at={(1,1)}, anchor=north east, legend cell align=left, align=left, font=\scriptsize, draw=white!12!black}
]
\addplot [color=black, line width=1.0pt]
  table[row sep=crcr]{%
11	19.2030828492053\\
11.33	19.0916492645253\\
11.66	18.9726652133048\\
11.99	18.8674909626484\\
12.32	18.7555449443259\\
12.65	18.6548487142774\\
12.98	18.5537708713584\\
13.31	18.4512919526866\\
13.64	18.3626585264504\\
13.97	18.275236781263\\
14.3	18.1795353066469\\
14.63	18.0927917143647\\
14.96	18.0025687347144\\
15.29	17.9218026799362\\
15.62	17.8426003354422\\
15.95	17.7622392239892\\
16.28	17.6791948571998\\
16.61	17.6083109532499\\
16.94	17.5350390843153\\
17.27	17.4665905440609\\
17.6	17.3935265782215\\
17.93	17.3193811543189\\
18.26	17.2531343588175\\
18.59	17.186245455551\\
18.92	17.1256067507667\\
19.25	17.0572749385363\\
19.58	16.9966929115253\\
19.91	16.9346221504223\\
20.24	16.8696005113565\\
20.57	16.8149948368734\\
20.9	16.7604625797426\\
21.23	16.6920635855345\\
21.56	16.6411544294435\\
21.89	16.5880119922155\\
22.22	16.5321136123568\\
22.55	16.4776713010111\\
22.88	16.4300873090251\\
23.21	16.3697267595748\\
23.54	16.3197493050863\\
23.87	16.2792583019989\\
24.2	16.2283504030691\\
24.53	16.1847011095544\\
24.86	16.1307226272987\\
25.19	16.0862348156576\\
25.52	16.0403233581892\\
25.85	15.9919821003645\\
26.18	15.9508519453608\\
26.51	15.9023912268606\\
26.84	15.8555970350664\\
27.17	15.8159842707196\\
27.5	15.7751569820039\\
27.83	15.7325943940677\\
28.16	15.6946220561751\\
28.49	15.6451511121886\\
28.82	15.6062525325889\\
29.15	15.5707721936092\\
29.48	15.5325007489838\\
29.81	15.4944748040444\\
30.14	15.4562371243338\\
30.47	15.4183508274189\\
30.8	15.3779951865477\\
31.13	15.3457354433807\\
31.46	15.3095554063538\\
31.79	15.2728854221365\\
32.12	15.2377471054397\\
32.45	15.2025154414147\\
32.78	15.1669159440707\\
33.11	15.1324813838645\\
33.44	15.0986921061972\\
33.77	15.0650110244437\\
34.1	15.0328727755954\\
34.43	15.0026493758535\\
34.76	14.967574894676\\
35.09	14.9371333972667\\
35.42	14.902413566796\\
35.75	14.8743307324679\\
36.08	14.840788909429\\
36.41	14.8122327128763\\
36.74	14.7839031406406\\
37.07	14.7555315274154\\
37.4	14.7229326010474\\
37.73	14.692879656142\\
38.06	14.6656451411265\\
38.39	14.6346971963199\\
38.72	14.6064608619786\\
39.05	14.5783519052157\\
39.38	14.5493142209968\\
39.71	14.5190946969248\\
40.04	14.492701212195\\
40.37	14.4673569250027\\
40.7	14.4374990305462\\
41.03	14.4132093766376\\
41.36	14.3866812477083\\
41.69	14.360445748875\\
42.02	14.3343228351778\\
42.35	14.3112511261268\\
42.68	14.286580263221\\
43.01	14.2626821786343\\
43.34	14.2346413754628\\
43.67	14.2141694034735\\
44	14.1766209741189\\
};
\addlegendentry{plug-in charging}
\addplot [color=blue, dashed, line width=1.0pt]
  table[row sep=crcr]{%
11	12.238730192722\\
11.33	12.3197186176398\\
11.66	12.4028774253131\\
11.99	12.490100356382\\
12.32	12.5815541328564\\
12.65	12.6706104534341\\
12.98	12.7634772454527\\
13.31	12.8549835961737\\
13.64	12.9506001668276\\
13.97	13.043605493029\\
14.3	13.1398679021787\\
14.63	13.2357557194027\\
14.96	13.3281869437394\\
15.29	13.4274613151466\\
15.62	13.5207025369824\\
15.95	13.6149850517454\\
16.28	13.7107417705914\\
16.61	13.8051663835002\\
16.94	13.90061552329\\
17.27	13.9951466087841\\
17.6	14.0890055218387\\
17.93	14.1829909180855\\
18.26	14.275426922725\\
18.59	14.3679630026382\\
18.92	14.4605603776769\\
19.25	14.554484963596\\
19.58	14.6434757434526\\
19.91	14.7358035093967\\
20.24	14.8257122145719\\
20.57	14.9161952384857\\
20.9	15.0059337974017\\
21.23	15.0953066972621\\
21.56	15.1832590849763\\
21.89	15.2714167236438\\
22.22	15.3594602211087\\
22.55	15.4456854924314\\
22.88	15.5325117272046\\
23.21	15.6204488037353\\
23.54	15.7033128344378\\
23.87	15.7903906603599\\
24.2	15.8741398905524\\
24.53	15.9570731454621\\
24.86	16.0420020976408\\
25.19	16.1239488758743\\
25.52	16.2059717468452\\
25.85	16.288742821554\\
26.18	16.3717502911186\\
26.51	16.4524133768976\\
26.84	16.5331202067068\\
27.17	16.6137676296066\\
27.5	16.691298791239\\
27.83	16.7716166579973\\
28.16	16.8505225890842\\
28.49	16.9291214151462\\
28.82	17.0029442142122\\
29.15	17.0828759215098\\
29.48	17.1586632588669\\
29.81	17.2367674086019\\
30.14	17.3131517565073\\
30.47	17.3877625749549\\
30.8	17.4620289286314\\
31.13	17.5388868260054\\
31.46	17.6121977581151\\
31.79	17.6859383042823\\
32.12	17.7586369995084\\
32.45	17.8317099307249\\
32.78	17.9042682597007\\
33.11	17.9763727813772\\
33.44	18.0472489258307\\
33.77	18.1198578266923\\
34.1	18.1883156465753\\
34.43	18.2591916928024\\
34.76	18.3306267993882\\
35.09	18.4008591661437\\
35.42	18.4704871061235\\
35.75	18.538258171048\\
36.08	18.6079355447243\\
36.41	18.6768353078138\\
36.74	18.7433144263538\\
37.07	18.81150298987\\
37.4	18.8781858770115\\
37.73	18.9444531055814\\
38.06	19.0109676146979\\
38.39	19.0764238001801\\
38.72	19.1412527210664\\
39.05	19.2076990441675\\
39.38	19.2735594797231\\
39.71	19.3386451671371\\
40.04	19.4015588270534\\
40.37	19.465658177723\\
40.7	19.529458403452\\
41.03	19.5926260633036\\
41.36	19.655517533312\\
41.69	19.7171396607438\\
42.02	19.7821124656067\\
42.35	19.8430441038964\\
42.68	19.9054285550939\\
43.01	19.9625581247372\\
43.34	20.0236330654813\\
43.67	20.0868758892553\\
44	20.1462622459275\\
};
\addlegendentry{battery swapping}
\end{axis}
\end{tikzpicture}%
\vspace*{-0.3in}
\caption{Vehicle searching time $t_m$ (minute) as a function of $s$.}
\label{tm_s}
\end{subfigure}
\hfill
\begin{subfigure}[b]{0.32\linewidth}
\centering
%
%
\begin{tikzpicture}

\begin{axis}[%
width=1.694in,
height=1.03in,
at={(1.358in,0.0in)},
scale only axis,
xmin=11,
xmax=44,
xtick={11, 22, 33, 44},
xlabel style={font=\color{white!15!black}},
xlabel={Charging speed},
ymin=0,
ymax=25,
ylabel style={font=\color{white!15!black}},
ylabel={Waiting time},
axis background/.style={fill=white},
legend style={at={(1,1)}, anchor=north east, legend cell align=left, align=left, font=\scriptsize, draw=white!12!black}
]
\addplot [color=black, line width=1.0pt]
  table[row sep=crcr]{%
11	1.49372947752543\\
11.33	1.49515338294839\\
11.66	1.49957338289551\\
11.99	1.50367051142496\\
12.32	1.50919307343669\\
12.65	1.51309622164317\\
12.98	1.51665743564217\\
13.31	1.5217935982959\\
13.64	1.52455267925154\\
13.97	1.527660541035\\
14.3	1.53254796736987\\
14.63	1.53727484861273\\
14.96	1.54277126012544\\
15.29	1.54574762559105\\
15.62	1.55104715811144\\
15.95	1.55526817632418\\
16.28	1.56020006795036\\
16.61	1.56404531016655\\
16.94	1.56711230506838\\
17.27	1.57041087964265\\
17.6	1.57499442592389\\
17.93	1.58043103879757\\
18.26	1.58291630374242\\
18.59	1.58902317464416\\
18.92	1.58964105428119\\
19.25	1.59699613632409\\
19.58	1.5984306323609\\
19.91	1.60302617583269\\
20.24	1.60682846116449\\
20.57	1.61016358566777\\
20.9	1.6118435017871\\
21.23	1.61772462330075\\
21.56	1.62090413072803\\
21.89	1.62471078353628\\
22.22	1.62756107411636\\
22.55	1.63303806346866\\
22.88	1.63532386164789\\
23.21	1.64070168743745\\
23.54	1.6423189987348\\
23.87	1.64490815668319\\
24.2	1.64944781010274\\
24.53	1.65242151669644\\
24.86	1.65614759779422\\
25.19	1.65892431514386\\
25.52	1.65971531573794\\
25.85	1.66361717938641\\
26.18	1.66907380170785\\
26.51	1.67114542204235\\
26.84	1.67438710935378\\
27.17	1.67741723406311\\
27.5	1.68074128925486\\
27.83	1.68336016092582\\
28.16	1.68697792304764\\
28.49	1.68796610282232\\
28.82	1.69374767804942\\
29.15	1.69487185956857\\
29.48	1.69838207799708\\
29.81	1.70079129072481\\
30.14	1.70343276221839\\
30.47	1.70619565707191\\
30.8	1.71073392231643\\
31.13	1.71328924928016\\
31.46	1.71487717625536\\
31.79	1.71828756597399\\
32.12	1.72098997711396\\
32.45	1.72097995631349\\
32.78	1.72523319366376\\
33.11	1.7297285146614\\
33.44	1.73105400954794\\
33.77	1.73317194485897\\
34.1	1.73378710136126\\
34.43	1.7369405526406\\
34.76	1.73850517895267\\
35.09	1.74316607550681\\
35.42	1.74689006274866\\
35.75	1.7476587490271\\
36.08	1.74898058634652\\
36.41	1.75206692368498\\
36.74	1.75425368077591\\
37.07	1.75615139603261\\
37.4	1.7595314213378\\
37.73	1.76123992971266\\
38.06	1.76399814544026\\
38.39	1.76681211356414\\
38.72	1.76795435186108\\
39.05	1.77115892516626\\
39.38	1.77368845838202\\
39.71	1.77460496866581\\
40.04	1.77735168949398\\
40.37	1.77983452097688\\
40.7	1.78157460530873\\
41.03	1.78240182928034\\
41.36	1.78520272148714\\
41.69	1.786200561528\\
42.02	1.78962115475963\\
42.35	1.79233895378487\\
42.68	1.79219774556847\\
43.01	1.79398570237008\\
43.34	1.79857106014669\\
43.67	1.79678848390889\\
44	1.80412897078916\\
};
\addlegendentry{plug-in charging}
\addplot [color=blue, dashed, line width=1.0pt]
  table[row sep=crcr]{%
11	23.494692573419\\
11.33	22.5068545138904\\
11.66	21.5871336626423\\
11.99	20.7454924649101\\
12.32	19.9790167709633\\
12.65	19.2384426814278\\
12.98	18.5637929025112\\
13.31	17.9174288164014\\
13.64	17.3329077406052\\
13.97	16.764226871229\\
14.3	16.2467498842031\\
14.63	15.7554451518763\\
14.96	15.2704092474639\\
15.29	14.8513545841858\\
15.62	14.4207687293432\\
15.95	14.018975515729\\
16.28	13.6466764833895\\
16.61	13.2870226591527\\
16.94	12.9517340852762\\
17.27	12.6293435078059\\
17.6	12.3203425785473\\
17.93	12.0279551404314\\
18.26	11.7430459585647\\
18.59	11.4730972550815\\
18.92	11.2171089377503\\
19.25	10.9801554750398\\
19.58	10.7327611862141\\
19.91	10.5123593171976\\
20.24	10.2923008890467\\
20.57	10.0854307184802\\
20.9	9.88548237727954\\
21.23	9.69365014288733\\
21.56	9.50522363353484\\
21.89	9.32652862410386\\
22.22	9.1558117759507\\
22.55	8.98606906742592\\
22.88	8.82645000224287\\
23.21	8.6783616420528\\
23.54	8.51848405170344\\
23.87	8.38098677193091\\
24.2	8.2379384569872\\
24.53	8.09841210849545\\
24.86	7.97191208042561\\
25.19	7.84095808425721\\
25.52	7.71595110175852\\
25.85	7.59882221166136\\
26.18	7.48755726502025\\
26.51	7.37366445641943\\
26.84	7.26475265368664\\
27.17	7.16027502435973\\
27.5	7.0507353851697\\
27.83	6.95404723316568\\
28.16	6.85724770629922\\
28.49	6.76355612231287\\
28.82	6.66002613422201\\
29.15	6.57783354302266\\
29.48	6.48751615595798\\
29.81	6.40715522640507\\
30.14	6.32540612137421\\
30.47	6.24217804853139\\
30.8	6.16129742448755\\
31.13	6.09029642821516\\
31.46	6.01303942454261\\
31.79	5.93983348929973\\
32.12	5.86684515889877\\
32.45	5.797551293359\\
32.78	5.72964879952044\\
33.11	5.66322586498328\\
33.44	5.59639348063021\\
33.77	5.53608707471089\\
34.1	5.46848051870891\\
34.43	5.40881517369686\\
34.76	5.35261240604055\\
35.09	5.2958280476376\\
35.42	5.23979516417997\\
35.75	5.18177060954653\\
36.08	5.12988492646152\\
36.41	5.07823887380696\\
36.74	5.02338631978769\\
37.07	4.97397469963905\\
37.4	4.92324771468476\\
37.73	4.87344702446746\\
38.06	4.82586634253245\\
38.39	4.77785359424857\\
38.72	4.73026326434804\\
39.05	4.68740428947455\\
39.38	4.64491886465391\\
39.71	4.60244025262024\\
40.04	4.55736255481876\\
40.37	4.51595901097624\\
40.7	4.47540153553069\\
41.03	4.43506213298313\\
41.36	4.39556814342957\\
41.69	4.35514207681485\\
42.02	4.32187239701567\\
42.35	4.28276895141768\\
42.68	4.24740438140002\\
43.01	4.20432192222724\\
43.34	4.16914794086797\\
43.67	4.1386963095698\\
44	4.10295366093938\\
};
\addlegendentry{battery swapping}
\end{axis}
\end{tikzpicture}%
\vspace*{-0.3in}
\caption{Vehicle waiting time $t_w$ (minute) as a function of $s$.}
\label{tw_s}
\end{subfigure}
\hfill
\begin{subfigure}[b]{0.32\linewidth}
\centering
%
%
\begin{tikzpicture}

\begin{axis}[%
width=1.694in,
height=1.03in,
at={(1.358in,0.0in)},
scale only axis,
xmin=11,
xmax=44,
xtick={11, 22, 33, 44},
xlabel style={font=\color{white!15!black}},
xlabel={Charging speed},
ymin=0,
ymax=140,
ylabel style={font=\color{white!15!black}},
ylabel={Charging time},
axis background/.style={fill=white},
legend style={at={(1,1)}, anchor=north east, legend cell align=left, align=left, font=\scriptsize, draw=white!12!black}
]
\addplot [color=black, line width=1.0pt]
  table[row sep=crcr]{%
11	122.727272727273\\
11.33	119.152691968226\\
11.66	115.780445969125\\
11.99	112.593828190158\\
12.32	109.577922077922\\
12.65	106.719367588933\\
12.98	104.006163328197\\
13.31	101.427498121713\\
13.64	98.9736070381232\\
13.97	96.6356478167502\\
14.3	94.4055944055944\\
14.63	92.2761449077238\\
14.96	90.2406417112299\\
15.29	88.2930019620667\\
15.62	86.4276568501921\\
15.95	84.6394984326019\\
16.28	82.9238329238329\\
16.61	81.2763395544852\\
16.94	79.6930342384888\\
17.27	78.1702374059062\\
17.6	76.7045454545455\\
17.93	75.292805354155\\
18.26	73.9320920043811\\
18.59	72.6196880043034\\
18.92	71.353065539112\\
19.25	70.1298701298701\\
19.58	68.9479060265577\\
19.91	67.8051230537418\\
20.24	66.699604743083\\
20.57	65.6295576081672\\
20.9	64.5933014354067\\
21.23	63.5892604804522\\
21.56	62.6159554730983\\
21.89	61.6719963453632\\
22.22	60.7560756075608\\
22.55	59.8669623059867\\
22.88	59.0034965034965\\
23.21	58.1645842309349\\
23.54	57.3491928632116\\
23.87	56.5563468789275\\
24.2	55.7851239669421\\
24.53	55.0346514472075\\
24.86	54.3041029766694\\
25.19	53.5926955140929\\
25.52	52.8996865203762\\
25.85	52.2243713733075\\
26.18	51.5660809778457\\
26.51	50.9241795548849\\
26.84	50.2980625931446\\
27.17	49.6871549503128\\
27.5	49.0909090909091\\
27.83	48.5088034495149\\
28.16	47.9403409090909\\
28.49	47.3850473850474\\
28.82	46.8424705065926\\
29.15	46.3121783876501\\
29.48	45.7937584803256\\
29.81	45.2868165045287\\
30.14	44.7909754479098\\
30.47	44.3058746307844\\
30.8	43.8311688311688\\
31.13	43.3665274654674\\
31.46	42.9116338207247\\
31.79	42.4661843346964\\
32.12	42.0298879202989\\
32.45	41.6024653312789\\
32.78	41.1836485661989\\
33.11	40.773180308064\\
33.44	40.3708133971292\\
33.77	39.9763103346165\\
34.1	39.5894428152493\\
34.43	39.2099912866686\\
34.76	38.8377445339471\\
35.09	38.4724992875463\\
35.42	38.1140598531903\\
35.75	37.7622377622378\\
36.08	37.4168514412417\\
36.41	37.0777258994782\\
36.74	36.7446924333152\\
37.07	36.4175883463717\\
37.4	36.096256684492\\
37.73	35.7805459846276\\
38.06	35.470310036784\\
38.39	35.1654076582443\\
38.72	34.8657024793388\\
39.05	34.5710627400768\\
39.38	34.2813610970036\\
39.71	33.9964744396877\\
40.04	33.7162837162837\\
40.37	33.4406737676492\\
40.7	33.1695331695332\\
41.03	32.9027540823787\\
41.36	32.6402321083172\\
41.69	32.3818661549532\\
42.02	32.1275583055688\\
42.35	31.8772136953955\\
42.68	31.630740393627\\
43.01	31.3880492908626\\
43.34	31.1490539916936\\
43.67	30.9136707121594\\
44	30.6818181818182\\
};
\addlegendentry{plug-in charging}
\addplot [color=blue, dashed, line width=1.0pt]
  table[row sep=crcr]{%
11	2\\
11.33	2\\
11.66	2\\
11.99	2\\
12.32	2\\
12.65	2\\
12.98	2\\
13.31	2\\
13.64	2\\
13.97	2\\
14.3	2\\
14.63	2\\
14.96	2\\
15.29	2\\
15.62	2\\
15.95	2\\
16.28	2\\
16.61	2\\
16.94	2\\
17.27	2\\
17.6	2\\
17.93	2\\
18.26	2\\
18.59	2\\
18.92	2\\
19.25	2\\
19.58	2\\
19.91	2\\
20.24	2\\
20.57	2\\
20.9	2\\
21.23	2\\
21.56	2\\
21.89	2\\
22.22	2\\
22.55	2\\
22.88	2\\
23.21	2\\
23.54	2\\
23.87	2\\
24.2	2\\
24.53	2\\
24.86	2\\
25.19	2\\
25.52	2\\
25.85	2\\
26.18	2\\
26.51	2\\
26.84	2\\
27.17	2\\
27.5	2\\
27.83	2\\
28.16	2\\
28.49	2\\
28.82	2\\
29.15	2\\
29.48	2\\
29.81	2\\
30.14	2\\
30.47	2\\
30.8	2\\
31.13	2\\
31.46	2\\
31.79	2\\
32.12	2\\
32.45	2\\
32.78	2\\
33.11	2\\
33.44	2\\
33.77	2\\
34.1	2\\
34.43	2\\
34.76	2\\
35.09	2\\
35.42	2\\
35.75	2\\
36.08	2\\
36.41	2\\
36.74	2\\
37.07	2\\
37.4	2\\
37.73	2\\
38.06	2\\
38.39	2\\
38.72	2\\
39.05	2\\
39.38	2\\
39.71	2\\
40.04	2\\
40.37	2\\
40.7	2\\
41.03	2\\
41.36	2\\
41.69	2\\
42.02	2\\
42.35	2\\
42.68	2\\
43.01	2\\
43.34	2\\
43.67	2\\
44	2\\
};
\addlegendentry{battery swapping}
\end{axis}
\end{tikzpicture}%
\vspace*{-0.3in}
\caption{Vehicle charging time $t_s$ (minute) as a function of $s$.}
\label{ts_s}
\end{subfigure}
\caption{Market outcomes on electrified AMoD services under different charging speed. Black lines present the results when charging vehicles; blue dashed lines show the results when swapping batteries.}
\label{AMoD_service_s}
\end{figure}

Figure \ref{AMoD_service_s} shows the market outcomes of the AMoD system under different charging speeds. As charging speed improves, the TNC lowers the ride fare (Figure \ref{pf_s}) and hires a larger fleet of vehicles (Figure \ref{N_s}). This increases the vehicle idle time (Figure \ref{wv_s}) and the number of idle vehicles, resulting in a reduced passenger waiting time (Figure \ref{wc_s}). The total passenger travel cost drops with the charging speed (Figure \ref{c_s}), and therefore more passengers choose the AMoD services (Figure \ref{lambda_s}). This is the case for both plug-in charging and battery swapping. Figure \ref{gamma_s}-\ref{ts_s} illustrates the charging demand and the average time spent on charging for vehicles under the two charging strategies. Numerical results show that the charging demand (Figure \ref{gamma_s}) and the number of operating vehicles (Figure \ref{N1_s}) under two distinct charging strategies monotonically increase with the charging speed. For plug-in charging, as charging speed improves, the vehicle searching time reduces (Figure \ref{tm_s}), the vehicle waiting time slightly decreases (Figure \ref{tw_s}), and the average vehicle charging time sharply drops (Figure \ref{ts_s}). Under battery swapping, there is a small number of non-operating vehicles, which experience a little change with charging speed (Figure \ref{N2_s}). As the charging speed improves, the vehicle searching time rises with the charging speed (Figure \ref{tm_s}), the vehicle waiting time reduces as the charging speed increases (Figure \ref{tw_s}), and the vehicle charging time (i.e., the service time of battery swapping) remains constant under different charging speeds (Figure \ref{ts_s}).

Overall, we summarize the key findings as follows:
\begin{itemize}
    \item When deploying charging stations for AMoD services, improved charging speed results in a transformation of the infrastructure deployment from {\em sparsely distributed large stations} to {\em densely distributed small stations} (Figure \ref{optimal_K_s} and Figure \ref{optimal_Q_s}).
    \item When deploying battery swapping stations for AMoD services, improved charging speed reduces the need for charging infrastructure (Figure \ref{optimal_K_s}) and leads to a smaller number of battery swapping stations.
    \item Under the strategy of plug-in charging, improved charging speed brings a Pareto improvement: a higher charging speed always leads to a higher passenger surplus (Figure \ref{PS_s}), a higher platform profit (Figure \ref{profit_s}) and increased social welfare (Figure \ref{SW_s}).
    \item Under the strategy of battery swapping, passengers always benefit from the improved charging technology (Figure \ref{PS_s}) and the social welfare monotonically increases with the charging speed (Figure \ref{SW_s}). However, the TNC may get hurt when the charging speed is relatively high (Figure \ref{profit_s}).
\end{itemize} 

To explain the impacts of charging speed on charging infrastructure planning decisions, we first note that improving charging speed directly increases the service rate of each charger. This motivates the government to reduce the total supply of chargers (i.e., $KQ$) for cost-saving (Figure \ref{infrastructure_cost_s}). Specifically, for battery swapping stations, the government has to reduce the number of battery swapping stations $K$ as $Q$ is fixed (Figure \ref{optimal_K_s}). On the other hand, for plug-in charging, the government is flexible in determining the number of charging stations $K$ and the number of chargers per station $Q$. In this case, it faces a trade-off: on the one hand, it can build a large number of small charging stations (i.e., large $K$ and small $Q$), which leads to densely distributed charging infrastructures. This reduces the time vehicles search for/travel to nearby charging stations.  On the other hand, it can build a small number of large charging stations (i.e., small $K$ and large $Q$), which leads to sparsely distributed large charging centers. This increases the service capability of each charging station and results in a shorter waiting time at stations. The government should carefully control the knobs to achieve a balance between the vehicle searching time to stations and the vehicle waiting time at stations. As the charging speed improves, it helps to curb the vehicle waiting time at stations. In this case, the improved charging speed alters the trade-off incurred by infrastructure deployment: the marginal effect of allocating more chargers per station on reducing vehicle waiting time diminishes, and the vehicle searching time gradually becomes the bottleneck. Consequently, the government builds more stations and reduces the number of chargers at each station, which leads to the transformation of the infrastructure deployment from {\em sparsely distributed large stations} to {\em densely distributed small stations}.

{A counter-intuitive result is that increased charging speed may not always lead to a higher TNC profit: in certain regimes, the TNC profit reduces as charging speed improves.} This suggests that, unlike plug-in charging, where all stakeholders benefit from the increased charging speed, the TNC platform may not necessarily welcome fast charging technology when battery swapping stations are deployed. We comment that such a difference is due to the limited flexibility of battery swapping stations for AMoD services. In the deployment of charging stations, the government is flexible to jointly control the number of charging stations $K$ and the number of chargers per station $Q$. The improved charging speed promotes the government to lower the total supply of charging infrastructures for cost-saving. In this case, although the total supply of charging infrastructures reduces, the government builds fewer charging stations and allocates more chargers to each station to balance the vehicle searching time and the vehicle waiting time, both the vehicle searching time (Figure \ref{tm_s}), the vehicle waiting time (Figure \ref{tw_s}), and the vehicle charging time (Figure \ref{ts_s}) reduce as charging speed increases. This leads to a more efficient use of the TNC fleet, which benefits the TNC platform. On the other hand, for battery swapping, the government can only control the total number of stations $K$ to adjust the supply of charging infrastructures\footnote{As mentioned in Section \ref{social_welfare_maximization}, the battery swapping station typically has a fixed number of chargers $Q$, which is difficult to adjust in short term because of its closed structure, and because of the need of standardization for the convenience of mass production.}. As charging speed increases, the government has to decrease the number of battery swapping stations to reduce infrastructure costs, which leads to a sparser distribution of battery swapping stations and increases the vehicle searching time to stations. Although improving charging speed increases each battery charger's service rate and helps curb the vehicle waiting time at battery swapping stations, the reduced supply of charging infrastructures also leads to a longer vehicle searching time to stations. When the negative effect of reduced charging infrastructure supply dominates, the vehicle searching time significantly increases (Figure \ref{tm_s}), the number of non-operating vehicles increases (Figure \ref{N2_s}), and the fleet utilization rate drops. To maintain the operation of AMoD services, the TNC needs to recruit a larger fleet of vehicles, which leads to increased vehicle operating costs and reduced profit.

\subsection{Impacts of battery capacity} \label{impacts_battery_capacity}

This section investigates the impacts of battery capacity on the charging infrastructure planning and AMoD market. To this end, we will solve the bi-level optimization under distinct values of the model parameter $C$ under different charging strategies and compare the corresponding results.

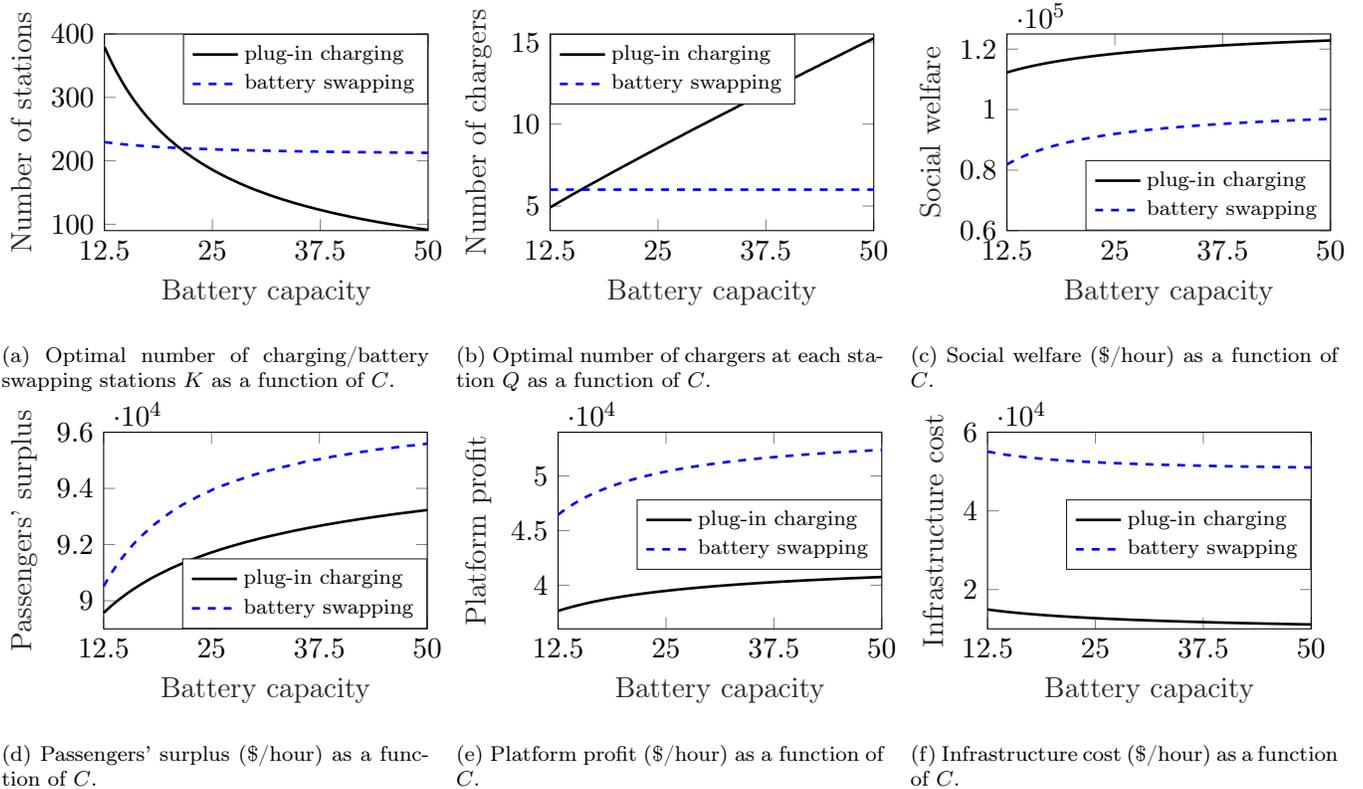
\begin{figure}[h]
\centering
\begin{subfigure}[b]{0.32\linewidth}
\centering
%
%
\begin{tikzpicture}

\begin{axis}[%
width=1.694in,
height=1.03in,
at={(1.358in,0.0in)},
scale only axis,
xmin=12.5,
xmax=50,
xtick={12.5,   25, 37.5,   50},
xlabel style={font=\color{white!15!black}},
xlabel={Battery capacity},
ymin=90,
ymax=400,
ylabel style={font=\color{white!15!black}},
ylabel={Number of stations},
axis background/.style={fill=white},
legend style={at={(1,1)}, anchor=north east, legend cell align=left, align=left, font=\scriptsize, draw=white!12!black}
]
\addplot [color=black, line width=1.0pt]
  table[row sep=crcr]{%
12.5	379.531953950936\\
12.875	368.167052468782\\
13.25	357.394766898598\\
13.625	347.195817234188\\
14	337.72138946224\\
14.375	328.690846771197\\
14.75	319.947020297796\\
15.125	311.669555217373\\
15.5	303.920083103988\\
15.875	296.545774993682\\
16.25	289.734811368542\\
16.625	282.845448724584\\
17	276.526573861139\\
17.375	270.297548968177\\
17.75	264.600517227373\\
18.125	258.975926225842\\
18.5	253.615441846501\\
18.875	248.30233756843\\
19.25	243.293431712792\\
19.625	238.513838098682\\
20	233.819733112096\\
20.375	229.721711147557\\
20.75	225.204589489747\\
21.125	221.15214617664\\
21.5	217.305470643316\\
21.875	213.406941644926\\
22.25	209.765330175014\\
22.625	206.223047964022\\
23	202.701821674108\\
23.375	199.356496848628\\
23.75	196.072007206155\\
24.125	192.90652971618\\
24.5	189.978736683739\\
24.875	186.971075151652\\
25.25	184.098408028511\\
25.625	181.376688920775\\
26	178.624457739957\\
26.375	176.08274816951\\
26.75	173.451341125581\\
27.125	171.077556285162\\
27.5	168.723609885462\\
27.875	166.378287919223\\
28.25	163.937113082262\\
28.625	161.882511303838\\
29	159.651210733414\\
29.375	157.476151161613\\
29.75	155.499056857749\\
30.125	153.538472356313\\
30.5	151.684807947061\\
30.875	149.791750760156\\
31.25	147.893607045603\\
31.625	146.066618694871\\
32	144.321121396671\\
32.375	142.578624323841\\
32.75	140.935709241389\\
33.125	139.254598124938\\
33.5	137.742631275014\\
33.875	136.140779121169\\
34.25	134.545477445639\\
34.625	133.041366749587\\
35	131.676432689813\\
35.375	130.233682753989\\
35.75	128.784690868843\\
36.125	127.491368934898\\
36.5	126.131109374586\\
36.875	124.752032695405\\
37.25	123.413724477938\\
37.625	122.300474215353\\
38	121.022197141067\\
38.375	119.775056691138\\
38.75	118.495682606806\\
39.125	117.301165330436\\
39.5	116.208542683538\\
39.875	115.136327190921\\
40.25	114.105908500477\\
40.625	112.979096339967\\
41	111.778077083176\\
41.375	110.827919209887\\
41.75	109.810471708742\\
42.125	108.746636995508\\
42.5	107.74545731471\\
42.875	106.848354891046\\
43.25	105.840239016028\\
43.625	104.985722764993\\
44	103.972358425303\\
44.375	103.024983423553\\
44.75	102.240785327942\\
45.125	101.276815440359\\
45.5	100.492525677621\\
45.875	99.6424017960074\\
46.25	98.8248323998594\\
46.625	97.9610454709144\\
47	97.1783488461864\\
47.375	96.3697123084225\\
47.75	95.576543992864\\
48.125	94.7597671938761\\
48.5	94.0716773468464\\
48.875	93.384238347596\\
49.25	92.6185805360878\\
49.625	91.9014677145605\\
50	91.1689982068743\\
};
\addlegendentry{plug-in charging}
\addplot [color=blue, dashed, line width=1.0pt]
  table[row sep=crcr]{%
12.5	229.537957906723\\
12.875	228.864383650944\\
13.25	228.273022174834\\
13.625	227.642822253256\\
14	227.064334601163\\
14.375	226.514387130726\\
14.75	226.016235351562\\
15.125	225.499725341797\\
15.5	225.079536438079\\
15.875	224.681091308594\\
16.25	224.241638171952\\
16.625	223.828125\\
17	223.449707029795\\
17.375	223.082733899355\\
17.75	222.73864634335\\
18.125	222.463989257812\\
18.5	222.110366774724\\
18.875	221.78029975621\\
19.25	221.514892531559\\
19.625	221.258544921875\\
20	220.961001518253\\
20.375	220.731976628304\\
20.75	220.445251464844\\
21.125	220.188903436065\\
21.5	220.025634765625\\
21.875	219.793701171876\\
22.25	219.535732269287\\
22.625	219.358825334348\\
23	219.110083582927\\
23.375	218.963813775918\\
23.75	218.773651123229\\
24.125	218.627929687131\\
24.5	218.392944335915\\
24.875	218.225097656159\\
25.25	218.116760253906\\
25.625	217.896270751953\\
26	217.748643457889\\
26.375	217.618560791016\\
26.75	217.425536737028\\
27.125	217.309570405621\\
27.5	217.163084260756\\
27.875	217.064619064331\\
28.25	216.864007709955\\
28.625	216.760218329727\\
29	216.638189554214\\
29.375	216.50390625\\
29.75	216.404700325847\\
30.125	216.297709941501\\
30.5	216.169738769713\\
30.875	216.046142578136\\
31.25	215.969848632812\\
31.625	215.84091186524\\
32	215.742492674326\\
32.375	215.664672851562\\
32.75	215.557861328125\\
33.125	215.435791388154\\
33.5	215.360832213628\\
33.875	215.240573883057\\
34.25	215.198492619675\\
34.625	215.136709809303\\
35	215.00244140625\\
35.375	214.937591546914\\
35.75	214.900207519554\\
36.125	214.762139321829\\
36.5	214.649002207437\\
36.875	214.648795116227\\
37.25	214.48669433603\\
37.625	214.439392085478\\
38	214.391708374023\\
38.375	214.294433593795\\
38.75	214.24330471491\\
39.125	214.21203021855\\
39.5	214.083862211555\\
39.875	214.0625\\
40.25	213.957203552309\\
40.625	213.892406225205\\
41	213.854980655014\\
41.375	213.769615441561\\
41.75	213.703149184585\\
42.125	213.623064756393\\
42.5	213.653159141904\\
42.875	213.549789786612\\
43.25	213.500213623047\\
43.625	213.446044921878\\
44	213.4765625\\
44.375	213.367748283963\\
44.75	213.335990905762\\
45.125	213.2143020659\\
45.5	213.232421875\\
45.875	213.125610351562\\
46.25	213.097953802162\\
46.625	213.040161551911\\
47	212.99438476567\\
47.375	212.915045022965\\
47.75	212.939548492432\\
48.125	212.841796875\\
48.5	212.792968750028\\
48.875	212.774658203125\\
49.25	212.744137598202\\
49.625	212.677002698183\\
50	212.6708984375\\
};
\addlegendentry{battery swapping}
\end{axis}
\end{tikzpicture}%
\vspace*{-0.3in}
\caption{Optimal number of charging/battery swapping stations $K$ as a function of $C$.}
\label{optimal_K_C}
\end{subfigure}
\hfill
\begin{subfigure}[b]{0.32\linewidth}
\centering
%
%
\begin{tikzpicture}

\begin{axis}[%
width=1.694in,
height=1.03in,
at={(1.358in,0.0in)},
scale only axis,
xmin=12.5,
xmax=50,
xtick={12.5,   25, 37.5,   50},
xlabel style={font=\color{white!15!black}},
xlabel={Battery capacity},
ymin=3.5,
ymax=15.5,
ylabel style={font=\color{white!15!black}},
ylabel={Number of chargers},
axis background/.style={fill=white},
legend style={at={(0,1)}, anchor=north west, legend cell align=left, align=left, font=\scriptsize, draw=white!12!black}
]
\addplot [color=black, line width=1.0pt]
  table[row sep=crcr]{%
12.5	4.91558865146944\\
12.875	5.02980689437378\\
13.25	5.14472674539727\\
13.625	5.25944141093362\\
14	5.37117368779145\\
14.375	5.48271426680217\\
14.75	5.59652240418143\\
15.125	5.7116905355551\\
15.5	5.82261645090117\\
15.875	5.93399863195265\\
16.25	6.03977693045022\\
16.625	6.15411618699982\\
17	6.26139122104463\\
17.375	6.3739829391017\\
17.75	6.47944518581844\\
18.125	6.58899104489949\\
18.5	6.6957938877891\\
18.875	6.8070119673982\\
19.25	6.91531953975527\\
19.625	7.0238483924036\\
20	7.13479411347044\\
20.375	7.23452918413321\\
20.75	7.34654919535083\\
21.125	7.45303083876651\\
21.5	7.55622960184179\\
21.875	7.6636757732604\\
22.25	7.76766583350997\\
22.625	7.87353525679391\\
23	7.98186814589125\\
23.375	8.08658223498536\\
23.75	8.19454788529572\\
24.125	8.29888506032888\\
24.5	8.40140018770341\\
24.875	8.50864914016548\\
25.25	8.6132669524587\\
25.625	8.71504575798975\\
26	8.82277043843409\\
26.375	8.92485347238528\\
26.75	9.03123812314949\\
27.125	9.13112623599197\\
27.5	9.2316557532507\\
27.875	9.33592116065821\\
28.25	9.44855898610361\\
28.625	9.54286443187495\\
29	9.65084471651469\\
29.375	9.75700638313066\\
29.75	9.85505686613076\\
30.125	9.95703881970985\\
30.5	10.0555872389538\\
30.875	10.1579932551236\\
31.25	10.2619177696821\\
31.625	10.3646856674814\\
32	10.4660845346854\\
32.375	10.5666468266429\\
32.75	10.66699519886\\
33.125	10.7704221087228\\
33.5	10.8669335872055\\
33.875	10.9701088947793\\
34.25	11.0761360944461\\
34.625	11.1779874978037\\
35	11.2710198885078\\
35.375	11.371807266173\\
35.75	11.4754876989255\\
36.125	11.5715413229365\\
36.5	11.6722513266328\\
36.875	11.7767765158299\\
37.25	11.8800447783648\\
37.625	11.9689886297276\\
38	12.073041794456\\
38.375	12.1738842692214\\
38.75	12.2811340702145\\
39.125	12.3825706098516\\
39.5	12.4785207668199\\
39.875	12.5723911564641\\
40.25	12.6662763105488\\
40.625	12.7698680782764\\
41	12.8815254564218\\
41.375	12.9721974029843\\
41.75	13.0716388115423\\
42.125	13.176012611331\\
42.5	13.2767435224463\\
42.875	13.3660427741283\\
43.25	13.4739131125743\\
43.625	13.5619444257016\\
44	13.6711676767059\\
44.375	13.7744021042146\\
44.75	13.8611275323196\\
45.125	13.9673338499482\\
45.5	14.0614719506091\\
45.875	14.1593783305716\\
46.25	14.2551415792003\\
46.625	14.3594658119138\\
47	14.4543607181447\\
47.375	14.5544849358444\\
47.75	14.6539108125188\\
48.125	14.7572253253391\\
48.5	14.8493760218262\\
48.875	14.9389764150255\\
49.25	15.0384017277331\\
49.625	15.1390991355117\\
50	15.2393999084055\\
};
\addlegendentry{plug-in charging}
\addplot [color=blue, dashed, line width=1.0pt]
  table[row sep=crcr]{%
12.5	6\\
12.875	6\\
13.25	6\\
13.625	6\\
14	6\\
14.375	6\\
14.75	6\\
15.125	6\\
15.5	6\\
15.875	6\\
16.25	6\\
16.625	6\\
17	6\\
17.375	6\\
17.75	6\\
18.125	6\\
18.5	6\\
18.875	6\\
19.25	6\\
19.625	6\\
20	6\\
20.375	6\\
20.75	6\\
21.125	6\\
21.5	6\\
21.875	6\\
22.25	6\\
22.625	6\\
23	6\\
23.375	6\\
23.75	6\\
24.125	6\\
24.5	6\\
24.875	6\\
25.25	6\\
25.625	6\\
26	6\\
26.375	6\\
26.75	6\\
27.125	6\\
27.5	6\\
27.875	6\\
28.25	6\\
28.625	6\\
29	6\\
29.375	6\\
29.75	6\\
30.125	6\\
30.5	6\\
30.875	6\\
31.25	6\\
31.625	6\\
32	6\\
32.375	6\\
32.75	6\\
33.125	6\\
33.5	6\\
33.875	6\\
34.25	6\\
34.625	6\\
35	6\\
35.375	6\\
35.75	6\\
36.125	6\\
36.5	6\\
36.875	6\\
37.25	6\\
37.625	6\\
38	6\\
38.375	6\\
38.75	6\\
39.125	6\\
39.5	6\\
39.875	6\\
40.25	6\\
40.625	6\\
41	6\\
41.375	6\\
41.75	6\\
42.125	6\\
42.5	6\\
42.875	6\\
43.25	6\\
43.625	6\\
44	6\\
44.375	6\\
44.75	6\\
45.125	6\\
45.5	6\\
45.875	6\\
46.25	6\\
46.625	6\\
47	6\\
47.375	6\\
47.75	6\\
48.125	6\\
48.5	6\\
48.875	6\\
49.25	6\\
49.625	6\\
50	6\\
};
\addlegendentry{battery swapping}
\end{axis}
\end{tikzpicture}%
\vspace*{-0.3in}
\caption{Optimal number of chargers at each station $Q$ as a function of $C$.}
\label{optimal_Q_C}
\end{subfigure}
\hfill
\begin{subfigure}[b]{0.32\linewidth}
\centering
%
%
\begin{tikzpicture}

\begin{axis}[%
width=1.694in,
height=1.03in,
at={(1.358in,0.0in)},
scale only axis,
xmin=12.5,
xmax=50,
xtick={12.5,   25, 37.5,   50},
xlabel style={font=\color{white!15!black}},
xlabel={Battery capacity},
ymin=60000,
ymax=125000,
ylabel style={font=\color{white!15!black}},
ylabel={Social welfare},
axis background/.style={fill=white},
legend style={at={(1,0)}, anchor=south east, legend cell align=left, align=left, font=\scriptsize, draw=white!12!black}
]
\addplot [color=black, line width=1.0pt]
  table[row sep=crcr]{%
12.5	112286.526647929\\
12.875	112597.795759581\\
13.25	112895.684416447\\
13.625	113181.122046835\\
14	113454.953357958\\
14.375	113717.941512164\\
14.75	113970.786229765\\
15.125	114214.12310082\\
15.5	114448.534412383\\
15.875	114674.545542141\\
16.25	114892.655295777\\
16.625	115103.308599853\\
17	115306.919473413\\
17.375	115503.875854097\\
17.75	115694.53119897\\
18.125	115879.213164407\\
18.5	116058.229006935\\
18.875	116231.861470567\\
19.25	116400.376075914\\
19.625	116564.018839166\\
20	116723.021566254\\
20.375	116877.593231517\\
20.75	117027.947868961\\
21.125	117174.265687212\\
21.5	117316.724579061\\
21.875	117455.4939826\\
22.25	117590.725531708\\
22.625	117722.575079561\\
23	117851.17388625\\
23.375	117976.659368701\\
23.75	118099.151074232\\
24.125	118218.767372663\\
24.5	118335.622717764\\
24.875	118449.817037712\\
25.25	118561.449736617\\
25.625	118670.618773127\\
26	118777.413893998\\
26.375	118881.914916503\\
26.75	118984.20651446\\
27.125	119084.367015278\\
27.5	119182.466462037\\
27.875	119278.576656782\\
28.25	119372.760988954\\
28.625	119465.08819891\\
29	119555.614602586\\
29.375	119644.402112637\\
29.75	119731.496915812\\
30.125	119816.963002508\\
30.5	119900.843211968\\
30.875	119983.18902748\\
31.25	120064.04724366\\
31.625	120143.462719089\\
32	120221.474727587\\
32.375	120298.124186079\\
32.75	120373.455310842\\
33.125	120447.500892832\\
33.5	120520.299342378\\
33.875	120591.884858875\\
34.25	120662.291132465\\
34.625	120731.548436746\\
35	120799.690383817\\
35.375	120866.744681453\\
35.75	120932.742255311\\
36.125	120997.705129896\\
36.5	121061.668097321\\
36.875	121124.653209184\\
37.25	121186.680913524\\
37.625	121247.781407469\\
38	121307.975106553\\
38.375	121367.28300615\\
38.75	121425.729456423\\
39.125	121483.332724774\\
39.5	121540.116030355\\
39.875	121596.093376083\\
40.25	121651.29066368\\
40.625	121705.720538982\\
41	121759.402166894\\
41.375	121812.352192074\\
41.75	121864.591064879\\
42.125	121916.125884841\\
42.5	121966.980069069\\
42.875	122017.166185932\\
43.25	122066.698648941\\
43.625	122115.592752705\\
44	122163.859710525\\
44.375	122211.51359304\\
44.75	122258.570712884\\
45.125	122305.035361261\\
45.5	122350.932635308\\
45.875	122396.260860291\\
46.25	122441.038905079\\
46.625	122485.278340642\\
47	122528.984861915\\
47.375	122572.172847446\\
47.75	122614.851526482\\
48.125	122657.033582096\\
48.5	122698.724752354\\
48.875	122739.937390763\\
49.25	122780.675037258\\
49.625	122820.955919748\\
50	122860.780855432\\
};
\addlegendentry{plug-in charging}
\addplot [color=blue, dashed, line width=1.0pt]
  table[row sep=crcr]{%
12.5	81869.4669424178\\
12.875	82463.8676243452\\
13.25	83023.8753776753\\
13.625	83552.4040805363\\
14	84052.0331390207\\
14.375	84525.0824961631\\
14.75	84973.609920222\\
15.125	85399.4817540285\\
15.5	85804.3744533088\\
15.875	86189.8060711067\\
16.25	86557.1468781782\\
16.625	86907.6414794176\\
17	87242.4237136421\\
17.375	87562.5339062995\\
17.75	87868.9089047162\\
18.125	88162.4166754727\\
18.5	88443.8584721637\\
18.875	88713.9603235244\\
19.25	88973.3922407183\\
19.625	89222.7797715036\\
20	89462.6891854509\\
20.375	89693.6560368877\\
20.75	89916.1730588212\\
21.125	90130.6831189521\\
21.5	90337.6318178601\\
21.875	90537.4024695469\\
22.25	90730.355909066\\
22.625	90916.8414057153\\
23	91097.1798338302\\
23.375	91271.6753649721\\
23.75	91440.5961965313\\
24.125	91604.2044934054\\
24.5	91762.7634810833\\
24.875	91916.4984371856\\
25.25	92065.6020546963\\
25.625	92210.3176557457\\
26	92350.8164379066\\
26.375	92487.2808787247\\
26.75	92619.8795691146\\
27.125	92748.7796438049\\
27.5	92874.1408585982\\
27.875	92996.0886533077\\
28.25	93114.7743690542\\
28.625	93230.3258160565\\
29	93342.8645925976\\
29.375	93452.5039531027\\
29.75	93559.3603189911\\
30.125	93663.5311787265\\
30.5	93765.1214632882\\
30.875	93864.2318568358\\
31.25	93960.9345500299\\
31.625	94055.3351849726\\
32	94147.4962288071\\
32.375	94237.5207502104\\
32.75	94325.4571257693\\
33.125	94411.388076499\\
33.5	94495.3857045273\\
33.875	94577.5128736615\\
34.25	94657.8235492755\\
34.625	94736.3816825363\\
35	94813.2536035982\\
35.375	94888.4766892725\\
35.75	94962.1163370091\\
36.125	95034.2236922735\\
36.5	95104.8277507728\\
36.875	95173.9880492196\\
37.25	95241.7478676663\\
37.625	95308.1559132463\\
38	95373.2400053555\\
38.375	95437.0452156582\\
38.75	95499.6128800798\\
39.125	95560.9614637826\\
39.5	95621.1451968819\\
39.875	95680.1910526352\\
40.25	95738.1291006949\\
40.625	95794.9956277328\\
41	95850.8118802621\\
41.375	95905.6101846245\\
41.75	95959.425868167\\
42.125	96012.265133705\\
42.5	96064.1700476538\\
42.875	96115.1749164652\\
43.25	96165.2791292191\\
43.625	96214.5202130638\\
44	96262.9162816255\\
44.375	96310.4949740149\\
44.75	96357.267296477\\
45.125	96403.260788669\\
45.5	96448.4885597423\\
45.875	96492.9764663824\\
46.25	96536.7409726199\\
46.625	96579.7975033975\\
47	96622.160794845\\
47.375	96663.8493611881\\
47.75	96704.8832437677\\
48.125	96745.2735653855\\
48.5	96785.0403783722\\
48.875	96824.1889132689\\
49.25	96862.7391379624\\
49.625	96900.7062378859\\
50	96938.0987801086\\
};
\addlegendentry{battery swapping}
\end{axis}
\end{tikzpicture}%
\vspace*{-0.3in}
\caption{Social welfare (\$/hour) as a function of $C$.}
\label{SW_C}
\end{subfigure}
\begin{subfigure}[b]{0.32\linewidth}
\centering
%
%
\begin{tikzpicture}

\begin{axis}[%
width=1.694in,
height=1.03in,
at={(1.358in,0.0in)},
scale only axis,
xmin=12.5,
xmax=50,
xtick={12.5,   25, 37.5,   50},
xlabel style={font=\color{white!15!black}},
xlabel={Battery capacity},
ymin=89000,
ymax=96000,
ylabel style={font=\color{white!15!black}},
ylabel={Passengers' surplus},
axis background/.style={fill=white},
legend style={at={(1,0)}, anchor=south east, legend cell align=left, align=left, font=\scriptsize, draw=white!12!black}
]
\addplot [color=black, line width=1.0pt]
  table[row sep=crcr]{%
12.5	89569.3492682826\\
12.875	89677.0118603993\\
13.25	89781.0349468484\\
13.625	89880.3831652041\\
14	89975.2471491001\\
14.375	90065.2045705828\\
14.75	90151.6053082942\\
15.125	90238.9316369338\\
15.5	90319.3548051106\\
15.875	90398.606246715\\
16.25	90472.7779908093\\
16.625	90547.3393112259\\
17	90616.5001015896\\
17.375	90686.6115034994\\
17.75	90752.4741867196\\
18.125	90817.5816502726\\
18.5	90877.9186151492\\
18.875	90937.5089312531\\
19.25	90994.7489476545\\
19.625	91052.475778044\\
20	91108.6889009802\\
20.375	91164.3442819316\\
20.75	91213.1967998495\\
21.125	91265.7868272042\\
21.5	91315.439212729\\
21.875	91361.9696830352\\
22.25	91408.3596304816\\
22.625	91455.4468633001\\
23	91500.5037280635\\
23.375	91542.7961904276\\
23.75	91586.3316265609\\
24.125	91625.1555854988\\
24.5	91667.872785382\\
24.875	91707.3752914519\\
25.25	91745.1703445896\\
25.625	91782.282983794\\
26	91820.1902820265\\
26.375	91857.6997731825\\
26.75	91891.0209277877\\
27.125	91926.2557508016\\
27.5	91959.2727153135\\
27.875	91992.9233575761\\
28.25	92026.3074220255\\
28.625	92057.4041283874\\
29	92089.6171079361\\
29.375	92119.2533412077\\
29.75	92148.2899325734\\
30.125	92179.0727298685\\
30.5	92209.444721243\\
30.875	92238.2125522294\\
31.25	92264.9612575537\\
31.625	92291.749956195\\
32	92319.218772608\\
32.375	92343.1323545002\\
32.75	92370.3166681013\\
33.125	92395.1656150898\\
33.5	92421.8007242955\\
33.875	92446.1495257294\\
34.25	92470.813636711\\
34.625	92495.1247776513\\
35	92518.5358800733\\
35.375	92541.2606096558\\
35.75	92563.6735493696\\
36.125	92587.9533850634\\
36.5	92609.385522991\\
36.875	92630.4447220338\\
37.25	92651.0465195422\\
37.625	92673.8821097147\\
38	92695.3165381732\\
38.375	92714.2776726909\\
38.75	92733.9339104637\\
39.125	92753.1513705774\\
39.5	92773.8012298724\\
39.875	92792.6222402948\\
40.25	92812.9930600652\\
40.625	92831.6865241909\\
41	92848.79721063\\
41.375	92867.8634557855\\
41.75	92886.6513825067\\
42.125	92903.6878571856\\
42.5	92921.3554805883\\
42.875	92937.5625094499\\
43.25	92956.678073594\\
43.625	92972.6137400648\\
44	92989.0149718189\\
44.375	93005.0275675217\\
44.75	93021.8338339232\\
45.125	93035.5390475378\\
45.5	93054.890567516\\
45.875	93069.8171320116\\
46.25	93084.7884603088\\
46.625	93100.1306264367\\
47	93114.971096416\\
47.375	93129.7095037708\\
47.75	93144.0723229072\\
48.125	93157.5304551983\\
48.5	93174.2498195483\\
48.875	93188.3169866392\\
49.25	93200.3131287764\\
49.625	93216.449116796\\
50	93229.6998976844\\
};
\addlegendentry{plug-in charging}
\addplot [color=blue, dashed, line width=1.0pt]
  table[row sep=crcr]{%
12.5	90513.2933738779\\
12.875	90716.7070696337\\
13.25	90915.9038564736\\
13.625	91090.4367254624\\
14	91258.3006677284\\
14.375	91416.581370832\\
14.75	91570.7019034887\\
15.125	91708.859984208\\
15.5	91853.0404996501\\
15.875	91990.429020572\\
16.25	92110.2522609412\\
16.625	92225.5084781505\\
17	92338.5234095121\\
17.375	92445.5316255049\\
17.75	92549.1625462146\\
18.125	92658.4948231656\\
18.5	92746.598338139\\
18.875	92832.7885127869\\
19.25	92925.0375849136\\
19.625	93013.4320312859\\
20	93088.9851072971\\
20.375	93172.2924059604\\
20.75	93240.2698087827\\
21.125	93309.4596346632\\
21.5	93391.7456902135\\
21.875	93457.3646208095\\
22.25	93514.4043399433\\
22.625	93582.9286311892\\
23	93634.7480564521\\
23.375	93702.3862585572\\
23.75	93758.8160906118\\
24.125	93820.5893176963\\
24.5	93863.007190566\\
24.875	93915.2624557642\\
25.25	93976.0427836084\\
25.625	94013.563294185\\
26	94062.3301773334\\
26.375	94112.155617\\
26.75	94148.143406004\\
27.125	94196.4582386629\\
27.5	94237.1569845072\\
27.875	94284.9345329567\\
28.25	94311.8896838576\\
28.625	94355.1811375678\\
29	94393.4272980262\\
29.375	94427.8024873333\\
29.75	94467.2036396715\\
30.125	94503.6853444979\\
30.5	94534.8423697327\\
30.875	94565.4648662443\\
31.25	94603.6017264246\\
31.625	94630.6547795456\\
32	94662.1814022984\\
32.375	94696.3991932041\\
32.75	94724.0552649127\\
33.125	94747.7649141399\\
33.5	94779.2292085032\\
33.875	94801.2052037339\\
34.25	94836.808773887\\
34.625	94867.7751728972\\
35	94884.2640660421\\
35.375	94912.8398333487\\
35.75	94945.6987566676\\
36.125	94958.8778638022\\
36.5	94975.8994072815\\
36.875	95013.2797343019\\
37.25	95019.582749522\\
37.625	95046.6502182284\\
38	95072.9253393233\\
38.375	95089.2152587205\\
38.75	95113.4792711522\\
39.125	95140.7968536353\\
39.5	95149.3240085163\\
39.875	95177.2510629806\\
40.25	95188.8415551926\\
40.625	95207.4430049461\\
41	95230.6039353045\\
41.375	95244.2220097564\\
41.75	95260.8528965315\\
42.125	95274.3928566939\\
42.5	95308.094626696\\
42.875	95316.2751381513\\
43.25	95334.0538413704\\
43.625	95350.5032877101\\
44	95382.3908893157\\
44.375	95387.6802229764\\
44.75	95406.9963319377\\
45.125	95409.0057898592\\
45.5	95436.8472332636\\
45.875	95440.8297546637\\
46.25	95459.2840691826\\
46.625	95471.6945595284\\
47	95485.9839525825\\
47.375	95493.6025451616\\
47.75	95520.3742060341\\
48.125	95523.8283180752\\
48.5	95536.1373557832\\
48.875	95553.8426645341\\
49.25	95568.9300914699\\
49.625	95576.8208239903\\
50	95595.8695995834\\
};
\addlegendentry{battery swapping}
\end{axis}
\end{tikzpicture}%
\vspace*{-0.3in}
\caption{Passengers' surplus (\$/hour) as a function of $C$.} 
\label{PS_C}
\end{subfigure}
\hfill
\begin{subfigure}[b]{0.32\linewidth}
\centering
%
%
\begin{tikzpicture}

\begin{axis}[%
width=1.694in,
height=1.03in,
at={(1.358in,0.0in)},
scale only axis,
xmin=12.5,
xmax=50,
xtick={12.5,   25, 37.5,   50},
xlabel style={font=\color{white!15!black}},
xlabel={Battery capacity},
ymin=36000,
ymax=54000,
ylabel style={font=\color{white!15!black}},
ylabel={Platform profit},
axis background/.style={fill=white},
legend style={at={(1,0.3)}, anchor=south east, legend cell align=left, align=left, font=\scriptsize, draw=white!12!black}
]
\addplot [color=black, line width=1.0pt]
  table[row sep=crcr]{%
12.5	37642.1611053367\\
12.875	37735.2573294918\\
13.25	37824.2367970244\\
13.625	37909.1873525461\\
14	37991.3881359296\\
14.375	38069.6809012591\\
14.75	38143.9062594527\\
15.125	38216.4718539323\\
15.5	38286.0602123972\\
15.875	38353.5570804169\\
16.25	38419.3463421842\\
16.625	38481.2793239489\\
17	38541.9478675011\\
17.375	38600.2400754307\\
17.75	38657.7733923622\\
18.125	38712.751984107\\
18.5	38765.5641947034\\
18.875	38815.9284062034\\
19.25	38865.2417060002\\
19.625	38913.8233674856\\
20	38960.3779086446\\
20.375	39008.676337794\\
20.75	39050.5638347536\\
21.125	39094.5089841186\\
21.5	39137.3656056694\\
21.875	39177.3771678031\\
22.25	39217.4618072698\\
22.625	39256.7637275268\\
23	39294.1838664643\\
23.375	39330.7648250338\\
23.75	39366.5910638061\\
24.125	39400.8847271765\\
24.5	39436.4490846573\\
24.875	39469.4119688593\\
25.25	39501.7892630057\\
25.625	39533.9849363505\\
26	39564.9242946065\\
26.375	39596.3169547426\\
26.75	39625.0285025499\\
27.125	39655.1573651553\\
27.5	39683.9800179903\\
27.875	39712.0099300788\\
28.25	39738.2094306831\\
28.625	39766.2669446342\\
29	39792.1498433834\\
29.375	39817.1152680272\\
29.75	39842.8233669412\\
30.125	39868.1985092061\\
30.5	39893.6370438104\\
30.875	39917.6452264091\\
31.25	39940.4622594156\\
31.625	39963.1894771678\\
32	39986.0724084041\\
32.375	40007.6155776478\\
32.75	40030.0229133469\\
33.125	40050.9816968315\\
33.5	40073.2188276223\\
33.875	40093.5687089801\\
34.25	40113.4296483952\\
34.625	40133.5015328361\\
35	40154.176037301\\
35.375	40173.4227905349\\
35.75	40192.0057929439\\
36.125	40211.9248964164\\
36.5	40230.154644147\\
36.875	40247.6229587441\\
37.25	40264.9189784835\\
37.625	40284.4031801052\\
38	40301.5069215073\\
38.375	40318.026761438\\
38.75	40333.8864646458\\
39.125	40350.1010527729\\
39.5	40367.2005057504\\
39.875	40383.7826698911\\
40.25	40400.6733294805\\
40.625	40415.8592617052\\
41	40429.582119598\\
41.375	40445.9419023099\\
41.75	40461.1622735867\\
42.125	40475.2145115965\\
42.5	40489.6950082895\\
42.875	40504.7211310335\\
43.25	40518.6780498757\\
43.625	40533.4633136881\\
44	40546.2331049047\\
44.375	40559.3464131663\\
44.75	40574.1173944416\\
45.125	40586.0330538435\\
45.5	40600.6247162858\\
45.875	40613.439446651\\
46.25	40626.3462639756\\
46.625	40638.4939809168\\
47	40651.221031231\\
47.375	40663.3555521918\\
47.75	40675.3404150965\\
48.125	40686.633016951\\
48.5	40699.7206122232\\
48.875	40712.1398778022\\
49.25	40723.0452809141\\
49.625	40734.95024639\\
50	40745.967541134\\
};
\addlegendentry{plug-in charging}
\addplot [color=blue, dashed, line width=1.0pt]
  table[row sep=crcr]{%
12.5	46445.2834661534\\
12.875	46674.612630938\\
13.25	46893.4968431618\\
13.625	47096.2446958554\\
14	47289.1727755714\\
14.375	47471.9540367053\\
14.75	47646.8045011083\\
15.125	47810.5558518517\\
15.5	47970.4226987977\\
15.875	48122.8389645972\\
16.25	48264.8877785055\\
16.625	48400.8830012671\\
17	48531.8299912807\\
17.375	48656.8584166398\\
17.75	48777.0214809057\\
18.125	48895.2792741821\\
18.5	49003.7481599585\\
18.875	49108.4437522279\\
19.25	49211.9288633788\\
19.625	49311.3985214679\\
20	49404.3444425346\\
20.375	49497.0380217202\\
20.75	49582.763601601\\
21.125	49666.5603089444\\
21.5	49752.0384713966\\
21.875	49830.5261299877\\
22.25	49904.5273137515\\
22.625	49980.0308547696\\
23	50048.8518372807\\
23.375	50120.604412635\\
23.75	50187.4563754945\\
24.125	50254.3183006204\\
24.5	50314.0629311368\\
24.875	50375.2594188996\\
25.25	50437.5817320254\\
25.625	50491.8593420294\\
26	50548.1606904666\\
26.375	50603.5798515685\\
26.75	50653.8649799973\\
27.125	50706.618302491\\
27.5	50756.1240966723\\
27.875	50806.6626957905\\
28.25	50850.2465355859\\
28.625	50897.5970776232\\
29	50942.6027875829\\
29.375	50985.6389657694\\
29.75	51029.284757523\\
30.125	51071.2962201888\\
30.5	51111.0163982867\\
30.875	51149.8412093443\\
31.25	51190.0964954803\\
31.625	51226.4992530847\\
32	51263.513068347\\
32.375	51300.6430413812\\
32.75	51335.2885796066\\
33.125	51368.2130955161\\
33.5	51402.7562272948\\
33.875	51434.0454018612\\
34.25	51468.6530041104\\
34.625	51501.4168638719\\
35	51529.5754750562\\
35.375	51560.6588271831\\
35.75	51592.4673850344\\
36.125	51618.2592657103\\
36.5	51644.6888732763\\
36.875	51676.4191428122\\
37.25	51698.9717587914\\
37.625	51726.9597955327\\
38	51754.3246757978\\
38.375	51778.4940194486\\
38.75	51804.526740506\\
39.125	51831.0518625993\\
39.5	51851.9481191388\\
39.875	51877.9399896546\\
40.25	51899.0163980564\\
40.625	51921.730116836\\
41	51945.403302161\\
41.375	51966.0958808427\\
41.75	51987.3287759358\\
42.125	52007.4078185456\\
42.5	52032.8336150149\\
42.875	52050.8493271007\\
43.25	52071.2765573799\\
43.625	52091.0677066044\\
44	52114.9003923098\\
44.375	52131.0743391896\\
44.75	52150.9087819221\\
45.125	52165.6874946258\\
45.5	52187.4225764788\\
45.875	52202.2931960937\\
46.25	52220.9658159562\\
46.625	52237.7417163276\\
47	52254.8291860234\\
47.375	52269.857621538\\
47.75	52290.0006759173\\
48.125	52303.4764973103\\
48.5	52319.2155225959\\
48.875	52336.2642174847\\
49.25	52352.4020700609\\
49.625	52366.3660614595\\
50	52383.2448055252\\
};
\addlegendentry{battery swapping}
\end{axis}
\end{tikzpicture}%
\vspace*{-0.3in}
\caption{Platform profit (\$/hour) as a function of $C$.}
\label{profit_C}
\end{subfigure}
\hfill
\begin{subfigure}[b]{0.32\linewidth}
\centering
%
%
\begin{tikzpicture}

\begin{axis}[%
width=1.694in,
height=1.03in,
at={(1.358in,0.0in)},
scale only axis,
xmin=12.5,
xmax=50,
xtick={12.5,   25, 37.5,   50},
xlabel style={font=\color{white!15!black}},
xlabel={Battery capacity},
ymin=10000,
ymax=60000,
ylabel style={font=\color{white!15!black}},
ylabel={Infrastructure cost},
axis background/.style={fill=white},
legend style={at={(1,0.3)}, anchor=south east, legend cell align=left, align=left, font=\scriptsize, draw=white!12!black}
]
\addplot [color=black, line width=1.0pt]
  table[row sep=crcr]{%
12.5	14924.98372569\\
12.875	14814.47343031\\
13.25	14709.5873274259\\
13.625	14608.4484709154\\
14	14511.6819270716\\
14.375	14416.9439596778\\
14.75	14324.7253379816\\
15.125	14241.2803900459\\
15.5	14156.8806051243\\
15.875	14077.6177849908\\
16.25	13999.4690372165\\
16.625	13925.3100353215\\
17	13851.5284956775\\
17.375	13782.9757248333\\
17.75	13715.7163801118\\
18.125	13651.120469973\\
18.5	13585.2538029179\\
18.875	13521.57586689\\
19.25	13459.6145777407\\
19.625	13402.2803063635\\
20	13346.0452433713\\
20.375	13295.4273882082\\
20.75	13235.8127656417\\
21.125	13186.0301241112\\
21.5	13136.0802393375\\
21.875	13083.8528682385\\
22.25	13035.0959060431\\
22.625	12989.6355112658\\
23	12943.5137082775\\
23.375	12896.9016467602\\
23.75	12853.7716161351\\
24.125	12807.272940012\\
24.5	12768.6991522754\\
24.875	12726.9702225994\\
25.25	12685.5098709778\\
25.625	12645.6491470178\\
26	12607.7006826353\\
26.375	12572.1018114223\\
26.75	12531.842915878\\
27.125	12497.0461006787\\
27.5	12460.7862712668\\
27.875	12426.3566308732\\
28.25	12391.7558637543\\
28.625	12358.5828741119\\
29	12326.1523487339\\
29.375	12291.9664965977\\
29.75	12259.6163837025\\
30.125	12230.3082365661\\
30.5	12202.238553085\\
30.875	12172.6687511585\\
31.25	12141.3762733092\\
31.625	12111.4767142736\\
32	12083.8164534253\\
32.375	12052.6237460689\\
32.75	12026.8842706066\\
33.125	11998.6464190891\\
33.5	11974.7202095401\\
33.875	11947.8333758346\\
34.25	11921.952152641\\
34.625	11897.0778737409\\
35	11873.0215335572\\
35.375	11847.9387187382\\
35.75	11822.9370870026\\
36.125	11802.1731515833\\
36.5	11777.8720698174\\
36.875	11753.4144715943\\
37.25	11729.2845845015\\
37.625	11710.5038823508\\
38	11688.848353128\\
38.375	11665.0214279788\\
38.75	11642.0909186861\\
39.125	11619.9196985759\\
39.5	11600.8857052673\\
39.875	11580.3115341032\\
40.25	11562.3757258659\\
40.625	11541.8252469141\\
41	11518.9771633345\\
41.375	11501.4531660212\\
41.75	11483.2225912141\\
42.125	11462.7764839412\\
42.5	11444.0704198087\\
42.875	11425.1174545517\\
43.25	11408.6574745284\\
43.625	11390.4843010476\\
44	11371.388366199\\
44.375	11352.8603876485\\
44.75	11337.3805154809\\
45.125	11316.5367401206\\
45.5	11304.5826484938\\
45.875	11286.995718372\\
46.25	11270.0958192059\\
46.625	11253.3462667114\\
47	11237.2072657318\\
47.375	11220.8922085168\\
47.75	11204.5612115217\\
48.125	11187.1298900536\\
48.5	11175.2456794178\\
48.875	11160.5194736789\\
49.25	11142.6833724327\\
49.625	11130.4434434381\\
50	11114.8865833861\\
};
\addlegendentry{plug-in charging}
\addplot [color=blue, dashed, line width=1.0pt]
  table[row sep=crcr]{%
12.5	55089.1098976135\\
12.875	54927.4520762265\\
13.25	54785.5253219601\\
13.625	54634.2773407814\\
14	54495.4403042792\\
14.375	54363.4529113742\\
14.75	54243.896484375\\
15.125	54119.9340820312\\
15.5	54019.088745139\\
15.875	53923.4619140625\\
16.25	53817.9931612685\\
16.625	53718.75\\
17	53627.9296871508\\
17.375	53539.8561358452\\
17.75	53457.2751224041\\
18.125	53391.357421875\\
18.5	53306.4880259339\\
18.875	53227.2719414905\\
19.25	53163.5742075741\\
19.625	53102.0507812501\\
20	53030.6403643808\\
20.375	52975.6743907928\\
20.75	52906.8603515625\\
21.125	52845.3368246555\\
21.5	52806.15234375\\
21.875	52750.4882812503\\
22.25	52688.5757446289\\
22.625	52646.1180802435\\
23	52586.4200599026\\
23.375	52551.3153062202\\
23.75	52505.6762695749\\
24.125	52470.7031249113\\
24.5	52414.3066406195\\
24.875	52374.0234374782\\
25.25	52348.0224609375\\
25.625	52295.1049804688\\
26	52259.6744298933\\
26.375	52228.4545898438\\
26.75	52182.1288168867\\
27.125	52154.296897349\\
27.5	52119.1402225813\\
27.875	52095.5085754395\\
28.25	52047.3618503893\\
28.625	52022.4523991345\\
29	51993.1654930115\\
29.375	51960.9375\\
29.75	51937.1280782034\\
30.125	51911.4503859602\\
30.5	51880.7373047312\\
30.875	51851.0742187527\\
31.25	51832.763671875\\
31.625	51801.8188476576\\
32	51778.1982418383\\
32.375	51759.5214843749\\
32.75	51733.88671875\\
33.125	51704.589933157\\
33.5	51686.5997312707\\
33.875	51657.7377319336\\
34.25	51647.6382287219\\
34.625	51632.8103542328\\
35	51600.5859375\\
35.375	51585.0219712593\\
35.75	51576.049804693\\
36.125	51542.9134372389\\
36.5	51515.760529785\\
36.875	51515.7108278945\\
37.25	51476.8066406472\\
37.625	51465.4541005148\\
38	51454.0100097656\\
38.375	51430.6640625109\\
38.75	51418.3931315783\\
39.125	51410.8872524521\\
39.5	51380.1269307733\\
39.875	51375\\
40.25	51349.7288525541\\
40.625	51334.1774940493\\
41	51325.1953572034\\
41.375	51304.7077059746\\
41.75	51288.7558043003\\
42.125	51269.5355415344\\
42.5	51276.758194057\\
42.875	51251.9495487868\\
43.25	51240.0512695312\\
43.625	51227.0507812506\\
44	51234.375\\
44.375	51208.2595881511\\
44.75	51200.6378173828\\
45.125	51171.432495816\\
45.5	51175.78125\\
45.875	51150.146484375\\
46.25	51143.5089125189\\
46.625	51129.6387724585\\
47	51118.6523437609\\
47.375	51099.6108055115\\
47.75	51105.4916381837\\
48.125	51082.03125\\
48.5	51070.3125000068\\
48.875	51065.91796875\\
49.25	51058.5930235684\\
49.625	51042.4806475639\\
50	51041.015625\\
};
\addlegendentry{battery swapping}
\end{axis}
\end{tikzpicture}%
\vspace*{-0.3in}
\caption{Infrastructure cost (\$/hour) as a function of $C$.}
\label{infrastructure_cost_C}
\end{subfigure}
\caption{Charging infrastructure planning under different battery capacities. Black lines present the results when deploying charging stations; blue dashed lines show the results when deploying battery swapping stations.}
\label{infrastructure_planning_C}
\end{figure}

Figure \ref{infrastructure_planning_C} presents the infrastructure planning of charging stations and battery swapping stations under different battery capacities. Numerical results show that for both charging strategies, increased battery capacity will lead to a larger passenger surplus (Figure \ref{PS_C}), larger platform profit (Figure \ref{profit_C}), smaller infrastructure cost (Figure \ref{infrastructure_cost_C}), and large social welfare (Figure \ref{SW_C}).  However, for plug-in charging, increased battery capacity leads to a smaller number of charging stations (Figure \ref{optimal_K_C}) but a larger number of chargers at each station (Figure \ref{optimal_Q_C}), whereas for battery swapping, increased battery capacity leads to a reduced number of battery swapping stations (Figure \ref{optimal_K_C}).

\begin{figure}[h]
\centering
\begin{subfigure}[b]{0.32\linewidth}
\centering
%
%
\begin{tikzpicture}

\begin{axis}[%
width=1.694in,
height=1.03in,
at={(1.358in,0.0in)},
scale only axis,
xmin=12.5,
xmax=50,
xtick={12.5,   25, 37.5,   50},
xlabel style={font=\color{white!15!black}},
xlabel={Battery capacity},
ymin=204,
ymax=218,
ylabel style={font=\color{white!15!black}},
ylabel={Passenger arrivals},
axis background/.style={fill=white},
legend style={at={(1,0)}, anchor=south east, legend cell align=left, align=left, font=\scriptsize, draw=white!12!black}
]
\addplot [color=black, line width=1.0pt]
  table[row sep=crcr]{%
12.5	205.211040584823\\
12.875	205.428676069012\\
13.25	205.638893544925\\
13.625	205.839607843997\\
14	206.031211671442\\
14.375	206.212859408043\\
14.75	206.387283167785\\
15.125	206.563533587223\\
15.5	206.725814129772\\
15.875	206.885695388592\\
16.25	207.035297495141\\
16.625	207.18565476271\\
17	207.32509412128\\
17.375	207.466423131835\\
17.75	207.599162970542\\
18.125	207.730357224577\\
18.5	207.851917857825\\
18.875	207.971954521796\\
19.25	208.087238393123\\
19.625	208.203484441019\\
20	208.316664649591\\
20.375	208.428704742296\\
20.75	208.527035899335\\
21.125	208.632875285437\\
21.5	208.732788579498\\
21.875	208.826407460336\\
22.25	208.919731744259\\
22.625	209.014446666426\\
23	209.105066122361\\
23.375	209.190115579678\\
23.75	209.277654356189\\
24.125	209.35571074413\\
24.5	209.441584993378\\
24.875	209.52098782065\\
25.25	209.596950512266\\
25.625	209.671534003881\\
26	209.74770666419\\
26.375	209.823072175221\\
26.75	209.890015831212\\
27.125	209.960797490729\\
27.5	210.027117599716\\
27.875	210.094704391116\\
28.25	210.161749615256\\
28.625	210.224195627874\\
29	210.288877651009\\
29.375	210.348380662969\\
29.75	210.40667504625\\
30.125	210.468470085825\\
30.5	210.529435349907\\
30.875	210.587175933189\\
31.25	210.640859804216\\
31.625	210.694620002331\\
32	210.749740984642\\
32.375	210.797724393705\\
32.75	210.852266822754\\
33.125	210.902120037994\\
33.5	210.955552978061\\
33.875	211.004395927418\\
34.25	211.053868060157\\
34.625	211.102628925603\\
35	211.149581517413\\
35.375	211.195154663826\\
35.75	211.240099757516\\
36.125	211.288785462579\\
36.5	211.331758302137\\
36.875	211.3739809199\\
37.25	211.415284116435\\
37.625	211.461062976853\\
38	211.504030303943\\
38.375	211.54203757996\\
38.75	211.581436096926\\
39.125	211.619953089114\\
39.5	211.661338739356\\
39.875	211.699057044209\\
40.25	211.739879060165\\
40.625	211.77733774428\\
41	211.811623117745\\
41.375	211.849825032126\\
41.75	211.887467346471\\
42.125	211.921598887991\\
42.5	211.956993216867\\
42.875	211.989459967424\\
43.25	212.027751394144\\
43.625	212.059671463865\\
44	212.092522634439\\
44.375	212.124593956325\\
44.75	212.158253394548\\
45.125	212.185700926399\\
45.5	212.224454604514\\
45.875	212.254345390662\\
46.25	212.284324590899\\
46.625	212.315045097738\\
47	212.344759802151\\
47.375	212.374268954773\\
47.75	212.403024961814\\
48.125	212.429968656603\\
48.5	212.463440080399\\
48.875	212.491600736953\\
49.25	212.515614611911\\
49.625	212.547914386202\\
50	212.574437709228\\
};
\addlegendentry{plug-in charging}
\addplot [color=blue, dashed, line width=1.0pt]
  table[row sep=crcr]{%
12.5	207.117002942979\\
12.875	207.527080789683\\
13.25	207.928436271176\\
13.625	208.279917398685\\
14	208.617810040845\\
14.375	208.936270471706\\
14.75	209.246228127022\\
15.125	209.523971990588\\
15.5	209.81371100458\\
15.875	210.089694767625\\
16.25	210.33030899173\\
16.625	210.561677850813\\
17	210.788476634329\\
17.375	211.003156478754\\
17.75	211.211000881912\\
18.125	211.430216160159\\
18.5	211.606819265811\\
18.875	211.779545891703\\
19.25	211.96436952331\\
19.625	212.141426533784\\
20	212.292727908107\\
20.375	212.459521504607\\
20.75	212.595594200406\\
21.125	212.734067853157\\
21.5	212.898717628472\\
21.875	213.029990853435\\
22.25	213.144081882057\\
22.625	213.281120842182\\
23	213.384735352907\\
23.375	213.519957882944\\
23.75	213.632753450422\\
24.125	213.756209769644\\
24.5	213.840971537476\\
24.875	213.945377407666\\
25.25	214.066797547671\\
25.625	214.141741745574\\
26	214.239138167003\\
26.375	214.338635296836\\
26.75	214.410491387373\\
27.125	214.506949509586\\
27.5	214.588192603547\\
27.875	214.683554902692\\
28.25	214.737350930523\\
28.625	214.82374200358\\
29	214.900056332326\\
29.375	214.968639918418\\
29.75	215.047243108811\\
30.125	215.120014601402\\
30.5	215.182158994321\\
30.875	215.243232077561\\
31.25	215.319284548152\\
31.625	215.373228894979\\
32	215.436088593418\\
32.375	215.504307969528\\
32.75	215.559440699758\\
33.125	215.606702864787\\
33.5	215.669418192119\\
33.875	215.713218028945\\
34.25	215.784173072905\\
34.625	215.8458809877\\
35	215.878736895362\\
35.375	215.935673699593\\
35.75	216.00113912712\\
36.125	216.027394446775\\
36.5	216.061303238416\\
36.875	216.135763429095\\
37.25	216.148318045413\\
37.625	216.202229729969\\
38	216.254559435006\\
38.375	216.287000656789\\
38.75	216.33531952884\\
39.125	216.389715380277\\
39.5	216.406694155805\\
39.875	216.462298111025\\
40.25	216.485374027023\\
40.625	216.5224067735\\
41	216.568514139464\\
41.375	216.595622806306\\
41.75	216.628727523925\\
42.125	216.655678465257\\
42.5	216.722756632761\\
42.875	216.739037744168\\
43.25	216.774420222161\\
43.625	216.807155732064\\
44	216.870610017351\\
44.375	216.881134917229\\
44.75	216.919569487903\\
45.125	216.923567725802\\
45.5	216.978961852948\\
45.875	216.986885245694\\
46.25	217.023599748437\\
46.625	217.04828913687\\
47	217.076715369614\\
47.375	217.091870764469\\
47.75	217.145124169944\\
48.125	217.151994704522\\
48.5	217.17647794194\\
48.875	217.211693162328\\
49.25	217.241700163781\\
49.625	217.2573933485\\
50	217.295276388434\\
};
\addlegendentry{battery swapping}
\end{axis}
\end{tikzpicture}%
\vspace*{-0.3in}
\caption{Passenger arrival rate $\lambda$ (per minute) as a function of $C$.}
\label{lambda_C}
\end{subfigure}
\hfill
\begin{subfigure}[b]{0.32\linewidth}
\centering
%
%
\begin{tikzpicture}

\begin{axis}[%
width=1.694in,
height=1.03in,
at={(1.358in,0.0in)},
scale only axis,
xmin=12.5,
xmax=50,
xtick={12.5,   25, 37.5,   50},
xlabel style={font=\color{white!15!black}},
xlabel={Battery capacity},
ymin=23.25,
ymax=23.75,
ylabel style={font=\color{white!15!black}},
ylabel={Passenger cost},
axis background/.style={fill=white},
legend style={at={(1,1)}, anchor=north east, legend cell align=left, align=left, font=\scriptsize, draw=white!12!black}
]
\addplot [color=black, line width=1.0pt]
  table[row sep=crcr]{%
12.5	23.7443444413029\\
12.875	23.7356050204246\\
13.25	23.7271698222007\\
13.625	23.7191217529655\\
14	23.7114442657264\\
14.375	23.704170465892\\
14.75	23.6971902823758\\
15.125	23.690141312788\\
15.5	23.6836548829138\\
15.875	23.6772679358535\\
16.25	23.6712948356973\\
16.625	23.6652947108554\\
17	23.6597330704215\\
17.375	23.6540988018839\\
17.75	23.6488094623162\\
18.125	23.6435840961157\\
18.5	23.6387445436989\\
18.875	23.633967656075\\
19.25	23.6293817687439\\
19.625	23.6247594516709\\
20	23.6202608198479\\
20.375	23.6158092301403\\
20.75	23.6119037381975\\
21.125	23.6077015104864\\
21.5	23.6037359721505\\
21.875	23.6000214911615\\
22.25	23.5963198851983\\
22.625	23.5925643310283\\
23	23.5889723065879\\
23.375	23.5856020818733\\
23.75	23.5821342372176\\
24.125	23.5790429117686\\
24.5	23.5756429216265\\
24.875	23.5725000389467\\
25.25	23.5694941186574\\
25.625	23.5665435320321\\
26	23.5635308504592\\
26.375	23.5605508612785\\
26.75	23.5579045174023\\
27.125	23.5551071092012\\
27.5	23.5524866397765\\
27.875	23.5498167282803\\
28.25	23.5471688161675\\
28.625	23.544703089385\\
29	23.5421496230144\\
29.375	23.5398011043619\\
29.75	23.5375007484419\\
30.125	23.5350627492864\\
30.5	23.5326579873993\\
30.875	23.5303808801132\\
31.25	23.5282641560621\\
31.625	23.5261448077812\\
32	23.5239722144278\\
32.375	23.5220812781411\\
32.75	23.5199322355485\\
33.125	23.517968300468\\
33.5	23.5158637115659\\
33.875	23.5139402429543\\
34.25	23.5119923198699\\
34.625	23.5100727204488\\
35	23.5082246067771\\
35.375	23.5064310695946\\
35.75	23.5046625193028\\
36.125	23.5027470813592\\
36.5	23.501056665193\\
36.875	23.4993959985922\\
37.25	23.4977717222873\\
37.625	23.4959717012021\\
38	23.4942824815438\\
38.375	23.4927884650136\\
38.75	23.491239962111\\
39.125	23.4897263044969\\
39.5	23.4881001314106\\
39.875	23.4866182568151\\
40.25	23.4850146552818\\
40.625	23.4835433680536\\
41	23.4821968829744\\
41.375	23.4806967665272\\
41.75	23.4792188125296\\
42.125	23.4778788636342\\
42.5	23.4764895016989\\
42.875	23.4752152032772\\
43.25	23.4737124683513\\
43.625	23.4724599226392\\
44	23.471170980289\\
44.375	23.469912772711\\
44.75	23.4685924061574\\
45.125	23.467515825592\\
45.5	23.4659959499274\\
45.875	23.4648238015175\\
46.25	23.4636483037067\\
46.625	23.4624438611787\\
47	23.4612789703001\\
47.375	23.4601222520642\\
47.75	23.4589951658351\\
48.125	23.4579392111026\\
48.5	23.4566275594404\\
48.875	23.4555241354893\\
49.25	23.4545832775687\\
49.625	23.4533178990279\\
50	23.4522789212031\\
};
\addlegendentry{plug-in charging}
\addplot [color=blue, dashed, line width=1.0pt]
  table[row sep=crcr]{%
12.5	23.6680339271376\\
12.875	23.6516814494833\\
13.25	23.6356992465643\\
13.625	23.6217212395392\\
14	23.6082995649847\\
14.375	23.5956640180806\\
14.75	23.5833790664653\\
15.125	23.5723819376371\\
15.5	23.5609209649587\\
15.875	23.5500146046578\\
16.25	23.5405143233017\\
16.625	23.5313863895252\\
17	23.5224456911585\\
17.375	23.5139890493656\\
17.75	23.5058075041424\\
18.125	23.4971845650695\\
18.5	23.4902424212195\\
18.875	23.483456643095\\
19.25	23.4761999743911\\
19.625	23.4692524587129\\
20	23.4633188270378\\
20.375	23.4567811099084\\
20.75	23.4514502401021\\
21.125	23.4460277930148\\
21.5	23.4395835793431\\
21.875	23.4344482182085\\
22.25	23.4299868393379\\
22.625	23.4246303475484\\
23	23.4205819472033\\
23.375	23.4153006588363\\
23.75	23.4108970945305\\
24.125	23.4060792171701\\
24.5	23.4027725323495\\
24.875	23.3987007750295\\
25.25	23.3939672391977\\
25.625	23.3910465041349\\
26	23.387251837128\\
26.375	23.3833765824769\\
26.75	23.3805786922888\\
27.125	23.3768239034003\\
27.5	23.3736623088266\\
27.875	23.3699523396147\\
28.25	23.3678599751742\\
28.625	23.3645006195564\\
29	23.3615338952799\\
29.375	23.3588683379587\\
29.75	23.3558140978556\\
30.125	23.3529871583015\\
30.5	23.350573581107\\
30.875	23.3482020902291\\
31.25	23.3452496068314\\
31.625	23.3431558429221\\
32	23.3407165098165\\
32.375	23.3380697562376\\
32.75	23.3359311652209\\
33.125	23.3340981789002\\
33.5	23.3316663034355\\
33.875	23.3299681980538\\
34.25	23.3272178090508\\
34.625	23.3248263784829\\
35	23.3235532757932\\
35.375	23.3213474078152\\
35.75	23.3188116261515\\
36.125	23.3177947868981\\
36.5	23.316481665804\\
36.875	23.3135986964239\\
37.25	23.3131126722644\\
37.625	23.3110258263825\\
38	23.3090005665304\\
38.375	23.3077452020428\\
38.75	23.3058756720185\\
39.125	23.3037713655591\\
39.5	23.3031146168202\\
39.875	23.3009640774291\\
40.25	23.300071706906\\
40.625	23.29863975036\\
41	23.2968571431135\\
41.375	23.2958091897472\\
41.75	23.2945295689391\\
42.125	23.2934879159288\\
42.5	23.290895741864\\
42.875	23.2902666582428\\
43.25	23.2888996339194\\
43.625	23.28763501635\\
44	23.285184072025\\
44.375	23.2847775926872\\
44.75	23.2832933389888\\
45.125	23.2831389469402\\
45.5	23.2810001068696\\
45.875	23.2806942055775\\
46.25	23.2792768577493\\
46.625	23.2783238289482\\
47	23.2772266493402\\
47.375	23.2766417311905\\
47.75	23.2745866582481\\
48.125	23.2743215469613\\
48.5	23.2733768673396\\
48.875	23.2720182277666\\
49.25	23.27086064852\\
49.625	23.2702552976197\\
50	23.2687941187298\\
};
\addlegendentry{battery swapping}
\end{axis}
\end{tikzpicture}%
\vspace*{-0.3in}
\caption{Passenger travel cost $c$ (\$/trip) as a function of $C$.}
\label{c_C}
\end{subfigure}
\hfill
\begin{subfigure}[b]{0.32\linewidth}
\centering
%
%
\begin{tikzpicture}

\begin{axis}[%
width=1.694in,
height=1.03in,
at={(1.358in,0.0in)},
scale only axis,
xmin=12.5,
xmax=50,
xtick={12.5,   25, 37.5,   50},
xlabel style={font=\color{white!15!black}},
xlabel={Battery capacity},
ymin=13.98,
ymax=14.18,
ylabel style={font=\color{white!15!black}},
ylabel={Ride fare},
axis background/.style={fill=white},
legend style={at={(1,1)}, anchor=north east, legend cell align=left, align=left, font=\scriptsize, draw=white!12!black}
]
\addplot [color=black, line width=1.0pt]
  table[row sep=crcr]{%
12.5	14.1638657585028\\
12.875	14.1608358370609\\
13.25	14.1579114571544\\
13.625	14.1551210631407\\
14	14.1524591395533\\
14.375	14.1499371409689\\
14.75	14.1475167815928\\
15.125	14.1450724294299\\
15.5	14.1428231139882\\
15.875	14.1406083949804\\
16.25	14.1385368626164\\
16.625	14.1364561768558\\
17	14.134527302312\\
17.375	14.1325734245292\\
17.75	14.1307390147654\\
18.125	14.1289268527847\\
18.5	14.1272483478821\\
18.875	14.1255913443992\\
19.25	14.1240008579633\\
19.625	14.122397859312\\
20	14.1208374661886\\
20.375	14.1192935464957\\
20.75	14.1179388257891\\
21.125	14.1164813405912\\
21.5	14.1151059654005\\
21.875	14.1138174175319\\
22.25	14.1125335326432\\
22.625	14.1112306889618\\
23	14.1099849433103\\
23.375	14.1088157455711\\
23.75	14.1076128576083\\
24.125	14.1065407422996\\
24.5	14.1053612758349\\
24.875	14.1042710936125\\
25.25	14.1032283867707\\
25.625	14.102204994793\\
26	14.1011597927933\\
26.375	14.1001262076634\\
26.75	14.0992083182413\\
27.125	14.0982378865548\\
27.5	14.0973287922784\\
27.875	14.096402710579\\
28.25	14.0954842406128\\
28.625	14.0946288983009\\
29	14.0937431127325\\
29.375	14.0929283654575\\
29.75	14.0921302744828\\
30.125	14.0912847014314\\
30.5	14.0904506134422\\
30.875	14.0896605872814\\
31.25	14.0889263099902\\
31.625	14.0881909754097\\
32	14.0874375180553\\
32.375	14.0867816480769\\
32.75	14.0860358681294\\
33.125	14.0853544895163\\
33.5	14.0846246865744\\
33.875	14.0839573327719\\
34.25	14.083281782201\\
34.625	14.0826159934647\\
35	14.081974600117\\
35.375	14.081352561943\\
35.75	14.0807389227423\\
36.125	14.0800744875379\\
36.5	14.0794882407342\\
36.875	14.0789119822427\\
37.25	14.0783485061391\\
37.625	14.0777240699286\\
38	14.0771380913598\\
38.375	14.0766198459036\\
38.75	14.0760828372008\\
39.125	14.075557710107\\
39.5	14.0749937721232\\
39.875	14.0744795962457\\
40.25	14.0739231812033\\
40.625	14.0734129009985\\
41	14.0729456936994\\
41.375	14.0724253931383\\
41.75	14.0719126338884\\
42.125	14.0714478675709\\
42.5	14.0709659759403\\
42.875	14.0705237126156\\
43.25	14.0700026316782\\
43.625	14.0695680735906\\
44	14.069121138761\\
44.375	14.0686845677594\\
44.75	14.0682263652943\\
45.125	14.0678530834688\\
45.5	14.0673257846797\\
45.875	14.066919219344\\
46.25	14.0665113893446\\
46.625	14.0660936171985\\
47	14.0656895235529\\
47.375	14.0652882245356\\
47.75	14.0648973944386\\
48.125	14.0645309733611\\
48.5	14.0640759412938\\
48.875	14.063693292539\\
49.25	14.0633667276037\\
49.625	14.0629279932021\\
50	14.0625675222135\\
};
\addlegendentry{plug-in charging}
\addplot [color=blue, dashed, line width=1.0pt]
  table[row sep=crcr]{%
12.5	14.1374061697601\\
12.875	14.1317349819071\\
13.25	14.1261920653884\\
13.625	14.1213442592364\\
14	14.1166889162777\\
14.375	14.1123056471691\\
14.75	14.1080446710476\\
15.125	14.1042297198651\\
15.5	14.1002544540643\\
15.875	14.0964715088687\\
16.25	14.0931755485257\\
16.625	14.0900091168554\\
17	14.0869076424826\\
17.375	14.0839742407027\\
17.75	14.0811355929419\\
18.125	14.0781446035257\\
18.5	14.0757369680111\\
18.875	14.0733825776399\\
19.25	14.0708653044846\\
19.625	14.0684551564801\\
20	14.0663973322452\\
20.375	14.0641290460426\\
20.75	14.0622797325836\\
21.125	14.0603993516065\\
21.5	14.0581634584991\\
21.875	14.0563821750823\\
22.25	14.0548344183201\\
22.625	14.0529765981124\\
23	14.0515728858898\\
23.375	14.049740018182\\
23.75	14.0482131031752\\
24.125	14.0465418996673\\
24.5	14.0453952468677\\
24.875	14.0439822924873\\
25.25	14.0423403602251\\
25.625	14.0413277463654\\
26	14.0400114035535\\
26.375	14.0386674183979\\
26.75	14.0376974267562\\
27.125	14.0363952090382\\
27.5	14.0352982973547\\
27.875	14.0340118085222\\
28.25	14.0332856881852\\
28.625	14.0321210167\\
29	14.0310918636053\\
29.375	14.0301676223243\\
29.75	14.0291078884639\\
30.125	14.0281280762304\\
30.5	14.0272909046091\\
30.875	14.0264679317815\\
31.25	14.0254444862646\\
31.625	14.0247178514783\\
32	14.0238724372583\\
32.375	14.0229539076127\\
32.75	14.022212548182\\
33.125	14.0215777128447\\
33.5	14.0207340883346\\
33.875	14.0201448002157\\
34.25	14.0191912336165\\
34.625	14.0183622401625\\
35	14.0179201374429\\
35.375	14.0171556123602\\
35.75	14.0162762071879\\
36.125	14.0159238521425\\
36.5	14.015468511555\\
36.875	14.0144688895916\\
37.25	14.0143000980584\\
37.625	14.0135761703588\\
38	14.0128742196995\\
38.375	14.0124393798954\\
38.75	14.0117901538903\\
39.125	14.0110608998562\\
39.5	14.0108332549425\\
39.875	14.0100877183902\\
40.25	14.0097782661525\\
40.625	14.0092818274929\\
41	14.0086636540699\\
41.375	14.0083006086259\\
41.75	14.0078563060601\\
42.125	14.0074953268779\\
42.5	14.0065970333911\\
42.875	14.0063786482007\\
43.25	14.0059047040393\\
43.625	14.0054663072382\\
44	14.0046157204\\
44.375	14.0044752803963\\
44.75	14.0039605816067\\
45.125	14.0039078463242\\
45.5	14.003166057352\\
45.875	14.0030599496881\\
46.25	14.002568108861\\
46.625	14.0022381855651\\
47	14.0018575630397\\
47.375	14.0016546547043\\
47.75	14.0009426100771\\
48.125	14.000850487328\\
48.5	14.000523122985\\
48.875	14.000051840795\\
49.25	13.999650426964\\
49.625	13.9994408105208\\
50	13.9989344938301\\
};
\addlegendentry{battery swapping}
\end{axis}
\end{tikzpicture}%
\vspace*{-0.3in}
\caption{Per-trip ride fare $p_f$ (\$/trip) as a function of $C$.}
\label{pf_C}
\end{subfigure}
\begin{subfigure}[b]{0.32\linewidth}
\centering
%
%
\begin{tikzpicture}

\begin{axis}[%
width=1.694in,
height=1.03in,
at={(1.358in,0.0in)},
scale only axis,
xmin=12.5,
xmax=50,
xtick={12.5,   25, 37.5,   50},
xlabel style={font=\color{white!15!black}},
xlabel={Battery capacity},
ymin=3.58,
ymax=3.74,
ylabel style={font=\color{white!15!black}},
ylabel={Waiting time},
axis background/.style={fill=white},
legend style={at={(1,1)}, anchor=north east, legend cell align=left, align=left, font=\scriptsize, draw=white!12!black}
]
\addplot [color=black, line width=1.0pt]
  table[row sep=crcr]{%
12.5	3.71336383054265\\
12.875	3.71115084626499\\
13.25	3.70901487017297\\
13.625	3.70697701155998\\
14	3.70503299464072\\
14.375	3.70319121121051\\
14.75	3.70142383751281\\
15.125	3.69963910207676\\
15.5	3.69799680966106\\
15.875	3.69637966700509\\
16.25	3.69486743142671\\
16.625	3.69334826899207\\
17	3.69194022019745\\
17.375	3.690513712153\\
17.75	3.68917459207394\\
18.125	3.68785164470192\\
18.5	3.68662643248715\\
18.875	3.68541717506814\\
19.25	3.68425616696923\\
19.625	3.6830858885112\\
20	3.68194703630206\\
20.375	3.6808200324204\\
20.75	3.67983136139863\\
21.125	3.6787675077113\\
21.5	3.67776356850776\\
21.875	3.67682328435255\\
22.25	3.6758861831609\\
22.625	3.67493552018081\\
23	3.67402610979754\\
23.375	3.67317299856678\\
23.75	3.67229510837568\\
24.125	3.67151246878645\\
24.5	3.67065180069441\\
24.875	3.66985618036207\\
25.25	3.66909524491733\\
25.625	3.66834827024773\\
26	3.66758568126588\\
26.375	3.66683126109114\\
26.75	3.66616131750427\\
27.125	3.66545318707226\\
27.5	3.66478986337134\\
27.875	3.66411396034932\\
28.25	3.66344363393593\\
28.625	3.66281945390856\\
29	3.66217306600075\\
29.375	3.66157858097071\\
29.75	3.6609963077361\\
30.125	3.66037908831588\\
30.5	3.65977029998338\\
30.875	3.65919391195031\\
31.25	3.65865807987285\\
31.625	3.65812164045407\\
32	3.65757158774128\\
32.375	3.65709287986982\\
32.75	3.65654897961981\\
33.125	3.65605186470996\\
33.5	3.65551900193467\\
33.875	3.65503213572965\\
34.25	3.65453896808873\\
34.625	3.65405299495504\\
35	3.65358527389926\\
35.375	3.6531312045161\\
35.75	3.65268356455832\\
36.125	3.65219867977571\\
36.5	3.65177070715457\\
36.875	3.65135039393394\\
37.25	3.6509392310652\\
37.625	3.65048357801298\\
38	3.65005596518758\\
38.375	3.64967775934494\\
38.75	3.64928570732955\\
39.125	3.6489025559651\\
39.5	3.64849083693308\\
39.875	3.6481157599106\\
40.25	3.64770987367383\\
40.625	3.64733739033142\\
41	3.64699658499029\\
41.375	3.64661681139106\\
41.75	3.64624270489967\\
42.125	3.64590348684622\\
42.5	3.64555175416998\\
42.875	3.64522925994639\\
43.25	3.64484877390432\\
43.625	3.64453172443744\\
44	3.64420536493334\\
44.375	3.64388690114401\\
44.75	3.64355272901673\\
45.125	3.6432801326059\\
45.5	3.6428954128867\\
45.875	3.64259867526104\\
46.25	3.64230112959772\\
46.625	3.64199621859696\\
47	3.64170133594853\\
47.375	3.64140853780178\\
47.75	3.64112316720792\\
48.125	3.64085590610137\\
48.5	3.64052388300255\\
48.875	3.64024451277144\\
49.25	3.64000641471513\\
49.625	3.63968601001002\\
50	3.63942302286419\\
};
\addlegendentry{plug-in charging}
\addplot [color=blue, dashed, line width=1.0pt]
  table[row sep=crcr]{%
12.5	3.69404176642536\\
12.875	3.68990173161867\\
13.25	3.68585549657979\\
13.625	3.68231665903206\\
14	3.67891885608796\\
14.375	3.67572029880291\\
14.75	3.67261023078206\\
15.125	3.66982644099688\\
15.5	3.66692500422261\\
15.875	3.66416399061595\\
16.25	3.66175921502947\\
16.625	3.6594485552984\\
17	3.65718529018446\\
17.375	3.65504449948172\\
17.75	3.65297360899243\\
18.125	3.65079068276891\\
18.5	3.64903312139861\\
18.875	3.64731552924617\\
19.25	3.64547855422733\\
19.625	3.64371988458635\\
20	3.64221763364055\\
20.375	3.64056281545187\\
20.75	3.63921337500716\\
21.125	3.63784048116598\\
21.5	3.63620934916434\\
21.875	3.63490931904116\\
22.25	3.63378000814643\\
22.625	3.63242393388993\\
23	3.63139886097423\\
23.375	3.63006226381951\\
23.75	3.62894728347103\\
24.125	3.62772764244292\\
24.5	3.62689042072937\\
24.875	3.62585987695433\\
25.25	3.62466158099716\\
25.625	3.62392199913549\\
26	3.62296140836222\\
26.375	3.62198029615465\\
26.75	3.62127180834598\\
27.125	3.62032119936517\\
27.5	3.61952093467903\\
27.875	3.61858160119863\\
28.25	3.61805204922054\\
28.625	3.61720139645599\\
29	3.61645039987386\\
29.375	3.61577547117614\\
29.75	3.61500240674099\\
30.125	3.61428646591903\\
30.5	3.61367545600693\\
30.875	3.61307525521226\\
31.25	3.61232756611119\\
31.625	3.61179767110223\\
32	3.61117987308457\\
32.375	3.61051001884684\\
32.75	3.60996845621663\\
33.125	3.60950405661065\\
33.5	3.60888845546546\\
33.875	3.60845868133259\\
34.25	3.60776223854044\\
34.625	3.60715664275987\\
35	3.60683454974818\\
35.375	3.60627588971126\\
35.75	3.60563388331922\\
36.125	3.60537633130062\\
36.5	3.60504385823604\\
36.875	3.60431387861718\\
37.25	3.60419092023487\\
37.625	3.6036626573735\\
38	3.60314974683368\\
38.375	3.60283171401061\\
38.75	3.60235872795667\\
39.125	3.60182576190035\\
39.5	3.60165944258826\\
39.875	3.60111486784454\\
40.25	3.60088893052461\\
40.625	3.60052632669268\\
41	3.60007499575337\\
41.375	3.59980952756639\\
41.75	3.59948576080582\\
42.125	3.59922193374065\\
42.5	3.59856539088098\\
42.875	3.59840620544268\\
43.25	3.59806005034112\\
43.625	3.59773980973326\\
44	3.59711951613371\\
44.375	3.59701640011276\\
44.75	3.59664060363647\\
45.125	3.59660120178915\\
45.5	3.59605970911534\\
45.875	3.59598226972457\\
46.25	3.59562354608071\\
46.625	3.59538203231903\\
47	3.59510429701568\\
47.375	3.59495623119622\\
47.75	3.59443567758567\\
48.125	3.59436862776483\\
48.5	3.59412935827695\\
48.875	3.59378542130685\\
49.25	3.59349233393643\\
49.625	3.59333894848795\\
50	3.59296884686035\\
};
\addlegendentry{battery swapping}
\end{axis}
\end{tikzpicture}%
\vspace*{-0.3in}
\caption{Passenger waiting time $w^c$ (minute) as a function of $C$.}
\label{wc_C}
\end{subfigure}
\hfill
\begin{subfigure}[b]{0.32\linewidth}
\centering
%
%
\begin{tikzpicture}

\begin{axis}[%
width=1.694in,
height=1.03in,
at={(1.358in,0.0in)},
scale only axis,
xmin=12.5,
xmax=50,
xtick={12.5,   25, 37.5,   50},
xlabel style={font=\color{white!15!black}},
xlabel={Battery capacity},
ymin=18.04,
ymax=18.22,
ylabel style={font=\color{white!15!black}},
ylabel={Vehicle idle time},
axis background/.style={fill=white},
legend style={at={(1,0)}, anchor=south east, legend cell align=left, align=left, font=\scriptsize, draw=white!12!black}
]
\addplot [color=black, line width=1.0pt]
  table[row sep=crcr]{%
12.5	18.0501639040625\\
12.875	18.052551653088\\
13.25	18.0548743694223\\
13.625	18.0571059752054\\
14	18.0592496584537\\
14.375	18.0612939059527\\
14.75	18.0632673180141\\
15.125	18.0652718863251\\
15.5	18.0671272361641\\
15.875	18.0689647064954\\
16.25	18.0706910832691\\
16.625	18.0724350339406\\
17	18.0740586168641\\
17.375	18.0757119640205\\
17.75	18.0772707360239\\
18.125	18.0788177438549\\
18.5	18.080256067541\\
18.875	18.0816806246876\\
19.25	18.083054593596\\
19.625	18.0844452723497\\
20	18.0858027842627\\
20.375	18.0871516923455\\
20.75	18.0883384532787\\
21.125	18.089620293865\\
21.5	18.0908340720816\\
21.875	18.0919736499045\\
22.25	18.0931134816545\\
22.625	18.0942725084675\\
23	18.0953857957331\\
23.375	18.0964318690419\\
23.75	18.0975119053799\\
24.125	18.0984778850697\\
24.5	18.0995419498094\\
24.875	18.1005285028722\\
25.25	18.1014742726049\\
25.625	18.1024053159492\\
26	18.1033571826885\\
26.375	18.1043020322506\\
26.75	18.1051428793786\\
27.125	18.1060330990207\\
27.5	18.106868645262\\
27.875	18.1077223799236\\
28.25	18.1085707943852\\
28.625	18.109362182142\\
29	18.1101833310922\\
29.375	18.110939816274\\
29.75	18.1116819463418\\
30.125	18.1124709951263\\
30.5	18.1132506049253\\
30.875	18.1139893618545\\
31.25	18.1146776731389\\
31.625	18.1153674170359\\
32	18.1160770261736\\
32.375	18.1166952689263\\
32.75	18.1173975243171\\
33.125	18.1180409496887\\
33.5	18.1187329889432\\
33.875	18.1193650869523\\
34.25	18.1200072901432\\
34.625	18.1206408694771\\
35	18.1212502118189\\
35.375	18.1218439870705\\
35.75	18.1224292685994\\
36.125	18.1230647172659\\
36.5	18.1236268019251\\
36.875	18.1241784407552\\
37.25	18.124719269982\\
37.625	18.1253194337348\\
38	18.1258834862052\\
38.375	18.1263830383327\\
38.75	18.1269019403849\\
39.125	18.127408980398\\
39.5	18.1279552264931\\
39.875	18.1284525177865\\
40.25	18.1289912824992\\
40.625	18.1294870393208\\
41	18.1299404148175\\
41.375	18.130446893247\\
41.75	18.1309458934577\\
42.125	18.1313992221849\\
42.5	18.1318698086418\\
42.875	18.1323007759303\\
43.25	18.1328113279656\\
43.625	18.1332364294806\\
44	18.1336752719512\\
44.375	18.1341029311087\\
44.75	18.1345519102619\\
45.125	18.1349195603459\\
45.5	18.1354378925854\\
45.875	18.1358384083956\\
46.25	18.1362400268722\\
46.625	18.1366522897622\\
47	18.1370512144974\\
47.375	18.137447530456\\
47.75	18.1378347543359\\
48.125	18.1381968462721\\
48.5	18.1386474701047\\
48.875	18.1390274636893\\
49.25	18.1393505960525\\
49.625	18.1397873562193\\
50	18.1401453739075\\
};
\addlegendentry{plug-in charging}
\addplot [color=blue, dashed, line width=1.0pt]
  table[row sep=crcr]{%
12.5	18.0716381596374\\
12.875	18.0764232112981\\
13.25	18.0811638315274\\
13.625	18.0853624159663\\
14	18.0894381601805\\
14.375	18.0933142717636\\
14.75	18.097123941847\\
15.125	18.1005640653719\\
15.5	18.1041840997323\\
15.875	18.1076596697168\\
16.25	18.1107089864722\\
16.625	18.1136617377809\\
17	18.1165745092215\\
17.375	18.1193489378774\\
17.75	18.122047917099\\
18.125	18.1249141220334\\
18.5	18.1272368114675\\
18.875	18.1295153353181\\
19.25	18.1319669694201\\
19.625	18.1343263879153\\
20	18.1363536330998\\
20.375	18.1385940597912\\
20.75	18.1404300415967\\
21.125	18.1423078004279\\
21.5	18.1445446837646\\
21.875	18.1463368248103\\
22.25	18.1478983150475\\
22.625	18.1497817390034\\
23	18.151211767412\\
23.375	18.1530772048207\\
23.75	18.1546433466436\\
24.125	18.1563601843829\\
24.5	18.1575433862512\\
24.875	18.1590003823961\\
25.25	18.1607022065991\\
25.625	18.1617571727683\\
26	18.163128259702\\
26.375	18.1645335616141\\
26.75	18.1655519763389\\
27.125	18.1669200396307\\
27.5	18.1680731873502\\
27.875	18.1694323330844\\
28.25	18.1701982769753\\
28.625	18.1714348503162\\
29	18.1725270019143\\
29.375	18.1735119568554\\
29.75	18.1746400112401\\
30.125	18.1756904088249\\
30.5	18.1765864206669\\
30.875	18.1774667428522\\
31.25	18.1785692789497\\
31.625	18.1793492016258\\
32	18.1802637882463\\
32.375	18.1812531219365\\
32.75	18.182057047721\\
33.125	18.1827493321172\\
33.5	18.183663822657\\
33.875	18.1843024861254\\
34.25	18.1853420166912\\
34.625	18.1862479741078\\
35	18.1867278392392\\
35.375	18.1875661387076\\
35.75	18.1885293084167\\
36.125	18.1889171212198\\
36.5	18.1894170942371\\
36.875	18.1905167125048\\
37.25	18.19070124879\\
37.625	18.1914974883206\\
38	18.1922737304198\\
38.375	18.192756476895\\
38.75	18.1934697204515\\
39.125	18.1942795220326\\
39.5	18.1945323172235\\
39.875	18.195360479512\\
40.25	18.1957041332599\\
40.625	18.1962565288446\\
41	18.1969442899535\\
41.375	18.1973504301957\\
41.75	18.1978428314712\\
42.125	18.1982467170141\\
42.5	18.1992530866885\\
42.875	18.1994960892168\\
43.25	18.2000269203589\\
43.625	18.2005185976124\\
44	18.2014690615075\\
44.375	18.2016293037845\\
44.75	18.2022074094424\\
45.125	18.2022707306613\\
45.5	18.2031045301984\\
45.875	18.2032238206926\\
46.25	18.2037760486163\\
46.625	18.2041507771126\\
47	18.2045793587625\\
47.375	18.2048079931832\\
47.75	18.2056153889147\\
48.125	18.205718580889\\
48.5	18.2060899728471\\
48.875	18.2066226953791\\
49.25	18.2070774612266\\
49.625	18.2073166041246\\
50	18.2078928671463\\
};
\addlegendentry{battery swapping}
\end{axis}
\end{tikzpicture}%
\vspace*{-0.3in}
\caption{Vehicle idle time $w^v$ (minute) as a function of $C$.}
\label{wv_C}
\end{subfigure}
\hfill
\begin{subfigure}[b]{0.32\linewidth}
\centering
%
%
\begin{tikzpicture}

\begin{axis}[%
width=1.694in,
height=1.03in,
at={(1.358in,0.0in)},
scale only axis,
xmin=12.5,
xmax=50,
xtick={12.5,   25, 37.5,   50},
xlabel style={font=\color{white!15!black}},
xlabel={Battery capacity},
ymin=8400,
ymax=9100,
ytick={8400,8600,8800,9000},
yticklabels={{8.4K},{8.6K},{8.8K},{9.0K}},
ylabel style={font=\color{white!15!black}},
ylabel={Number of vehicles},
axis background/.style={fill=white},
legend style={at={(1,0.3)}, anchor=south east, legend cell align=left, align=left, font=\scriptsize, draw=white!12!black}
]
\addplot [color=black, line width=1.0pt]
  table[row sep=crcr]{%
12.5	8973.39259522367\\
12.875	8976.95461295442\\
13.25	8980.45386690526\\
13.625	8983.79233863989\\
14	8986.90458789245\\
14.375	8989.82894168293\\
14.75	8992.69691057465\\
15.125	8995.75222099769\\
15.5	8998.3821220185\\
15.875	9001.0410209948\\
16.25	9003.35524716282\\
16.625	9005.953313808\\
17	9008.1499963466\\
17.375	9010.58430407647\\
17.75	9012.68753477599\\
18.125	9014.8882664021\\
18.5	9016.80386624642\\
18.875	9018.81092184519\\
19.25	9020.67622241704\\
19.625	9022.63117552297\\
20	9024.58214337639\\
20.375	9026.36854546038\\
20.75	9027.96749680255\\
21.125	9029.7631357382\\
21.5	9031.36784973516\\
21.875	9032.8798808479\\
22.25	9034.37390358416\\
22.625	9035.97970024254\\
23	9037.52759753485\\
23.375	9038.88400100847\\
23.75	9040.39918604346\\
24.125	9041.59623282627\\
24.5	9043.05437195297\\
24.875	9044.39685844584\\
25.25	9045.62566868916\\
25.625	9046.80540428563\\
26	9048.13666806861\\
26.375	9049.40274350676\\
26.75	9050.4726017804\\
27.125	9051.61864924862\\
27.5	9052.65322914966\\
27.875	9053.79527069284\\
28.25	9055.03274128834\\
28.625	9055.9461609709\\
29	9057.0998078311\\
29.375	9058.08515219086\\
29.75	9058.96860509974\\
30.125	9060.02793372793\\
30.5	9061.0463888054\\
30.875	9062.01598680121\\
31.25	9062.8848098786\\
31.625	9063.7626088587\\
32	9064.69026584584\\
32.375	9065.39118321111\\
32.75	9066.32378425326\\
33.125	9067.14474651706\\
33.5	9068.03986243857\\
33.875	9068.85591921416\\
34.25	9069.73171439654\\
34.625	9070.56229058108\\
35	9071.27354731411\\
35.375	9072.01766057891\\
35.75	9072.77720212316\\
36.125	9073.6139639446\\
36.5	9074.30992686252\\
36.875	9075.02229497242\\
37.25	9075.70555811663\\
37.625	9076.44253133106\\
38	9077.2113611737\\
38.375	9077.80027339121\\
38.75	9078.4935126483\\
39.125	9079.12458891308\\
39.5	9079.82403294218\\
39.875	9080.39579964071\\
40.25	9081.08371331546\\
40.625	9081.73535337163\\
41	9082.34292621141\\
41.375	9082.95038765573\\
41.75	9083.60752766065\\
42.125	9084.18680654798\\
42.5	9084.79353763955\\
42.875	9085.23572436368\\
43.25	9086.00393178576\\
43.625	9086.43787855054\\
44	9087.04459246167\\
44.375	9087.59431478259\\
44.75	9088.10534922345\\
45.125	9088.53043447498\\
45.5	9089.27689236396\\
45.875	9089.75032116443\\
46.25	9090.22149085366\\
46.625	9090.77486631768\\
47	9091.24614151914\\
47.375	9091.74704722552\\
47.75	9092.22479108926\\
48.125	9092.66796254914\\
48.5	9093.28041901206\\
48.875	9093.70350160595\\
49.25	9094.04321276858\\
49.625	9094.68161453963\\
50	9095.124215465\\
};
\addlegendentry{plug-in charging}
\addplot [color=blue, dashed, line width=1.0pt]
  table[row sep=crcr]{%
12.5	8469.58765290131\\
12.875	8472.69690989629\\
13.25	8476.09338517809\\
13.625	8478.33390735121\\
14	8480.60920006556\\
14.375	8482.68382596035\\
14.75	8484.89674569413\\
15.125	8486.41020814336\\
15.5	8488.70157688809\\
15.875	8490.86964189907\\
16.25	8492.15459692197\\
16.625	8493.42573459146\\
17	8494.82292007015\\
17.375	8496.07168684036\\
17.75	8497.33505491463\\
18.125	8499.22035553336\\
18.5	8499.87261896715\\
18.875	8500.59888321883\\
19.25	8501.93367465652\\
19.625	8503.1882613991\\
20	8503.73877090104\\
20.375	8504.98284778858\\
20.75	8505.33481527926\\
21.125	8505.91766900518\\
21.5	8507.53526570382\\
21.875	8508.14707339065\\
22.25	8508.29867463202\\
22.625	8509.35635481686\\
23	8509.38719191307\\
23.375	8510.60778876179\\
23.75	8511.16772855836\\
24.125	8512.19218155716\\
24.5	8511.98898861554\\
24.875	8512.54884110253\\
25.25	8513.7788833292\\
25.625	8513.50270600826\\
26	8514.07489434429\\
26.375	8514.79630621951\\
26.75	8514.64521009482\\
27.125	8515.4071831996\\
27.5	8515.71588718353\\
27.875	8516.57398506744\\
28.25	8516.0702579092\\
28.625	8516.7444121249\\
29	8517.13088897896\\
29.375	8517.30867984707\\
29.75	8517.88321391851\\
30.125	8518.3106214294\\
30.5	8518.42310576622\\
30.875	8518.54685260241\\
31.25	8519.23142220902\\
31.625	8519.20205760226\\
32	8519.52248486607\\
32.375	8520.06792961866\\
32.75	8520.20549728812\\
33.125	8520.11260447903\\
33.5	8520.58683526653\\
33.875	8520.4485876836\\
34.25	8521.27799114\\
34.625	8521.82453784976\\
35	8521.41328360319\\
35.375	8521.86150440064\\
35.75	8522.63350892739\\
36.125	8522.08982483595\\
36.5	8521.83782943994\\
36.875	8523.00667790806\\
37.25	8522.077153582\\
37.625	8522.59455592555\\
38	8523.08413435803\\
38.375	8522.91605123205\\
38.75	8523.31666712273\\
39.125	8523.95087995354\\
39.5	8523.32268872465\\
39.875	8524.04397881544\\
40.25	8523.66948361461\\
40.625	8523.79593854503\\
41	8524.2546957716\\
41.375	8524.08115212183\\
41.75	8524.13259124592\\
42.125	8523.99213624393\\
42.5	8525.24962168813\\
42.875	8524.77899073538\\
43.25	8524.98261998082\\
43.625	8525.11257160531\\
44	8526.31478928847\\
44.375	8525.71458141846\\
44.75	8526.08948060122\\
45.125	8525.29667029787\\
45.5	8526.28529342083\\
45.875	8525.65662356273\\
46.25	8526.03214062308\\
46.625	8526.00906077437\\
47	8526.12759301101\\
47.375	8525.8032997031\\
47.75	8526.80293650044\\
48.125	8526.21954228329\\
48.5	8526.25487527133\\
48.875	8526.67115872604\\
49.25	8526.92046968772\\
49.625	8526.68942530538\\
50	8527.23317500649\\
};
\addlegendentry{battery swapping}
\end{axis}
\end{tikzpicture}%
\vspace*{-0.3in}
\caption{Total number of vehicles $N$ as a function of $C$.}
\label{N_C}
\end{subfigure}
\caption{Market outcomes on electrified AMoD services under different battery capacities. Black lines present the results when charging vehicles; blue dashed lines show the results when swapping batteries.}
\label{AMoD_service_C}
\end{figure}
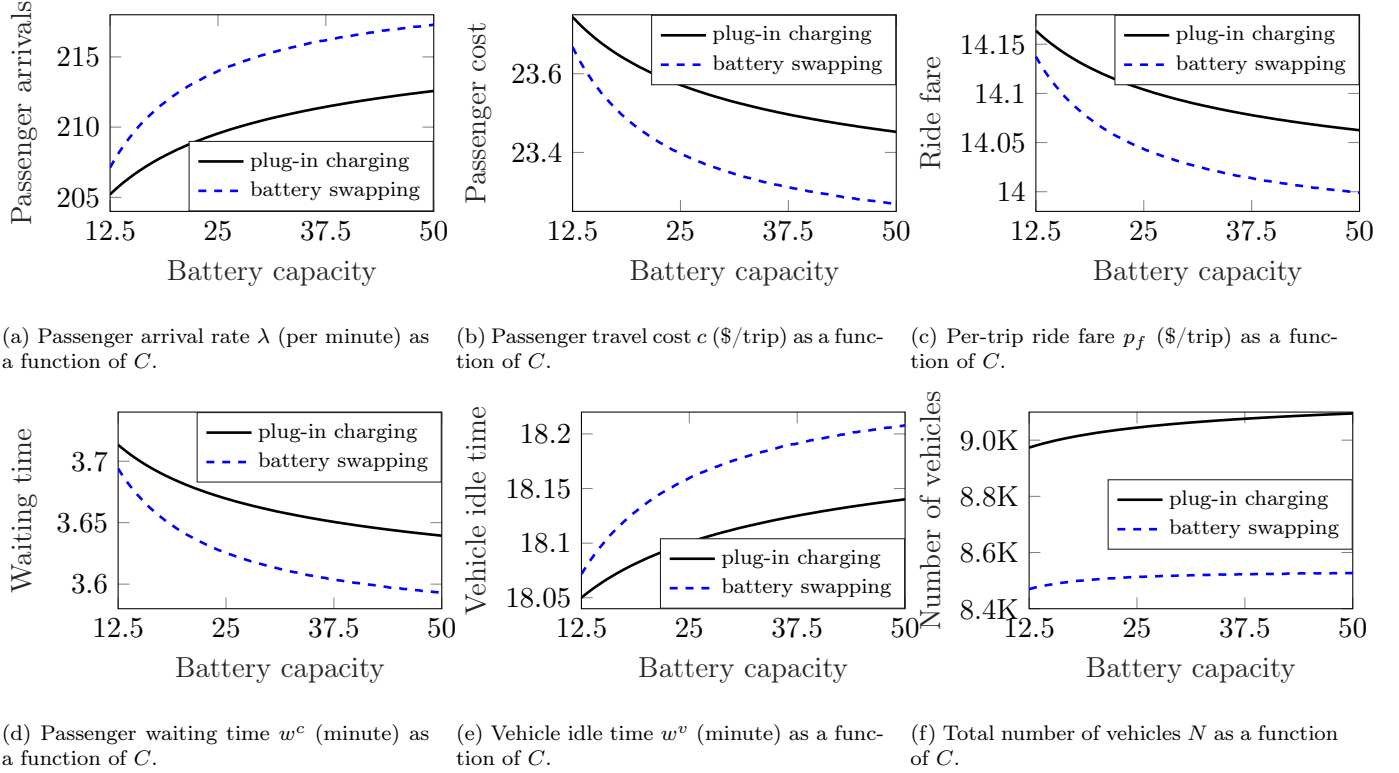

Figure \ref{AMoD_service_C} presents the impacts of battery capacity on the AMoD system. Under both plug-in charging and battery swapping, as the battery capacity increases, the TNC reduces the ride fare (Figure \ref{pf_C}) and recruits more autonomous vehicles (Figure \ref{N_C}). The vehicle idle time increases (Figure \ref{wv_C}), which lowers the passenger waiting time and improves the service quality (Figure \ref{wc_C}). Therefore, the total passenger cost reduces with battery capacity (Figure \ref{c_C}), attracting more passengers to choose AMoD services (Figure \ref{lambda_C}).

\begin{figure}[h]
\centering
\begin{subfigure}[b]{0.32\linewidth}
\centering
%
%
\begin{tikzpicture}

\begin{axis}[%
width=1.694in,
height=1.03in,
at={(1.358in,0.0in)},
scale only axis,
xmin=12.5,
xmax=50,
xtick={12.5,   25, 37.5,   50},
xlabel style={font=\color{white!15!black}},
xlabel={Battery capacity},
ymin=400,
ymax=1800,
ytick={500,1000,1500},
yticklabels={{0.5K},{1K},{1.5K}},
ylabel style={font=\color{white!15!black}},
ylabel={Charging arrivals},
axis background/.style={fill=white},
legend style={at={(1,1)}, anchor=north east, legend cell align=left, align=left, font=\scriptsize, draw=white!12!black}
]
\addplot [color=black, line width=1.0pt]
  table[row sep=crcr]{%
12.5	1593.96135309245\\
12.875	1548.31280724541\\
13.25	1505.23411374706\\
13.625	1464.49316012705\\
14	1425.88711960857\\
14.375	1389.2591758405\\
14.75	1354.48412388083\\
15.125	1321.46523619088\\
15.5	1289.96732215052\\
15.875	1259.96130968657\\
16.25	1231.27981780163\\
16.625	1203.94135923426\\
17	1177.7420403041\\
17.375	1152.71134853599\\
17.75	1128.68496777657\\
18.125	1105.66812750688\\
18.5	1083.54234297545\\
18.875	1062.30964256591\\
19.25	1041.88310412984\\
19.625	1022.24844441521\\
20	1003.34987112651\\
20.375	985.121302722816\\
20.75	967.532304419586\\
21.125	950.58927891927\\
21.5	934.212921400548\\
21.875	918.388076502349\\
22.25	903.093146595156\\
22.625	888.318077746772\\
23	874.018616705891\\
23.375	860.155654720844\\
23.75	846.748066226625\\
24.125	833.722347011526\\
24.5	821.121952741431\\
24.875	808.890331617102\\
25.25	797.009696137203\\
25.625	785.471091270526\\
26	774.281109324924\\
26.375	763.401730620637\\
26.75	752.809241307277\\
27.125	742.515476572261\\
27.5	732.492077652466\\
27.875	722.748592248023\\
28.25	713.273428672219\\
28.625	704.014882718336\\
29	695.018615529937\\
29.375	686.237061466535\\
29.75	677.667002463154\\
30.125	669.325478444161\\
30.5	661.185128273352\\
30.875	653.238850663577\\
31.25	645.475167905085\\
31.625	637.896072318525\\
32	630.498289339558\\
32.375	623.253997117875\\
32.75	616.193362554682\\
33.125	609.284278113824\\
33.5	602.534712682346\\
33.875	595.92922207209\\
34.25	589.472923097054\\
34.625	583.152723801716\\
35	576.958565661213\\
35.375	570.898634342819\\
35.75	564.967069978329\\
36.125	559.163374621539\\
36.5	553.469509749428\\
36.875	547.892775673625\\
37.25	542.426201516779\\
37.625	537.071149808133\\
38	531.824841362666\\
38.375	526.66903224783\\
38.75	521.620101429771\\
39.125	516.663555197695\\
39.5	511.805080510847\\
39.875	507.02980849625\\
40.25	502.350703517453\\
40.625	497.756231723847\\
41	493.24348533238\\
41.375	488.811395701868\\
41.75	484.462305262692\\
42.125	480.186226192426\\
42.5	475.986998288897\\
42.875	471.851074675207\\
43.25	467.806503875539\\
43.625	463.811420460617\\
44	459.895122544039\\
44.375	456.041513036977\\
44.75	452.249618953908\\
45.125	448.516764493232\\
45.5	444.862523832842\\
45.875	441.253319163094\\
46.25	437.702412703126\\
46.625	434.213298657218\\
47	430.77514692886\\
47.375	427.393139619014\\
47.75	424.062993977642\\
48.125	420.783064024337\\
48.5	417.562135347741\\
48.875	414.380903366075\\
49.25	411.244208558173\\
49.625	408.169874544644\\
50	405.13189886037\\
};
\addlegendentry{plug-in charging}
\addplot [color=blue, dashed, line width=1.0pt]
  table[row sep=crcr]{%
12.5	1627.20996772176\\
12.875	1580.72359434407\\
13.25	1536.91402971182\\
13.625	1495.27125124447\\
14	1455.85504493985\\
14.375	1418.44712043028\\
14.75	1382.95815024101\\
15.125	1349.09396536791\\
15.5	1316.99906502176\\
15.875	1286.39265257583\\
16.25	1257.04349597172\\
16.625	1229.01064581464\\
17	1202.23004058335\\
17.375	1176.5774121926\\
17.75	1152.00689813548\\
18.125	1128.54395517887\\
18.5	1105.84498143923\\
18.875	1084.05557353745\\
19.25	1063.19941872182\\
19.625	1043.12635946248\\
20	1023.70633520778\\
20.375	1005.09289321798\\
20.75	987.031281957183\\
21.125	969.639113919246\\
21.5	952.984881395388\\
21.875	936.772800238368\\
22.25	921.048378110319\\
22.625	905.954501264296\\
23	891.227848986959\\
23.375	877.113318617362\\
23.75	863.366266474608\\
24.125	850.098785885855\\
24.5	837.097563383502\\
24.875	824.572297386327\\
25.25	812.491704181817\\
25.625	800.60085617868\\
26	789.142618055776\\
26.375	778.02510299212\\
26.75	767.128093341299\\
27.125	756.625039638091\\
27.5	746.362037856559\\
27.875	736.429232417197\\
28.25	726.625764570172\\
28.625	717.192768487041\\
29	707.975667721167\\
29.375	698.973658765246\\
29.75	690.235143203617\\
30.125	681.700154821277\\
30.5	673.345937137738\\
30.875	665.195242991092\\
31.25	657.289631562303\\
31.625	649.508406885439\\
32	641.939568177085\\
32.375	634.565041401053\\
32.75	627.324460454225\\
33.125	620.228152348797\\
33.5	613.337493410505\\
33.875	606.548833692307\\
34.25	599.987125872125\\
34.625	593.544573428699\\
35	587.16371988644\\
35.375	580.985482850073\\
35.75	574.961824768464\\
36.125	568.961679947549\\
36.5	563.106852191278\\
36.875	557.478126522913\\
37.25	551.8057728606\\
37.625	546.353325193723\\
38	541.006276020609\\
38.375	535.715848432995\\
38.75	530.568520771494\\
39.125	525.536328442504\\
39.5	520.510542478031\\
39.875	515.673373667875\\
40.25	510.850385724728\\
40.625	506.150635307126\\
41	501.559308886727\\
41.375	497.008509416472\\
41.75	492.554351422182\\
42.125	488.16669005719\\
42.5	483.947988433364\\
42.875	479.690668687407\\
43.25	475.550515351254\\
43.625	471.476728373078\\
44	467.540199677407\\
44.375	463.556912569818\\
44.75	459.701021792035\\
45.125	455.836991604737\\
45.5	452.145716102847\\
45.875	448.41642828902\\
46.25	444.808063495267\\
46.625	441.233821445381\\
47	437.724950804218\\
47.375	434.245490116748\\
47.75	430.89791938382\\
48.125	427.51073043261\\
48.5	424.211420391177\\
48.875	420.98442428557\\
49.25	417.797050528068\\
49.625	414.630784360755\\
50	411.555179296581\\
};
\addlegendentry{battery swapping}
\end{axis}
\end{tikzpicture}%
\vspace*{-0.3in}
\caption{Charging arrival rate $\gamma$ (per hour) as a function of $C$.}
\label{gamma_C}
\end{subfigure}
\hfill
\begin{subfigure}[b]{0.32\linewidth}
\centering
%
%
\begin{tikzpicture}

\begin{axis}[%
width=1.694in,
height=1.03in,
at={(1.358in,0.0in)},
scale only axis,
xmin=12.5,
xmax=50,
xtick={12.5,   25, 37.5,   50},
xlabel style={font=\color{white!15!black}},
xlabel={Battery capacity},
ymin=7800,
ymax=8300,
ytick={7800,8000,8200},
yticklabels={{7.8K},{8K},{8.2K}},
ylabel style={font=\color{white!15!black}},
ylabel={Number of vehicles},
axis background/.style={fill=white},
legend style={at={(1,0)}, anchor=south east, legend cell align=left, align=left, font=\scriptsize, draw=white!12!black}
]
\addplot [color=black, line width=1.0pt]
  table[row sep=crcr]{%
12.5	7811.05613474758\\
12.875	7819.37601072682\\
13.25	7827.41606724694\\
13.625	7835.09591493714\\
14	7842.4302764213\\
14.375	7849.38631790511\\
14.75	7856.06819221001\\
15.125	7862.82251944741\\
15.5	7869.04375832397\\
15.875	7875.17524091787\\
16.25	7880.91423140873\\
16.625	7886.68423767352\\
17	7892.03679171431\\
17.375	7897.46368316491\\
17.75	7902.56218737703\\
18.125	7907.60283043478\\
18.5	7912.27453449933\\
18.875	7916.88873238757\\
19.25	7921.32156921884\\
19.625	7925.79263175829\\
20	7930.14667333513\\
20.375	7934.45803875932\\
20.75	7938.24261356721\\
21.125	7942.31720452391\\
21.5	7946.16454283113\\
21.875	7949.77010011797\\
22.25	7953.36519774553\\
22.625	7957.014451137\\
23	7960.50689374538\\
23.375	7963.78504232775\\
23.75	7967.15991363157\\
24.125	7970.16988855371\\
24.5	7973.48172011053\\
24.875	7976.54460554339\\
25.25	7979.47527513789\\
25.625	7982.35330510169\\
26	7985.29295625784\\
26.375	7988.20214880128\\
26.75	7990.78664059861\\
27.125	7993.51962228542\\
27.5	7996.08069880895\\
27.875	7998.69120153659\\
28.25	8001.28116363325\\
28.625	8003.69376023703\\
29	8006.1930962739\\
29.375	8008.49259262153\\
29.75	8010.74564154621\\
30.125	8013.13450886905\\
30.5	8015.49159317563\\
30.875	8017.72414441728\\
31.25	8019.80017862902\\
31.625	8021.87930913075\\
32	8024.01158372282\\
32.375	8025.86790081423\\
32.75	8027.97792724336\\
33.125	8029.90689308003\\
33.5	8031.97488296168\\
33.875	8033.86518627258\\
34.25	8035.78026243523\\
34.625	8037.66797031912\\
35	8039.48557915373\\
35.375	8041.25027440847\\
35.75	8042.99063315492\\
36.125	8044.87615932872\\
36.5	8046.54070364659\\
36.875	8048.17610541869\\
37.25	8049.77615991923\\
37.625	8051.54977856145\\
38	8053.21465157847\\
38.375	8054.68748415917\\
38.75	8056.21446359944\\
39.125	8057.70726113165\\
39.5	8059.31154822167\\
39.875	8060.77359986792\\
40.25	8062.35609781812\\
40.625	8063.80850746809\\
41	8065.13782928139\\
41.375	8066.61928356469\\
41.75	8068.07905577988\\
42.125	8069.40287143991\\
42.5	8070.77578390468\\
42.875	8072.03502919886\\
43.25	8073.52064974997\\
43.625	8074.75900479709\\
44	8076.03376084989\\
44.375	8077.27813183968\\
44.75	8078.58417284418\\
45.125	8079.64949186437\\
45.5	8081.15351800271\\
45.875	8082.31373669104\\
46.25	8083.477390808\\
46.625	8084.66997531319\\
47	8085.82356389\\
47.375	8086.96922013028\\
47.75	8088.08584992952\\
48.125	8089.13198265982\\
48.5	8090.43174109213\\
48.875	8091.52545717593\\
49.25	8092.45795916966\\
49.625	8093.71264491697\\
50	8094.74284015385\\
};
\addlegendentry{plug-in charging}
\addplot [color=blue, dashed, line width=1.0pt]
  table[row sep=crcr]{%
12.5	7884.04954127284\\
12.875	7899.79329179501\\
13.25	7915.21580239797\\
13.625	7928.73305330031\\
14	7941.73727258698\\
14.375	7954.0021036326\\
14.75	7965.94828142103\\
15.125	7976.65943412718\\
15.5	7987.84068316103\\
15.875	7998.49781204615\\
16.25	8007.79400089965\\
16.625	8016.73798430986\\
17	8025.50982785815\\
17.375	8033.81719629674\\
17.75	8041.86340512612\\
18.125	8050.35449732979\\
18.5	8057.19837000025\\
18.875	8063.8939494926\\
19.25	8071.0617334781\\
19.625	8077.93105586793\\
20	8083.80376905957\\
20.375	8090.27944901301\\
20.75	8095.56441909247\\
21.125	8100.94484839487\\
21.5	8107.34369993584\\
21.875	8112.44751771208\\
22.25	8116.88436272502\\
22.625	8122.21551005718\\
23	8126.24779048598\\
23.375	8131.51033539433\\
23.75	8135.90232759006\\
24.125	8140.71026639655\\
24.5	8144.01232573944\\
24.875	8148.07980170495\\
25.25	8152.81185944876\\
25.625	8155.73367343458\\
26	8159.53102663729\\
26.375	8163.41140348819\\
26.75	8166.21460305856\\
27.125	8169.97793341271\\
27.5	8173.14798326373\\
27.875	8176.87023059513\\
28.25	8178.96997662594\\
28.625	8182.3433663448\\
29	8185.32338928497\\
29.375	8188.00231386562\\
29.75	8191.07259298727\\
30.125	8193.91638146198\\
30.5	8196.34478728175\\
30.875	8198.73137124947\\
31.25	8201.70465647596\\
31.625	8203.81329429233\\
32	8206.27163127368\\
32.375	8208.93905499368\\
32.75	8211.09571274971\\
33.125	8212.94515783361\\
33.5	8215.39858727556\\
33.875	8217.11209504039\\
34.25	8219.88900149512\\
34.625	8222.30447934873\\
35	8223.59013237929\\
35.375	8225.81944221248\\
35.75	8228.38264350936\\
36.125	8229.41095783883\\
36.5	8230.73887955632\\
36.875	8233.65529251419\\
37.25	8234.14686845611\\
37.625	8236.25856541933\\
38	8238.30902056344\\
38.375	8239.58050806716\\
38.75	8241.47322010227\\
39.125	8243.60537948733\\
39.5	8244.27091866648\\
39.875	8246.4506036188\\
40.25	8247.35519856869\\
40.625	8248.80711420836\\
41	8250.61485984804\\
41.375	8251.67808820295\\
41.75	8252.9758150058\\
42.125	8254.03291833814\\
42.5	8256.66424215339\\
42.875	8257.30268341583\\
43.25	8258.69071464696\\
43.625	8259.97504364178\\
44	8262.46414563312\\
44.375	8262.87751907064\\
44.75	8264.38570898212\\
45.125	8264.54322591558\\
45.5	8266.71710411128\\
45.875	8267.0280598934\\
46.25	8268.46884829529\\
46.625	8269.43841320118\\
47	8270.55418461769\\
47.375	8271.14909115572\\
47.75	8273.24031971641\\
48.125	8273.50993577542\\
48.5	8274.47144314894\\
48.875	8275.85415720107\\
49.25	8277.03255951002\\
49.625	8277.64911022966\\
50	8279.13727678258\\
};
\addlegendentry{battery swapping}
\end{axis}
\end{tikzpicture}%
\vspace*{-0.3in}
\caption{Number of operating vehicles $N_1$ as a function of $C$.}
\label{N1_C}
\end{subfigure}
\hfill
\begin{subfigure}[b]{0.32\linewidth}
\centering
%
%
\begin{tikzpicture}

\begin{axis}[%
width=1.694in,
height=1.03in,
at={(1.358in,0.0in)},
scale only axis,
xmin=12.5,
xmax=50,
xtick={12.5,   25, 37.5,   50},
xlabel style={font=\color{white!15!black}},
xlabel={Battery capacity},
ymin=200,
ymax=1200,
ytick={400,800,1200},
yticklabels={{0.4K},{0.8K},{1.2K}},
ylabel style={font=\color{white!15!black}},
ylabel={Number of vehicles},
axis background/.style={fill=white},
legend style={at={(1,0.3)}, anchor=south east, legend cell align=left, align=left, font=\scriptsize, draw=white!12!black}
]
\addplot [color=black, line width=1.0pt]
  table[row sep=crcr]{%
12.5	1162.33646047609\\
12.875	1157.5786022276\\
13.25	1153.03779965831\\
13.625	1148.69642370275\\
14	1144.47431147116\\
14.375	1140.44262377782\\
14.75	1136.62871836464\\
15.125	1132.92970155028\\
15.5	1129.33836369453\\
15.875	1125.86578007693\\
16.25	1122.44101575409\\
16.625	1119.26907613448\\
17	1116.1132046323\\
17.375	1113.12062091156\\
17.75	1110.12534739896\\
18.125	1107.28543596732\\
18.5	1104.52933174709\\
18.875	1101.92218945762\\
19.25	1099.3546531982\\
19.625	1096.83854376468\\
20	1094.43547004126\\
20.375	1091.91050670106\\
20.75	1089.72488323535\\
21.125	1087.44593121429\\
21.5	1085.20330690403\\
21.875	1083.10978072993\\
22.25	1081.00870583864\\
22.625	1078.96524910554\\
23	1077.02070378947\\
23.375	1075.09895868072\\
23.75	1073.23927241189\\
24.125	1071.42634427256\\
24.5	1069.57265184244\\
24.875	1067.85225290245\\
25.25	1066.15039355127\\
25.625	1064.45209918395\\
26	1062.84371181077\\
26.375	1061.20059470548\\
26.75	1059.68596118179\\
27.125	1058.09902696321\\
27.5	1056.57253034071\\
27.875	1055.10406915624\\
28.25	1053.75157765509\\
28.625	1052.25240073387\\
29	1050.90671155719\\
29.375	1049.59255956933\\
29.75	1048.22296355352\\
30.125	1046.89342485888\\
30.5	1045.55479562977\\
30.875	1044.29184238393\\
31.25	1043.08463124959\\
31.625	1041.88329972796\\
32	1040.67868212302\\
32.375	1039.52328239688\\
32.75	1038.3458570099\\
33.125	1037.23785343703\\
33.5	1036.06497947689\\
33.875	1034.99073294159\\
34.25	1033.95145196131\\
34.625	1032.89432026195\\
35	1031.78796816038\\
35.375	1030.76738617044\\
35.75	1029.78656896824\\
36.125	1028.73780461588\\
36.5	1027.76922321593\\
36.875	1026.84618955373\\
37.25	1025.92939819739\\
37.625	1024.89275276961\\
38	1023.99670959522\\
38.375	1023.11278923204\\
38.75	1022.27904904886\\
39.125	1021.41732778143\\
39.5	1020.51248472051\\
39.875	1019.62219977279\\
40.25	1018.72761549733\\
40.625	1017.92684590353\\
41	1017.20509693003\\
41.375	1016.33110409104\\
41.75	1015.52847188077\\
42.125	1014.78393510807\\
42.5	1014.01775373487\\
42.875	1013.20069516481\\
43.25	1012.48328203578\\
43.625	1011.67887375345\\
44	1011.01083161178\\
44.375	1010.31618294291\\
44.75	1009.52117637927\\
45.125	1008.88094261061\\
45.5	1008.12337436125\\
45.875	1007.43658447339\\
46.25	1006.74410004566\\
46.625	1006.10489100449\\
47	1005.42257762914\\
47.375	1004.77782709525\\
47.75	1004.13894115974\\
48.125	1003.53597988932\\
48.5	1002.84867791993\\
48.875	1002.17804443003\\
49.25	1001.58525359892\\
49.625	1000.96896962267\\
50	1000.38137531115\\
};
\addlegendentry{plug-in charging}
\addplot [color=blue, dashed, line width=1.0pt]
  table[row sep=crcr]{%
12.5	585.538111628479\\
12.875	572.903618101274\\
13.25	560.877582780118\\
13.625	549.600854050904\\
14	538.871927478583\\
14.375	528.681722327751\\
14.75	518.948464273101\\
15.125	509.750774016176\\
15.5	500.860893727061\\
15.875	492.371829852929\\
16.25	484.360596022324\\
16.625	476.687750281605\\
17	469.313092212002\\
17.375	462.254490543619\\
17.75	455.471649788509\\
18.125	448.865858203572\\
18.5	442.674248966904\\
18.875	436.704933726229\\
19.25	430.87194117842\\
19.625	425.257205531171\\
20	419.935001841471\\
20.375	414.703398775572\\
20.75	409.770396186792\\
21.125	404.972820610311\\
21.5	400.191565767981\\
21.875	395.699555678568\\
22.25	391.414311907\\
22.625	387.140844759676\\
23	383.139401427085\\
23.375	379.097453367467\\
23.75	375.265400968297\\
24.125	371.481915160605\\
24.5	367.976662876107\\
24.875	364.469039397589\\
25.25	360.967023880437\\
25.625	357.769032573673\\
26	354.543867707002\\
26.375	351.384902731324\\
26.75	348.430607036264\\
27.125	345.429249786897\\
27.5	342.567903919797\\
27.875	339.703754472315\\
28.25	337.100281283268\\
28.625	334.401045780102\\
29	331.807499693998\\
29.375	329.306365981448\\
29.75	326.810620931239\\
30.125	324.394239967418\\
30.5	322.078318484467\\
30.875	319.815481352945\\
31.25	317.526765733058\\
31.625	315.388763309939\\
32	313.25085359239\\
32.375	311.128874624972\\
32.75	309.109784538409\\
33.125	307.167446645426\\
33.5	305.188247990969\\
33.875	303.336492643211\\
34.25	301.388989644872\\
34.625	299.520058501032\\
35	297.823151223898\\
35.375	296.042062188158\\
35.75	294.250865418025\\
36.125	292.678866997119\\
36.5	291.098949883618\\
36.875	289.351385393877\\
37.25	287.930285125886\\
37.625	286.335990506216\\
38	284.775113794593\\
38.375	283.335543164888\\
38.75	281.843447020461\\
39.125	280.345500466209\\
39.5	279.051770058164\\
39.875	277.59337519664\\
40.25	276.314285045924\\
40.625	274.988824336675\\
41	273.639835923556\\
41.375	272.403063918876\\
41.75	271.156776240113\\
42.125	269.959217905792\\
42.5	268.585379534732\\
42.875	267.476307319544\\
43.25	266.291905333855\\
43.625	265.137527963524\\
44	263.850643655353\\
44.375	262.837062347819\\
44.75	261.703771619106\\
45.125	260.753444382294\\
45.5	259.568189309554\\
45.875	258.628563669334\\
46.25	257.563292327788\\
46.625	256.570647573193\\
47	255.573408393318\\
47.375	254.654208547387\\
47.75	253.562616784026\\
48.125	252.709606507874\\
48.5	251.783432122393\\
48.875	250.817001524973\\
49.25	249.887910177704\\
49.625	249.040315075716\\
50	248.095898223905\\
};
\addlegendentry{battery swapping}
\end{axis}
\end{tikzpicture}%
\vspace*{-0.3in}
\caption{Number of non-operating vehicles $N_2$ as a function of $C$.}
\label{N2_C}
\end{subfigure}
\begin{subfigure}[b]{0.32\linewidth}
\centering
%
%
\begin{tikzpicture}

\begin{axis}[%
width=1.694in,
height=1.03in,
at={(1.358in,0.0in)},
scale only axis,
xmin=12.5,
xmax=50,
xtick={12.5,   25, 37.5,   50},
xlabel style={font=\color{white!15!black}},
xlabel={Battery capacity},
ymin=100,
ymax=450,
ylabel style={font=\color{white!15!black}},
ylabel={Number of vehicles},
axis background/.style={fill=white},
legend style={at={(1,1)}, anchor=north east, legend cell align=left, align=left, font=\scriptsize, draw=white!12!black}
]
\addplot [color=black, line width=1.0pt]
  table[row sep=crcr]{%
12.5	308.184507567093\\
12.875	303.943902036698\\
13.25	299.907375666099\\
13.625	296.044702555779\\
14	292.255740945598\\
14.375	288.633448073068\\
14.75	285.227944891229\\
15.125	281.945872237797\\
15.5	278.712340209282\\
15.875	275.593219026378\\
16.25	272.466812096455\\
16.625	269.642239463815\\
17	266.771192492323\\
17.375	264.092886766313\\
17.75	261.357270171892\\
18.125	258.792855073036\\
18.5	256.280303606839\\
18.875	253.932277760813\\
19.25	251.600205739641\\
19.625	249.319849813436\\
20	247.154776345906\\
20.375	244.819427503967\\
20.75	242.847727990657\\
21.125	240.771191222695\\
21.5	238.708414789428\\
21.875	236.798613078749\\
22.25	234.867475789424\\
22.625	233.000626288948\\
23	231.232591330585\\
23.375	229.466367530388\\
23.75	227.773714817982\\
24.125	226.102397269916\\
24.5	224.39457557173\\
24.875	222.822796539096\\
25.25	221.256367053846\\
25.625	219.68310577476\\
26	218.215397039736\\
26.375	216.696514084845\\
26.75	215.304595426532\\
27.125	213.828790788784\\
27.5	212.408646383998\\
27.875	211.055229595719\\
28.25	209.833388062746\\
28.625	208.419842021508\\
29	207.189392751129\\
29.375	205.979482214602\\
29.75	204.696137689905\\
30.125	203.463231669072\\
30.5	202.213069851792\\
30.875	201.041282634683\\
31.25	199.922662662295\\
31.625	198.806980351649\\
32	197.686109265018\\
32.375	196.605223037435\\
32.75	195.507616900473\\
33.125	194.478854772048\\
33.5	193.377113543199\\
33.875	192.379040786938\\
34.25	191.419641114975\\
34.625	190.434727296901\\
35	189.38596358933\\
35.375	188.4319454665\\
35.75	187.520263173119\\
36.125	186.532930764037\\
36.5	185.626416988959\\
36.875	184.768931539508\\
37.25	183.914562729524\\
37.625	182.925791137431\\
38	182.093016922409\\
38.375	181.264091688852\\
38.75	180.492950540043\\
39.125	179.685839690607\\
39.5	178.830978951048\\
39.875	177.985447750686\\
40.25	177.137347437198\\
40.625	176.390359897782\\
41	175.727701909359\\
41.375	174.893601529971\\
41.75	174.138701442464\\
42.125	173.443875740132\\
42.5	172.724042711155\\
42.875	171.94051268119\\
43.25	171.276601991479\\
43.625	170.503581372628\\
44	169.885786364729\\
44.375	169.235040624286\\
44.75	168.470284792189\\
45.125	167.87299972164\\
45.5	167.153752012974\\
45.875	166.503390076286\\
46.25	165.845268606235\\
46.625	165.247007518114\\
47	164.597437664995\\
47.375	163.988900570922\\
47.75	163.384894761687\\
48.125	162.818386979947\\
48.5	162.161911386791\\
48.875	161.51770352098\\
49.25	160.95627896253\\
49.625	160.375089976297\\
50	159.819599126201\\
};
\addlegendentry{plug-in charging}
\addplot [color=blue, dashed, line width=1.0pt]
  table[row sep=crcr]{%
12.5	404.551231809878\\
12.875	393.57184115581\\
13.25	383.159386777282\\
13.625	373.293296968604\\
14	363.915758794883\\
14.375	354.995169012927\\
14.75	346.494550371346\\
15.125	338.396903045128\\
15.5	330.654661282683\\
15.875	323.256657007655\\
16.25	316.190898988066\\
16.625	309.425075257509\\
17	302.938778282998\\
17.375	296.718562468166\\
17.75	290.746489950044\\
18.125	285.000621624117\\
18.5	279.49048382477\\
18.875	274.18724784408\\
19.25	269.073206757876\\
19.625	264.146026945895\\
20	259.402862883865\\
20.375	254.818388458394\\
20.75	250.401964406729\\
21.125	246.132861334336\\
21.5	241.995092189501\\
21.875	238.003771912747\\
22.25	234.146150553987\\
22.625	230.401880817986\\
23	226.785224144295\\
23.375	223.26812234812\\
23.75	219.86432013451\\
24.125	216.557766473168\\
24.5	213.360472690603\\
24.875	210.248826520149\\
25.25	207.219969267022\\
25.625	204.290579443727\\
26	201.435011663655\\
26.375	198.65652376431\\
26.75	195.961069419941\\
27.125	193.329653161951\\
27.5	190.771603466341\\
27.875	188.275447258845\\
28.25	185.854996885311\\
28.625	183.48615531638\\
29	181.179064009924\\
29.375	178.930810822737\\
29.75	176.734327463073\\
30.125	174.592112864348\\
30.5	172.503528918913\\
30.875	170.464153145647\\
31.25	168.467995450248\\
31.625	166.523328066624\\
32	164.62033775552\\
32.375	162.758555265147\\
32.75	160.941289799709\\
33.125	159.165792387104\\
33.5	157.424869644607\\
33.875	155.725914550323\\
34.25	154.056315911395\\
34.625	152.42396905102\\
35	150.832422501061\\
35.375	149.267851034717\\
35.75	147.733089945724\\
36.125	146.238372564037\\
36.5	144.771663223699\\
36.875	143.324618065272\\
37.25	141.919884959752\\
37.625	140.533057699303\\
38	139.173163615691\\
38.375	137.843483006427\\
38.75	136.535328205103\\
39.125	135.250227015197\\
39.5	133.996900695532\\
39.875	132.758274696368\\
40.25	131.548971790084\\
40.625	130.358481227575\\
41	129.18729337695\\
41.375	128.040695717967\\
41.75	126.912935033713\\
42.125	125.805972765452\\
42.5	124.709982526928\\
42.875	123.642815222989\\
43.25	122.589899346791\\
43.625	121.555158532988\\
44	120.531634582786\\
44.375	119.535214008275\\
44.75	118.549736534268\\
45.125	117.586804396152\\
45.5	116.629654580628\\
45.875	115.696676568313\\
46.25	114.773123392063\\
46.625	113.866308691868\\
47	112.972935660103\\
47.375	112.095797263759\\
47.75	111.225257893884\\
48.125	110.376279317845\\
48.5	109.537016607195\\
48.875	108.708440609014\\
49.25	107.893121431718\\
49.625	107.092354385587\\
50	106.299501385258\\
};
\addlegendentry{battery swapping}
\end{axis}
\end{tikzpicture}%
\vspace*{-0.3in}
\caption{Number of searching vehicles $N_2^m$ as a function of $C$.}
\label{N2m_C}
\end{subfigure}
\hfill
\begin{subfigure}[b]{0.32\linewidth}
\centering
%
%
\begin{tikzpicture}

\begin{axis}[%
width=1.694in,
height=1.03in,
at={(1.358in,0.0in)},
scale only axis,
xmin=12.5,
xmax=50,
xtick={12.5,   25, 37.5,   50},
xlabel style={font=\color{white!15!black}},
xlabel={Battery capacity},
ymin=0,
ymax=140,
ylabel style={font=\color{white!15!black}},
ylabel={Number of vehicles},
axis background/.style={fill=white},
legend style={at={(1,0.3)}, anchor=south east, legend cell align=left, align=left, font=\scriptsize, draw=white!12!black}
]
\addplot [color=black, line width=1.0pt]
  table[row sep=crcr]{%
12.5	39.0580791685371\\
12.875	38.1313068292586\\
13.25	37.225114608865\\
13.625	36.3632040534227\\
14	35.5741292951957\\
14.375	34.8300580939989\\
14.75	34.0927396316841\\
15.125	33.3272144193795\\
15.5	32.6694714852517\\
15.875	32.0135968620575\\
16.25	31.4529611417796\\
16.625	30.8099008732715\\
17	30.2759568375803\\
17.375	29.6857472029017\\
17.75	29.1888790347639\\
18.125	28.6647931690132\\
18.5	28.2044822065575\\
18.875	27.7178638291512\\
19.25	27.271502956312\\
19.625	26.8158235201783\\
20	26.3580718645784\\
20.375	25.9700842571127\\
20.75	25.5741650157925\\
21.125	25.171164289203\\
21.5	24.8121635191245\\
21.875	24.4582014629952\\
22.25	24.1212182052131\\
22.625	23.7656747202892\\
23	23.4160503765236\\
23.375	23.1087463235259\\
23.75	22.7728341578125\\
24.125	22.4968352077437\\
24.5	22.1899372730515\\
24.875	21.8916245870923\\
25.25	21.6192381011588\\
25.625	21.3632187534376\\
26	21.0747712163425\\
26.375	20.8109633202982\\
26.75	20.5685255519522\\
27.125	20.3311874553162\\
27.5	20.1102965976868\\
27.875	19.8690528322394\\
28.25	19.6010566836534\\
28.625	19.4151307109484\\
29	19.1725067455497\\
29.375	18.9588813310213\\
29.75	18.7752808203933\\
30.125	18.5626007208358\\
30.5	18.3630089096342\\
30.875	18.1653523486038\\
31.25	17.9806459813589\\
31.625	17.7982758142153\\
32	17.6129940862167\\
32.375	17.4606347675369\\
32.75	17.2791782321848\\
33.125	17.1095649709651\\
33.5	16.9414301440203\\
33.875	16.7756849763375\\
34.25	16.5998629160316\\
34.625	16.4370131709149\\
35	16.302240101584\\
35.375	16.1542920271457\\
35.75	16.0019659518138\\
36.125	15.8503639708088\\
36.5	15.7121973511233\\
36.875	15.5685538020217\\
37.25	15.4312715655976\\
37.625	15.3046066150048\\
38	15.1578029181265\\
38.375	15.038165667759\\
38.75	14.8996877195735\\
39.125	14.7762409363304\\
39.5	14.6510233985216\\
39.875	14.544376912593\\
40.25	14.4241664728757\\
40.625	14.2995668056058\\
41	14.1735491677202\\
41.375	14.0686685861465\\
41.75	13.9501740635046\\
42.125	13.8373184351928\\
42.5	13.725407180519\\
42.875	13.6441213913704\\
43.25	13.5081271645159\\
43.625	13.4300243883273\\
44	13.3138246677815\\
44.375	13.2103274588848\\
44.75	13.1257368885053\\
45.125	13.0358020716352\\
45.5	12.9186973048718\\
45.875	12.8310846268185\\
46.25	12.745971040901\\
46.625	12.6453587180343\\
47	12.5620165509258\\
47.375	12.4718815013397\\
47.75	12.3855388440008\\
48.125	12.3009356114581\\
48.5	12.2055298091008\\
48.875	12.133977871993\\
49.25	12.0653589845262\\
49.625	11.9626513804532\\
50	11.8828921523729\\
};
\addlegendentry{plug-in charging}
\addplot [color=blue, dashed, line width=1.0pt]
  table[row sep=crcr]{%
12.5	126.746547561209\\
12.875	126.640990467329\\
13.25	126.487728345775\\
13.625	126.465182040817\\
14	126.427667185704\\
14.375	126.404982633814\\
14.75	126.355308893722\\
15.125	126.384072125451\\
15.5	126.30626361032\\
15.875	126.235417759413\\
16.25	126.268247168534\\
16.625	126.295653496941\\
17	126.299979242892\\
17.375	126.316681002366\\
17.75	126.324929900616\\
18.125	126.247104740159\\
18.5	126.322265760826\\
18.875	126.382500097567\\
19.25	126.358753796483\\
19.625	126.340299936527\\
20	126.40859445068\\
20.375	126.381913876579\\
20.75	126.467389048156\\
21.125	126.518655478667\\
21.5	126.4303108653\\
21.875	126.470023757875\\
22.25	126.566548749336\\
22.625	126.540480566214\\
23	126.646582316558\\
23.375	126.592220398769\\
23.75	126.622205284634\\
24.125	126.587522491242\\
24.5	126.71293807272\\
24.875	126.734469631229\\
25.25	126.663997807355\\
25.625	126.79175792399\\
26	126.804102108155\\
26.375	126.794208867277\\
26.75	126.898601171613\\
27.125	126.878761970342\\
27.5	126.917565858238\\
27.875	126.880666132897\\
28.25	127.024425578951\\
28.625	127.008464847487\\
29	127.029246760035\\
29.375	127.07643319987\\
29.75	127.068455361379\\
30.125	127.078788609027\\
30.5	127.129924994297\\
30.875	127.178153440929\\
31.25	127.1491158974\\
31.625	127.2151550138\\
32	127.232530230967\\
32.375	127.218151313123\\
32.75	127.257679390226\\
33.125	127.327382513362\\
33.5	127.318795232678\\
33.875	127.392283636478\\
34.25	127.333102871072\\
34.625	127.311270335723\\
35	127.418604726622\\
35.375	127.408028391772\\
35.75	127.352381313351\\
36.125	127.475105101497\\
36.5	127.557058253543\\
36.875	127.444163111175\\
37.25	127.616874404114\\
37.625	127.591155300455\\
38	127.568407644882\\
38.375	127.634865210695\\
38.75	127.622501456308\\
39.125	127.577395836261\\
39.5	127.704517946698\\
39.875	127.645988044677\\
40.25	127.736967065016\\
40.625	127.75865526553\\
41	127.733898917048\\
41.375	127.795417887027\\
41.75	127.825362825661\\
42.125	127.881022138434\\
42.5	127.743797393358\\
42.875	127.843803140308\\
43.25	127.850322142023\\
43.625	127.866478484766\\
44	127.734335749986\\
44.375	127.849951253884\\
44.75	127.830667691771\\
45.125	127.972073599317\\
45.5	127.867010858831\\
45.875	127.98467282472\\
46.25	127.963233485883\\
46.625	127.996544833146\\
47	128.009641039741\\
47.375	128.083561613069\\
47.75	127.974094910682\\
48.125	128.082969508942\\
48.5	128.106034835492\\
48.875	128.075746773107\\
49.25	128.068220395051\\
49.625	128.12693454477\\
50	128.077890862094\\
};
\addlegendentry{battery swapping}
\end{axis}
\end{tikzpicture}%
\vspace*{-0.3in}
\caption{Number of waiting vehicles $N_2^w$ as a function of $C$.}
\label{N2w_C}
\end{subfigure}
\hfill
\begin{subfigure}[b]{0.32\linewidth}
\centering
%
%
\begin{tikzpicture}

\begin{axis}[%
width=1.694in,
height=1.03in,
at={(1.358in,0.0in)},
scale only axis,
xmin=12.5,
xmax=50,
xtick={12.5,   25, 37.5,   50},
xlabel style={font=\color{white!15!black}},
xlabel={Battery capacity},
ymin=0,
ymax=900,
ylabel style={font=\color{white!15!black}},
ylabel={Number of vehicles},
axis background/.style={fill=white},
legend style={at={(1,0.3)}, anchor=south east, legend cell align=left, align=left, font=\scriptsize, draw=white!12!black}
]
\addplot [color=black, line width=1.0pt]
  table[row sep=crcr]{%
12.5	815.093873740458\\
12.875	815.503393361647\\
13.25	815.905309383351\\
13.625	816.288517093545\\
14	816.644441230363\\
14.375	816.979117610749\\
14.75	817.30803384173\\
15.125	817.656614893105\\
15.5	817.956551999992\\
15.875	818.258964188491\\
16.25	818.521242515858\\
16.625	818.81693579739\\
17	819.066055302394\\
17.375	819.341986942342\\
17.75	819.579198192306\\
18.125	819.827787725269\\
18.5	820.044545933692\\
18.875	820.272047867657\\
19.25	820.482944502245\\
19.625	820.702870431071\\
20	820.922621830779\\
20.375	821.120994939984\\
20.75	821.302990228898\\
21.125	821.503575702391\\
21.5	821.682728595482\\
21.875	821.852966188181\\
22.25	822.020011844\\
22.625	822.198948096302\\
23	822.372062082361\\
23.375	822.523844826807\\
23.75	822.692723436095\\
24.125	822.827111794898\\
24.5	822.988138997661\\
24.875	823.137831776267\\
25.25	823.27478839627\\
25.625	823.40577465575\\
26	823.553543554692\\
26.375	823.693117300336\\
26.75	823.812840203305\\
27.125	823.939048719105\\
27.5	824.053587359024\\
27.875	824.179786728285\\
28.25	824.317132908689\\
28.625	824.417428001415\\
29	824.544812060516\\
29.375	824.654196023705\\
29.75	824.751545043224\\
30.125	824.86759246897\\
30.5	824.978716868341\\
30.875	825.085207400643\\
31.25	825.181322605932\\
31.625	825.278043562092\\
32	825.379578771785\\
32.375	825.457424591913\\
32.75	825.559061877239\\
33.125	825.649433694017\\
33.5	825.74643578967\\
33.875	825.836007178311\\
34.25	825.931947930305\\
34.625	826.022579794135\\
35	826.099764469465\\
35.375	826.181148676796\\
35.75	826.264339843306\\
36.125	826.354509881035\\
36.5	826.43060887585\\
36.875	826.508704212201\\
37.25	826.583563902273\\
37.625	826.662355017177\\
38	826.745889754689\\
38.375	826.810531875429\\
38.75	826.88641078924\\
39.125	826.955247154492\\
39.5	827.030482370937\\
39.875	827.092375109509\\
40.25	827.16610158726\\
40.625	827.236919200143\\
41	827.303845852947\\
41.375	827.368833974924\\
41.75	827.439596374802\\
42.125	827.502740932743\\
42.5	827.568303843196\\
42.875	827.616061092253\\
43.25	827.698552879789\\
43.625	827.745267992498\\
44	827.81122057927\\
44.375	827.87081485974\\
44.75	827.925154698576\\
45.125	827.972140817335\\
45.5	828.050925043404\\
45.875	828.102109770284\\
46.25	828.152860398528\\
46.625	828.212524768341\\
47	828.263123413218\\
47.375	828.317045022987\\
47.75	828.368507554053\\
48.125	828.416657297913\\
48.5	828.481236724041\\
48.875	828.526363037057\\
49.25	828.563615651865\\
49.625	828.631228265916\\
50	828.678884032576\\
};
\addlegendentry{plug-in charging}
\addplot [color=blue, dashed, line width=1.0pt]
  table[row sep=crcr]{%
12.5	54.240332257392\\
12.875	52.6907864781357\\
13.25	51.2304676570606\\
13.625	49.8423750414824\\
14	48.5285014979951\\
14.375	47.2815706810095\\
14.75	46.0986050080337\\
15.125	44.9697988455972\\
15.5	43.8999688340588\\
15.875	42.8797550858608\\
16.25	41.9014498657239\\
16.625	40.9670215271546\\
17	40.0743346861116\\
17.375	39.2192470730866\\
17.75	38.4002299378492\\
18.125	37.6181318392955\\
18.5	36.8614993813078\\
18.875	36.1351857845818\\
19.25	35.4399806240607\\
19.625	34.7708786487492\\
20	34.1235445069261\\
20.375	33.5030964405993\\
20.75	32.9010427319061\\
21.125	32.3213037973082\\
21.5	31.7661627131796\\
21.875	31.2257600079456\\
22.25	30.7016126036773\\
22.625	30.1984833754765\\
23	29.707594966232\\
23.375	29.2371106205787\\
23.75	28.7788755491536\\
24.125	28.3366261961952\\
24.5	27.9032521127834\\
24.875	27.4857432462109\\
25.25	27.0830568060606\\
25.625	26.686695205956\\
26	26.3047539351925\\
26.375	25.9341700997373\\
26.75	25.57093644471\\
27.125	25.220834654603\\
27.5	24.8787345952186\\
27.875	24.5476410805732\\
28.25	24.2208588190057\\
28.625	23.9064256162347\\
29	23.5991889240389\\
29.375	23.2991219588415\\
29.75	23.0078381067872\\
30.125	22.7233384940426\\
30.5	22.4448645712579\\
30.875	22.1731747663697\\
31.25	21.9096543854101\\
31.625	21.6502802295146\\
32	21.3979856059028\\
32.375	21.1521680467018\\
32.75	20.9108153484742\\
33.125	20.6742717449599\\
33.5	20.4445831136835\\
33.875	20.2182944564102\\
34.25	19.9995708624042\\
34.625	19.78481911429\\
35	19.5721239962147\\
35.375	19.3661827616691\\
35.75	19.1653941589488\\
36.125	18.965389331585\\
36.5	18.7702284063759\\
36.875	18.5826042174304\\
37.25	18.39352576202\\
37.625	18.2117775064574\\
38	18.0335425340203\\
38.375	17.8571949477665\\
38.75	17.6856173590498\\
39.125	17.5178776147501\\
39.5	17.3503514159344\\
39.875	17.1891124555958\\
40.25	17.0283461908243\\
40.625	16.8716878435709\\
41	16.7186436295576\\
41.375	16.5669503138824\\
41.75	16.4184783807394\\
42.125	16.2722230019063\\
42.5	16.1315996144455\\
42.875	15.9896889562469\\
43.25	15.8516838450418\\
43.625	15.7158909457693\\
44	15.5846733225802\\
44.375	15.4518970856606\\
44.75	15.3233673930678\\
45.125	15.1945663868246\\
45.5	15.0715238700949\\
45.875	14.9472142763007\\
46.25	14.8269354498422\\
46.625	14.7077940481794\\
47	14.5908316934739\\
47.375	14.4748496705583\\
47.75	14.3632639794607\\
48.125	14.250357681087\\
48.5	14.1403806797059\\
48.875	14.0328141428523\\
49.25	13.9265683509356\\
49.625	13.8210261453585\\
50	13.7185059765527\\
};
\addlegendentry{battery swapping}
\end{axis}
\end{tikzpicture}%
\vspace*{-0.3in}
\caption{Number of charging vehicles $N_2^s$ as a function of $C$.}
\label{N2s_C}
\end{subfigure}
\begin{subfigure}[b]{0.32\linewidth}
\centering
%
%
\begin{tikzpicture}

\begin{axis}[%
width=1.694in,
height=1.03in,
at={(1.358in,0.0in)},
scale only axis,
xmin=12.5,
xmax=50,
xtick={12.5,   25, 37.5,   50},
xlabel style={font=\color{white!15!black}},
xlabel={Battery capacity},
ymin=10,
ymax=24,
ylabel style={font=\color{white!15!black}},
ylabel={Searching time},
axis background/.style={fill=white},
legend style={at={(0,1)}, anchor=north west, legend cell align=left, align=left, font=\scriptsize, draw=white!12!black}
]
\addplot [color=black, line width=1.0pt]
  table[row sep=crcr]{%
12.5	11.600701872822\\
12.875	11.7783913152837\\
13.25	11.9545806035258\\
13.625	12.1288938978764\\
14	12.2978489780801\\
14.375	12.4656415343859\\
14.75	12.6348300373135\\
15.125	12.8015114366764\\
15.5	12.9636930528428\\
15.875	13.1238896103057\\
16.25	13.2772489968816\\
16.625	13.4379753994986\\
17	13.5906429436844\\
17.375	13.7463496183182\\
17.75	13.8935457262311\\
18.125	14.0436093960623\\
18.5	14.1912481003603\\
18.875	14.3422746581192\\
19.25	14.4891612931821\\
19.625	14.6336158010627\\
20	14.7797762350882\\
20.375	14.9110222361835\\
20.75	15.0598213753497\\
21.125	15.197174840627\\
21.5	15.3310926869795\\
21.875	15.4704935182035\\
22.25	15.6042027342311\\
22.625	15.7376484026954\\
23	15.8737528179033\\
23.375	16.0063843982884\\
23.75	16.1398926483302\\
24.125	16.2717766710017\\
24.5	16.3966807723913\\
24.875	16.5280350991688\\
25.25	16.6564874775945\\
25.625	16.780994861523\\
26	16.9097807820722\\
26.375	17.0313877000519\\
26.75	17.1600918489775\\
27.125	17.2787340494962\\
27.5	17.3988486317616\\
27.875	17.5210493822692\\
28.25	17.6510196197852\\
28.625	17.7626792107087\\
29	17.8863749650635\\
29.375	18.0094746070179\\
29.75	18.1236038006175\\
30.125	18.2389499478209\\
30.5	18.3500560921442\\
30.875	18.465645369726\\
31.25	18.5837664346858\\
31.625	18.6996273197661\\
32	18.8123691316047\\
32.375	18.9269758987443\\
32.75	19.0369739871831\\
33.125	19.1515384615635\\
33.5	19.2563624441482\\
33.875	19.3693177305206\\
34.25	19.4838100562042\\
34.625	19.5936384611643\\
35	19.6949287031336\\
35.375	19.803720043935\\
35.75	19.914816965914\\
36.125	20.015573898089\\
36.5	20.1232133354191\\
36.875	20.2341340944681\\
37.25	20.3435485470922\\
37.625	20.4359282232286\\
38	20.5435702991055\\
38.375	20.6502468066385\\
38.75	20.7614258014952\\
39.125	20.8668683381609\\
39.5	20.9647366656719\\
39.875	21.0621282735099\\
40.25	21.1570139582031\\
40.625	21.2622583492608\\
41	21.3761811926548\\
41.375	21.4676175393392\\
41.75	21.5668421940122\\
42.125	21.6720763253164\\
42.5	21.7725328631336\\
42.875	21.8637432752963\\
43.25	21.9676213014407\\
43.625	22.0568412744083\\
44	22.1640688979207\\
44.375	22.265741488832\\
44.75	22.3509687214608\\
45.125	22.457086960125\\
45.5	22.5445493460963\\
45.875	22.6405172963347\\
46.25	22.7339759333774\\
46.625	22.8339861578352\\
47	22.9257568137528\\
47.375	23.0217406929514\\
47.75	23.1170694564734\\
48.125	23.2164838702534\\
48.5	23.3012379704991\\
48.875	23.3868456112164\\
49.25	23.4833136534875\\
49.625	23.5747564891033\\
50	23.6692691307356\\
};
\addlegendentry{plug-in charging}
\addplot [color=blue, dashed, line width=1.0pt]
  table[row sep=crcr]{%
12.5	14.9169894421045\\
12.875	14.9389245240864\\
13.25	14.958262311489\\
13.625	14.9789530156989\\
14	14.9980216805136\\
14.375	15.016217266044\\
14.75	15.0327564277037\\
15.125	15.0499629409954\\
15.5	15.064004374698\\
15.875	15.077355566065\\
16.25	15.0921221103958\\
16.625	15.1060567120999\\
17	15.1188425537598\\
17.375	15.1312727608064\\
17.75	15.1429556760789\\
18.125	15.1523006427665\\
18.5	15.1643578539021\\
18.875	15.1756379213675\\
19.25	15.184726516199\\
19.625	15.1935203947109\\
20	15.2037466583352\\
20.375	15.2116320896001\\
20.75	15.2215214847218\\
21.125	15.2303794969332\\
21.5	15.2360292538008\\
21.875	15.2440659155893\\
22.25	15.253019675321\\
22.625	15.2591690088045\\
23	15.2678279343769\\
23.375	15.272926606569\\
23.75	15.2795629390722\\
24.125	15.2846542120985\\
24.5	15.2928749543753\\
24.875	15.2987550408799\\
25.25	15.3025539732022\\
25.625	15.310294352081\\
26	15.3154834415052\\
26.375	15.3200602140203\\
26.75	15.3268590568556\\
27.125	15.3309480681019\\
27.5	15.3361178990996\\
27.875	15.3395959017703\\
28.25	15.3466892544267\\
28.625	15.3503629745281\\
29	15.3546856710291\\
29.375	15.3594466897785\\
29.75	15.3629668848319\\
30.125	15.3667660154886\\
30.5	15.3713138585665\\
30.875	15.3757100588222\\
31.25	15.378425646224\\
31.625	15.3830182613167\\
32	15.3865266373586\\
32.375	15.3893024020794\\
32.75	15.3931147224685\\
33.125	15.3974751179235\\
33.5	15.4001545318127\\
33.875	15.4044560866458\\
34.25	15.4059621550175\\
34.625	15.4081741329572\\
35	15.4129845621489\\
35.375	15.4153095498157\\
35.75	15.4166503146761\\
36.125	15.4216051152182\\
36.5	15.4256687866966\\
36.875	15.4256762279674\\
37.25	15.4315041929369\\
37.625	15.4332060832034\\
38	15.4349222681168\\
38.375	15.4384250616773\\
38.75	15.4402671315556\\
39.125	15.4413942133397\\
39.5	15.4460157587899\\
39.875	15.446786451519\\
40.25	15.4505869584645\\
40.625	15.4529271091569\\
41	15.454279215397\\
41.375	15.457364607495\\
41.75	15.4597682063638\\
42.125	15.4626657649313\\
42.5	15.4615767199247\\
42.875	15.4653183762757\\
43.25	15.4671138467268\\
43.625	15.4690763574828\\
44	15.4679706257495\\
44.375	15.4719143346099\\
44.75	15.4730658729619\\
45.125	15.477480752345\\
45.5	15.4768231249723\\
45.875	15.480700875715\\
46.25	15.4817054111184\\
46.625	15.4838051605657\\
47	15.4854689620788\\
47.375	15.4883539124778\\
47.75	15.4874627456454\\
48.125	15.4910187923684\\
48.5	15.4927959986822\\
48.875	15.4934626087647\\
49.25	15.494573927032\\
49.625	15.4970192892012\\
50	15.4972416919076\\
};
\addlegendentry{battery swapping}
\end{axis}
\end{tikzpicture}%
\vspace*{-0.3in}
\caption{Vehicle searching time $t_m$ (minute) as a function of $C$.}
\label{tm_C}
\end{subfigure}
\hfill
\begin{subfigure}[b]{0.32\linewidth}
\centering
%
%
\begin{tikzpicture}

\begin{axis}[%
width=1.694in,
height=1.03in,
at={(1.358in,0.0in)},
scale only axis,
xmin=12.5,
xmax=50,
xtick={12.5,   25, 37.5,   50},
xlabel style={font=\color{white!15!black}},
xlabel={Battery capacity},
ymin=0,
ymax=20,
ylabel style={font=\color{white!15!black}},
ylabel={Waiting time},
axis background/.style={fill=white},
legend style={at={(0,1)}, anchor=north west, legend cell align=left, align=left, font=\scriptsize, draw=white!12!black}
]
\addplot [color=black, line width=1.0pt]
  table[row sep=crcr]{%
12.5	1.47022683176453\\
12.875	1.47765903572538\\
13.25	1.48382690515292\\
13.625	1.48979339925089\\
14	1.49692617904963\\
14.375	1.5042574646848\\
14.75	1.51021657754108\\
15.125	1.51319369620854\\
15.5	1.51954879434254\\
15.875	1.52450380575677\\
16.25	1.53269601371052\\
16.625	1.53545190404626\\
17	1.54240686677515\\
17.375	1.54517853444946\\
17.75	1.55165771857122\\
18.125	1.55551882825717\\
18.5	1.56179307930543\\
18.875	1.56552455433999\\
19.25	1.57051224930394\\
19.625	1.57393187536825\\
20	1.57620423083236\\
20.375	1.58173927527501\\
20.75	1.58594177573022\\
21.125	1.58877223933055\\
21.5	1.59356585318431\\
21.875	1.59789976081637\\
22.25	1.60257344191936\\
22.625	1.60521385181563\\
23	1.60747493901974\\
23.375	1.61194636319811\\
23.75	1.61366775309889\\
24.125	1.61901635155038\\
24.5	1.62143544200473\\
24.875	1.6238264000506\\
25.25	1.62752635552156\\
25.625	1.63187816770304\\
26	1.63311006526173\\
26.375	1.6356496836898\\
26.75	1.639341635836\\
27.125	1.64288999463065\\
27.5	1.64727760568858\\
27.875	1.64945761599666\\
28.25	1.64882547665975\\
28.625	1.65466366017572\\
29	1.65513610575145\\
29.375	1.65763836396462\\
29.75	1.66234573194355\\
30.125	1.66399767993153\\
30.5	1.66637222687584\\
30.875	1.66848793486343\\
31.25	1.6713869293887\\
31.625	1.67409174502595\\
32	1.67610231945271\\
32.375	1.68091675448024\\
32.75	1.68250870089352\\
33.125	1.68488492996388\\
33.5	1.68701617889542\\
33.875	1.68902792697501\\
34.25	1.68963108556202\\
34.625	1.69118783125195\\
35	1.69532868443346\\
35.375	1.69777516238848\\
35.75	1.69942286573562\\
36.125	1.70079422475089\\
36.5	1.7033130542172\\
36.875	1.70491977553971\\
37.25	1.70691661159959\\
37.625	1.70978537429974\\
38	1.71008968433537\\
38.375	1.71320105192924\\
38.75	1.71385508480978\\
39.125	1.71596089420433\\
39.5	1.71757068732862\\
39.875	1.72112684527113\\
40.25	1.72280039086772\\
40.625	1.72368310762274\\
41	1.72412403884087\\
41.375	1.72688305262758\\
41.75	1.72771015354109\\
42.125	1.72899400446123\\
42.5	1.73014060004073\\
42.875	1.73496962796096\\
43.25	1.73252749407389\\
43.625	1.73734717980723\\
44	1.73698184848742\\
44.375	1.73804275460516\\
44.75	1.74139276254555\\
45.125	1.74385482598814\\
45.5	1.74238511172849\\
45.875	1.74472359566446\\
46.25	1.74721052536844\\
46.625	1.74734750277885\\
47	1.74968541808662\\
47.375	1.75087716837814\\
47.75	1.75241023431351\\
48.125	1.7540062796938\\
48.5	1.7538270991362\\
48.875	1.75693104196072\\
49.25	1.76032032550598\\
49.625	1.75848127848221\\
50	1.75985532402647\\
};
\addlegendentry{plug-in charging}
\addplot [color=blue, dashed, line width=1.0pt]
  table[row sep=crcr]{%
12.5	4.67351663554517\\
12.875	4.80695009249402\\
13.25	4.93798843268385\\
13.625	5.07460496958917\\
14	5.21045007709242\\
14.375	5.34690285509423\\
14.75	5.48195802765406\\
15.125	5.62084222610769\\
15.5	5.75427577581002\\
15.875	5.88787960689814\\
16.25	6.02691542050067\\
16.625	6.16572300298802\\
17	6.30328514407817\\
17.375	6.44156583460005\\
17.75	6.57938403520357\\
18.125	6.71203478575099\\
18.5	6.85388646045598\\
18.875	6.99498272132818\\
19.25	7.13085907900843\\
19.625	7.26701796712126\\
20	7.40887831421044\\
20.375	7.54449154278338\\
20.75	7.6877435209987\\
21.125	7.82880890400243\\
21.5	7.96006190655475\\
21.875	8.10036480941977\\
22.25	8.24494467982287\\
22.625	8.38058514348931\\
23	8.5262090358048\\
23.375	8.65969432079697\\
23.75	8.79966314655803\\
24.125	8.93455146105147\\
24.5	9.08230607389804\\
24.875	9.22183318791642\\
25.25	9.35374457280691\\
25.625	9.50224499103151\\
26	9.6411547829388\\
26.375	9.77815818895712\\
26.75	9.92522127189167\\
27.125	10.0614245093737\\
27.5	10.2028955992504\\
27.875	10.3375037720679\\
28.25	10.4888457117199\\
28.625	10.6254667164661\\
29	10.7655603901403\\
29.375	10.9082594120373\\
29.75	11.0456666785998\\
30.125	11.1848695685578\\
30.5	11.3281971108076\\
30.875	11.4713526394805\\
31.25	11.6066747252819\\
31.625	11.751825257243\\
32	11.8920100774223\\
32.375	12.0288521755537\\
32.75	12.1714698608834\\
33.125	12.3174720816372\\
33.5	12.455015054571\\
33.875	12.6016844705799\\
34.25	12.7335835100779\\
34.625	12.8695915388754\\
35	13.0204166651781\\
35.375	13.1577843666689\\
35.75	13.2898264713108\\
36.125	13.4429199287989\\
36.5	13.5914231294291\\
36.875	13.7165019089878\\
37.25	13.8762819108477\\
37.625	14.0119387308805\\
38	14.1479032646219\\
38.375	14.2950632038274\\
38.75	14.432349051248\\
39.125	14.5653941238682\\
39.5	14.7206837354793\\
39.875	14.8519579907829\\
40.25	15.0028623606263\\
40.625	15.1447390978388\\
41	15.2804140990507\\
41.375	15.4277541087257\\
41.75	15.5709146562099\\
42.125	15.7177076694994\\
42.5	15.8377098919523\\
42.875	15.9907804948716\\
43.25	16.1308190841836\\
43.625	16.2722532150413\\
44	16.3923019887642\\
44.375	16.548123579276\\
44.75	16.684409426821\\
45.125	16.8444521997394\\
45.5	16.9680268512923\\
45.875	17.1248863445605\\
46.25	17.2609146264605\\
46.625	17.4052674947525\\
47	17.5465859286138\\
47.375	17.6973943810402\\
47.75	17.8196397550979\\
48.125	17.9761059161249\\
48.5	18.1191776568432\\
48.875	18.2537508826542\\
49.25	18.3919278845814\\
49.625	18.5408714515454\\
50	18.6722797775503\\
};
\addlegendentry{battery swapping}
\end{axis}
\end{tikzpicture}%
\vspace*{-0.3in}
\caption{Vehicle waiting time $t_w$ (minute) as a function of $C$.}
\label{tw_C}
\end{subfigure}
\hfill
\begin{subfigure}[b]{0.32\linewidth}
\centering
%
%
\begin{tikzpicture}

\begin{axis}[%
width=1.694in,
height=1.03in,
at={(1.358in,0.0in)},
scale only axis,
xmin=12.5,
xmax=50,
xtick={12.5,   25, 37.5,   50},
xlabel style={font=\color{white!15!black}},
xlabel={Battery capacity},
ymin=0,
ymax=140,
ylabel style={font=\color{white!15!black}},
ylabel={Charging time},
axis background/.style={fill=white},
legend style={at={(0,1)}, anchor=north west, legend cell align=left, align=left, font=\scriptsize, draw=white!12!black}
]
\addplot [color=black, line width=1.0pt]
  table[row sep=crcr]{%
12.5	30.6818181818182\\
12.875	31.6022727272727\\
13.25	32.5227272727273\\
13.625	33.4431818181818\\
14	34.3636363636364\\
14.375	35.2840909090909\\
14.75	36.2045454545455\\
15.125	37.125\\
15.5	38.0454545454545\\
15.875	38.9659090909091\\
16.25	39.8863636363636\\
16.625	40.8068181818182\\
17	41.7272727272727\\
17.375	42.6477272727273\\
17.75	43.5681818181818\\
18.125	44.4886363636364\\
18.5	45.4090909090909\\
18.875	46.3295454545455\\
19.25	47.25\\
19.625	48.1704545454545\\
20	49.0909090909091\\
20.375	50.0113636363636\\
20.75	50.9318181818182\\
21.125	51.8522727272727\\
21.5	52.7727272727273\\
21.875	53.6931818181818\\
22.25	54.6136363636364\\
22.625	55.5340909090909\\
23	56.4545454545455\\
23.375	57.375\\
23.75	58.2954545454545\\
24.125	59.2159090909091\\
24.5	60.1363636363636\\
24.875	61.0568181818182\\
25.25	61.9772727272727\\
25.625	62.8977272727273\\
26	63.8181818181818\\
26.375	64.7386363636364\\
26.75	65.6590909090909\\
27.125	66.5795454545455\\
27.5	67.5\\
27.875	68.4204545454545\\
28.25	69.3409090909091\\
28.625	70.2613636363636\\
29	71.1818181818182\\
29.375	72.1022727272727\\
29.75	73.0227272727273\\
30.125	73.9431818181818\\
30.5	74.8636363636364\\
30.875	75.7840909090909\\
31.25	76.7045454545455\\
31.625	77.625\\
32	78.5454545454545\\
32.375	79.4659090909091\\
32.75	80.3863636363636\\
33.125	81.3068181818182\\
33.5	82.2272727272727\\
33.875	83.1477272727273\\
34.25	84.0681818181818\\
34.625	84.9886363636364\\
35	85.9090909090909\\
35.375	86.8295454545455\\
35.75	87.75\\
36.125	88.6704545454546\\
36.5	89.5909090909091\\
36.875	90.5113636363636\\
37.25	91.4318181818182\\
37.625	92.3522727272727\\
38	93.2727272727273\\
38.375	94.1931818181818\\
38.75	95.1136363636364\\
39.125	96.0340909090909\\
39.5	96.9545454545455\\
39.875	97.875\\
40.25	98.7954545454545\\
40.625	99.7159090909091\\
41	100.636363636364\\
41.375	101.556818181818\\
41.75	102.477272727273\\
42.125	103.397727272727\\
42.5	104.318181818182\\
42.875	105.238636363636\\
43.25	106.159090909091\\
43.625	107.079545454545\\
44	108\\
44.375	108.920454545455\\
44.75	109.840909090909\\
45.125	110.761363636364\\
45.5	111.681818181818\\
45.875	112.602272727273\\
46.25	113.522727272727\\
46.625	114.443181818182\\
47	115.363636363636\\
47.375	116.284090909091\\
47.75	117.204545454545\\
48.125	118.125\\
48.5	119.045454545455\\
48.875	119.965909090909\\
49.25	120.886363636364\\
49.625	121.806818181818\\
50	122.727272727273\\
};
\addlegendentry{plug-in charging}
\addplot [color=blue, dashed, line width=1.0pt]
  table[row sep=crcr]{%
12.5	2\\
12.875	2\\
13.25	2\\
13.625	2\\
14	2\\
14.375	2\\
14.75	2\\
15.125	2\\
15.5	2\\
15.875	2\\
16.25	2\\
16.625	2\\
17	2\\
17.375	2\\
17.75	2\\
18.125	2\\
18.5	2\\
18.875	2\\
19.25	2\\
19.625	2\\
20	2\\
20.375	2\\
20.75	2\\
21.125	2\\
21.5	2\\
21.875	2\\
22.25	2\\
22.625	2\\
23	2\\
23.375	2\\
23.75	2\\
24.125	2\\
24.5	2\\
24.875	2\\
25.25	2\\
25.625	2\\
26	2\\
26.375	2\\
26.75	2\\
27.125	2\\
27.5	2\\
27.875	2\\
28.25	2\\
28.625	2\\
29	2\\
29.375	2\\
29.75	2\\
30.125	2\\
30.5	2\\
30.875	2\\
31.25	2\\
31.625	2\\
32	2\\
32.375	2\\
32.75	2\\
33.125	2\\
33.5	2\\
33.875	2\\
34.25	2\\
34.625	2\\
35	2\\
35.375	2\\
35.75	2\\
36.125	2\\
36.5	2\\
36.875	2\\
37.25	2\\
37.625	2\\
38	2\\
38.375	2\\
38.75	2\\
39.125	2\\
39.5	2\\
39.875	2\\
40.25	2\\
40.625	2\\
41	2\\
41.375	2\\
41.75	2\\
42.125	2\\
42.5	2\\
42.875	2\\
43.25	2\\
43.625	2\\
44	2\\
44.375	2\\
44.75	2\\
45.125	2\\
45.5	2\\
45.875	2\\
46.25	2\\
46.625	2\\
47	2\\
47.375	2\\
47.75	2\\
48.125	2\\
48.5	2\\
48.875	2\\
49.25	2\\
49.625	2\\
50	2\\
};
\addlegendentry{battery swapping}
\end{axis}
\end{tikzpicture}%
\vspace*{-0.3in}
\caption{Vehicle charging time $t_s$ (minute) as a function of $C$.}
\label{ts_C}
\end{subfigure}
\caption{Market outcomes on vehicle charging under different battery capacities. Black lines present the results when charging vehicles; blue dashed lines show the results when swapping batteries.}
\label{vehicle_charging_C}
\end{figure}
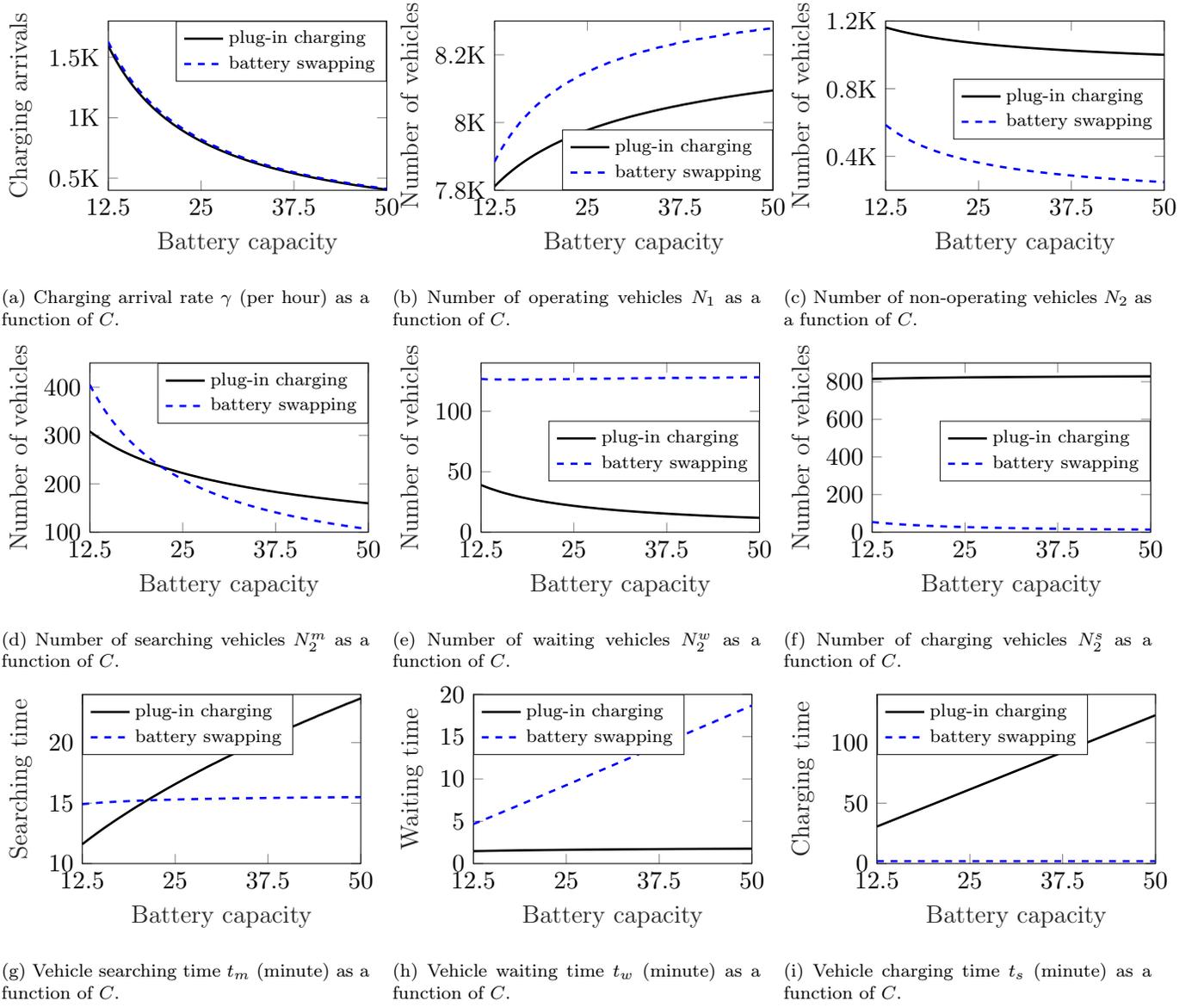

Figure \ref{vehicle_charging_C} shows the charging demand and the average time spent
on charging for TNC vehicles under the two charging strategies. Numerical results show that under both strategies of vehicle charging and battery swapping, as the battery capacity increases, the charging demand largely drops (Figure \ref{gamma_C}), and the number of operating vehicles increases (Figure \ref{N1_C}). The number of non-operating vehicles under battery swapping significantly decreases, while that under plug-in charging has a relatively small reduction (Figure \ref{N2_C}). More specifically, under plug-in charging, the number of searching vehicles reduces (Figure \ref{N2m_C}), the number of waiting vehicles decreases (Figure \ref{N2w_C}), and the number of charging vehicles remains almost steady (Figure \ref{N2s_C}). As for distinct time periods in the vehicle's charging process, we find that the vehicle searching time to stations increases with the battery capacity (Figure \ref{tm_C}). Vehicles also have a longer waiting time at charging stations (Figure \ref{tw_C}), and the vehicle charging time linearly increases with the battery capacity (Figure \ref{ts_C}). On the other hand, under battery swapping, the number of searching vehicles declines (Figure \ref{N2m_C}), whereas the number of waiting vehicles at stations slightly increases (Figure \ref{N2w_C}). The number of charging vehicles keeps decreasing with battery capacity and becomes close to zero (Figure \ref{N2s_C}). Correspondingly, the vehicle searching time has a slight increase (Figure \ref{tm_C}), the vehicle waiting time at battery swapping stations rises significantly (Figure \ref{tw_C}), and the vehicle charging time (i.e., the service time of battery swapping) still remains constant (Figure \ref{ts_C}).

The key findings regarding the impacts of battery capacity are outlined below:
\begin{itemize}
    \item When deploying charging stations for AMoD systems, the increased battery capacity results in a transformation of infrastructure deployment from {\em densely distributed small stations} to {\em sparsely distributed large stations} (Figure \ref{optimal_K_C} and Figure \ref{optimal_Q_C}).
    \item When deploying battery swapping stations for AMoD systems, the enlarged battery capacity lowers the supply of charging infrastructure (Figure \ref{optimal_K_C} and Figure \ref{optimal_Q_C}).
    \item Under both strategies of vehicle charging and battery swapping, improved battery capacity leads to a Pareto improvement: passengers enjoy higher surplus (Figure \ref{PS_C}), the TNC earns increased profit (Figure \ref{profit_C}), and the government produces higher social welfare (Figure \ref{SW_C}).
\end{itemize}

Analogous to the impacts of charging speed, the improved battery capacity reduces the need for charging infrastructures under both charging strategies. This is intuitive because increasing battery capacity extends vehicles' operating hours and directly reduces the demand for recharging electricity. The marginal benefit of offering charging infrastructures diminishes as the charging demand drops. Therefore, the government cuts down the total supply of charging infrastructures to save infrastructure costs. When planning battery swapping stations, due to its standard configuration, the government would offer fewer battery swapping stations for cost-saving as battery capacity enlarges (Figure \ref{optimal_K_C}). On the other hand, when deploying charging stations, although the total number of chargers declines (Figure \ref{infrastructure_cost_C}), the improved battery capacity transforms the infrastructure deployment from {\em densely distributed small stations} to {\em sparsely distributed large stations} (Figure \ref{optimal_K_C} and Figure \ref{optimal_Q_C}), which is contrary to the transformation caused by improved charging speed (Figure \ref{optimal_K_s} and Figure \ref{optimal_Q_s}). Similarly, this can be explained as a consequence of the trade-off between the vehicle searching time and the vehicle waiting time. As battery capacity increases, the charging demand largely drops (Figure \ref{gamma_C}), which lowers the arrival rate of vehicles at charging stations and alleviates the congestion at stations. However, a larger battery capacity also requires a longer time to charge the depleted battery (Figure \ref{ts_C}). This raises the service/occupied time of each electric vehicle charger and extends vehicles' waiting time for vacant chargers. In this case, expanding battery capacity breaks the balance achieved by infrastructure deployment: the vehicle waiting time at stations increases with battery capacity and becomes the bottleneck. To achieve the optimal trade-off, the government reduces the number of charging stations and deploys more chargers at each station, which leads to the transformation of infrastructure deployment from {\em densely distributed small stations} to {\em sparsely distributed large stations}. Under such transformation, although the vehicle searching time increases owing to the reduced density of charging stations (Figure \ref{tm_C}), the vehicle waiting time remains at a relatively low level (Figure \ref{tw_C}). Overall, given the reduced charging demand, both the number of searching vehicles (Figure \ref{N2m_C}) and waiting vehicles (Figure \ref{N2w_C}) decrease as battery capacity. This improves the fleet utilization rate and brings higher profit for the TNC.

\subsection{Comparison of plug-in charging and battery swapping} \label{comparison_strategies}

This section compares two distinct charging strategies under different infrastructure costs, charging speeds, and battery capacities. In particular, we fix the per charger cost ratio of battery swapping stations and charging stations as $r=\frac{\phi_{bs}}{\phi_{vc}}=5$, and decrease the per charger operating cost of battery swapping stations $\phi_{bs}$ from $40\$/\text{hour}$ to $20\$/\text{hour}$ to reflect the trend of reduced infrastructure costs over time. The bi-level optimization for vehicle charging and battery swapping will be solved and compared under distinct combinations of model parameters $\phi_{bs}$, $s$, and $C$.

\begin{figure}[ht!]
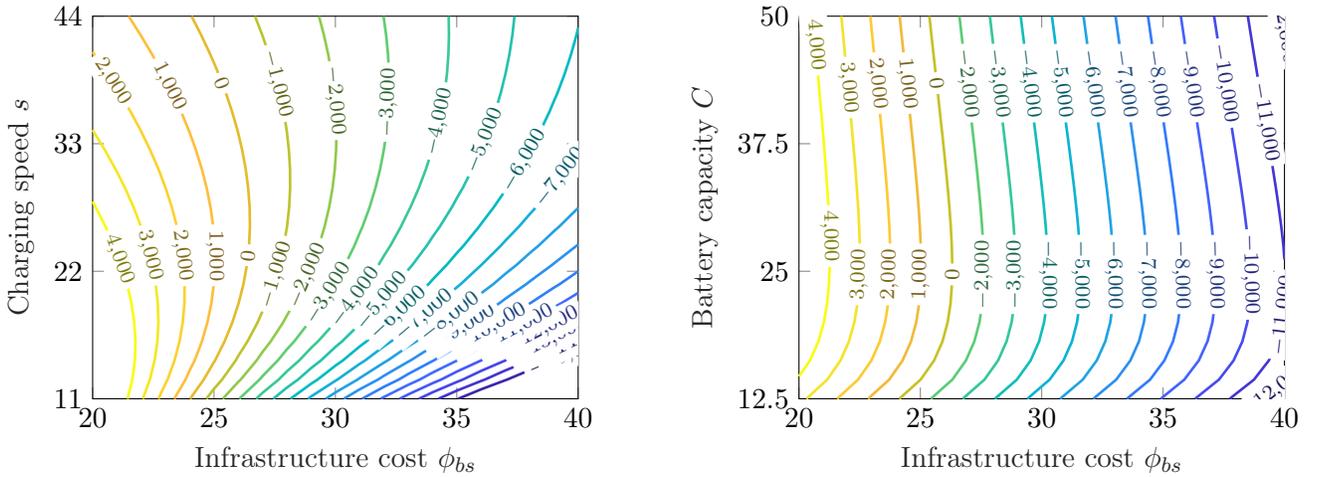

\centering
\begin{subfigure}[b]{0.48\textwidth}
\centering
\include{figures/SW_phi_s_fixed_ratio}
\vspace*{-0.3in}
\caption{The social welfare difference between battery swapping and vehicle charging under distinct $r$ and $s$.}
\label{SW_phi_s}
\end{subfigure}
\hfill
\begin{subfigure}[b]{0.48\textwidth}
\centering
\include{figures/SW_phi_C_fixed_ratio}
\vspace*{-0.3in}
\caption{The social welfare difference between battery swapping and vehicle charging under distinct $r$ and $C$.}
\label{SW_phi_C}
\end{subfigure}
\caption{The social welfare difference between battery swapping and vehicle charging under different cost ratios, charging speeds, and battery capacities.}
\label{comparison_strategies_fig}
\end{figure}

Figure \ref{comparison_strategies_fig} presents the difference in social welfare between battery swapping and plug-in charging under distinct technology conditions. Specifically, Figure \ref{SW_phi_s} shows the contour plot of social welfare difference between battery swapping and plug-in charging under distinct combinations of infrastructure costs and charging speeds. Figure \ref{SW_phi_C} illustrates the contour plot of social welfare differences under distinct infrastructure costs and battery capacities. Overall, several key findings are summarized below:

\begin{itemize}
    \item As the infrastructure cost reduces, the advantage of battery swapping over plug-in charging strengthens. There exists a break-even cost ratio below which battery swapping is superior to plug-in charging in promoting societal benefits.
    \item Under a fixed infrastructure cost, when the charging speed is either very small or very large, and advantage of battery swapping diminishes. Similarly, the advantage of battery swapping initially increases and then diminishes as battery capacity increases.
\end{itemize}

As the technology evolves and infrastructure develops, the marginal cost of supplying charging infrastructures reduces, which enables the government to deploy more charging stations and battery swapping stations for electrified AMoD services. The increased charging infrastructure supply improves the charging capability and reduces the vehicle searching time and vehicle waiting time for charging. However, the overlong charging time is still the major bottleneck in plug-in charging, which restricts the benefits of improving the charging supply. As the infrastructure costs reduce, the advantage of battery swapping over plug-in charging gradually emerges and becomes more and more significant. For example, in the numerical example with the charging speed of $22\text{kW}$ (Level 2 charging), the battery capacity of $25\text{kWh}$ (100 miles vehicle range), and the fixed cost ratio of 5, when the infrastructure cost of battery swapping $\phi_{bs}$ is below 25\text{\$/hour} (the corresponding infrastructure cost for plug-in charging $\phi_{vc}=5\text{\$/hour}$), battery swapping for electric AMoD services produces higher social welfare than plug-in charging. This provides technical guidance for the government on how to plan charging infrastructures according to distinct technological configurations.

{An interesting result is that the advantage of battery swapping initially strengthens as the charging speed/battery capacity improves but eventually diminishes as the charging speed/battery capacity further increases. We argue that this is again due to the difference in the flexibility of infrastructure deployment when planning charging stations and battery swapping stations. When the charging speed/battery capacity is relatively low, the government must supply a large amount of charging infrastructure (either charging stations or battery swapping stations) to accommodate the AMoD services. Given the high supply of charging infrastructure, the costs of charging infrastructure dominate social welfare, and the government expends much more infrastructure costs when deploying battery swapping stations than deploying charging stations (see Figure \ref{infrastructure_cost_s} and Figure \ref{infrastructure_cost_C}). As charging speed/battery capacity increases, the marginal benefit of supplying charging infrastructures reduces, and therefore the government lowers the charging infrastructure supply. The difference in infrastructure costs between the two charging strategies decreases, and the advantages of battery swapping over plug-in charging become increasingly significant. However, as mentioned before, the reduced charging infrastructure supply also negatively impacts and suppresses the TNC profit. When deploying charging stations, the government has the flexibility to determine the size of charging stations, by which the negative impacts of reduced charging infrastructure supply can be largely mitigated. On the other hand, when deploying battery swapping stations, the limited flexibility of infrastructure deployment magnifies the negative effects of reduced infrastructure supply. As charging speed/battery capacity increases, such negative effects gradually become the dominant factor and make battery swapping less advantageous compared to vehicle charging.}

\section{Conclusion} 

This paper investigates two prevalent charging solutions, plug-in charging and battery swapping, and the corresponding charging infrastructure planning for charging electric AMoD systems. An economic equilibrium model is developed to understand and reveal the operational mechanism of the TNC's electric AMoD system. The incentives of passengers, the operating/charging shift of TNC vehicles, and the balance of energy are captured, and vehicle charging processes at charging stations and battery swapping stations are characterized by queueing theory. A bi-level optimization framework is established to capture the intimate correlations between the government's charging infrastructure planning and the TNC's operation of the electric AMoD system. The upper-level social welfare maximization for the government and the lower-level profit maximization for the TNC are formulated, in which the structural difference of charging infrastructures and the flexibility of infrastructure deployment are considered. The government's optimal infrastructure deployment, the corresponding TNC's profit-maximizing operational strategies, and the market outcomes are obtained by solving the bi-level optimization problem.

We use the established bi-level optimization framework to investigate the impacts of involved charging technologies on the government's charging infrastructure planning and the TNC's operation of electric AMoD systems. We find that the evolution of charging technologies will reduce the need for charging infrastructures under both plug-in charging and battery swapping. In the planning of charging stations, the improved charging speed leads to a transformation of infrastructure deployment from {\em sparsely distributed large stations} to {\em densely distributed small stations}, while the enlarged battery capacity transforms the infrastructure deployment from {\em densely distributed small stations} to {\em sparsely distributed large stations}. The underlying reason is that the evolving charging technologies alter the trade-off between vehicle searching time and vehicle waiting time and promote the government to adjust the infrastructure deployment. We also show that the evolution of charging technologies has distinct impacts on the TNC's operation of electric AMoD systems under distinct charging solutions. The improved charging speed always leads to increased TNC profit under plug-in charging, while the TNC may get hurt when the charging speed is relatively high under battery swapping.
This demonstrates that the limited flexibility of infrastructure deployment of battery swapping stations restricts the potential of battery swapping for electric AMoD systems. We further compare plug-in charging and battery swapping for charging electric AMoD systems under distinct infrastructure costs. We find that there exists a break-even cost ratio at which battery swapping is equivalent to plug-in charging in promoting social welfare. 

This paper delivers a comprehensive economic analysis of the operation of the electric AMoD system and the corresponding charging infrastructure planning. Under the proposed bi-level optimization framework, we evaluate and compare two prevalent charging solutions plug-in charging and battery swapping separately, and investigate the impacts of evolving charging technologies on the charging infrastructure planning and the market outcomes at the aggregate level. One future extension is considering the synergy between plug-in charging and battery swapping for charging electric AMoD systems. Another future direction would be extending the current economic equilibrium model to a network equilibrium model and optimizing the zonal level charging infrastructure deployment to account for the spatial heterogeneity of AMoD services. We hope this work could stimulate more discussion in electrified, automated, and shared mobility for efficient and sustainable future transportation.

\section*{Acknowledgments}  {This research was supported by the Hong Kong Research Grants Council under project 26200420, project 16202922, and the National Science Foundation of China under project 72201225.}

\bibliographystyle{unsrt}
\bibliography{reference}

\section*{Appendix}

\subsection*{\bf{Appendix A: Proof of Lemma 2}}
{
To prove Lemma \ref{lemma_N_gamma}, we first note that given $\delta$, $s$, $t_s$, $C$, $K$, $Q$ and $V$, $P_V$ defined by (\ref{blocking_prob_charging}) and (\ref{blocking_prob_swapping}) and $t_w$ defined by (\ref{vehicle_waiting_time_charging}) and (\ref{vehicle_waiting_time_swapping}) are monotonically increasing functions of $\gamma$ due to the nature of M/M/Q/V queueing system and the mixed queueing system. Without loss of generality, we denote $P_V = \hat{P}_V(\gamma)$ and $t_w = \hat{t}_w(\gamma)$ defined by (\ref{blocking_prob_charging}) and (\ref{vehicle_waiting_time_charging}) under plug-in charging, and denote $P_V = \tilde{P}_V(\gamma)$ and $t_w = \tilde{t}_w(\gamma)$ defined by (\ref{blocking_prob_swapping}) and (\ref{vehicle_waiting_time_swapping}) under battery swapping. Below we show that in equilibrium (\ref{constraints_charging}), $N_1$ and $N$ can be uniquely determined as a function of $\gamma$ from (\ref{conservation_N_charging})-(\ref{energy_balance_charging}), i.e., (i) of Lemma \ref{lemma_N_gamma}. The proof of (ii) is analogous and is thus omitted.

First, combining (\ref{N2m_charging}), (\ref{blocking_prob_charging}) and (\ref{vehicle_searching_time_charging}) derives $N_2^m$ as a function of $\gamma$ (denoted as $\hat{N}_2^m(\gamma)$):
\begin{equation} \label{N2m_matching_charging}
    \hat{N}_2^m(\gamma) = \frac{\gamma B}{\sqrt{K (1-\hat{P}_V(\gamma))}} .
\end{equation}
Substituting (\ref{N2m_matching_charging}), (\ref{N2w_charging}) and (\ref{N2s_charging}) into (\ref{conservation_N2_charging}) further determines $N_2$ as a function of $\gamma$ (denoted as $\hat{N}_2(\gamma)$):
\begin{equation} \label{N2_function_gamma}
    \hat{N}_2(\gamma) = \frac{\gamma B}{\sqrt{K (1-\hat{P}_V(\gamma))}} + \gamma \hat{t}_w(\gamma) + \gamma \hat{t}_s ,
\end{equation}
where $\hat{t}_s = \frac{(1-\delta)C}{s}$ is a constant. Besides, by (\ref{energy_balance_charging}), we have:
\begin{equation} \label{N_energy_charging}
    N = \frac{\gamma(1-\delta)C}{l} + N_2^w + N_2^s .
\end{equation}
Substituting (\ref{N2w_charging}) and (\ref{N2s_charging}) into (\ref{N_energy_charging}) determines $N$ as a function of $\gamma$ (denoted as $\hat{N}(\gamma)$):
\begin{equation} \label{N_function_gamma}
    \hat{N}(\gamma) = \frac{\gamma(1-\delta)C}{l} + \gamma \hat{t}_w(\gamma) + \gamma \hat{t}_s.
\end{equation}
Finally, subtracting (\ref{N2_function_gamma}) from (\ref{N_function_gamma}) derives $N_1$ as a function of $\gamma$:
\begin{equation}
    \hat{N}_1(\gamma) = \frac{\gamma(1-\delta)C}{l} - \frac{\gamma B}{\sqrt{K (1-\hat{P}_V(\gamma))}} .
\end{equation}
This completes the proof.
}

\subsection*{\bf{Appendix B: Proof of Lemma 3}}
Given $\delta$, $s$, $t_s$, $C$, $K$, $Q$ and $V$, note that $\hat{P}_V(\gamma)$ and $\tilde{P}_V(\gamma)$ monotonically increase with $\gamma$ according to the nature of the queueing models. Besides, $\forall \gamma \geq 0$, $\hat{P}_V(\gamma) \geq 0$ and $\tilde{P}_V(\gamma) \geq 0$, and $\lim_{\gamma \to \infty} \hat{P}_V(\gamma)=1$ and $\lim_{\gamma \to \infty} \tilde{P}_V(\gamma)=1$. Below we show that if $K > \left(\frac{Bl}{(1-\delta)C}\right)^2$, there exist a unique $\hat{\gamma}_0 > 0$ and $\hat{\gamma}_* \in (0,\hat{\gamma}_0)$ such that $\hat{N}_1(\hat{\gamma}_0)=0$, and $\forall \gamma \in (0,\hat{\gamma}_0)$, $\hat{N}_1(\gamma)>0$ and $\hat{N}_1(\gamma)<\hat{N}_1(\hat{\gamma}_{*})$, i.e., (i) of Lemma \ref{lemma_property_N1}. The proof of (ii) is similar and is thus omitted.

We first rewrite $\hat{N}_1(\gamma)$ as:
\begin{equation} \label{N1_rewrite}
    \hat{N}_1(\gamma) = \gamma \left(\frac{(1-\delta)C}{l}-\frac{B}{\sqrt{K(1-\hat{P}_V(\gamma))}} \right) .
\end{equation}
Let $\hat{B}(\gamma)$ be the term inside the bracket in (\ref{N1_rewrite}), i.e., $\hat{B}(\gamma)=\frac{(1-\delta)C}{l}-\frac{B}{\sqrt{K(1-\hat{P}_V(\gamma))}}$. To guarantee the existence of positive $\hat{N}_1(\gamma)$, we need the maximum of $\hat{B}(\gamma)$ to be greater than zero, i.e., $\max_{\gamma \geq 0} \hat{B}(\gamma)>0$. $\hat{B}(\gamma)$ is a decreasing function of $\gamma$ since $\hat{P}_V(\gamma)$ monotonically increases with $\gamma$. The maximum of $\hat{B}(\gamma)$ is achieved at $\gamma=0$. Therefore, it suffices to show that $\hat{B}(0)=\frac{(1-\delta)C}{l}-\frac{B}{\sqrt{K(1-\hat{P}_V(0))}}>0$, which is equivalent to $K > \left(\frac{Bl}{(1-\delta)C}\right)^2$.

When $K > \left(\frac{Bl}{(1-\delta)C}\right)^2$, there exist $\gamma>0$ such that $\hat{N}_1(\gamma)>0$. Let $\hat{\gamma}_0=\hat{B}^{-1}(0)$, where $\hat{B}^{-1}(\cdot)$ is the inverse function of $\hat{B}(\gamma)$. We have $\hat{N}_1(\gamma)>0$ when $0<\gamma<\hat{\gamma}_0$ due to the monotonicity of $\hat{B}(\gamma)$. Besides, since $\hat{N}_1(0)=0$ and $\hat{N}_1(\hat{\gamma}_0)=0$, by the continuity of $\hat{N}_1(\gamma)$, there must exist a $\hat{\gamma}_{*} \in (0,\hat{\gamma}_0)$ such that $\hat{N}_1(\hat{\gamma}_{*}) > \hat{N}_1(\gamma)$, $\forall \gamma \in (0, \hat{\gamma}_0)$. This completes the proof.

\subsection*{\bf{Appendix C: Proof of Proposition 1}}

{
To guarantee the existence of market equilibrium, we need to show that there exists strictly positive $p_f$, $w^c$, $w^v$, $\lambda$, $\gamma$, $t_m$, $t_w$, $P_V$, $N$, $N_1$, $N_2$, $N_2^m$, $N_2^w$, $N_2^s$ satisfying (\ref{constraints_charging}) and (\ref{constraints_swapping}), respectively. We have shown in Lemma \ref{lemma_wc_gamma} that there exists strictly positive $w^c$ and $w^v$ satisfying (\ref{passenger_waiting_time_charging})-(\ref{conservation_N1_charging})/(\ref{passenger_waiting_time_swapping})-(\ref{conservation_N1_swapping}) if $N_1 \geq \left(\sqrt[3]{2}+\sqrt[3]{\frac{1}{4}}\right)\left(A\lambda \right)^{\frac{2}{3}}+\lambda\tau$, and $w^c=w^c(\lambda,N_1)$ and $w^v=w^v(\lambda,N_1)$ can uniquely derived from (\ref{passenger_waiting_time_charging})-(\ref{conservation_N1_charging})/(\ref{passenger_waiting_time_swapping})-(\ref{conservation_N1_swapping}). Besides, we have shown in Lemma \ref{lemma_N_gamma} that $P_V=\hat{P}_V(\gamma)$, $t_w = \hat{t}_w(\gamma)$, $N_1=\hat{N}_1(\gamma)$ ,$N_2=\hat{N}_2(\gamma)$, $N=\hat{N}(\gamma)$ can be derived from (\ref{conservation_N_charging})-(\ref{energy_balance_charging}), and $P_V = \tilde{P}_V(\gamma)$, $t_w=\tilde{t}_w(\gamma)$, $N_1=\tilde{N}_1(\gamma)$, $N_2=\tilde{N}_2(\gamma)$, $N=\tilde{N}(\gamma)$ can be derived from (\ref{conservation_N_swapping})-(\ref{energy_balance_swapping}), and $\hat{N}_1(\gamma)=\tilde{N}_1(\gamma)$. In this case, the market equilibrium (\ref{constraints_charging}) and (\ref{constraints_swapping}) are equivalent to the following:
\begin{subnumcases}{\label{equivalent_equilibrium_charging}}
    p_f = F_p^{-1}\left(\frac{\lambda}{\lambda_0}\right) - \alpha w^c\left(\lambda,\hat{N}_1(\gamma)\right) \\
    w^v = w^v\left(\lambda, \hat{N}_1(\gamma)\right) \\
    N = \frac{\gamma(1-\delta)C}{l} + \gamma \hat{t}_w(\gamma) + \gamma \hat{t}_s \\
    N_1 = \hat{N}_1(\gamma) = \frac{\gamma (1-\delta) C}{l} - \frac{\gamma B}{\sqrt{K(1-\hat{P}_V(\gamma))}} \\
    N_2 = \frac{\gamma B}{\sqrt{K(1-\hat{P}_V(\gamma))}} + \gamma \hat{t}_w(\gamma) + \gamma \hat{t}_s \\
    N_2^m = \frac{\gamma B}{\sqrt{K(1-\hat{P}_V(\gamma))}} \\
    N_2^w = \gamma \hat{t}_w(\gamma) \\
    N_2^s = \gamma \hat{t}_s  \\
    P_V = \hat{P}_V(\gamma) \\
    t_w = \hat{t}_w(\gamma)  \\
    \hat{N}_1(\gamma) \geq \left(\sqrt[3]{2}+\sqrt[3]{\frac{1}{4}}\right)\left(A\lambda \right)^{\frac{2}{3}}+\lambda\tau
\end{subnumcases}
and
\begin{subnumcases}{\label{equivalent_equilibrium_swapping}}
    p_f = F_p^{-1}\left(\frac{\lambda}{\lambda_0}\right) - \alpha w^c\left(\lambda,\tilde{N}_1(\gamma)\right) \\
    w^v = w^v\left(\lambda, \tilde{N}_1(\gamma)\right) \\
    N = \frac{\gamma(1-\delta)C}{l} + \gamma \tilde{t}_w(\gamma) + \gamma \tilde{t}_s \\
    N_1 = \tilde{N}_1(\gamma) = \frac{\gamma (1-\delta) C}{l} - \frac{\gamma B}{\sqrt{K(1-\tilde{P}_V(\gamma))}} \\
    N_2 = \frac{\gamma B}{\sqrt{K(1-\tilde{P}_V(\gamma))}} + \gamma \tilde{t}_w(\gamma) + \gamma \tilde{t}_s \\
    N_2^m = \frac{\gamma B}{\sqrt{K(1-\tilde{P}_V(\gamma))}} \\
    N_2^w = \gamma \tilde{t}_w(\gamma) \\
    N_2^s = \gamma \tilde{t}_s  \\
    P_V = \tilde{P}_V(\gamma) \\
    t_w = \tilde{t}_w(\gamma)  \\
    \tilde{N}_1(\gamma) \geq \left(\sqrt[3]{2}+\sqrt[3]{\frac{1}{4}}\right)\left(A\lambda \right)^{\frac{2}{3}}+\lambda\tau
\end{subnumcases}

where $F_p^{-1}(\cdot)$ is inverse function of $F_p(\cdot)$. For the profit maximization problem (\ref{optimalpricing_charging}) and (\ref{optimalpricing_swapping}), we can equivalently treat $\lambda$ and $\gamma$ as decision variables/free variables. To prove Proposition \ref{prop_equilibrium_charging_swapping}, it is equivalent to show that there exists $\lambda \in (0,\lambda_0)$ and $\gamma \in (0,+\infty)$ such that
\begin{subnumcases}{\label{condition_charging}}
    p_f = F_p^{-1}\left(\frac{\lambda}{\lambda_0}\right) - \alpha w^c\left(\lambda,\hat{N}_1(\gamma)\right)>0 \label{condition_pf_charging} \\
    N_1 = \frac{\gamma (1-\delta) C}{l} - \frac{\gamma B}{\sqrt{K(1-\hat{P}_V(\gamma))}} > 0 \label{condition_N1_charging} \\
    N = \frac{\gamma(1-\delta)C}{l} + \gamma \hat{t}_w(\gamma) + \gamma \hat{t}_s > 0 \label{condition_N_charging} \\
    N_2 = \frac{\gamma B}{\sqrt{K(1-\hat{P}_V(\gamma))}} + \gamma \hat{t}_w(\gamma) + \gamma \hat{t}_s > 0 \label{condition_N2_charging} \\
    N_2^m = \frac{\gamma B}{\sqrt{K(1-\hat{P}_V(\gamma))}} > 0 \label{condtion_N2m_charging} \\
    N_2^w = \gamma \hat{t}_w(\gamma)>0 \label{condition_N2w_charging} \\
    N_2^s = \gamma \hat{t}_s>0 \label{condition_N2s_charging} \\
    P_V = \hat{P}_V(\gamma) > 0 \label{condition_PV_charging} \\
    t_w = \hat{t}_w(\gamma)>0 \label{condition_tw_charging} \\
    \frac{\gamma (1-\delta) C}{l} - \frac{\gamma B}{\sqrt{K(1-\hat{P}_V(\gamma))}} \geq \left(\sqrt[3]{2}+\sqrt[3]{\frac{1}{4}}\right)\left(A\lambda \right)^{\frac{2}{3}}+\lambda\tau \label{condition_wc_charging}
\end{subnumcases}
and
\begin{subnumcases}{\label{condition_swapping}}
    p_f = F_p^{-1}\left(\frac{\lambda}{\lambda_0}\right) - \alpha w^c\left(\lambda,\tilde{N}_1(\gamma)\right)>0 \label{condition_pf_swapping} \\
    N_1 = \frac{\gamma (1-\delta) C}{l} - \frac{\gamma B}{\sqrt{K(1-\tilde{P}_V(\gamma))}} > 0 \label{condition_N1_swapping} \\
    N = \frac{\gamma(1-\delta)C}{l} + \gamma \tilde{t}_w(\gamma) + \gamma \tilde{t}_s > 0 \label{condition_N_swapping} \\
    N_2 = \frac{\gamma B}{\sqrt{K(1-\tilde{P}_V(\gamma))}} + \gamma \tilde{t}_w(\gamma) + \gamma \tilde{t}_s > 0 \label{condition_N2_swapping} \\
    N_2^m = \frac{\gamma B}{\sqrt{K(1-\tilde{P}_V(\gamma))}} > 0 \label{condtion_N2m_swapping} \\
    N_2^w = \gamma \tilde{t}_w(\gamma)>0 \label{condition_N2w_swapping} \\
    N_2^s = \gamma \tilde{t}_s>0 \label{condition_N2s_sawpping} \\
    P_V = \tilde{P}_V(\gamma) > 0 \label{condition_PV_swapping}
    t_w = \tilde{t}_w(\gamma)>0 \label{condition_tw_swapping} \\
    \frac{\gamma (1-\delta) C}{l} - \frac{\gamma B}{\sqrt{K(1-\tilde{P}_V(\gamma))}} \geq \left(\sqrt[3]{2}+\sqrt[3]{\frac{1}{4}}\right)\left(A\lambda \right)^{\frac{2}{3}}+\lambda\tau \label{condition_wc_swapping}
\end{subnumcases}
Note that (\ref{condition_charging}) and (\ref{condition_swapping}) have the same form. Below we show the existence of market equilibrium (\ref{constraints_charging}) by justifying (\ref{condition_charging}), i.e., (i) of Proposition \ref{prop_equilibrium_charging_swapping}. The proof of (ii) is analogous and is thereby omitted. 

First, the positivity of $P_V$ and $t_w$, i.e., condition (\ref{condition_PV_charging}) and (\ref{condition_tw_charging}) is naturally satisfied due to the nature of queueing models. Further, it is obvious that $N$, $N_2$, $N_2^m$, $N_2^w$, $N_2^s>0$, i.e., conditions (\ref{condition_N_charging})-(\ref{condition_N2s_charging}) given $t_w>0$ and $P_V \in [0,1)$. Therefore, $\forall \gamma \in (0,\infty)$, $N$, $N_2$, $N_2^m$, $N_2^w$, $N_2^s$, $t_w>0$, i.e., conditions (\ref{condition_N_charging})-(\ref{condition_tw_charging}).

Next, we show that there exists $\gamma>0$ such that $N_1>0$, i.e., condition (\ref{condition_N1_charging}). Based on the result of Lemma (\ref{lemma_property_N1}), $\hat{N}_1(\gamma)>0$ when $\gamma \in (0,\hat{\gamma}_0)$.

Next, we show that there exists $\gamma\in\left(0,\hat{\gamma}_0 \right)$ and $\lambda\in(0,\lambda_0)$, such that $\frac{\gamma (1-\delta) C}{l} - \frac{\gamma B}{\sqrt{K(1-\hat{P}_V(\gamma))}} \geq \left(\sqrt[3]{2}+\sqrt[3]{\frac{1}{4}}\right)\left(A\lambda \right)^{\frac{2}{3}}+\lambda\tau$, i.e., condition (\ref{condition_wc_charging}). Note that the LHS of (\ref{condition_wc_charging}) is $\hat{N}_1(\gamma)=\frac{\gamma (1-\delta) C}{l} - \frac{\gamma B}{\sqrt{K(1-\hat{P}_V(\gamma))}}$. Based on the result of Lemma \ref{lemma_property_N1}, it initially increases with $\gamma$ and then decreases within $(0,\hat{\gamma}_0)$, and the the maximum of $\hat{N}_1(\gamma)$ is achieved at $\hat{\gamma}_{*}$. The RHS of (\ref{condition_wc_charging}) is a continuous and strictly increasing function of $\lambda$. By continuity, there must exist $0<\gamma_1<\hat{\gamma}_{*}<\gamma_2<\hat{\gamma}_0$ and $0<\lambda_1<\lambda_0$ such that $\hat{N}_1(\gamma)\geq \left(\sqrt[3]{2}+\sqrt[3]{\frac{1}{4}}\right)\left(A\lambda \right)^{\frac{2}{3}}+\lambda\tau$ when $\gamma \in \left(\gamma_1,\gamma_2\right)$ and $\lambda \in (0,\lambda_1)$.

Next, we show that there exists $\gamma \in \left(\gamma_1,\gamma_2 \right)$ and $\lambda \in (0,\lambda_1)$ such that $p_f>0$, i.e., condition (\ref{condition_pf_charging}). Note that $w^c(\lambda,N_1)$ is derived from (\ref{passenger_waiting_time_charging})-(\ref{conservation_N1_charging}). Substituting (\ref{conservation_N1_charging}) into (\ref{passenger_waiting_time_charging}), we have:
\begin{equation} \label{wc_function}
    w^c = \frac{A}{\sqrt{N_1-\lambda\mu -\lambda w^c}} .
\end{equation}
Plugging (\ref{wc_function}) into (\ref{condition_pf_charging}) leads to:
\begin{equation}
    p_f = F_p^{-1}\left(\frac{\lambda}{\lambda_0}\right) - \frac{\alpha A}{\sqrt{N_1 - \lambda\tau - \lambda w^c}} .
\end{equation}
By assumption, $F_p(\cdot)$ is a strictly decreasing function, thus $F_p^{-1}(\frac{\lambda}{\lambda_0})$ monotonically decreases with $\lambda$. Therefore, $p_f$ is a decreasing function with respect to $\lambda$ and an increasing function with respect to $N_1$. The maximum of $p_f$ is
\begin{equation}
    {p_f^{\text{max}} = \lim_{\lambda \to 0} F_p^{-1}\left(\frac{\lambda}{\lambda_0}\right)-\frac{\alpha A}{\sqrt{{N_1^{\text{max}}}}-\lambda\tau-\lambda w^c}} = F_p^{-1}(0) - \frac{\alpha A}{\sqrt{{N_1^{\text{max}}}}} ,
\end{equation}
where ${N_1^{\text{max}}}$ is the maximum of $N_1$. Based on the result of Lemma \ref{lemma_property_N1}, the maximum of $\hat{N}_1(\gamma)$ is $\hat{N}_1^\text{max}=\hat{N}_1(\hat{\gamma}_{*})$. To guarantee the existence of positive $p_f$, we need that
\begin{equation}
    p_f^{\text{max}} = F_p^{-1}(0) - \frac{\alpha A}{\sqrt{{\hat{N}_1^{\text{max}}}}} > 0,
\end{equation}
which is equivalent to $F_p\left(\frac{\alpha A}{\sqrt{{\hat{N}_1^{\text{max}}}}}\right)>0$. Under this assumption, there must exist $0<\lambda_2<\lambda_1$ and $(\gamma_3,\gamma_4) \subset (\gamma_1,\gamma_2)$ such that $p_f>0$ when $\lambda \in (0,\lambda_2)$ and $\gamma \in \left(\gamma_3,\gamma_4\right)$.

Finally, we take the intersection of all the derived conditions to show the feasibility of (\ref{condition_charging}). Given the assumptions that $K > \left(\frac{Bl}{(1-\delta)C}\right)^{2}$ and $F_p\left(\frac{\alpha A}{\sqrt{{\hat{N}_1^{\text{max}}}}}\right)>0$, (\ref{condition_charging}) hold when $\lambda \in (0,\lambda_2)$ and $\gamma\in\left(\gamma_3,\gamma_4 \right)$. This completes the proof.
}

\subsection*{\bf{Appendix D: Calibration of model parameters $A$ and $B$}}
{
Both the passenger waiting time function (\ref{func_passenger_waiting_time}) and the vehicle searching time function (\ref{func_vehicle_searching_time}) follow the same 'square root law'. Here we calibrate the model parameters of the vehicle searching time function (\ref{func_vehicle_searching_time}) based on realistic geographic and traffic data of New York City to determine the value of the parameter $B$. The calibration of the passenger waiting time function (\ref{func_passenger_waiting_time}) follows the same approach, and the value of parameter $A$ is taken to be the same as that of $B$.

Consider the numerical simulation for New York City, where a fleet of vehicles and a number of charging/battery swapping stations are randomly generated and uniformly distributed across the city. We assume that vehicles travel to the nearest charging/battery swapping stations for energy top-up and calculate the average travel time to the charging/battery swapping stations. We vary the number of charging/battery swapping stations and calculate vehicles' average travel time to charging/battery swapping stations. The scatter plots in Figure \ref{fig:calibration} show the calculated average travel time under different supplies of charging infrastructures, which graphically indicates that the vehicle searching time follows the 'square root law'. Therefore, we fit the vehicle searching time function (\ref{func_vehicle_searching_time}) based on the simulation data and calibrate that $B=230$. The solid line in Figure \ref{fig:calibration} shows the fitted vehicle searching time function $t_m=\frac{230}{\sqrt{K}}$ with $R^2=0.9943$. It indicates that the proposed searching time function (\ref{func_vehicle_searching_time}) is a good fit for the vehicle searching time for charging under the nearest-neighbor matching.
\begin{figure}[h]
    \centering
    \includegraphics[width=0.67\linewidth]{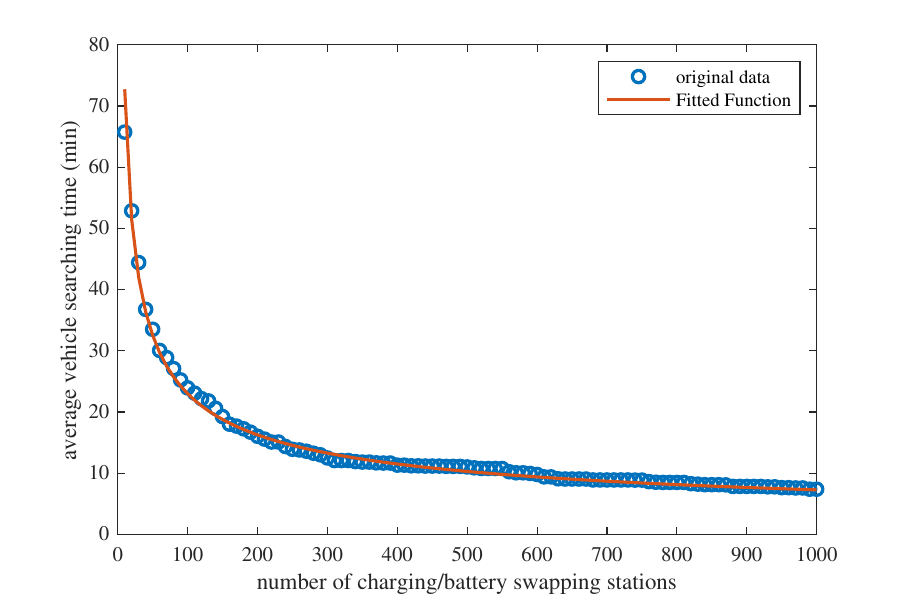}
    \caption{(1) scatter plots: simulated average vehicle searching time under distinct supplies of charging/battery swapping stations; (2) solid line: fitted vehicle searching time function $t_m=\frac{230}{\sqrt{K}}$ with $R^2=0.9943$.}
    \label{fig:calibration}
\end{figure}
}

\end{document}